%% file: sequential2022_article.tex
\documentclass[final,onefignum,onetabnum,letterpaper]{siamonline190516}

\input{sequential2022_shared}

\usepackage{lineno}

\ifpdf
\hypersetup{
	pdftitle={Xiao2022sequential},
 	pdfauthor={Xiao, Glaubitz, Gelb, and Song}
}
\fi

\begin{document}

\maketitle


\begin{abstract}
	A new algorithm is developed to jointly recover a temporal sequence of images from noisy and under-sampled Fourier data. Specifically we consider the case where each data set  is missing vital information that prevents its (individual) accurate recovery. Our new method is designed to restore the missing information in each individual image by ``borrowing'' it from  the other images in the sequence.  As a result, {\em all} of the individual reconstructions yield improved accuracy.   The use of high resolution Fourier edge detection methods is essential to our algorithm.   In particular, edge information is obtained directly from the Fourier data which leads to an accurate coupling term between data sets.  Moreover,  data loss is largely avoided as coarse reconstructions are not required to process inter- and intra-image information.  Numerical examples are provided to demonstrate the accuracy, efficiency and robustness of our new method.
\end{abstract}

\begin{keywords}
	Sequential image recovery, joint sparsity, edge detection, Fourier data
\end{keywords}

\begin{AMS}
	655F22, 65K10, 68U10
\end{AMS}

\begin{DOI}
    \href{https://doi.org/10.1007/s10915-022-01850-7}{10.1007/s10915-022-01850-7}
\end{DOI}

\input{1_introduction} 
\input{2_background}

\input{3_proposed}

\input{4_numerical} 
\input{5_summary}



\bibliographystyle{siamplain}
\bibliography{literature}

\end{document}

%% file: sequential2022_shared.tex
\usepackage{lipsum}
\usepackage{amsfonts}
\usepackage{epstopdf}
\ifpdf
  \DeclareGraphicsExtensions{.eps,.pdf,.png,.jpg}
\else
  \DeclareGraphicsExtensions{.eps}
\fi

\usepackage{datetime}
\newdateformat{monthyeardate}{%
	\monthname[\THEMONTH] \THEDAY, \THEYEAR}

\usepackage{academicons}
\usepackage{xcolor}

\usepackage{amsmath}
\allowdisplaybreaks
\usepackage{amssymb}
\usepackage{commath}
\usepackage{mathtools}
\usepackage{bbm}
\usepackage{bm}

\usepackage{color}
\usepackage{graphicx}
\usepackage[small]{caption}
\usepackage{subcaption}

\usepackage{relsize}
\usepackage{adjustbox}
\usepackage{algorithm}
\usepackage[noend]{algpseudocode}
\usepackage{booktabs}
\usepackage{tikz}
\usepackage{pifont}
\usetikzlibrary{bayesnet} 

\usepackage{enumitem}
\setlist[enumerate]{leftmargin=.5in}
\setlist[itemize]{leftmargin=.5in}

\newsiamremark{remark}{Remark}
\newsiamremark{example}{Example}
\newsiamremark{hypothesis}{Hypothesis}
\crefname{hypothesis}{Hypothesis}{Hypotheses}
\newsiamthm{claim}{Claim}

\headers{Sequential image recovery}{Y.\ Xiao, J.\ Glaubitz, A.\ Gelb, and G.\ Song}

\title{Sequential image recovery from noisy and under-sampled Fourier data\thanks{\monthyeardate\today 
\funding{This work is partially supported by the NSF grants DMS \#1521661 (GS), DMS \#1912685 (AG), DMS \#1939203 (GS), AFOSR grant \#FA9550-18-1-0316 (AG), and ONR MURI grant \#N00014-20-1-2595 (AG).}
\corresponding{Anne Gelb} 
}}
\author{ 
Yao Xiao\thanks{Department of Mathematics, Dartmouth College, Hanover, 03755, NH, USA (\email{yao.xiao.gr@dartmouth.edu}, \email{jan.glaubitz@dartmouth.edu}, \email{anne.e.gelb@dartmouth.edu})}
\and 
Jan Glaubitz\footnotemark[2]
\and 
Anne Gelb\footnotemark[2]
\and
Guohui Song\thanks{Department of Mathematics and Statistics, Old Dominion University, Norfolk, 23529, VA, USA (\email{gsong@odu.edu})}
}

\usepackage{amsopn}

\usepackage{comment}

\usepackage[normalem]{ulem}

\DeclareMathOperator{\diag}{diag}

\DeclareMathOperator*{\argmin}{arg\,min}

\newcommand{\intd}{\, \mathrm{d}} 

\newcommand{\N}{\mathbb{N}}

\DeclarePairedDelimiter\ceil{\lceil}{\rceil}



%% file: 1_introduction.tex
\section{Introduction}
\label{sec:intro1}

We are interested in jointly recovering  a temporal sequence of images from  noisy and under-sampled Fourier data.  We assume that each data collection is obstructed in some way, for example by some occlusion in the image scene or by the sensor not being able to access some Fourier bandwidths of data, either of which can prevent accurate individual image recovery.  It is intuitive to suggest that complementary information from other images in the sequence may be  ``borrowed'' to improve any of the individual recoveries.  Change, however, such as the insertion, deletion, translation or rotation of an object within the scene, may have occurred between the sequential data collections.   In this case  the ``borrowed'' information would be incompatible with the missing information we are seeking to replace, which may  lead to greater error in recovering each individual image.

From the above discussion it is evident that joint recovery of a temporal sequence of images requires incorporating (1) an accurate and robust individual image recovery process; (2) a coupling term that provides essential missing information from other data collections; and (3) a technique to determine change within the sequence of images so that the coupling term is appropriately used.   This investigation proposes a new joint recovery method  for a temporal sequence of images given noisy and under-sampled Fourier data that effectively incorporates these three critical components. To ensure accuracy, efficiency, and robustness, our new technique uses the Fourier concentration factor edge detection method, first developed in \cite{gelb1999detection}, to both improve the individual recovery process as well as to determine any internal changes to the image over time.  Moreover, our method avoids significant data loss by {\em directly} using the given Fourier data to determine the edges. Specifically, no initial image reconstruction is needed to determine change in the underlying image. This is especially pertinent in our data acquisition model, in which bands of Fourier data might be unobservable.  Furthermore, we employ the variance based joint sparsity (VBJS) method, \cite{adcock2019jointsparsity,gelb2019reducing,scarnati2019accelerated}, for each individual image reconstruction. The VBJS method is a weighted $\ell_1$ regularization method which uses the concentration factor edge detection method to determine the spatially varying weights in the image recovery. 
It was shown in \cite{gelb2019reducing} to be robust to noise and missing (even false) information.

As will be demonstrated in our numerical results, combining these three components yields a joint recovery process that reduces the effects of data loss while simultaneously enhancing each individual reconstruction. Our method is furthermore computationally efficient and robust to both occlusion and missing information.  

The rest of this paper is organized as follows. 
Section \ref{sec:introduction} provides necessary background information for the problem set up.  Our new algorithm is introduced in Section \ref{sec:l2-reg}.   A series of numerical experiments  in Section \ref{sec:test} demonstrates the efficacy and robustness of our new method.  
Finally, some concluding thoughts are offered in Section \ref{sec:summary}.

%% file: 2_background.tex
\section{Background Information}
\label{sec:introduction}

Let  $f_1,\dots,f_J$ be a sequence of two-dimensional images taken at times $t_1,\dots,t_J$ of the same scene in $[0,1]^{2}$. 
The corresponding  $(2N+1)^2$ Fourier samples are given by
\begin{equation}\label{eq:fourcoeff}
	\hat{f}^j_{k,l} = \int^1_0\int^1_0f_j(x,y)e^{-i2\pi (kx+ly)}\intd x \intd y, 
	\quad -N \le k,l \le N, \quad j = 1,\dots,J.
\end{equation}
As noted in the introduction, we  are interested in applications where the data are collected in the Fourier domain. 
We seek to simultaneously recover each ${\mathbf f}_j \in \mathbb R^{(2N+1) \times (2N+1)}$, the values of $f_j$ at uniform grid points 
\begin{equation}\label{eq:grid} 
	x_\mu = \mu\Delta x, \quad 
	y_\nu = \nu\Delta y, \quad 
	\Delta x=\Delta y=\frac{1}{2N}, \quad 
	0 \le \mu,\nu \le 2N.
\end{equation}
We further assume that each of the $J$ data sets  contains additive independent and identically distributed (iid) zero-mean complex Gaussian noise given by
\[
	\boldsymbol{\eta}_j= ( \eta_{k,l} )_{k,l=-N}^{N} \quad
\text{s.t.} \quad \eta_{k,l}\sim\mathbb{C}\mathcal{N}\left(0, \sigma^2_j\right), \quad
j=1,\dots,J,
\]
where $\sigma^2_j$ is the variance.
Finally, for modeling  purposes, we will employ the discrete Fourier transform which we denote by $F$. 
Hence we obtain the data model relating the Fourier coefficients in \eqref{eq:fourcoeff} to each underlying image as
\begin{equation}
\label{eq:forwardmodel}
	\mathbf{b}_j = F\mathbf{f}_j+\boldsymbol{\eta}_j, \quad j = 1,\dots,J.
\end{equation} 
We will also consider the case where different symmetric bands of Fourier data are missing from the acquired data. 
That is, for each $j$ there is a band $\mathcal{K}_j$  for which ${\bf b}_j(k,l)$ is not available. 
In each of our numerical experiments we chose $J = 6$ with the zeroed out data corresponding to the symmetric intervals given by  
\begin{equation}
\label{eq:bandzero}
\mathcal{K}_j = [\pm(10 + J(j-1)),\pm(35+J(j-1)].
\end{equation}
The forward model in \eqref{eq:forwardmodel} is accordingly modified as
\begin{equation}
\label{eq:forwardmodel_j}
\mathbf{b}_j = F_j{\mathbf{f}_j}+\boldsymbol{\eta}_j, \quad j = 1,\dots,J, 
\end{equation}
where $F_j$ contains only zeroes in the columns and rows corresponding to \eqref{eq:bandzero}.

\subsection{Edge Detection}
\label{sec:edge_detection}

Images are often assumed to be sparse in their corresponding edge domain.
Compressive sensing (CS) algorithms, see e.\,g.\ the classical papers \cite{candes2006robust,candes2006stable,candes2006near,donoho2006compressed},  are designed to exploit such sparsity.
In this investigation we use the weighted $\ell_1$ regularization CS algorithm, see e.\,g.\ \cite{candes2008enhancing,chartrand2008iteratively,daubechies2008iteratively,Liu_2012,langer2017automated}, realized by the VBJS method, \cite{adcock2019jointsparsity,gelb2019reducing,scarnati2019accelerated}, as the main reconstruction tool for each individual image. 
The weights are constructed so that the  $\ell_1$ penalty is spatially dependent, with a heavier penalty where one can reasonably assume that the underlying image is close to zero in the sparse domain. 
By contrast, a small penalty is applied where the image is presumably supported in the sparse domain.  In this way we see that edge detection is critical to produce the proper spatially varying weights for the penalty term in each individual recovery.

\begin{figure}[h!]
\centering
\begin{subfigure}[b]{0.22\textwidth}
\includegraphics[width=\textwidth,height=\textwidth]{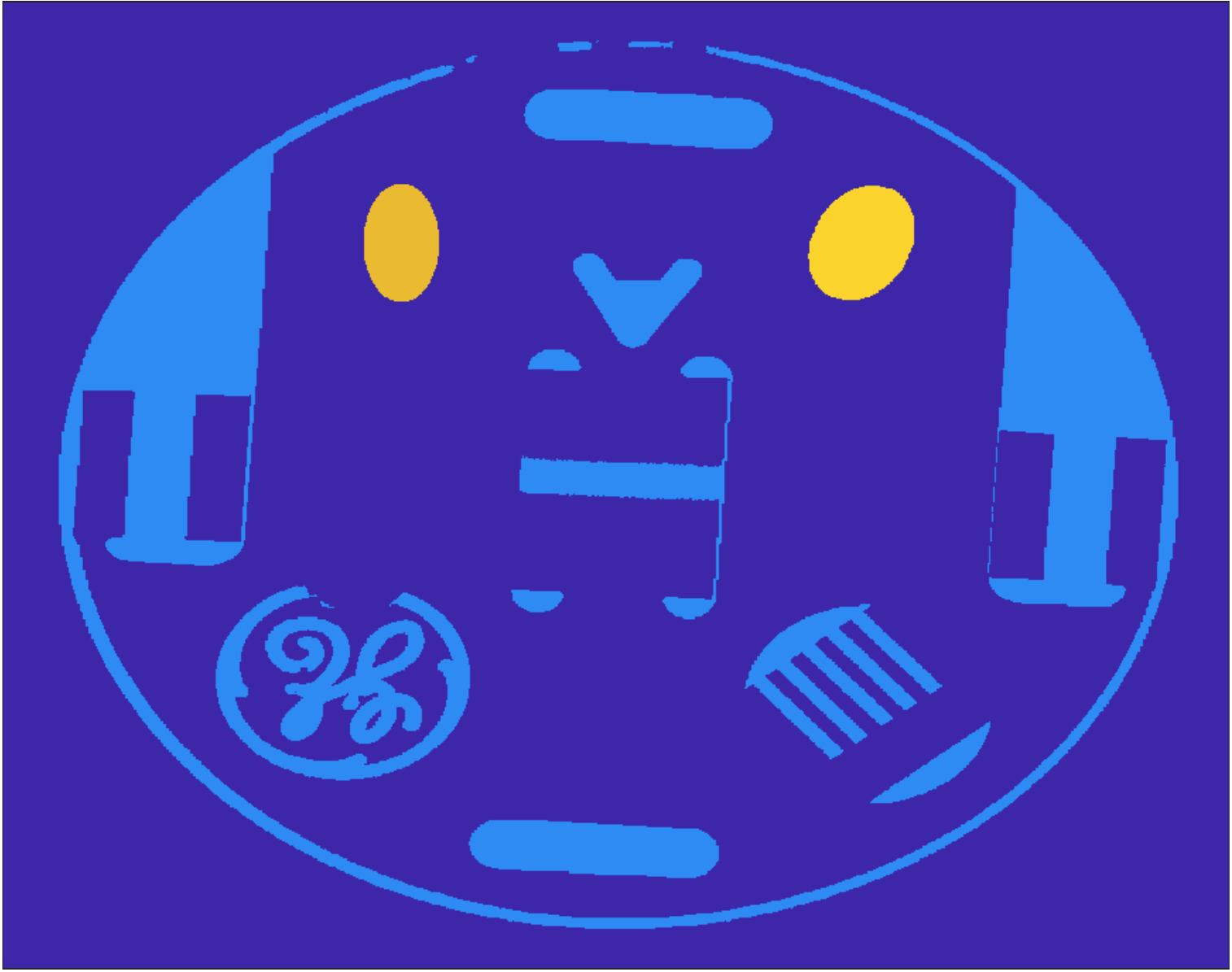}
\end{subfigure}
~
\begin{subfigure}[b]{0.22\textwidth}
\includegraphics[width=\textwidth,height=\textwidth]{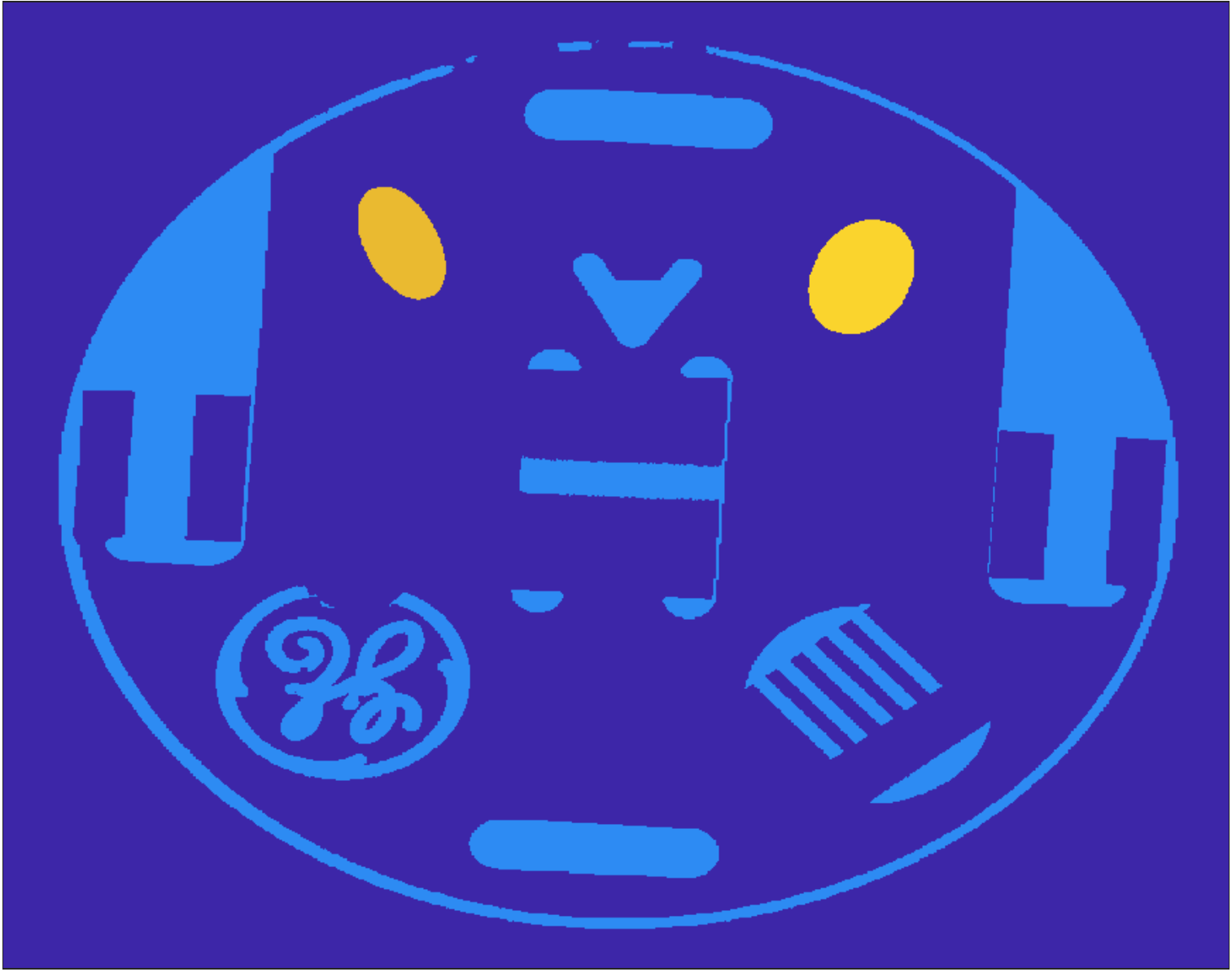}
\end{subfigure}
~
\begin{subfigure}[b]{0.22\textwidth}
\includegraphics[width=\textwidth,height=\textwidth]{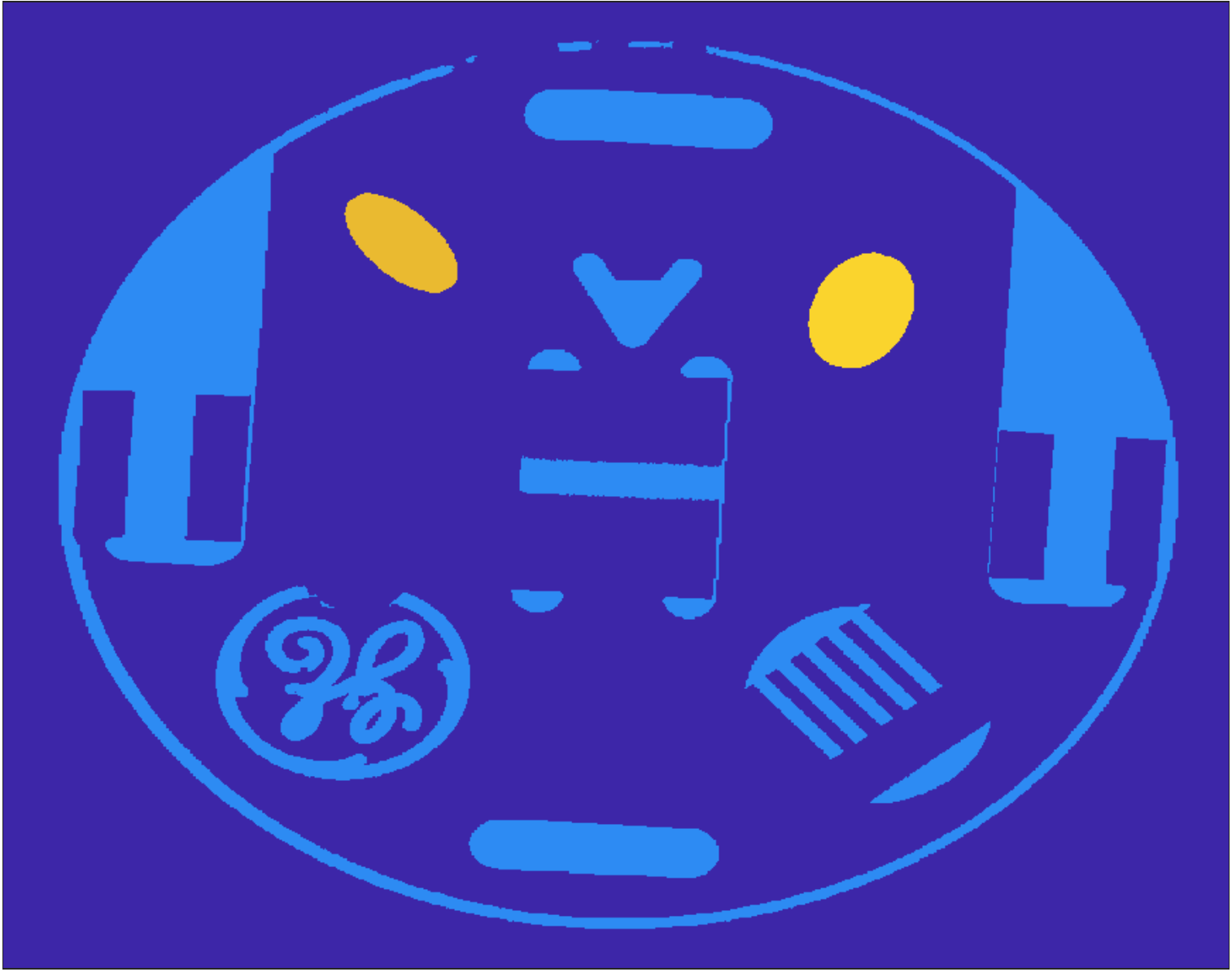}
\end{subfigure}
~
\begin{subfigure}[b]{0.22\textwidth}
\includegraphics[width=\textwidth,height=\textwidth]{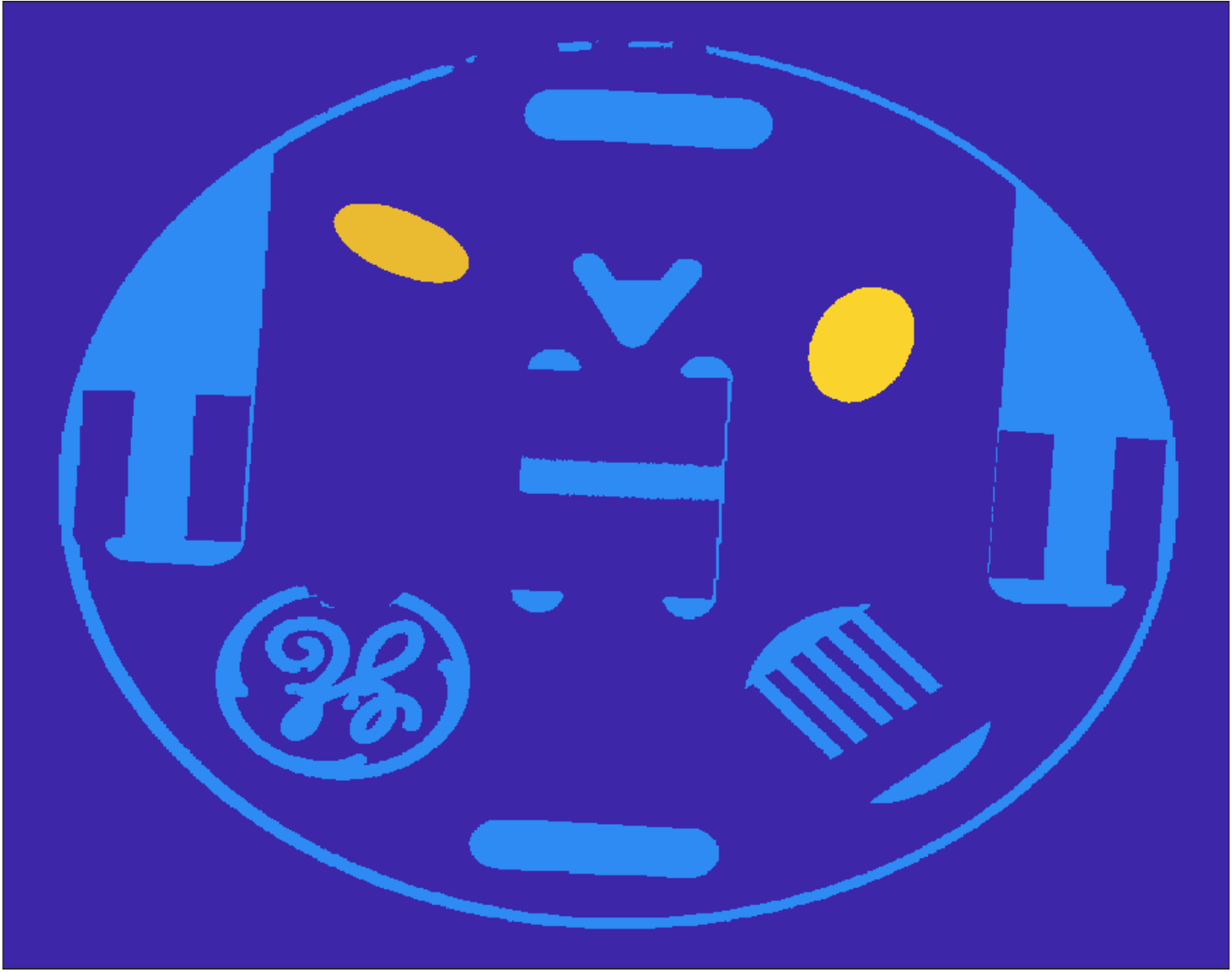}
\end{subfigure}
\\ 
\begin{subfigure}[b]{0.22\textwidth}
\includegraphics[width=\textwidth,height=\textwidth]{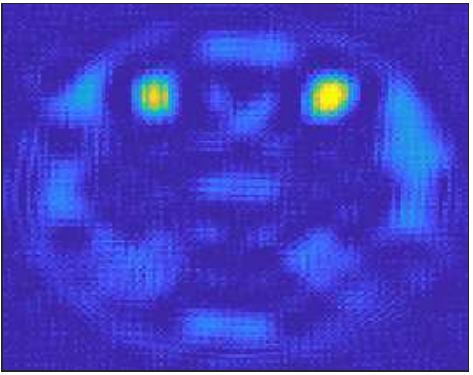}
\end{subfigure}
~
\begin{subfigure}[b]{0.22\textwidth}
\includegraphics[width=\textwidth,height=\textwidth]{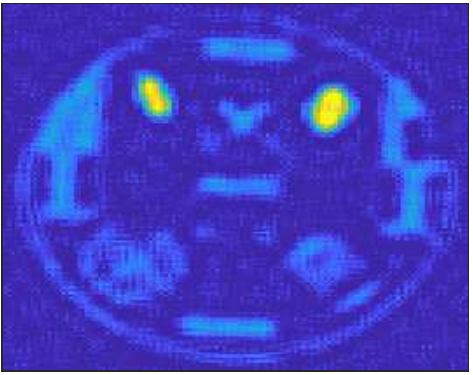}
\end{subfigure}
~
\begin{subfigure}[b]{0.22\textwidth}
\includegraphics[width=\textwidth,height=\textwidth]{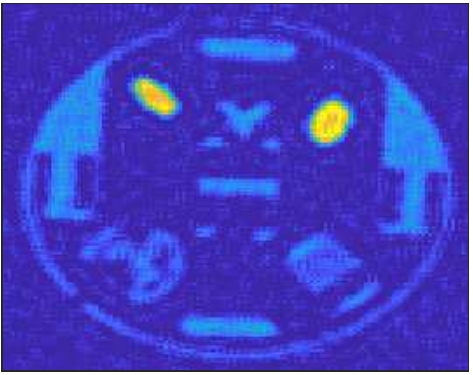}
\end{subfigure}
~
\begin{subfigure}[b]{0.22\textwidth}
\includegraphics[width=\textwidth,height=\textwidth]{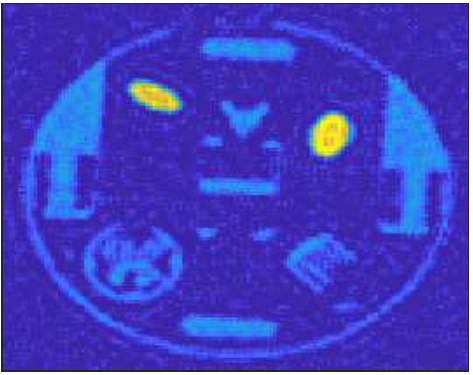}
\end{subfigure}
\\
\begin{subfigure}[b]{0.22\textwidth}
\includegraphics[width=\textwidth,height=\textwidth]{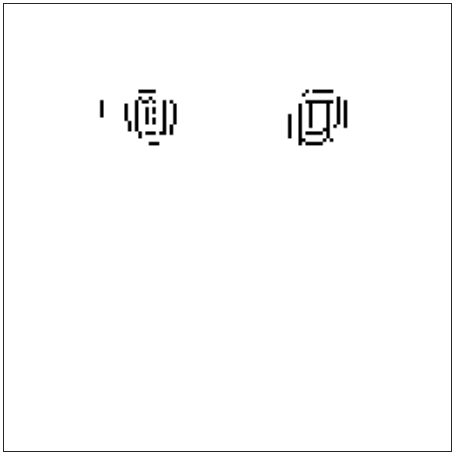}
\end{subfigure}
~
\begin{subfigure}[b]{0.22\textwidth}
\includegraphics[width=\textwidth,height=\textwidth]{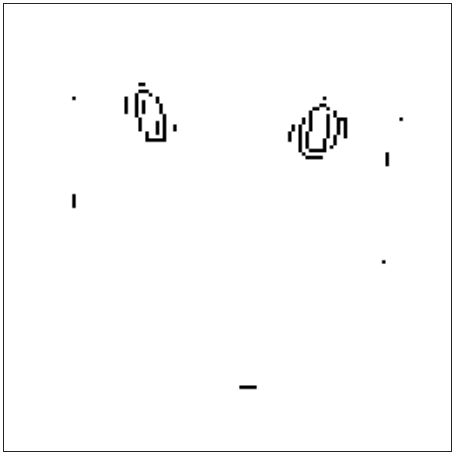}
\end{subfigure}
~
\begin{subfigure}[b]{0.22\textwidth}
\includegraphics[width=\textwidth,height=\textwidth]{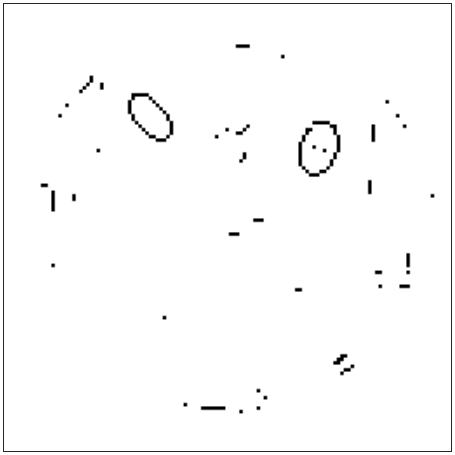}
\end{subfigure}
~
\begin{subfigure}[b]{0.22\textwidth}
\includegraphics[width=\textwidth,height=\textwidth]{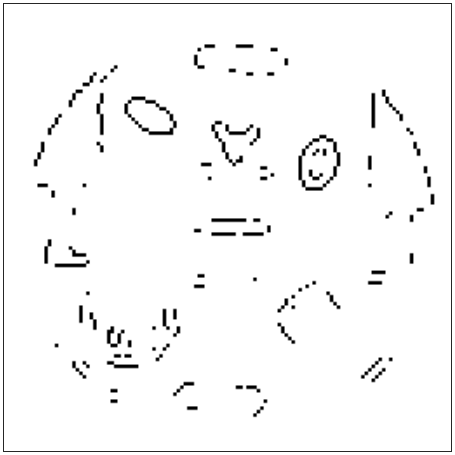}
\end{subfigure}
\caption{True underlying sequence of images. Objects that translate or rotate  between the sequential data collections are added to the underlying structure (top row); Individual CS recovery for each sequential image using \eqref{eq:l1_reg} (middle row); Edge masks determined using Sobel's method on the CS reconstruction (bottom row). In each case we are given under-sampled noisy Fourier data in \eqref{eq:forwardmodel_j} with SNR $= 2$dB.  There is also an obstruction in each of the images (see Figure \ref{fig:rec_GE}).}
\label{fig:sobel-edges}
\end{figure}

When data are acquired in the image domain, classical edge detection methods such as Canny’s operator, \cite{canny1986computational}, or Sobel’s method, \cite{sobel19683x3,jin2009edge, gao2010improved}, are commonly used to generate edge maps, which can then subsequently be used to detect change between images.  As is illustrated in Figure \ref{fig:sobel-edges} for  Sobel's method,  these methods may not perform well when the data are acquired indirectly, especially when occlusions are present in the underlying images.
The substantial information loss visible in the middle row for the standard CS recovery, \eqref{eq:l1_reg} is mainly due to the missing Fourier data in each band $\mathcal{K}_j$ given in \eqref{eq:bandzero}. This makes accurate edge detection difficult for the Sobel method (bottom row). As expected, losing information in the interval $\mathcal{K}_1$ is most problematic since it causes significant blur in the image reconstruction.

We will return to this example in Section \ref{sec:test} and demonstrate how using the concentration factor edge detection method described in Section \ref{sec:concfac} which works directly on the given Fourier data prevents information loss in the construction of the edge maps. As already stated this is needed both for the spatially dependent weights in each individual image recovery as well as to calculate the difference between images necessary to determine what inter-image information is to be ``borrowed''.

\subsection{Concentration factor edge detection}
\label{sec:concfac}

The concentration factor edge detection method determines the edges of piecewise smooth functions from finitely sampled Fourier data, \cite{gelb1999detection,gelb2000detection}. The method can be thought of as a  bandpass filter approximation that ``concentrates'' at the singular support of the underlying image. 
We briefly describe it below. 

For ease of presentation, we first consider  the method for one-dimensional data.  Let $f \colon\allowbreak  [0,1]\to\mathbb{R}$  be a piecewise analytic function with $\mathcal M$ distinct jump discontinuities located at $\{\xi_\mu\}_{\mu = 1}^\mathcal{M}$.  Suppose we are given the first $2N+1$  continuous Fourier coefficients,
\[\hat{f}_{k}=\int^1_0f(x)e^{-i2\pi kx}dx, \quad k = -N,\dots,N.\]
The corresponding jump function of $f$ is defined as\footnote{We will use jump function and edge function interchangeably throughout our exposition.} 
\begin{equation}\label{eq:1Djump}
	[f](x)=f(x^+)-f(x^-),
\end{equation} 
where $f(x^+)$ and $f(x^-)$ denote the right- and left-hand side limit of $f$ at $x$, respectively.
An equivalent formulation is given by
\begin{equation}
\label{eq:jumpfunction2}
[f](x) = \sum_{\mu = 1}^\mathcal{M}[f](\xi_\mu) I_{\xi_\mu}(x),
\end{equation}
where the indicator function $I_{\xi_\mu}(x)$ has value $1$ at $x = {\xi_\mu}$ and $0$
everywhere else.
Given $\{\hat{f}_k\}_{k = -N}^N$ we can now define the {\em concentration factor} edge detection method as
\begin{equation}
        S^\sigma_Nf(x)= i\sum_{1\leq\abs{k}\leq N}\text{sgn}(k)\sigma\left(\frac{\abs{k}}{N}\right)\hat{f}_ke^{i2\pi k x}.
        \label{eq:concsum1D}
   \end{equation}
It was shown in \cite{gelb1999detection,gelb2000detection} that  $S_N^\sigma f(x) \rightarrow [f](x)$ for admissible concentration factors $\sigma\left(\frac{\abs{k}}{N}\right)$.

\subsubsection{Edge detection for images}
\label{sec:2Dedgedetection}
Defining a jump in two dimensions is not as straightforward. 
The number of discontinuities is infinite, which would make the sum analogous to \eqref{eq:jumpfunction2} undefined.
Specifically, if we consider the parameterized function $g_{\boldsymbol{u}}(t) = f(\boldsymbol{x} + t \boldsymbol{u})$ and then use the jump height $[g_{\boldsymbol{u}}](0)$ to define the jump height $[f](\boldsymbol{x})$, we could obtain {\em different} values depending on the direction of $\boldsymbol{u}$. To circumvent this issue, a reasonable alternative for imaging problems is to consider only the normal direction with respect to the edge curve, yielding
\begin{definition}\label{def:2Djump}
Assume the discontinuities of $f$ form a finite number of closed and smooth edge curves $\Gamma_\mu$, $1\le \mu\le \mathcal{M}$. For a discontinuity point $(x,y)\in \Gamma_\mu$, we let ${\boldsymbol n}(x,y)$ be the normal direction at $(x,y)$ with respect to the edge curve $\Gamma_\mu$, and consider the corresponding one-dimensional  function $g(t) = f( (x,y) + t {\boldsymbol{n}}(x,y))$, $t\in \mathbb{R}$. We then define the jump function 
\begin{align}\label{eq:2Djump}
[f](x,y) = [g](0), \quad (x,y)\in \Gamma_\mu, 1\le \mu\le \mathcal{M}.
\end{align} 
\end{definition}

We seek to compute a corresponding matrix $G_j\in\mathbb R^{(2N+1) \times (2N+1)}$ 
for the edge function $[f_j](x,y)$ defined in Definition \ref{def:2Djump} for each underlying image $\mathbf f_j$, $j=1,\dots,J$. 

\begin{remark}
\label{rem:2Dedge}
In some applications it is sufficient to apply a  line-by-line approach of \eqref{eq:concsum1D} to two-dimensional functions  to obtain an approximation to $[f](x,y)$, \cite{ArchibaldGelb}.  Here, however, the noise and obstruction in \eqref{eq:forwardmodel_j} cause too much loss of accuracy when using the  line-by-line approach. 
In particular, in the extreme case, where edges in the underlying image are parallel with the Cartesian coordinates, an edge may be missed completely since in one direction the image is smooth. 
A more detailed discussion of this issue is provided in \cite{martinez2014edge} where the problem was mitigated by first rotating the image some small amount to ensure that the edges did not ``line up'' with the Cartesian coordinates.  A dimension by dimension  approach was then used in the rotated coordinate system.  
In this investigation we employ the procedure introduced in \cite{adcock2019jointsparsity}, which also employs rotation but does not use a  line-by-line approach.  It is briefly described below.
\end{remark}

\subsubsection{Concentration factor method using rotation}
\label{sec:concrotate}
Suppose we are given uniformly spaced Fourier samples \eqref{eq:fourcoeff} for a piecewise constant image with smooth edge curves.  Based on \eqref{eq:2Djump}, we seek to recover $[f](x,y)I_{\Gamma}(x,y)$ such that $I_\Gamma$ is the indicator function defined over a single edge curve.\footnote{As in the one-dimensional case, it is straightforward to show that the final approximation of the edge mask given in \eqref{eq:regularized_edge2} holds for multiple edge curves, see e.g. \cite{adcock2019jointsparsity,gelb2017detecting,viswanathan2012iterative}.} Following \cite{adcock2019jointsparsity} 
we begin by  parameterizing the curve $\Gamma$ as
\[x = u(s),\; y=v(s),\quad s\in[a,b],\]
where $u$ and $v$ are smooth functions defined over an arbitrary region $[a,b]$.  If we then let $\theta(s)$ denote the normal direction of each point $(u(s),v(s))$, we can parameterize the points around $\Gamma$ as
\begin{equation}\label{eq:parameterize}
    x = u(s)+r\cos\theta(s),\; y=v(s)+r\sin\theta(s),
\end{equation}
where $r \in [-\epsilon,\epsilon]$ and $\epsilon>0$ is sufficiently small. 

Since the indicator function $I_\Gamma$ has support of measure zero, we regularize it as  $h(\frac{r}{\epsilon})$ which is narrow in $[-\epsilon,\epsilon]$ and satisfies $h(0) = 1$.  While there are various choices for $h(\frac{r}{\epsilon})$, based on its simplicity and success shown in \cite{gelb2017detecting}, here we use
\begin{equation}\label{eq:gaussian}
h\left(\frac{r}{\epsilon}\right) =\mbox{exp}\biggl(-5\left(\frac{r}{\epsilon}\right)^2\biggr).
\end{equation}
We can now write the corresponding regularized edge function to $[f](x,y)I_\Gamma(x,y)$ as 
\begin{equation}
\label{eq:regedgefunc}  
H(x,y)= [f](u(s),v(s)) h\left(\frac{r}{\epsilon}\right).
\end{equation}

To relate the given Fourier data in \eqref{eq:fourcoeff} for each image in the sequence to the regularized edge function in \eqref{eq:regedgefunc},  we define the partial sum approximation of  $H(x,y)$ as
\begin{equation}\label{eq:regularized_edge}
    {H^j}(x,y) = \sum_{k = -N}^N\sum_{l= -N}^N\hat H^j(k,l)e^{2\pi i(kx+ly)},\quad j = 1,\dots,J.
\end{equation}
From \eqref{eq:regedgefunc} we are able to approximate the coefficients $\hat{H}^j(k,l)$ in \eqref{eq:regularized_edge} (after discarding higher order terms) as
\begin{equation}\label{eq:parameterize_Jhat}
    \tilde{H}^j(k,l)\approx\int^b_a[f_j](u(s),v(s))e^{-2\pi i(ku+lv)}(v^\prime\cos\theta-u^\prime\sin\theta)\epsilon\hat h(\epsilon(k\cos\theta+l\sin\theta))ds,
\end{equation} 
where $\hat{h}(\cdot)$ are the Fourier coefficients of $h(\frac{r}{\epsilon})$.  Substituting \eqref{eq:parameterize} into \eqref{eq:fourcoeff} and again eliminating higher order terms we obtain
\begin{equation}\label{parameterize_fhat}
    \hat f^j(k,l)\approx\int^b_a \frac{[f_j](u(s),v(s))(v^\prime\cos\theta-u^\prime\sin\theta)}{2\pi i(k\cos\theta+l\sin\theta)}e^{-2\pi i(ku+lv)}ds.
\end{equation}
There is no linear relationship between $\hat{f}^j(k,l)$ and {$\tilde{H}^j(k,l)$} since $\theta = \theta(s)$.  However,  we can {\em fix} $\theta = \theta_m$ for each $j$, yielding
\[\tilde H^j_{\theta_m}(k,l)=2\pi i(k\cos\theta_m+l\sin\theta_m)\epsilon\hat h(\epsilon(k\cos\theta_m+l\sin\theta_m))\hat f^j(k,l).\]
Hence to construct $M$ edge masks for {\em each} image in the sequence from given Fourier data we define $M$ rotation angles
\begin{equation}
    \label{eq:rotation-angles}
    \theta_m=\frac{\pi (m-1)}{M-1}, \quad m = 1,\dots, M,
\end{equation}
and compute $M$ Fourier partial sum approximations, 
\begin{equation}\label{eq:regularized_edge2}
    {H^j_{\theta_m}}(x,y) = \sum_{k = -N}^N\sum_{l= -N}^N\tilde H^j_{\theta_m}(k,l)e^{2\pi i(kx+ly)},\quad m = 1,\dots, M, \quad j = 1,\dots,J.
\end{equation}
Finally, the temporal sequence {\em edge masks} {${\tilde G}_j \in \mathbb{R}^{(2N+1) \times (2N+1)}$} are generated by averaging the results in \eqref{eq:regularized_edge2}, yielding 
\begin{equation}\label{eq:edge_mask}
    \tilde G_j(x_\mu,y_\nu)  = \frac{1}{M}\sum^{M}_{m=1}H^j_{\theta_m}(x_\mu,y_\nu),\quad \mu, \nu  = 0, 1, \dots, 2N, \quad j = 1,\dots, J.\end{equation}

We note that \eqref{eq:edge_mask} is used to recover {\em inter}-signal information, in particular to determine differences (or change) in the sequential images.  The approximations $H^j_{\theta_m}(x,y)$ in \eqref{eq:regularized_edge2} are also used to extract {\em intra-}signal information, specifically to improve the accuracy of each of the individual reconstructions. 
In this regard we employ a weighted $\ell_1$ regularization which we now describe.

\subsection{Individual image recovery} 
\label{sec:preliminary}
At the core of our new image recovery method is the import of ``missing'' information from other parts of the time sequenced data acquisitions, which will be incorporated into each individual image recovery. The initial recovery only uses {\em intra-}image information, and is designed to exploit the presumed sparsity of the image in the edge domain.  For this purpose we use the weighted $\ell_1$ regularization method, see e.g.~\cite{candes2008enhancing,chartrand2008iteratively, daubechies2008iteratively,Liu_2012, xie2014reweighted}, as realized by the VBJS approach, \cite{adcock2019jointsparsity, gelb2019reducing, scarnati2019accelerated}.\footnote{To be clear, the methods cited here are primarily {\em re-weighted} $\ell_1$ regularization methods (and often referred to as such), since the weights are iteratively adapted. By contrast the VBJS method does not iteratively adapt the weights, so it is a weighted $\ell_1$ regularization scheme.}

\subsubsection{Compressive Sensing}
\label{sec:SMV}
As already discussed, since the underlying image is sparse in the edge domain, it is appropriate to use CS algorithms, \cite{candes2006robust,candes2006stable,candes2006near,donoho2006compressed}, for each individual recovery.
Using the forward model given in \eqref{eq:forwardmodel_j} for each $j = 1,\dots,J$, the standard CS approach is to solve the unconstrained minimization problem 
\begin{equation}
\label{eq:l1_reg}
    \tilde{\mathbf{f}}_j = \argmin_{\mathbf s}\left\{
    \norm{F_j\mathbf s-\mathbf b_j}_2^2
    + \mu\norm{\mathcal L\mathbf s}_1\right\}. 
\end{equation}
The first term measures the difference between the forward model and the given data and is often referred to as the fidelity term.  The second term is known as the penalty or regularization term and expresses the prior belief that the underlying signal is sparse in some domain, for example the edge or gradient domain. Thus the sparsifying transform matrix ${\mathcal L}$ is designed to promote such sparsity.  As the images considered in our investigation are piecewise constant, we choose $\mathcal L$ as the first order difference (TV) operator. 

The choice of regularization parameter $\mu$ in \eqref{eq:l1_reg} has been the subject of many investigations, and is inherently problem dependent, \cite{landi2008lagrange,wen2011parameter,yang2015tv,gong2019adaptive}.   As we are primarily concerned with being able to compare our new method to the ``best'' case scenario for \eqref{eq:l1_reg}, we use Algorithm \ref{algo:EM_l1_reg} to experiment with various choices for $\mu$ and then pick the best solution according to the mean squared error (MSE), (see e.g. \cite{vogel2002Computational} for precedent). The algorithm requires us to choose the number of experiments to conduct, $K_{max}$, and a way to sample regularization parameter $\mu$.  To this end we note that typically $\mu$ is chosen to reflect both the SNR given in \eqref{eq:SNR} and the regularity in the signal.  For instance, $\mu$ might be chosen to approximate $\frac{\sigma}{{\xi}}$,  where $\sigma$ is the the standard deviation of the complex noise $\boldsymbol\eta$ in \eqref{eq:forwardmodel_j} and $\xi$ is the standard deviation of data in the transformed domain under $\mathcal{L}$ is used in \cite{sanders2020effective}.  Since $\sigma$ and $\xi$ are explicitly known in our test problems, to obtain the ``best'' case scenario we simply pick samples from $\mathcal{N}\left(\frac{\sigma}{\xi},\sigma\right)$.   For convenience we use $K_{max} = 10$, and note that larger $K_{max}$ did not improve the results in any of our experiments.

\begin{algorithm}[h!]
\caption{Image recovery via standard $\ell_1$-regularization}
\label{algo:EM_l1_reg}
	\hspace*{\algorithmicindent} \textbf{Input:}  Measurements $\{\mathbf{b}_j\}_{j = 1}^J$, forward operators $\{F_j\}_{j = 1}^J$, variances $\sigma^2$ and $\xi^2$ for the model \eqref{eq:forwardmodel_j}, and the sparse transform operator $\mathcal{L}$.\\
	\hspace*{\algorithmicindent} \textbf{Choose:} Number of iterations $K_{max}$.\\
	\hspace*{\algorithmicindent} \textbf{Output:}  Reconstructions $\check{\mathbf{f}}_j$, $j = 1,\dots,J$.
\begin{algorithmic}[1]
\For{$j = 1$ to $J$}
\For{ $i = 1$ to $K_{max}$}

\State{Sample candidate regularization parameter $\mu_k$ from $\mathcal{N}(\frac{\sigma}{\xi},\sigma)$.} 

\State{Calculate  $\tilde{\mathbf{f}}_{j_i}$ from \eqref{eq:l1_reg}. }

\State{Compute $\tilde{\mathbf{b}}_{j_i}=F_j\tilde{\mathbf{f}}_{j_i}$. }

\State{Determine the MSE 
$E_{j_i}=\text{MSE}(\tilde{\mathbf{b}}_{j_i}-\mathbf{b}_j)$.}

\EndFor

\State{Choose $i^\ast=\argmin_i E_{j_i}$ and set $\check{\mathbf{f}}_j=\tilde{\mathbf{f}}_{i^\ast}$.}
\EndFor
\end{algorithmic}
\end{algorithm}

\subsection{Weighted $\ell_1$ regularization}
\label{sec:wl1-reg}

Assuming that the model \eqref{eq:forwardmodel_j} correctly describes how the data are acquired, any single measurement $\mathbf b_j$ should be sufficient to reconstruct the underlying image using \eqref{eq:l1_reg}. 
Indeed the main goal of compressive sensing is to reconstruct images from noisy and under-sampled data.  
As noted previously and particularly apparent in the first image in the second row of Figure \ref{fig:sobel-edges}, the accuracy of CS algorithms deteriorate for severely under-sampled data with low SNR, \cite{shchukina2017pitfalls,kang2019compressive}, a problem which may be further exacerbated when the acquisition method is obstructed in some way.  Such difficulties persist even for optimally chosen regularization parameters.  This is because as designed, the regularization term in  \eqref{eq:l1_reg} has {\em global} impact.  The solution would potentially be more accurate if the regularization term varied spatially, specifically to be more heavily penalized in true sparse regions (in the sparse domain) and much less so in regions of support in the sparse domain.  As already discussed, the (re-)weighted $\ell_1$ regularization method is designed to accomplish this task. In its most basic form, the  method is written as 
\begin{equation}
\label{eq:wl1_reg}
    \tilde{\mathbf{f}}_j = \argmin_{\mathbf s}\left\{
    \norm{F_j\mathbf s-\mathbf b_j}_2^2
    + \norm{W\mathcal L\mathbf s}_1
    \right\},
\end{equation}
where (in the re-weighted form) the entries ${\boldsymbol w}_j$ of the diagonal weighting matrix $W= \diag({\boldsymbol w}_j)$ are often constructed iteratively, yielding significant computational expense.  A more serious concern is that the noise and incompleteness of the acquired data may yield poor choices for the weights being fed into \eqref{eq:wl1_reg}.  The resulting reconstruction error will propagate with the penalty possibly being enhanced at edge locations while being simultaneously reduced in the smooth regions (which are sparse in the edge domain).  As a consequence, using \eqref{eq:wl1_reg} may yield worse results than the original (global) approximation in \eqref{eq:l1_reg}.   These issues are discussed in more detail in \cite{adcock2019jointsparsity,churchill2018edge, gelb2019reducing, scarnati2019accelerated}.

\subsubsection{Variance based joint sparsity (VBJS)}
\label{sec:VBJS}
The VBJS method was developed in \cite{adcock2019jointsparsity,gelb2019reducing,scarnati2019accelerated} to mitigate the issues inherent to the iterative weighted $\ell_1$ method.  Although originally designed for multiple measurement vectors (MMV),\footnote{To be clear, the typical MMV reconstruction process uses multiple observations at a {\em single} snapshot in time, that is, when there is no change in the underlying scene.  This is in contrast to the problem of interest in this investigation, which considers a sequence of observations over time during which the underlying scene changes.} the method can also be employed to a single measurement vector that is numerically processed multiple times to form MMVs.

Unlike the typical re-weighted $\ell_1$ regularization schemes that update the weights at each iteration, the VBJS procedure constructs the  weighting matrix in \eqref{eq:wl1_reg} once and then uses it for all subsequent iterations.  
Assuming an image is sparse in its edge domain,  to construct the weighting matrix $W$ we seek first to determine the edges of the image from the given data.  The VBJS method is based on the intuitive observation that because edge approximation methods produce smaller errors in sparse regions away from edges but larger errors near them, \cite{gelb2000detection,PAGelb},  it follows that multiple measurements in the edge domain yield correspondingly small variance in sparse regions while producing larger variance in regions of singular support.   The  values in the weighting matrix $W$ are thus accordingly constructed to be inversely proportional to the variance in the sparse domain.  Finally, we reiterate that the edge information needed to determine the weighting matrix $W$ can be computed directly from the given Fourier data, as is discussed in Section \ref{sec:edge_detection}, rather than by  first obtaining a coarse recovery in the physical domain.   As is illustrated in  Figure \ref{fig:sobel-edges}, having to recover the images before obtaining the edges can lead to significant information loss.

To illustrate how the VBJS method works, we first define ${G}_j \in \mathbb{R}^{(2N+1) \times (2N+1)}$ as the corresponding sparse domain matrix to each image ${\bf f}_j \in \mathbb{R}^{(2N+1) \times (2N + 1)}$, $j = 1,\dots,J$.  Observe that the non-zero entries of ${G}_j$ correspond to the singular support in the sparse domain. There are no inherent restrictions on how to approximate ${G}_j$ from $\mathbf{b}_j$,  and we employ the concentration factor edge detection method to approximate each ${G}_j$ in the edge domain {\em directly} from the given measurements ${\bf b}_j$ in \eqref{eq:forwardmodel_j}.
Specifically, we use \eqref{eq:edge_mask} directly and define
\begin{equation}
    \label{eq:sparse-edgeM}
    {\tilde G}_j^m := H^j_{\theta_m} \approx G_j, \quad j = 1,\dots,J, \quad m = 1,\dots,M.
\end{equation}

\begin{remark}
We note that  ${\tilde G}^m_j$ in \eqref{eq:sparse-edgeM} could be computed in multiple ways.  For example,  we could apply $M$ different regularizations for the indicator function in \eqref{eq:2Djump}.  For simplicity here we fix the regularization $h$ in \eqref{eq:gaussian} and use $M$ evenly spaced rotation angles $\{\theta_m\}_{m = 1}^M$, $\theta_m \in [0,\pi]$, as is done in \eqref{eq:regularized_edge2}.
\end{remark}

To construct $W$ in \eqref{eq:wl1_reg} for each $j = 1,\dots, J$, we require spatially varying weights ${\boldsymbol w}_j = \{ w_j(i)\}_{i=1}^{(2N+1)^2}$ that will effectively penalize the pixel-values that are zero in the sparse domain while allowing nonzero values to pass through.  Ideally, this means
\begin{equation}
\label{eq:weights}
	w_j(i) =
	\begin{cases}
		0,& i\in\mathcal{E}_j, \\
		const,& i\not\in\mathcal{E}_j,
	\end{cases}
\end{equation}
for some arbitrary $const > 0$.  Here $\mathcal{E}_j=\{i\in[1,(2N+1)^2]\mid {g}_j(i)\neq0\}$ is the set of indices for nonzero pixel values in the sparse domain and ${\bf g}_j = \text{vec}(G_j)$,  where $\text{vec}(\cdot)$ shapes the elements of input into a column vector. However since ${G}_j$ is not explicitly known, we first must determine an approximation to $\mathcal{E}_j$.
In this regard, for each image $\mathbf f_j$, we use the edge domain approximations in \eqref{eq:sparse-edgeM} to compute the matrices
\[\tilde{\mathbf{P}}_j=[\tilde{\bf g}^1_j,\tilde{\bf g}^2_j\dots, \tilde{\bf g}^M_j]\in\mathbb R^{(2N+1)^2\times M}, \quad j=1,\dots,J.\]
Here $\tilde{\mathbf g}_j^m=\text{vec}(\tilde G^m_j)\in\mathbb R^{(2N+1)^2}$. 
Based on the above discussion, the VBJS method uses the pointwise variance across the column vectors, given by (we have dropped the subscript $j$ for ease of notation)  
\begin{equation}
\label{eq:variance}
    v_i = \frac{1}{M}\sum^M_{m=1}\left(\tilde{\mathbf{P}}_{i,m}\right)^2-\left(\frac{1}{M}\sum^M_{m=1}\tilde{\mathbf{P}}_{i,m}\right)^2, \quad 
    i=1,\dots (2N+1)^2, 
\end{equation}
to determine the set $\mathcal{E}_j$ in \eqref{eq:weights}.  
Specifically, from \eqref{eq:variance} we  define a vector  ${\mathbf r}\in\mathbb R^{(2N+1)^2}$  element-wise  as
\begin{equation}\label{eq:matrix_T}
    r_i=\frac{\abs{{\tilde g}_i v_i}}{\max_i\abs{{\tilde g}_iv_i}},
    \quad i=1,\dots (2N+1)^2,
\end{equation}
and then replace \eqref{eq:weights} with 
\begin{equation}\label{eq:weight_l1}
    w_j(i) = \begin{cases}
    (1-r_i)/c, & r_i>\tau^w,\\
    1, & r_i\le\tau^w,
    \end{cases}
\quad\quad j = 1,\dots, J,
\end{equation}
where $c=\abs{\{i\mid r_i>\tau^w\}}$ ($\abs{\cdot}$  refers to the cardinality of a set) means the number of pixels that are considered to contain an edge location.  Observe that \eqref{eq:weight_l1} scales the weights according to the magnitude of the nonzero component in the sparse domain instead of using a constant value, as in \eqref{eq:weights}. This is to prevent false edges, which would typically yield small nonzero values, from having undue influence on the optimization.\footnote{Standard weighted $\ell_1$ regularization schemes typically scale the weights to be inversely proportional to the magnitudes of the components in the sparse domain of the solution at each iteration.  Numerical experiments in \cite{adcock2019jointsparsity} demonstrate that the VBJS approach is more robust, especially in low SNR environments.} Following \cite{scarnati2019accelerated}, and consistent with \eqref{eq:grid}, we choose $\tau^w = \frac{1}{2N+1}$ for each sequential image. 
Further details describing the construction of \eqref{eq:weight_l1} can be found in \cite{scarnati2019accelerated}, while  Algorithm \ref{algo:mmv_recovery} summarizes the VBJS recovery process for each image in the temporal sequence.

\begin{algorithm}[h!]
\caption{Variance Based Joint Sparsity (VBJS) method}
\label{algo:mmv_recovery}
	\hspace*{\algorithmicindent} \textbf{Input:} Measurements $\{\mathbf{b}_j\}_{j = 1}^J$ in \eqref{eq:forwardmodel_j} with forward operators $\{F_j\}_{j = 1}^J$, variances $\sigma^2$ and $\xi^2$, and sparse transform operator $\mathcal{L}$. Here we use TV.\\
	\hspace*{\algorithmicindent} \textbf{Output:} Image reconstructions $\{\mathbf{\tilde f}_j\}_{j = 1}^J$.
\begin{algorithmic}[1]

\For{$j=1$ to $J$}

\State{Compute $M$ edge maps $\tilde{G}^m_j$  for $m = 1,\dots, M$ rotation  angles $\{\theta_m\}_{m = 1}^M$ using \eqref{eq:edge_mask}.}

\State{Calculate the spatially adaptive weight  $\boldsymbol w_j$  in \eqref{eq:weight_l1}.}

\State{Determine  $\mathbf{\tilde f}_j$ from \eqref{eq:wl1_reg}.}

\EndFor
\end{algorithmic}
\end{algorithm}

%% file: 3_proposed.tex
\section{Incorporating inter-image information}
\label{sec:l2-reg}

As described above and demonstrated in \cite{adcock2019jointsparsity,gelb2019reducing,scarnati2019accelerated}, assuming there is enough information to generate the weighting matrix $W$, the weighted $\ell_1$ regularization method \eqref{eq:wl1_reg} yields more accurate image recovery than the standard CS approach in \eqref{eq:l1_reg}.   If intra-image information is inaccurate or insufficient, the reconstruction may not yield significant improvement, however.

The problems considered in this investigation not only result in under-determined systems of noisy data,  but also contain obstacles in each sequential image that further complicate the reconstruction process. The question then becomes how information from other images in the sequence might be used to improve {\em each} individual image reconstruction.  To this end, we note that supplementary data sets have been used to develop techniques such as joint sparse coding, \cite{song2016coupled,song2020multimodal}, which uses coupled dictionaries to recover multi-spectral (infrared) data.  Multi-channel data was jointly recovered via the penalized weighted least-square with a chosen reference channel in \cite{Rigie2015joint,kazantsev2018joint}, while a  cross modal regularization was employed in \cite{eslahi2019joint}.  Finally, we note that in \cite{kazantsev2014multimodal,chen2017joint}, edge features extracted from  supplementary data sets were specifically used to improve gradient regularization in each individual recovery.

These aforementioned techniques typically do not extract edge information as a pre-processing step for the purpose of joint recovery, however, with edge information used only to improve  an individual recovery, if at all.  By contrast here we develop a method that incorporates the edge information from each individual image, realized {\em directly} from the acquired data,  which allows us to accurately determine the common {\em inter-}image information needed to compensate for the missing information in each image recovery.  It also allows us to correctly identify structures that are {\em not} common across images. 

With this in mind, our new joint recovery method can be written as  
\begin{equation}\label{eq:optModel}
\tilde{\mathbf{f}} = \argmin_{\mathbf s} \left\{\norm{ F\mathbf{s}-\mathbf b}_2^2 
+ \norm{ W\boldsymbol{\mathcal L}\mathbf{s} }_1
+ \beta \norm{ \Phi\mathbf{s} }_2^2 \right\}, 
\end{equation}
where we have stacked the $\mathbf b_j$'s, $j=1,\dots,J$, into a long vector $\mathbf b\in\mathbb{C}^{J(2N+1)^2}$. 
Accordingly, $F\in\mathbb{C}^{J(2N+1)\times J(2N+1)}$ and $W, \boldsymbol{\mathcal{L}}\in\mathbb R^{J(2N+1)\times J(2N+1)}$ are respectively given by
\begin{align*}
    F = \diag(F_1,\dots,F_J), \quad
    W = \diag(W_1,\dots,W_J), \quad  
    \boldsymbol{\mathcal L} =\diag(\mathcal{L},\dots,\mathcal{L}).
\end{align*}
The inter-image information matrix $\Phi$ is defined as
\begin{equation}
\label{eq:Phi}
\Phi = \begin{pmatrix}
\tilde{C}_1&-\tilde{C}_1&&&\\
&\tilde{C}_2&-\tilde{C}_2&&\\
&&\ddots&\ddots&\\
&&&\tilde{C}_{J-1}&-\tilde{C}_{J-1}
\end{pmatrix}\in\mathbb R^{(J-1)(2N+1)^2\times J(2N+1)^2}.
\end{equation}
Here  {each {\em change mask}} $\tilde{C}_j \in\mathbb R^{(2N+1)^2\times (2N+1)^2}$, $j=1,\dots,J-1$, is a diagonal matrix that highlights the changed and unchanged regions along the diagonal entries and is given by 
\begin{align} 
	& \tilde{C}_{j}{(k,k)} = 
	\begin{cases} 
		{0} & \text{ if {$(k,k)$} lies in a region of change}, \\
		{1} & \text{ otherwise.} 
	\end{cases}
\label{eq:changemask}
\end{align}
By defining \eqref{eq:changemask} in this way we are able to ``borrow'' information from sequential data sets in regions that presumably remain the same, that is, where  we expect the difference between two sequential images to mainly stem from noise or missing data.

\subsection{Construction of the change mask}
\label{sec:changemask}
Broadly speaking, change detection requires two ingredients: (1)  determining the {\em unit of analysis}, which defines a means for comparison,  and (2) choosing that method for comparison, \cite{Chen2012object,HUSSAIN2013change,TEWKESBURY2015a}. The idea is then to identify the state change of an object in the scene by comparing the chosen unit of analysis. 

There are two general approaches for selecting a unit of analysis. At a fundamental level,  change is determined by calculating the difference at a single pixel or in a neighborhood of a pixel, \cite{inglada2007new,thonfeld2016robust,YE2016a}.  By contrast, an object-based approach compares adjacent frames to determine change in local features of that object, \cite{mcdermid2008object,TEWKESBURY2015a}. While pixel-based techniques typically require less pre-processing (e.g.~no image segmentation is needed),  object-based change detection methods are less sensitive to misalignment error, \cite{CHEN2014assessment}. 

Because change detection in local features is critical to the coupling term in \eqref{eq:optModel},  our change mask construction approach, as defined in \eqref{eq:changemask}, falls under the category of object-based techniques.
To this end we must  identify  all closed contour regions which requires several processing steps. To begin, we create a sequence of $J$ binary edge maps,
\begin{align} 
	& {\tilde{U}}_{j}(k,l) = 
	\begin{cases} 
		1 & \text{ if ${\tilde G}_j(k,l) > \tau_j^u$}, \\
		0 & \text{ otherwise,} 
	\end{cases}
\label{eq:ucurve}
\end{align}
{where $\tilde{G}_j$ is defined in \eqref{eq:edge_mask}.}  The parameter $\tau^u_j$ should be proportional to the maximum value of ${\tilde G}_j$ for each $j = 1,\dots, J$, and in general depends on the SNR. To demonstrate the robustness of our new method, in all of our numerical experiments we  fix the threshold as
\begin{equation}
	\tau^u_j = \frac{1}{2} \max_{k,l}({\tilde G}_j(k,l)), \quad k,l = 1,\dots,2N+1, 
\label{eq:tau_u}
\end{equation}
and note that it may be possible to improve performance with some additional tuning, specifically if some prior information regarding the magnitude of the non-zero entries in the sparse domain is given.  
Although we assume that all objects within the scene are rigid bodies with closed boundaries,
due to noise or low resolution, the binary edge masks defined in \eqref{eq:ucurve} might contain only parts of the edge curves or scattered subsets of the actual edge points.  Hence more processing is needed to identify the closed contour regions and subsequently construct the change masks in \eqref{eq:changemask}.  These tasks are accomplished using the following four steps:
\begin{enumerate}[label= C.\arabic* ] 
	\item \label{item:C1} 
	{\bf Clustering}: 
	All edge points belonging to the same object are clustered, and gaps in potential objects are bridged.  
	The edge points belonging to each individual object are then clustered into sequences that contain only edges of a single object. 
	
	\item \label{item:C2}
	{\bf Ordering}:
	The centroid for each cluster is computed and the edge points are ordered  in a counter-clockwise manner.  
	\item \label{item:C3}
	{\bf Filling closed regions}: 
	The ordered edge points are connected in each cluster via piecewise-linear interpolation to form closed edge curves.  The enclosed region for each single object is then filled with points inside the curve having the value of $1$ and all other points being asigned the value $0$.  
	\item \label{item:C4} 
	{\bf Comparison of filled single objects}:  Completion of steps \ref{item:C1}-\ref{item:C3} allows us to construct the change mask $\tilde{C}_j$ of two sequential images as the union of all filled single objects that are otherwise not accounted for in both images.\footnote{Note that the change mask, as defined in \eqref{eq:changemask}, assigns nonzero values only when there is {\em no} change in the underlying image, which is consistent with the proposed optimization model \eqref{eq:optModel}.}
\end{enumerate}
 
Below we provide more detailed descriptions for each step in the above procedure.

\subsubsection*{\ref{item:C1} Clustering} 
Given a binary edge mask $\tilde{U}_j$, defined in \eqref{eq:ucurve},
we seek to first determine the edge points that enclose unconnected regions belong to $i = 1, \dots, I$ distinguishable objects in the image and then store the corresponding set of edge points as their own (binary) {\em single object edge} masks, given by $\tilde U_j^{(i)}, i=1, \dots, I$, which satisfy the following properties for $i, i' = 1, \dots, I$:
\begin{enumerate}[label= (\roman*) ] 
\item \label{item:i} 
$\tilde U_j=\tilde U_j^{(1)}+\cdots+\tilde U_j^{(I)}$.

\item \label{item:ii} 
{$\sum_{k,l=1}^{2N+1}\left({\tilde U}_j^{(i)}(k,l)+{\tilde U}_j^{(i')}(k,l)\right)=\sum_{k,l=1}^{2N+1}{\tilde U}_j^{(i)}(k,l)+\sum_{k,l=1}^{2N+1}{\tilde U}_j^{(i')}(k,l)$}, $i \ne i'$. Moreover, we assume that none of the objects can  embedded inside another.  This is critical for Step \ref{item:C2} and is explained in more detail in Remark \ref{remark:overlap}.

\item \label{item:iii} 
For a given integer $d > 0$, the number of nonzero entries in ${\tilde U}_j^{(i)}$ is greater than $d$. 
\end{enumerate}

\begin{remark}
The first property ensures that the union  of all of the {single object edge} masks recovers the binary edge mask given by \eqref{eq:ucurve}.  The second guarantees that there is no overlap and no embedding of enclosed regions in the single object edge masks.  Finally, the third property ensures that the distance between the edge points belonging to a particular object is guaranteed to be no greater than the distance between any two objects, which prevents two objects from being inadvertently combined and recognized as a single object. 
\end{remark}

To generate closed form curves, any small gaps must first be bridged between the scattered edge points.
In this regard, we assume that there exists an integer $d\in\N$ such that a gap in a single edge curve is no larger than $d$ while the distance between any two distinct edge curves is always greater than $d$. 
Morphological dilation can then be used to connect the scattered edge points that belong to the same edge curve, followed by morphological erosion to ``thin'' the dilated boundary, see e.g.~\cite{Hsiao2005contour,papari2008adaptive,papari2011edge,gao2019extracting}.
To ensure that each of the $I$ closed contours is smooth and contains no sharp corners, we use the MATLAB function \texttt{imclose}. 

To proceed with the segmentation, we start with any edge point in any of the $I$ closed contours.  We then consider all points in the $8$-directional neighborhood and gradually locate every point belonging to that closed contour in the binary edge mask ${\tilde U}_{j}$  in \eqref{eq:ucurve}.  Once we run out of candidate edge points for the $8$-directional neighborhood, the resulting cluster is identified as the contour of a single object and then stored in the single object edge mask ${\tilde U}_j^{(i)}$. 
The same boundary tracking technique is sequentially performed on edge points that have not yet been assigned to any previously  constructed ${\tilde U}_j^{(\ell)}$, $\ell = 1,\dots, i$, to form ${\tilde U}_j^{(i+1)}$ and so on until all $I$ single object edge masks are generated.
This tracking step is accomplished using the MATLAB function \texttt{bwboundaries}.

\subsubsection*{\ref{item:C2} Ordering} 
After all edge points have been clustered such that $\tilde{U}_j = {\tilde U}_j^{(1)} + \dots + \tilde{U}_j^{(I)}$,
we proceed to order the edge points in a counter-clockwise manner for every single object edge mask  ${\tilde U}_j^{(i)}$. 
To this end, let us denote the sequence of edge points contained in  ${\tilde U}_j^{(i)}$ by
\begin{equation}\label{eq:unordered-edge-points}
	X = \left( \mathbf{x}_m \right)_{m=1}^\mathcal{M}, \quad 
	\mathbf{x}_m = (k_m,l_m), \quad 
	m = 1,\dots,\mathcal{M}. 
\end{equation}
Our goal is to determine an ordered sequence of $\mathcal{\tilde{M}} \leq \mathcal{M}$ points from $X$, 
\begin{equation}\label{eq:ordered-eq}
	X_{\text{ord}} = \left( \tilde{\mathbf{x}}_m \right)_{m=1}^{\mathcal{\tilde{M}}}, \quad 
	\tilde{\mathbf{x}}_m = (\tilde{k}_m,\tilde{l}_m), \quad 
	m=1,\dots,\mathcal{\tilde{M}}, 
\end{equation}
whose elements are ordered in a counter-clockwise manner. 
To do so, we assume that the underlying object $O$ is star-shaped. 
That is, there exists a point $\mathbf{x}_0$ in the object $O$ such that for all $\mathbf{x} \in O$ the line segment from $\mathbf{x}_0$ to $\mathbf{x}$ is in $O$.  
Furthermore, we assume that this point $\mathbf{x}_0$ can be approximated by the centroid or geometric center of $O$,\footnote{These are reasonable assumptions in  many  applications, such as the cars and tanks considered in our synthetic aperture radar (SAR) data example.} which is given by
\begin{equation}\label{eq:centroid} 
	\mathbf{x}_c = (k_c,l_c), \quad 
	k_{c} = \frac{1}{\mathcal{M}} \sum_{m=1}^\mathcal{M} k_m, \quad 
	l_{c}  = \frac{1}{\mathcal{M}} \sum_{m=1}^\mathcal{M} l_m.
\end{equation}
We then compute the angle $\varphi_m \in [0,2\pi)$ and radius $r_m \geq 0$ of $\mathbf{x}_m$ with respect to $\mathbf{x}_c$  for each $m = 1,\dots,\mathcal{M}$ as
\begin{equation}\label{eq:radii-angles}
	r_m = \norm{ \mathbf{x}_m - \mathbf{x}_c }_2, \quad 
	\varphi_m = 
	\left\{
	\begin{array}{c l}
		\arccos\left( \frac{k_m-k_c}{r_m} \right) & \text{if } l_m - l_c \geq 0 \text{ and } r_m \neq 0, \\ 
		- \arccos \left( \frac{k_m-k_c}{r_m} \right) & \text{if } l_m - l_c < 0, \\ 
		0 & \text{if } r_m = 0.
	\end{array} \right.
\end{equation}
Note that multiple edge points might have the same angle (but different radii). 
In this case, we only use the edge point with the largest radius.\footnote{Using the smallest or average radius would cause the reconstructed edge curve to sometimes lie inside the actual edge curve, resulting in an inaccurate change mask, $\tilde{C}_j$ in \eqref{eq:changemask}.} 
Thus, in general, the ordered sequence \eqref{eq:ordered-eq} might contain fewer elements than the unordered sequence of edge points \eqref{eq:unordered-edge-points}, i.\,e.\ ${\mathcal{\tilde{M}} < \mathcal M}$. 
Algorithm \ref{algo:ordering} describes how the ordered sequence of points is generated for each single object mask.

\begin{algorithm}[h!]
\caption{Generating an ordered sequence of edge points}
\label{algo:ordering}
	\hspace*{\algorithmicindent} \textbf{Input:} $X = \left( x_d \right)_{d=1}^\mathcal{D}$ as in \eqref{eq:unordered-edge-points}.
	\Comment{unordered edge points.} \\ 
	\hspace*{\algorithmicindent} \textbf{Output:} $X_{\text{ord}} = \left( \tilde{\mathbf{x}}_d \right)_{d=1}^{\mathcal{\tilde{D}}}$  in \eqref{eq:ordered-eq}.
	\Comment{ordered edge points.}
\begin{algorithmic}[1]
    \State{Compute the centroid $\mathbf{x}_c$ using \eqref{eq:centroid}.} 
    \State{Compute the vectors $\mathbf{r} = (r_1,\dots,r_\mathcal{D})^T$ and $\boldsymbol{\varphi} = (\varphi_1,\dots,\varphi_\mathcal{D})^T$ in \eqref{eq:radii-angles}.} 
	\State{Construct the matrix $A = [\boldsymbol{\varphi}, \mathbf{r}, X^T]$.} 
\State{Sort the rows of $A$ in ascending order based on the first column $\boldsymbol{\varphi}$.} 
	\While{$A$ is nonempty}
	\State{Find the top block submatrix $B$ of $A$ that has constant first column.} 
		\Comment{points with the same angle.} 
		\State{Sort the rows of $B$ in descending order based on the second column $\mathbf{r}$.} 
		\State{Store the last two columns of the first row of $B$, $\tilde{\mathbf{x}} = (\tilde{k},\tilde{l})$, in $X_{\text{ord}}$.} 
		\Comment{point with the largest radius.}
	\EndWhile \\
	\Return{$X_{\text{ord}}$} 
\end{algorithmic}
\end{algorithm}  

\begin{remark}
\label{remark:overlap}
Property \ref{item:ii} for the single object masks $\tilde{U}_j^{(i)}$ is necessary for accurate construction of  $X_{\text{ord}}$.  Specifically, the objects in the single object masks must neither overlap other objects nor be contained within  another object.  If this assumption does not hold, then \eqref{eq:radii-angles} would produce a sequence of ordered points that encloses a  larger region  ``covering'' the particular region of interest.  Note that Figure \ref{fig:sobel-edges} shows a violation of this assumption, as the internal objects of interest are embedded in the skull of the phantom.  Here, however, because it is readily apparent where the skull is, we can simply extract it before proceeding to Step \ref{item:C1}, that is to generate clusters of the edge points belonging to the interior individual distinguishable objects.   We note that extracting the skull in this way is common practice in MRI reconstruction since it is much larger in magnitude from the internal features,  see \cite{archibald2002method, archibald2003improving}.   The synthetic aperture radar (SAR) image example in Section \ref{subsec:golf_course} demonstrates that no initial extraction is needed when Property \ref{item:ii} holds for the single object masks.   Other change detection methods, such as those previously mentioned  and others specifically designed for remote sensing images, e.g.~\cite{afaq2021analysis, Ash_2014, LI2016a, inglada2007new}, may include strategies for overlapping and embedded images.\footnote{To the best of our knowledge, however, these procedures all use pixelated reconstructed images to detect the change, where as in our approach we use edge maps generated from the Fourier data to avoid data loss.}  But as this is not the primary focus of our paper, we assume that Property \ref{item:ii} holds. 
\end{remark}

\subsubsection*{\ref{item:C3} Filling closed regions} 
We now use piecewise linear interpolation to construct closed edge curves for each sequence of ordered edge points $X_{\text{ord}}$.  
In a second step we ``fill'' the resulting closed curves by assigning the value $1$ to pixels within the closed edge curve and $0$ to those outside.  In this way we obtain the (binary)
{\em filled single object masks}, ${\tilde Q}_j^{(i)}$ for each $\tilde U_j^{(i)}$, 
defined as
\begin{equation} \label{eq:edge-region-mask}
\begin{aligned}
	{\tilde Q}_j^{(i)}(k,l) & = \left\{ 
	\begin{array}{c l} 
		0 & \text{if $(k,l)$ lies outside of the reconstructed edge curve}, \\
		1 & \text{otherwise}.
	\end{array} \right.
\end{aligned}
\end{equation}
Note that using piecewise linear interpolation results in polygon shaped reconstructed edge curves. 
In our implementation we have combined the two steps (reconstructing and filling) by using the MATLAB function \texttt{inpolygon}.

\subsubsection*{\ref{item:C4} Comparison of filled single objects in two images} 

Now consider two consecutive data sets, ${\mathbf b}_j$, ${\mathbf b}_{j+1}$, from \eqref{eq:forwardmodel_j} for which we have constructed the filled single object masks ${\tilde Q}_j^{(i)}$ and ${\tilde Q}_{j+1}^{(i')}$, $1 \le i \le I$ and $1 \le i' \le I'$, respectively,  where $I$ and $I'$ correspond to the number of distinguishable objects in each image.  

To determine regions of change  we must identify filled single object masks that do not have a corresponding counterpart in the other image.  This is accomplished by using the {\em least-squares measure}, given by 
\begin{equation}\label{eq:diff-measure}
	\text{diff}(A,B) = \frac{\sum_{k,l=1}^{2N+1}\left(A(k,l)-B(k,l)\right)^2}{{\sum_{k,l=1}^{2N+1}A(k,l)+\sum_{k,l=1}^{2N+1}B(k,l)}}, \quad 
	A,B \in \{0,1\}^{(2N+1) \times (2N+1)}. 
\end{equation}
It is easy to verify that the difference measure \eqref{eq:diff-measure} provides us with a value between $0$ and $1$. 
In particular, we get 
\begin{equation}\label{eq:similarity-measure}
	\text{diff}(A,B) = \left\{ 
	 \begin{array}{c l} 
	 	0 & \text{if $A = B$}, \\ 
		1 & \text{if $A$ and $B$ are disjoint}.
	 \end{array} \right.
\end{equation}
To allow for small differences {between} the matrices, we say that two matrices $A,B$ are equal if their difference measure lies below a certain small threshold $\tau_{\rm diff} > 0$.\footnote{Comparing two matrices by simply checking the relation for all elements ($A = B \iff A(k,l) = B(k,l)$ for all $k,l$) is not appropriate since any single misidentified edge point, for example, due to noise, may cause the objects in two {single object edge} masks of the same object, say at corresponding times $j$ and $j+1$, to be determined as different.} 
Here we choose $\tau_{\rm diff} = 10^{-3}$ to be consistent with our discretization in \eqref{eq:grid} {and the amount of noise in the acquired data}, although our experiments indicate that this could be chosen smaller. 
Based on the calculation of $\text{diff}(A,B)$, we characterize the regions of change by pairwise comparing the filled single object masks {${\tilde Q}^{(i)}_j$ and ${\tilde Q}^{(i')}_{j+1}$ for ${i = 1,\dots,I}$ and ${i' = 1,\dots,I'}$.}
This procedure is described by Algorithm \ref{algo:change-index}. 

\begin{algorithm}[h!]
\caption{Determining the Change Index Sets}
\label{algo:change-index}
	\hspace*{\algorithmicindent} \textbf{Input:} 
	$\{{\tilde Q}^{(i)}_{j}\}_{i=1}^{I}$, $\{{\tilde Q}^{(i')}_{j+1}\}_{i'=1}^{I'}$, and threshold $\tau_{\rm diff}$. (We choose $\tau_{diff} = 10^{-3}$.) \\
	\hspace*{\algorithmicindent} \textbf{Output:} ${\tilde{I}_{{j}} \subset \{1, \dots, I\}}$ and ${\tilde{I}'_{{j+1}} \subset \{1, \dots, I'\}}$
	\Comment{new index sets}
\begin{algorithmic}[1]
    \State{Initiate $\tilde{I}_{{j}} = \{1,\dots,I\}$, $\tilde{I}'_{{j+1}} = \{1,\dots,I'\}$}
    \For{${i = 1,\dots,I}$ and ${i' = 1,\dots,I'}$} 
	\Comment{compare all filled single object masks}
		\If{ $\text{diff}({\tilde Q}^{(i)}_j,{\tilde Q}^{(i')}_{j+1}) <\tau_{diff}$}  
		\Comment{there is no change for this pair}
			\State{${\tilde I}_{{j}} = {\tilde I}_{{j}} \setminus \{ i \}$} 
			\Comment{remove the indices from the index sets}
			\State{$\tilde{I}'_{{j+1}} = \tilde{I}'_{{j+1}} \setminus \{ i' \}$} 
		\EndIf 
	\EndFor\\
	\Return{$\tilde{I}_{{j}}$, $\tilde{I}'_{{j+1}}$} 
\end{algorithmic}
\end{algorithm}  

Algorithm \ref{algo:change-index} returns two new index sets, ${\tilde{I}_{j} \subset \{1, \dots, I\}}$ and ${\tilde{I}'_{j+1} \subset \{1, \dots, I'\}}$, which contain only the indices of the filled single object masks corresponding to change, such as those due to translation, rotation, insertion or deletion of objects in the underlying scene.
We stress once again that the coupling term in \eqref{eq:optModel} assumes that any differences over the {\em unchanged} regions between each pair of adjacent images stem either from noise or other sources of error.  The coupling term therefore aims to minimize such difference over the unchanged regions via $\ell_2$-norm to enhance the recovery of the shared features.

With Steps \ref{item:C1} - \ref{item:C4} completed, we are now able to construct the change mask $\tilde{C}$ in \eqref{eq:changemask} that constitutes $\Phi$ in \eqref{eq:optModel} as 
$$	\tilde{C}_{j} = \diag{(\boldsymbol1-\text{vec}(R_{j}))}, $$
where 
\begin{equation} 
\label{eq:CM-edge-region-mask}
R_{j} = \ceil*{ \frac{1}{2} \left( \ceil*{\frac{1}{\mid\tilde{I}_{j}\mid}\sum_{i \in \tilde{I}_{j}} \tilde{Q}^{(i)}_{j}} + \ceil*{\frac{1}{\mid\tilde{I}'_{j+1}\mid}\sum_{i' \in \tilde{I}'_{j+1}} \tilde{Q}^{(i')}_{j+1}} \right) }.
\end{equation}
Here $\abs{\cdot}$ denotes cardinality and $\boldsymbol1\in\mathbb R^{(2N+1)^2}$ is a vector with value $1$ for every entry. We also use $\ceil*{\cdot}$ to denote the ceiling operator and note that the summation and ceiling operators are performed element by element.  From \eqref{eq:CM-edge-region-mask} we observe that $\tilde{C}_j$ is constructed by comparing the filled single object masks ${\tilde Q}^{(i)}_{j}$ and ${\tilde Q}^{(i')}_{j+1}$, with $i \in \tilde{I}_{j}$ and $i' \in \tilde{I}'_{j+1}$ respectively.

\subsection{Coupling term regularization parameter}
\label{subsec:CM-regularizer}

The coupled term in \eqref{eq:optModel} boosts the accuracy in each image recovery  by promoting the emphasis on the unchanged regions (inter-image information). As discussed previously,  each of the time sequenced data collections may be missing vital information for independent signal recovery, due to obstacles or occlusions in the scene,  and/or missing bands of Fourier data. 
To improve each reconstruction, information is ``borrowed'' from other images in the sequence.  It is important that this information comes only from common features that remain unchanged in the adjacent data sets.  This information is captured using Steps \ref{item:C1} - \ref{item:C4} and then inserted into  \eqref{eq:Phi}.  What remains is to choose the regularization parameter $\beta$ in \eqref{eq:optModel}, which represents how much information is borrowed from the neighboring time frames. 
Several factors must be taken into account, as we now describe. 

\begin{figure}[h!]
    \centering
    \begin{subfigure}[b]{.21\textwidth}
    \includegraphics[width=\textwidth]{plot/true_f3_GE}
    \caption{${\bf f}_3$ (ground truth)}
    \end{subfigure}
    ~
    \begin{subfigure}[b]{.21\textwidth}
    \includegraphics[width=\textwidth]{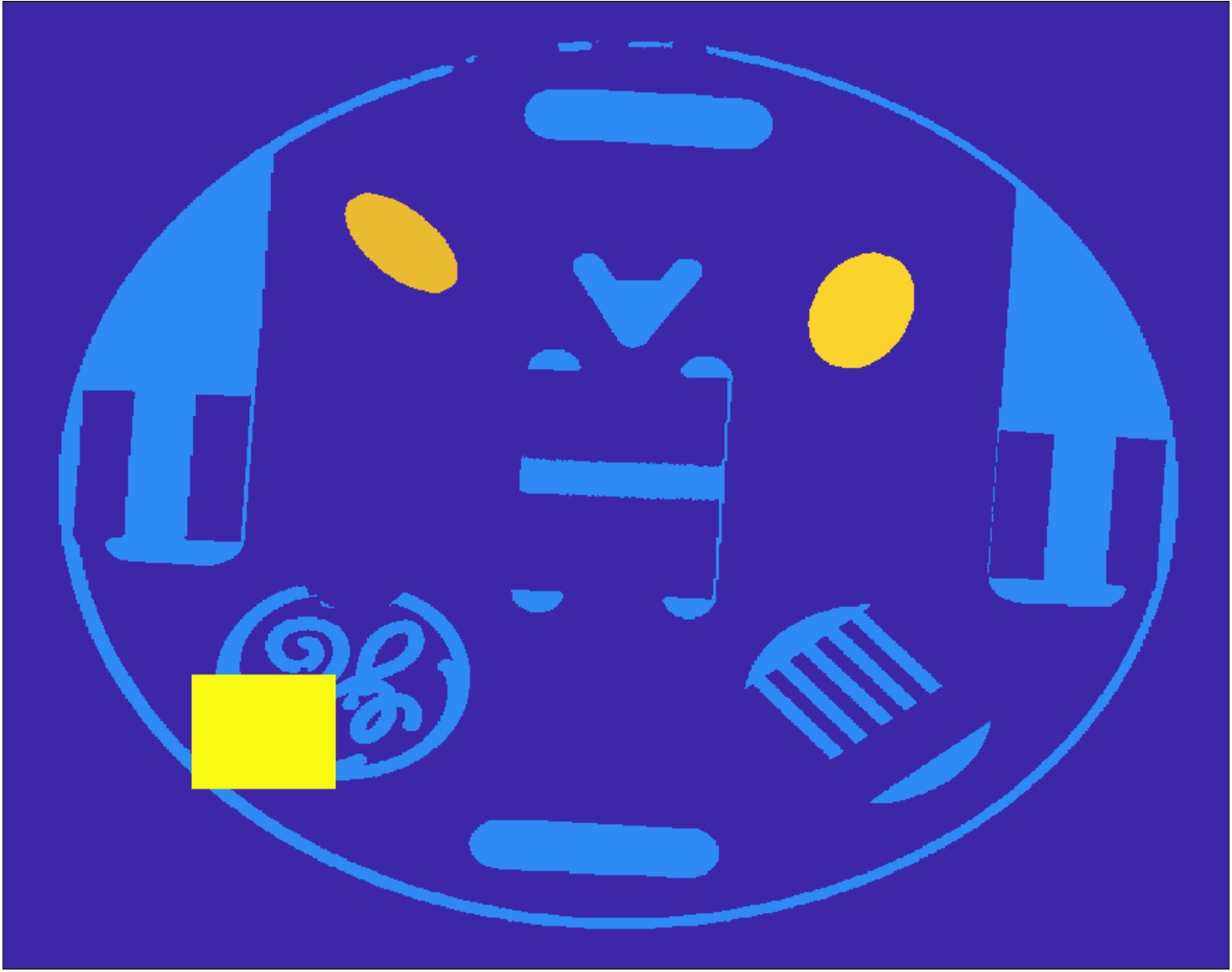}
    \caption{${\bf f}_3$ (with occlusion)}
    \end{subfigure}
    ~
    \begin{subfigure}[b]{.21\textwidth}
    \includegraphics[width=\textwidth]{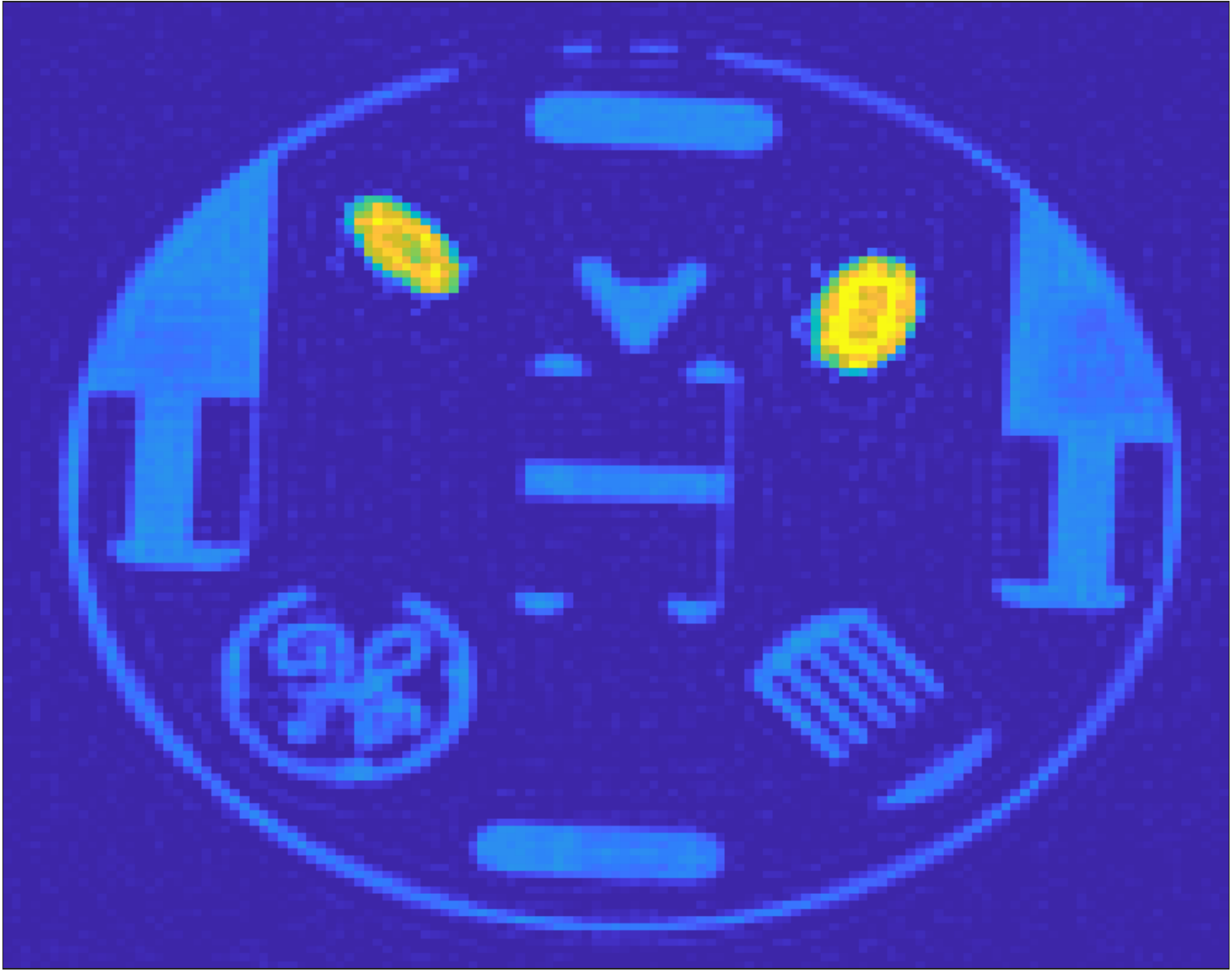}
    \caption{SNR $= 10$; $\beta=.06$.}
    \end{subfigure}
    ~
    \begin{subfigure}[b]{.21\textwidth}
    \includegraphics[width=\textwidth]{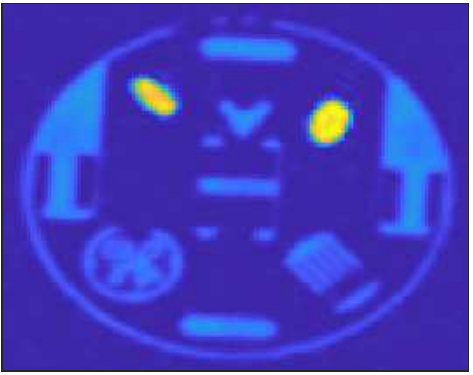}
    \caption{SNR $= 2$; $\beta=.5$.}
    \end{subfigure}
    
    \caption{Recovery ${\tilde{\mathbf f}}_3$ using \eqref{eq:optModel} for SNR $= 10,2$ and $\beta = .06,.5$ respectively.}
    \label{fig:beta_snr}
\end{figure}

First, as is the case for parameter selection of $\mu$ in the standard $\ell_1$ regularization approach given in \eqref{eq:l1_reg}, the choice of $\beta$ depends on the level of noise. Specifically, a  higher SNR value would suggest less regularization is needed and vice versa.  Figure \ref{fig:beta_snr} illustrates that this intuition seems to hold.  
Supplemental information becomes very important when data from any individual image is either missing, e.g.~in the case of missing bands of Fourier samples, physical obstacles, or blur. Similar to the low SNR case, in these situations it is also reasonable to choose $\beta$ to increase the influence of the borrowed information.  That is, we essentially want to more heavily penalize inter-image similarity.

However we must also consider the amount of change occurring between adjacent temporal frames. In particular, the coupled term in \eqref{eq:optModel} increases  accuracy by borrowing information from unchanged regions, that is, those that presumably share the same features in  consecutive frames.   In applications where significant change occurs between images, or equivalently when there are relatively few unchanged regions in adjacent frames, a small $\beta$ should be chosen so less emphasis is put on shared information.  In general using the coupled term in \eqref{eq:optModel} can be seen as most beneficial when objects move slowly compared to temporal resolution.
Finally, although we can heuristically determine the coupling parameter $\beta$, more work is needed to choose it robustly, which will be the subject of future work.  One possible option is to iteratively determine $\beta$ using some threshold, or, similarly to the approach described in Section \ref{sec:VBJS}, to choose $\beta$ to be spatially varying.  A probablistic approach such as sparse Bayesian learning may also be useful, see for example \cite{ji2008bayesian,tipping2001sparse,wipf2004sparse}. Here, in order to demonstrate the robsutness of our approach we simply choose $\beta = .5$.

We now have all of the necessary ingredients for computing \eqref{eq:optModel}, which are described in Algorithm \ref{algo:joint_recovery}.  

\begin{algorithm}[h!]
\caption{Recovering a temporal sequence of images using \eqref{eq:optModel}.}
\label{algo:joint_recovery}
	\hspace*{\algorithmicindent} \textbf{Input:} Measurements $\{\mathbf{b}_j\}_{j = 1}^J$, forward operators $\{F_j\}_{j = 1}^J$, variances $\sigma^2$ and $\xi^2$ for the model \eqref{eq:forwardmodel_j} and sparse transform $\mathcal{L}$, respectively. (We employ the $TV$ operator, but others may be used.) 
We also in advance choose parameters $\tau^w = \frac{1}{2N+1}$ {in \eqref{eq:weight_l1}} and $\beta = .5$ {in \eqref{eq:optModel}}. 
\\
	\hspace*{\algorithmicindent} \textbf{Output:} Reconstruction $\{\tilde{\mathbf{f}}_j\}_{j = 1}^J$.
\begin{algorithmic}[1]

\For{ $j=1,\dots,J$}
\State{Determine each weighting matrix $W_j$ from \eqref{eq:weight_l1} using the edge map ${\tilde G}_j$ in \eqref{eq:sparse-edgeM}.}
\State{{Define threshold parameter $\tau^u_j$  according to \eqref{eq:tau_u}.}}
\State{Form the binary edge masks  $\tilde{U}_j$ in \eqref{eq:ucurve} from the edge maps ${\tilde G}_j$ in \eqref{eq:sparse-edgeM}.}
\State{Form the single object edge masks $\{{\tilde U}^{(i)}_j\}_{i = 1}^{I}$ using Step \ref{item:C1}.}
\State{Form the filled single object masks $\{{\tilde Q}^{(i)}_j\}_{i = 1}^{I}$ using \eqref{eq:edge-region-mask}.}

\For{$j<J$}
\State{Construct change masks $\tilde{C}_j$ in \eqref{eq:changemask}.} 
\State{Construct inter-image information matrix $\Phi$ in \eqref{eq:Phi}.}
\EndFor

\EndFor

\State{Solve for the joint recovery $\{\tilde{\mathbf{f}}_j\}_{j = 1}^J$ by \eqref{eq:optModel} using the ADMM method, {\cite{boyd2011distributed}}}

\end{algorithmic}
\end{algorithm}

%% file: 4_numerical.tex
\section{Numerical Experiments}
\label{sec:test}

We consider two numerical experiments to validate the robustness and accuracy of our new joint recovery method as detailed in Algorithm \ref{algo:joint_recovery}. {Following \cite{gelb2019reducing}}, to solve the convex optimization problem in \eqref{eq:optModel} we employ the Alternating Direction Method of Multipliers (ADMM), see \cite{boyd2011distributed}.   In all experiments we assume we are given data as in \eqref{eq:forwardmodel_j}, that is, noisy Fourier data where bands are missing from each set in the temporal sequence.  To approximate the continuous Fourier samples in \eqref{eq:fourcoeff} we use 
highly resolved image data and then apply the discrete Fourier transform. In this way we can include noise and obstructions in the image domain into our calculated Fourier samples while also ensuring that the data fidelity errors in \eqref{eq:l1_reg}, such as aliasing, are properly accounted for.

Our new method will be tested for different signal to noise ratio (SNR) values calculated as
\begin{equation}
\label{eq:SNR}
	\text{SNR}_{{dB}_j} = 10\log_{10}\left(\frac{\bar{\mathbf{b}}_j}{\sigma_j}\right)^2,
\end{equation}
where $\bar{\mathbf{b}}_j$ is the mean over the data $\mathbf{b}_j$ given in \eqref{eq:forwardmodel}.

In our first example we generate a synthetic image so that we may analyze the efficacy of our method with respect to SNR and data resolution.  
 In the second experiment we use a pre-formed synthetic aperture radar (SAR) image given in \cite{SAR_Image_ref}.   As in the first case, we consider various levels of SNR and introduce deletions, insertions, rotations and translations of objects.  We also include occlusions (obstructions) in the physical domain. We note that in this second experiment, since we are starting from a given SAR image, we assume the Fourier data are acquired using a discrete Fourier transform, that is, we do not simulate the integral transform to obtain the data. 
We compare our new algorithm to individual recovery via the standard $\ell_1$ regularization \eqref{eq:l1_reg} using Algorithm \ref{algo:EM_l1_reg} and the VBJS method \eqref{eq:wl1_reg} using Algorithm \ref{algo:mmv_recovery}. 

\subsection{Experiment 1: Synthetic MRI Phantom}
\label{subsec:phantom}

\begin{figure}[h!]
    \centering
    \begin{subfigure}[b]{.19\textwidth}
        \includegraphics[width=\textwidth]{plot/true_f1_GE}
        \caption{$\mathbf{f}_1$}
    \end{subfigure}
    ~
    \begin{subfigure}[b]{.19\textwidth}
        \includegraphics[width=\textwidth]{plot/true_f2_GE}
        \caption{$\mathbf{f}_2$}
    \end{subfigure}
    ~ 
    \begin{subfigure}[b]{.19\textwidth}
        \includegraphics[width=\textwidth]{plot/true_f3_GE}
        \caption{$\mathbf{f}_3$}
    \end{subfigure}
    ~
    \begin{subfigure}[b]{.19\textwidth}
        \includegraphics[width=\textwidth]{plot/true_f4_GE}
        \caption{$\mathbf{f}_4$}
    \end{subfigure}
    \\ 
    \begin{subfigure}[b]{.19\textwidth}
        \includegraphics[width=\textwidth]{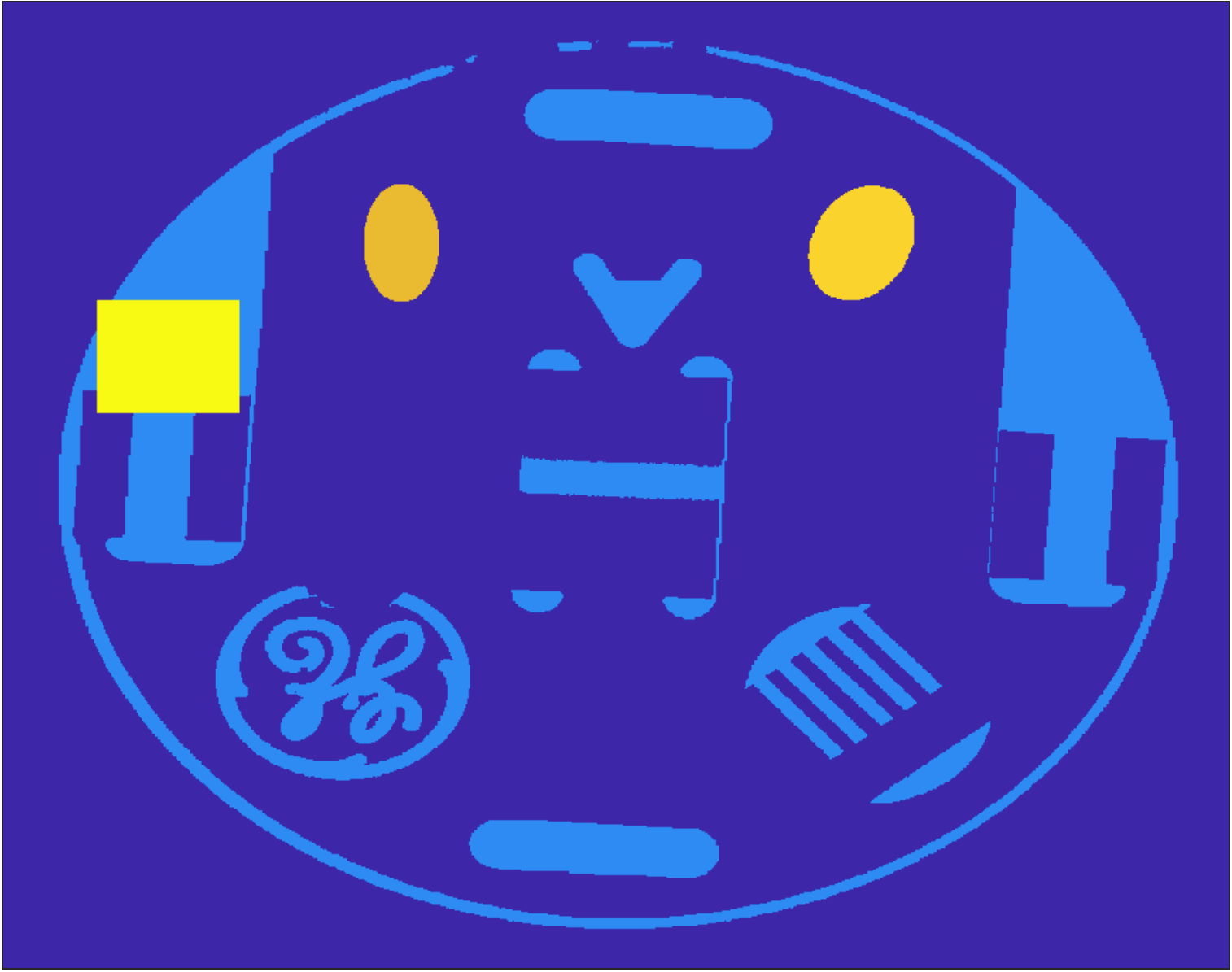}
        \caption{Obstacle$_1$}
    \end{subfigure}
    ~
    \begin{subfigure}[b]{.19\textwidth}
        \includegraphics[width=\textwidth]{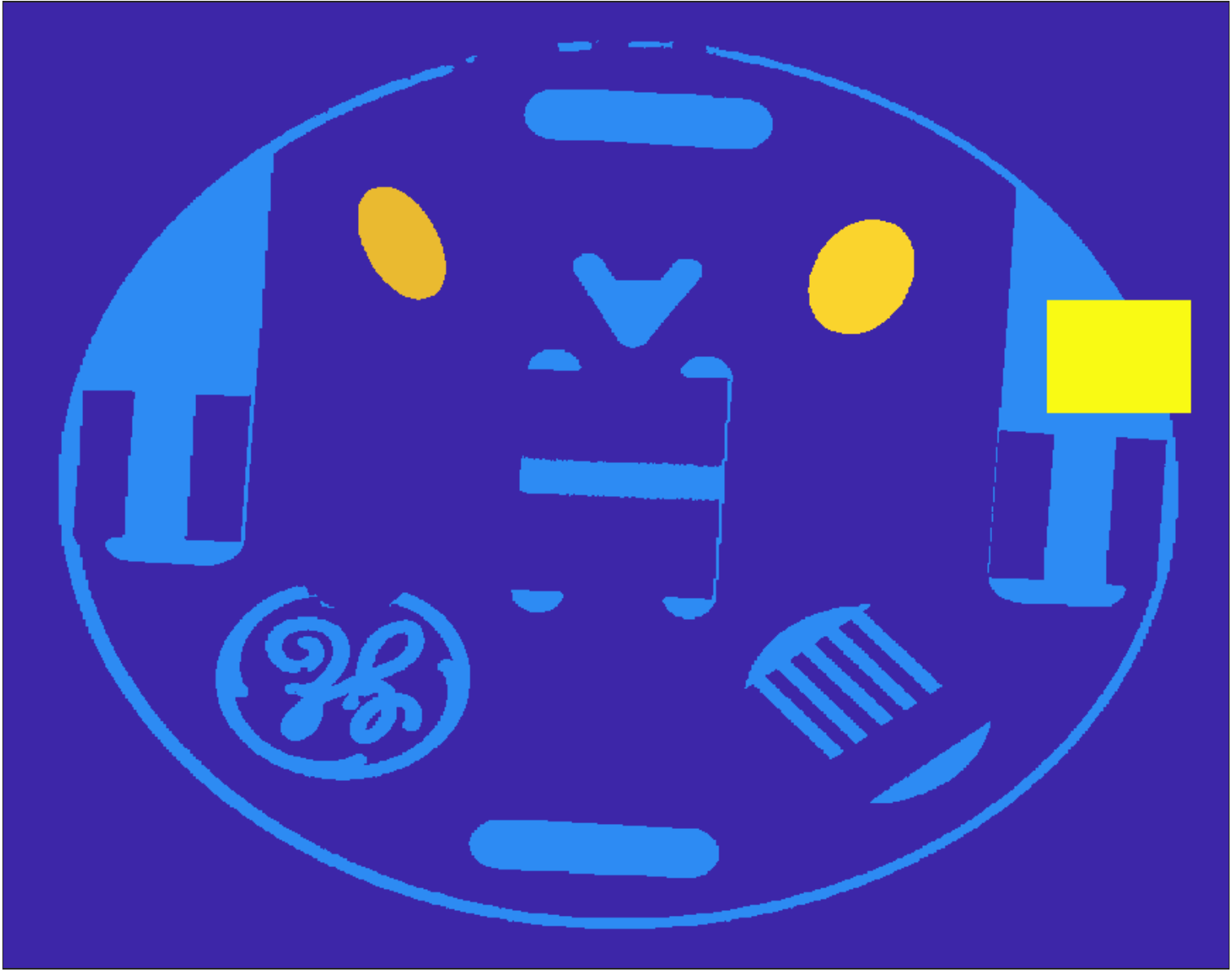}
        \caption{Obstacle$_2$}
    \end{subfigure}
    ~
    \begin{subfigure}[b]{.19\textwidth}
        \includegraphics[width=\textwidth]{plot/obstacle_f3_GE}
        \caption{Obstacle$_3$}
    \end{subfigure}
    ~
    \begin{subfigure}[b]{.19\textwidth}
        \includegraphics[width=\textwidth]{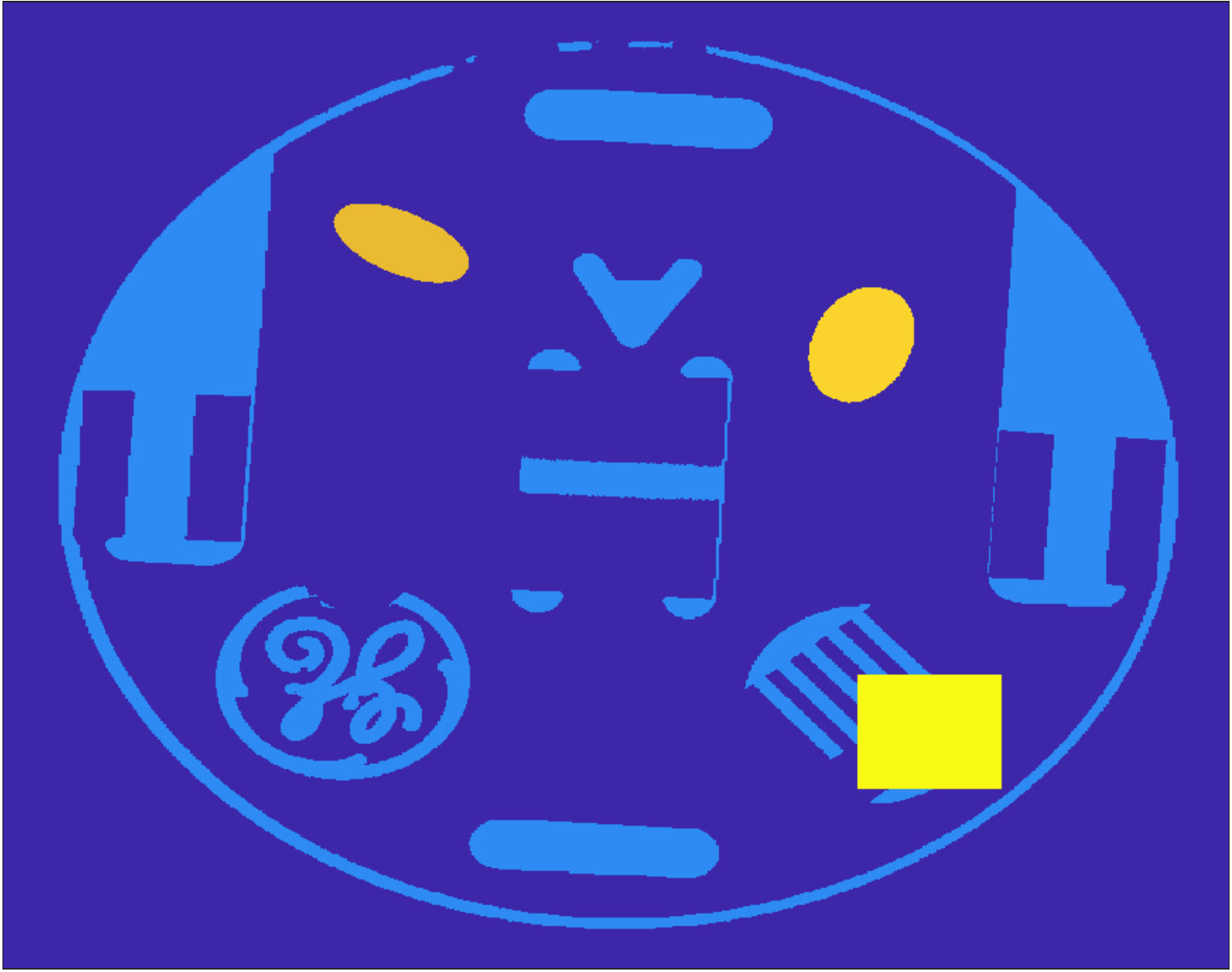}
        \caption{Obstacle$_4$}
    \end{subfigure} 
    \\
    \begin{subfigure}[b]{.19\textwidth}
        \includegraphics[width=\textwidth]{plot/rec_f1_std_GE}
        \caption{$\tilde{f}_{1}$}
    \end{subfigure}
    ~
    \begin{subfigure}[b]{.19\textwidth}
        \includegraphics[width=\textwidth]{plot/rec_f2_std_GE}
        \caption{$\tilde{f}_{2}$}
    \end{subfigure}
    ~
    \begin{subfigure}[b]{.19\textwidth}
        \includegraphics[width=\textwidth]{plot/rec_f3_std_GE}
        \caption{$\tilde{f}_{3}$}
    \end{subfigure}
    ~
    \begin{subfigure}[b]{.19\textwidth}
        \includegraphics[width=\textwidth]{plot/rec_f4_std_GE}
        \caption{$\tilde{f}_{4}$}
    \end{subfigure} 
    \\
    \begin{subfigure}[b]{.19\textwidth}
        \includegraphics[width=\textwidth]{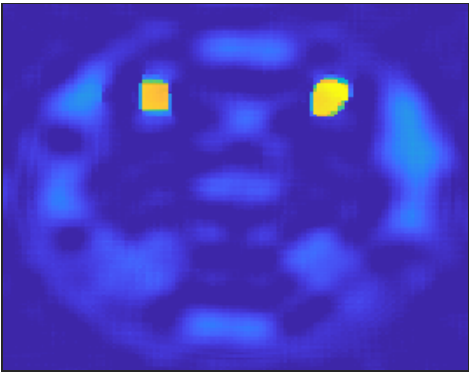}
        \caption{$\tilde{f}_1^\text{VBJS}$}
    \end{subfigure}
    ~
    \begin{subfigure}[b]{.19\textwidth}
        \includegraphics[width=\textwidth]{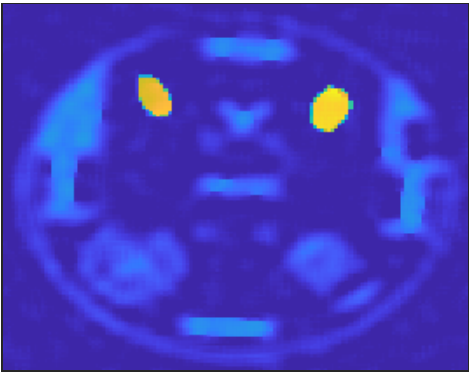}
        \caption{$\tilde{f}_2^\text{VBJS}$}
    \end{subfigure}
    ~
    \begin{subfigure}[b]{.19\textwidth}
        \includegraphics[width=\textwidth]{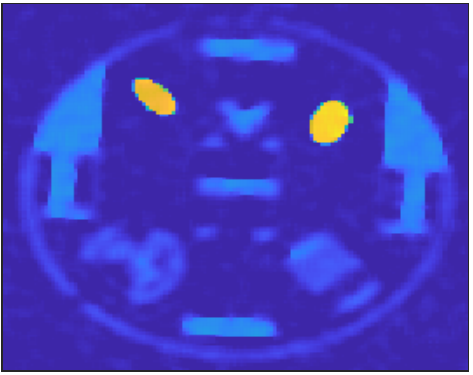}
        \caption{$\tilde{f}_3^\text{VBJS}$}
    \end{subfigure}
    ~
    \begin{subfigure}[b]{.19\textwidth}
        \includegraphics[width=\textwidth]{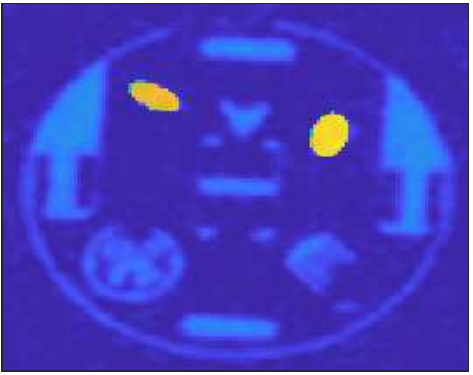}
        \caption{$\tilde{f}_4^\text{VBJS}$}
    \end{subfigure}   
    \\
    \begin{subfigure}[b]{.19\textwidth}
        \includegraphics[width=\textwidth]{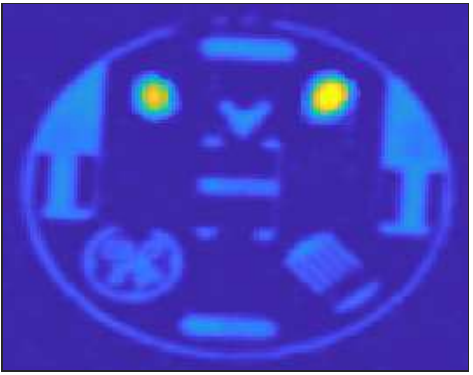}
        \caption{$\tilde{f}_1^\text{joint}$}
    \end{subfigure}
    ~
    \begin{subfigure}[b]{.19\textwidth}
        \includegraphics[width=\textwidth]{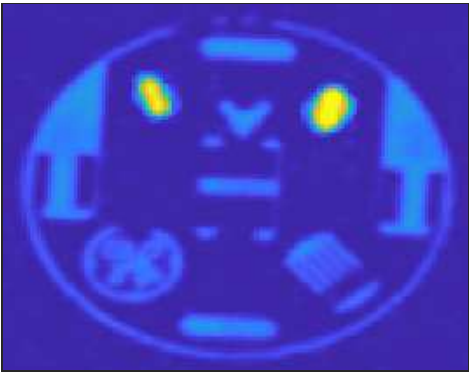}
        \caption{$\tilde{f}_2^\text{joint}$}
    \end{subfigure}
    ~
    \begin{subfigure}[b]{.19\textwidth}
        \includegraphics[width=\textwidth]{plot/rec_f3_joint_GE}
        \caption{$\tilde{f}_3^\text{joint}$}
    \end{subfigure}
    ~
    \begin{subfigure}[b]{.19\textwidth}
        \includegraphics[width=\textwidth]{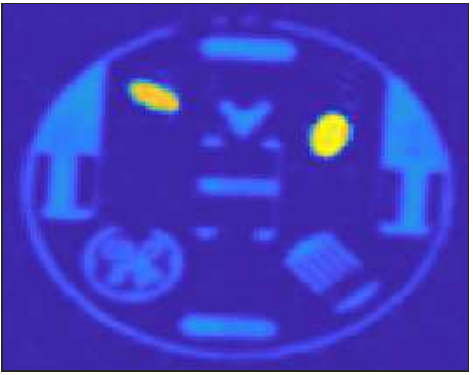}
        \caption{$\tilde{f}_4^\text{joint}$}
    \end{subfigure}
    \caption{(first row) Temporal sequence of MRI phantom images of hydroiodic acid, MnCl2, and rGNP-HI using a GE HTXT 1.5T clinical MRI scanner, \cite{Lalwanietal}.  (second {row})
    Physical images in scaled color with physical obstructions in yellow. (third {row})
    Individual reconstruction by standard $\ell_1$-regularization. (fourth {row}) 
    Individual reconstructions by VBJS $\ell_1$-regularization. (bottom {row}) 
    Joint recovery (last {row}). }
\label{fig:rec_GE}
\end{figure}

Figure \ref{fig:rec_GE} (first {row}) displays a sequence of ground truth images $\mathbf{f}_j$, $j=1,\dots,4$.\footnote{Although we incorporate information from $J = 6$ data sets, for better visualization we display only four.}   The image background is zero-valued and each of the static structures have magnitude $0.3$.  Each image also has two ellipses of magnitudes $0.8$ and $0.9$, respectively, that rotate and/or translate in time.

The third {row} in Figure \ref{fig:rec_GE} displays the individual recoveries using the standard $\ell_1$ reconstruction, \eqref{eq:l1_reg} as obtained using Algorithm \ref{algo:EM_l1_reg}, while the fourth {row} shows the individual recoveries using the VBJS method, \eqref{eq:wl1_reg} as obtained by Algorithm \ref{algo:mmv_recovery}. Finally, the results using our new joint recovery method \eqref{eq:optModel} as realized by Algorithm \ref{algo:joint_recovery} are displayed in the fifth {row}.  The parameters are all chosen as described in each of the given algorithms.  Observe that only the method provided by \eqref{eq:optModel} is able to recover information from the obscured region while still enhancing the recovery in the rest of the domain.  
In particular we see that when the edges of the static objects are not well detected, the performance of \eqref{eq:wl1_reg} is poor.  Indeed, falsely identified edges may cause the VBJS recovery method to yield worse results than those obtained using standard $\ell_1$ reconstruction, since the weighting matrix $W$  does not accurately penalize smooth regions.

\begin{figure}[h!]
\centering 
\begin{subfigure}[b]{0.22\textwidth}
    \includegraphics[width=\textwidth]{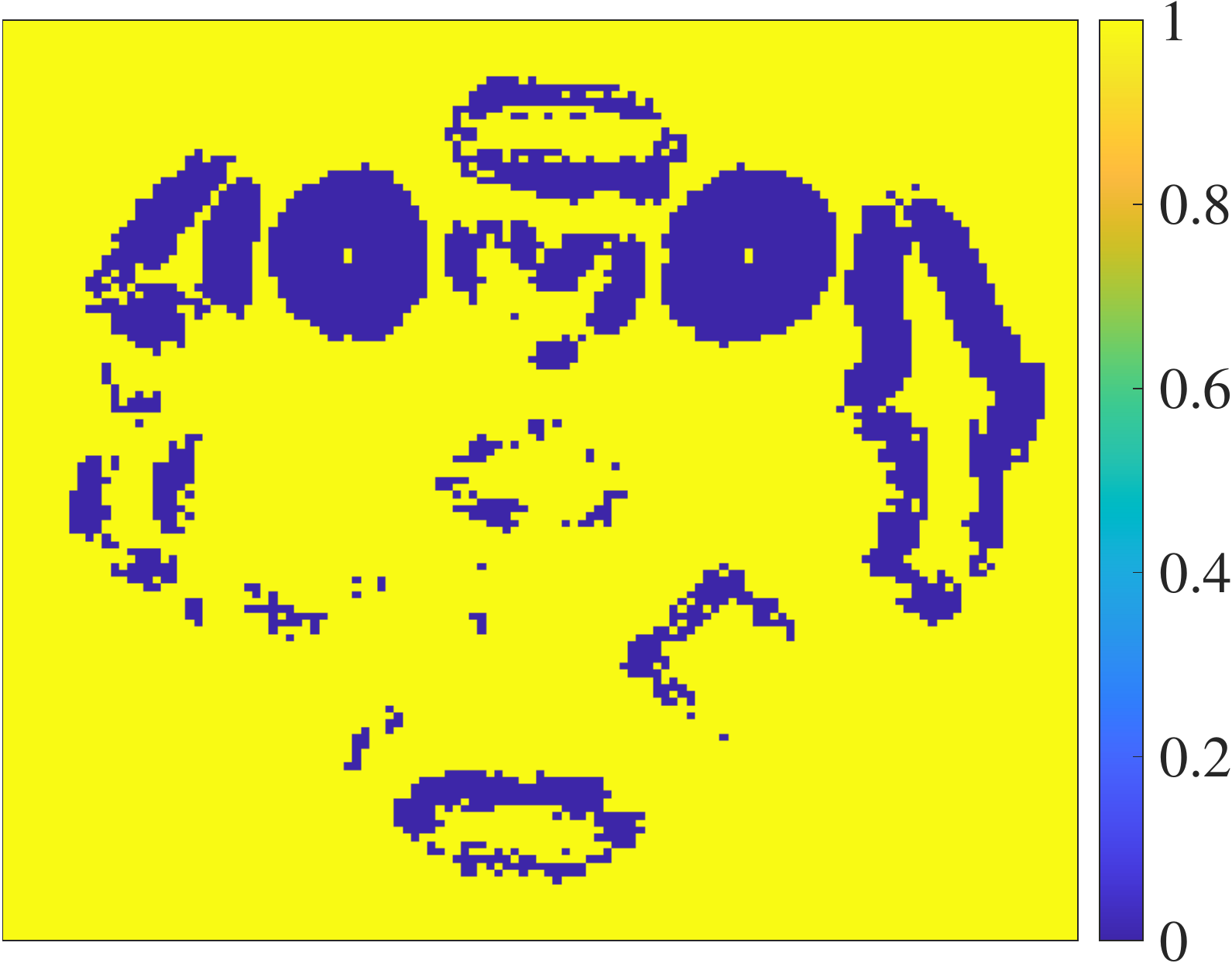}
    \caption{$W_1$}
\end{subfigure} 
~ 
\begin{subfigure}[b]{0.22\textwidth}
    \includegraphics[width=\textwidth]{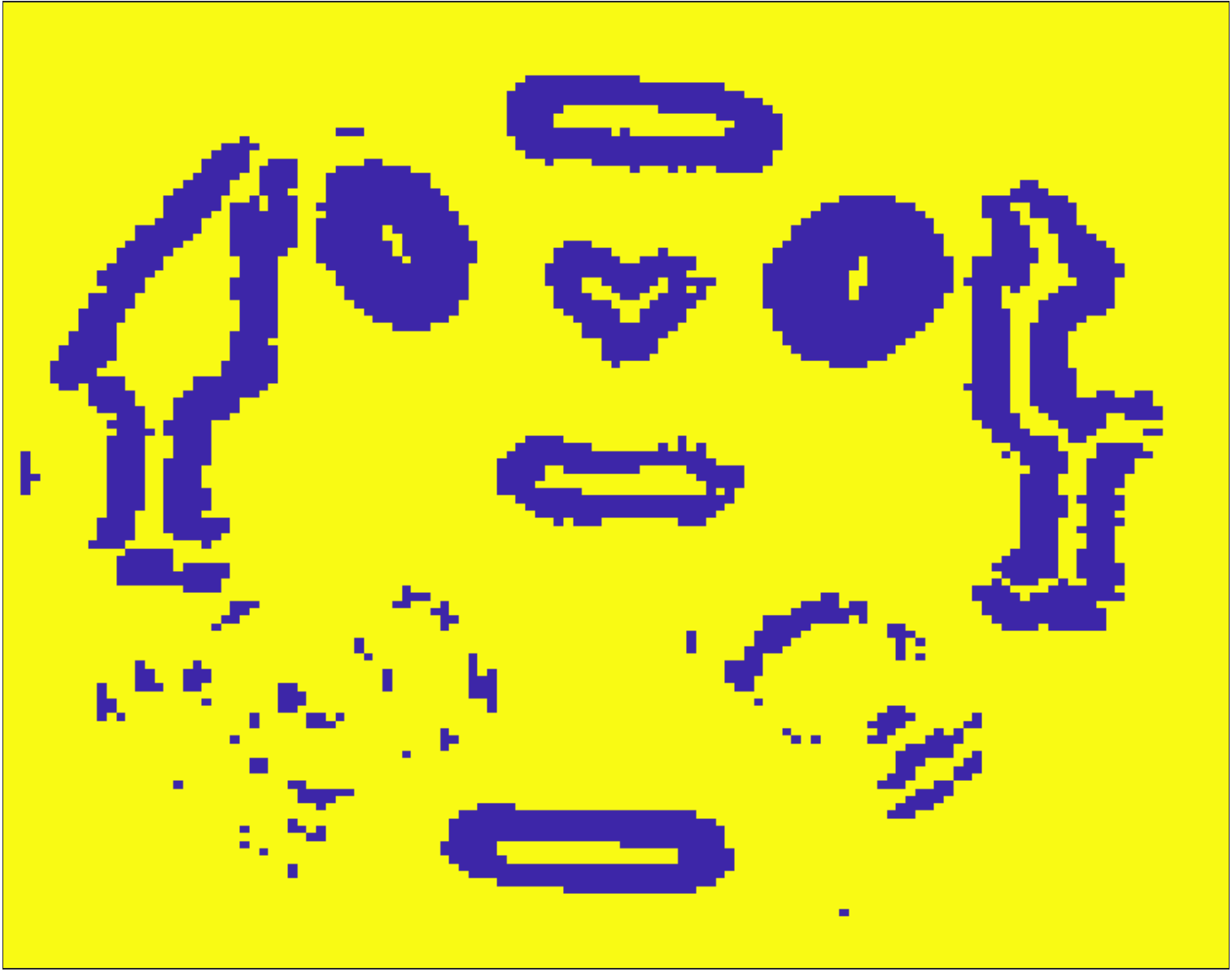}
    \caption{$W_2$}
\end{subfigure} 
~ 
\begin{subfigure}[b]{0.22\textwidth}
    \includegraphics[width=\textwidth]{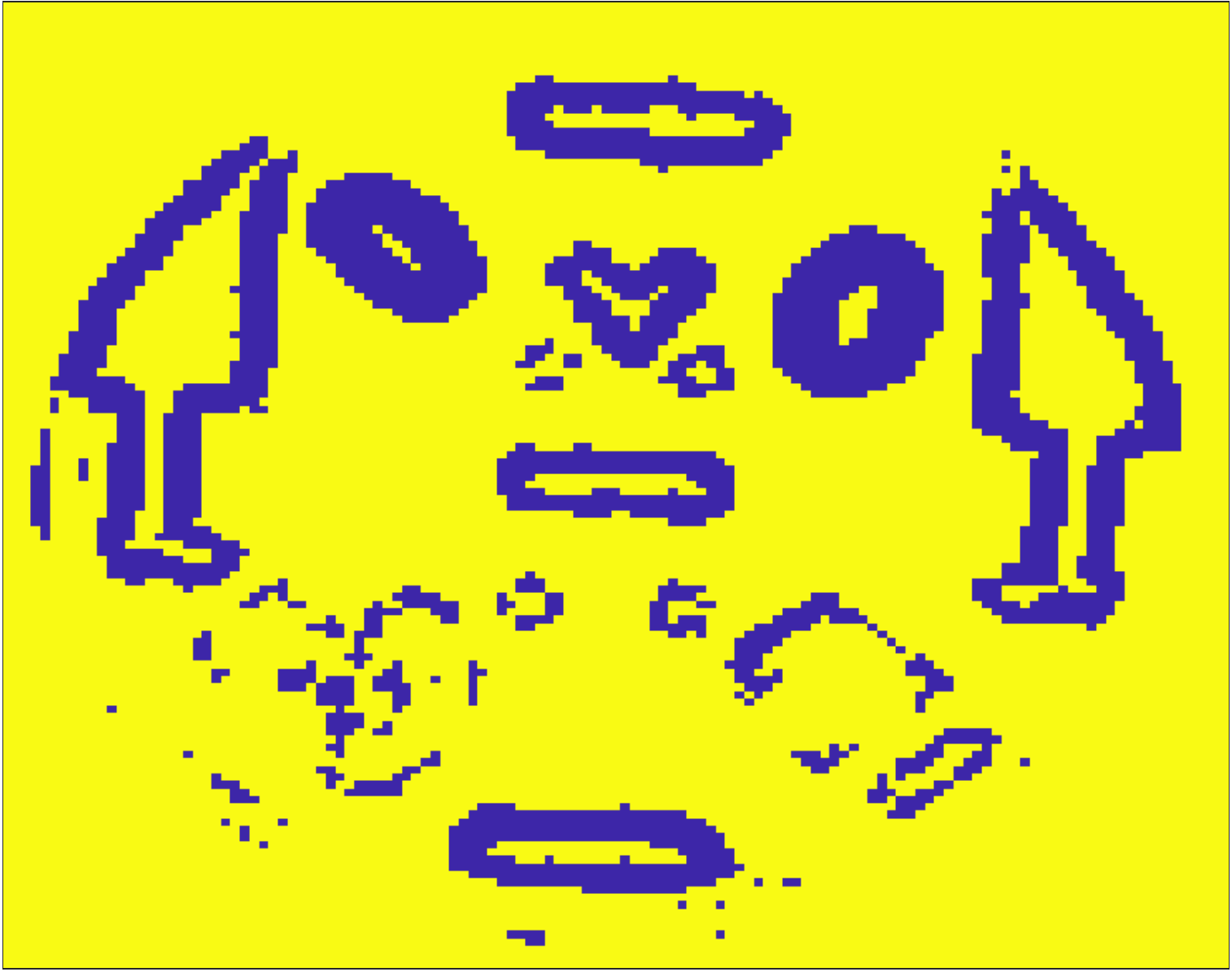}
    \caption{$W_3$}
\end{subfigure} 
~
\begin{subfigure}[b]{0.22\textwidth}
    \includegraphics[width=\textwidth]{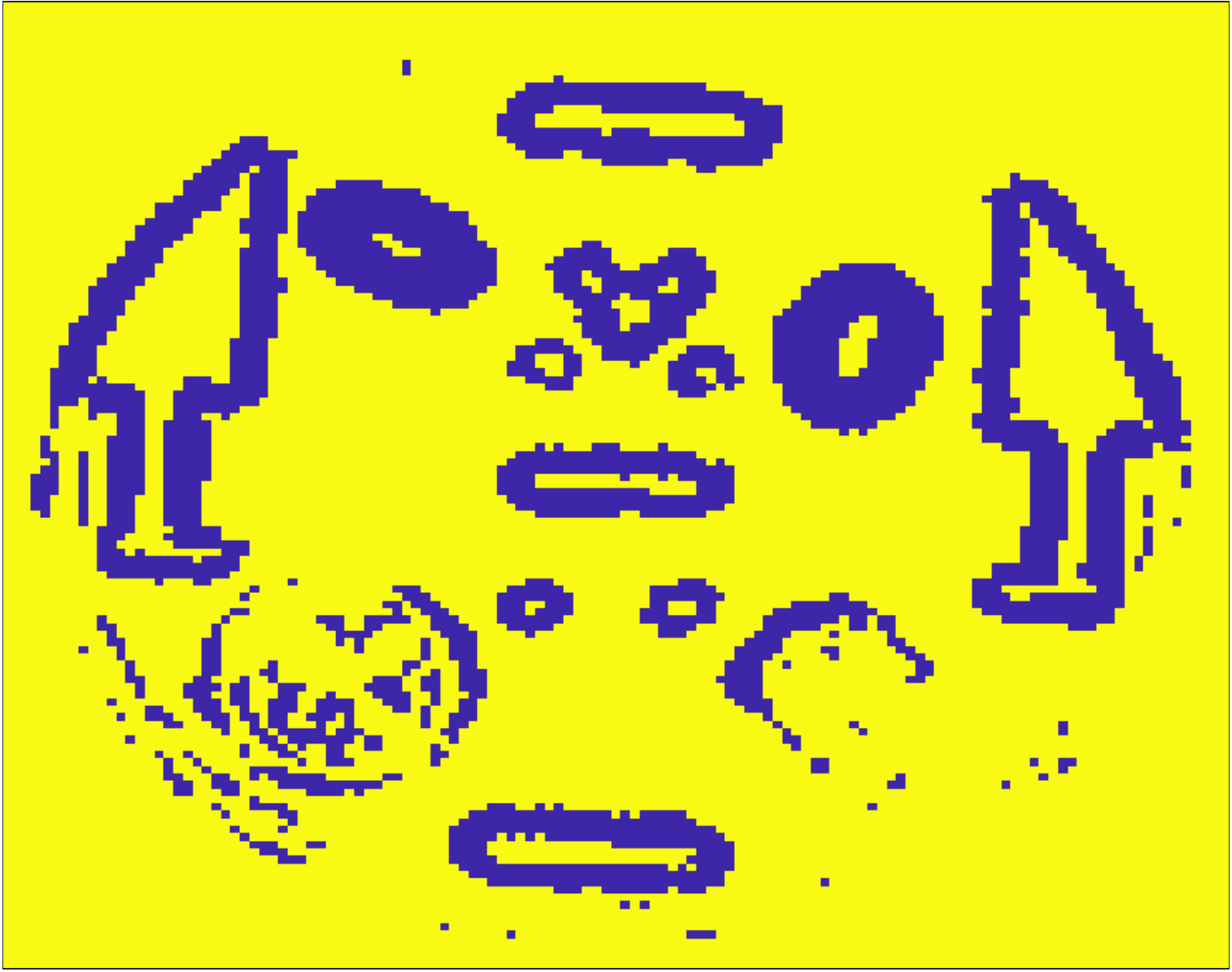}
    \caption{$W_4$}
\end{subfigure} 
\\ 
\begin{subfigure}[b]{0.22\textwidth}
    \includegraphics[width=\textwidth]{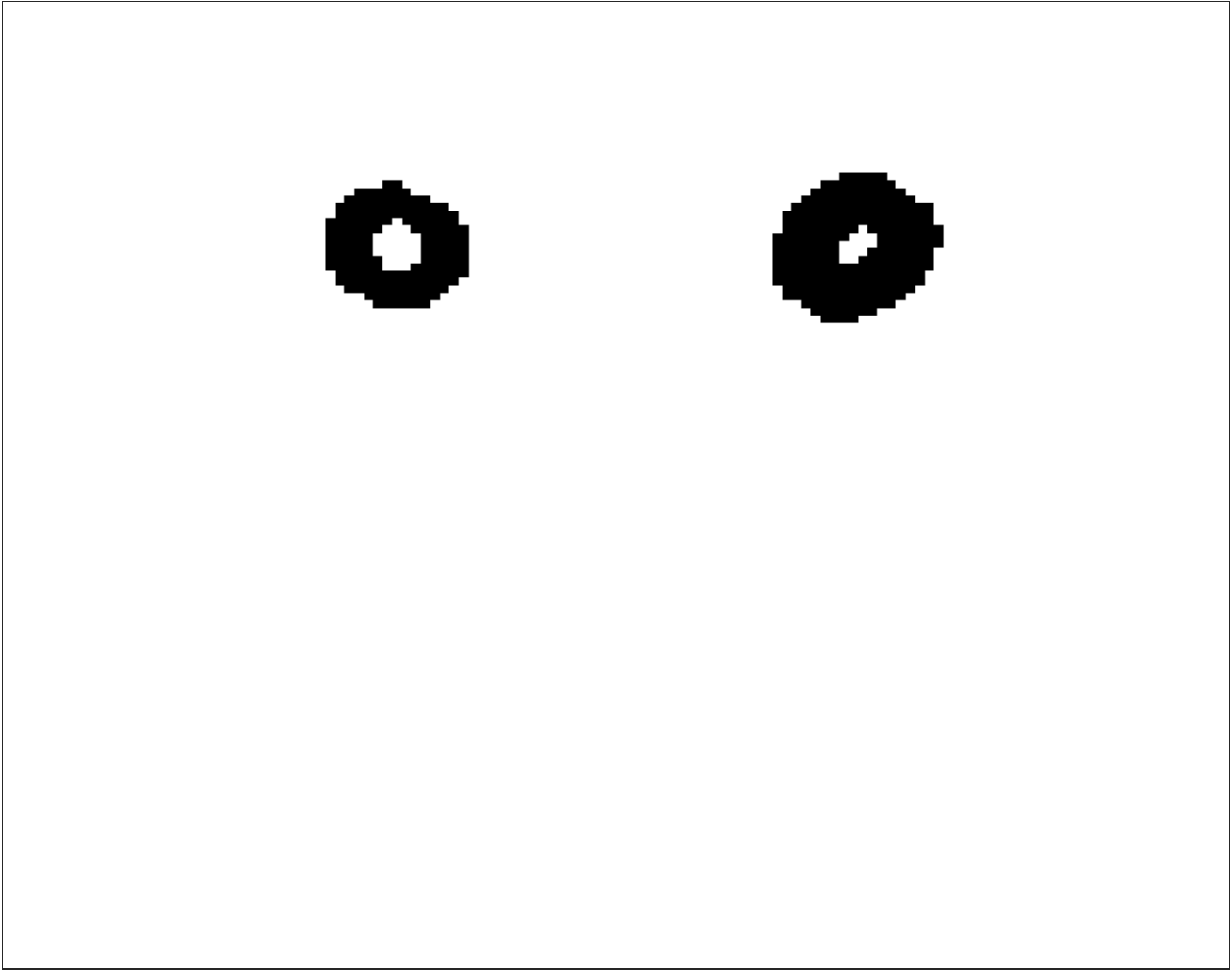}
    \caption{$\tilde{U}_1$}
\end{subfigure}
~ 
\begin{subfigure}[b]{0.22\textwidth}
    \includegraphics[width=\textwidth]{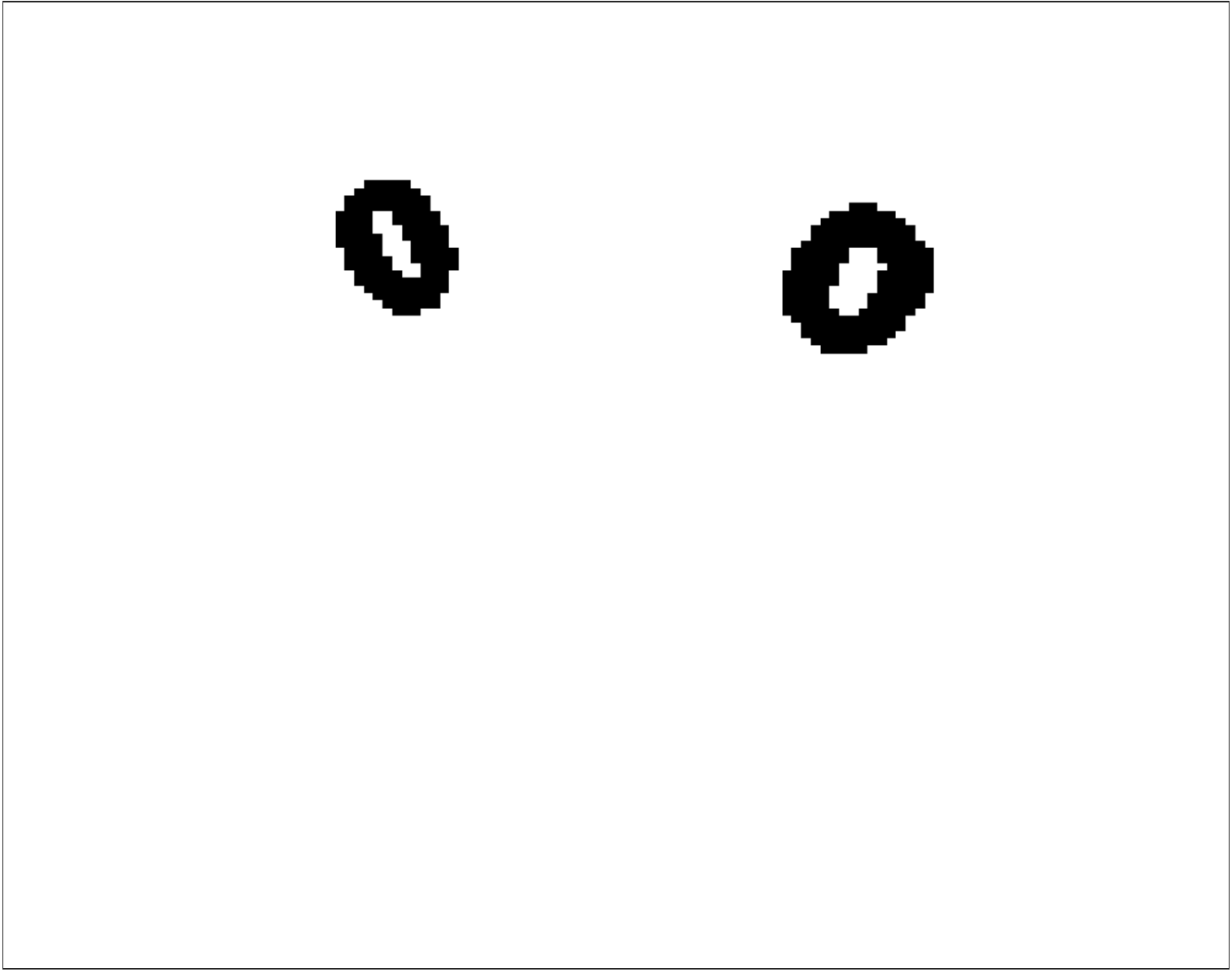}
    \caption{$\tilde{U}_2$}
\end{subfigure} 
~ 
\begin{subfigure}[b]{0.22\textwidth}
    \includegraphics[width=\textwidth]{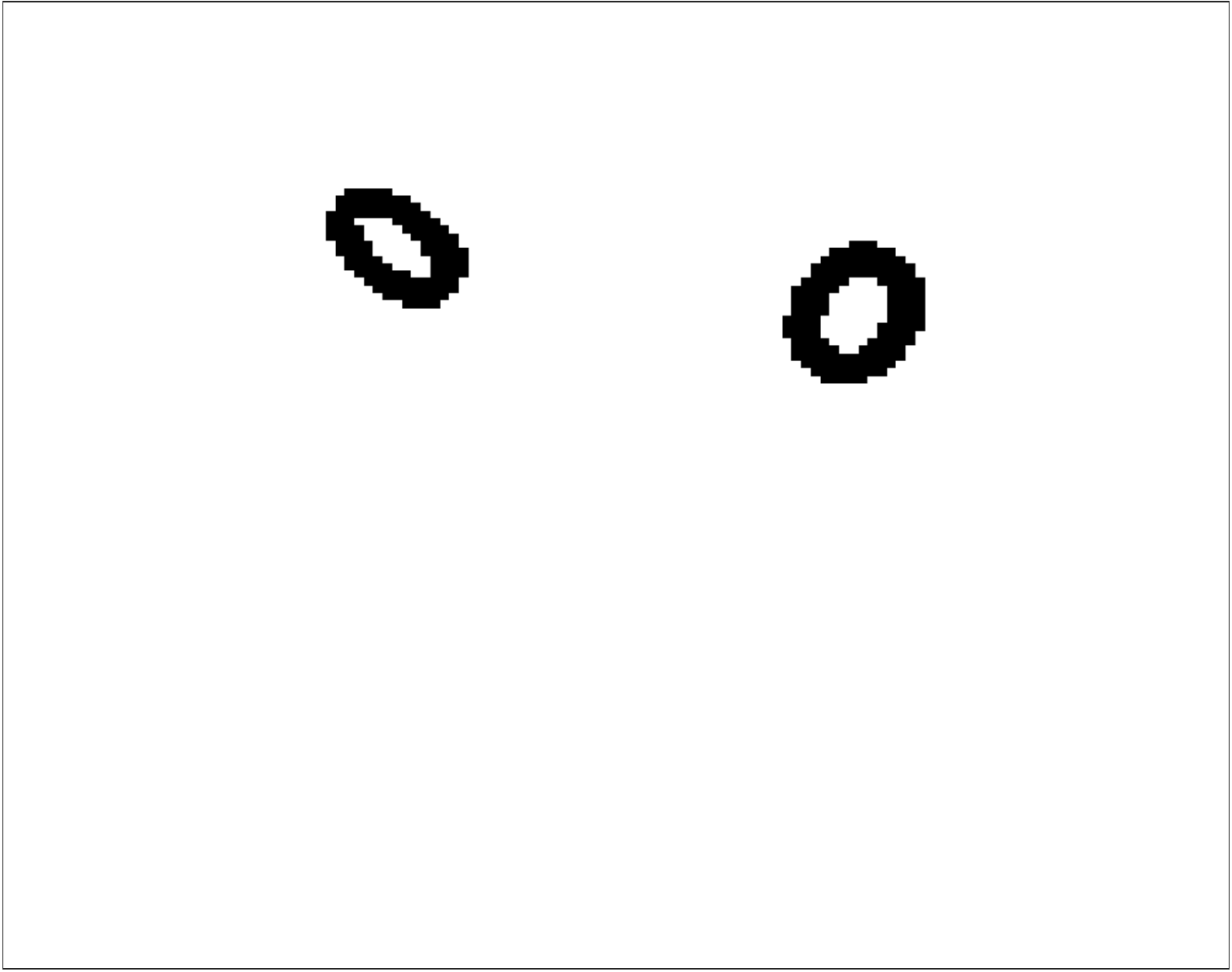}
    \caption{$\tilde{U}_3$}
\end{subfigure}
~
\begin{subfigure}[b]{0.22\textwidth}
    \includegraphics[width=\textwidth]{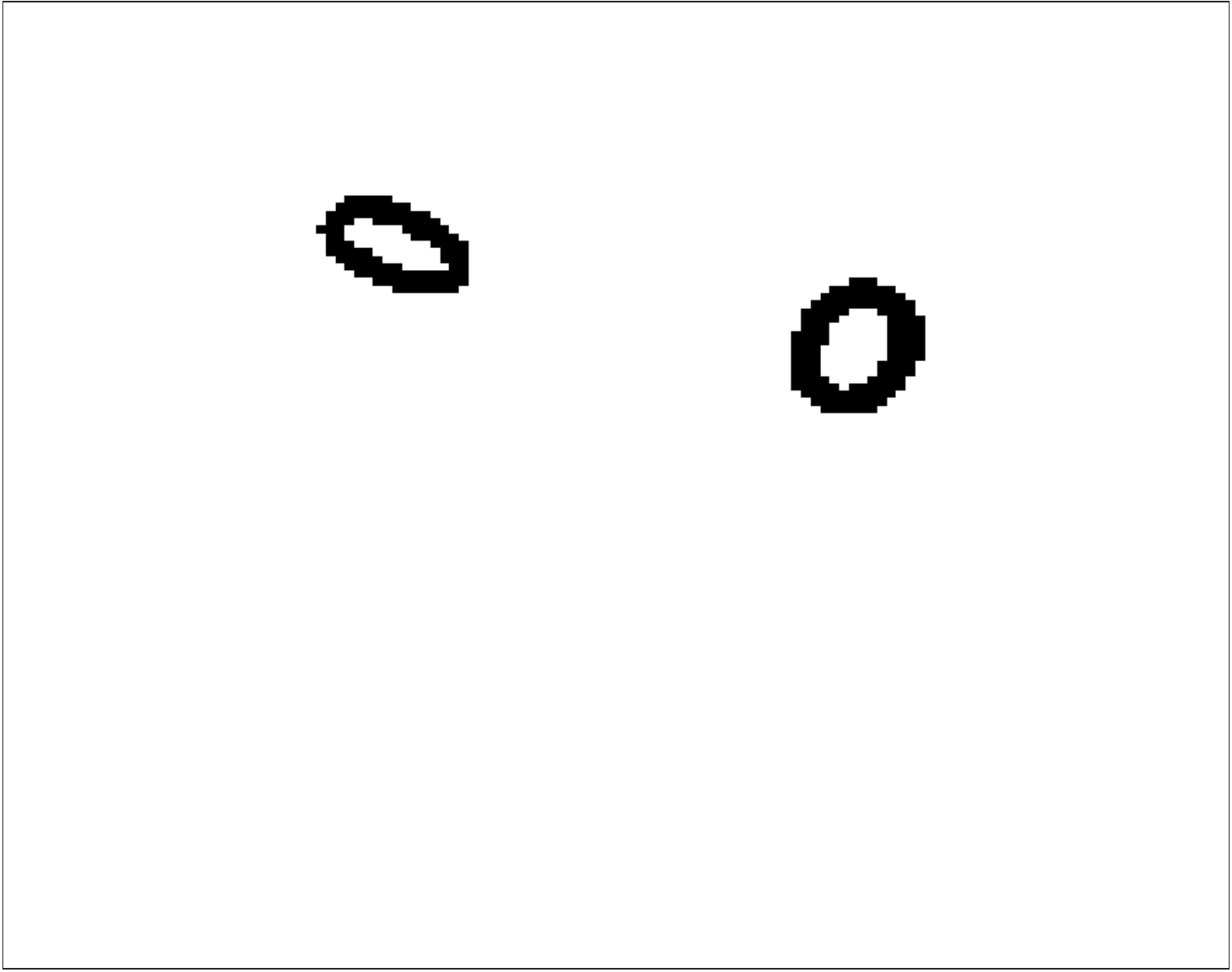}
    \caption{$\tilde{U}_4$}
\end{subfigure} 
\\ 
\begin{subfigure}[b]{0.22\textwidth}
    \includegraphics[width=\textwidth]{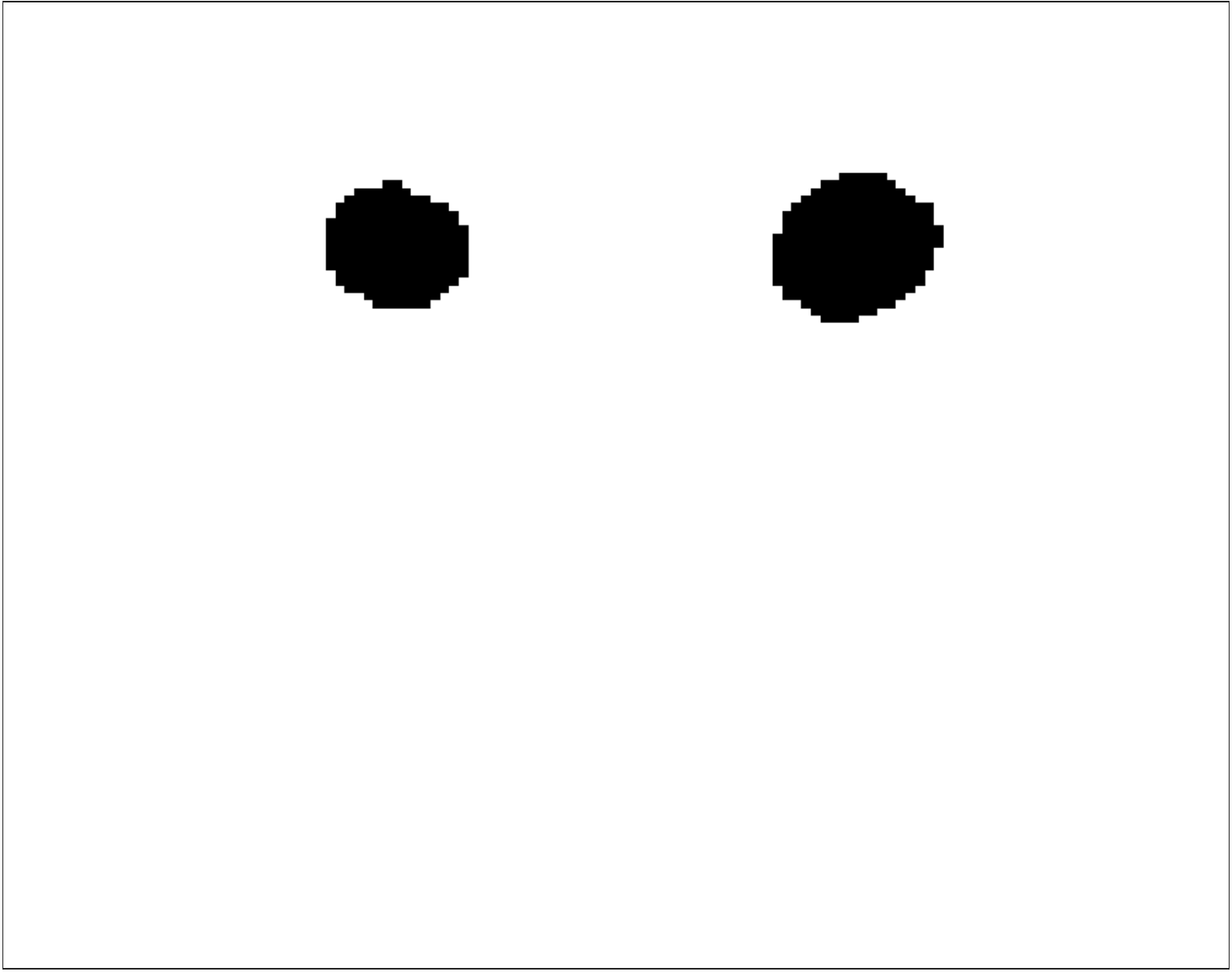}
    \caption{$\tilde{Q}_1$}
\end{subfigure}
~ 
\begin{subfigure}[b]{0.22\textwidth}
    \includegraphics[width=\textwidth]{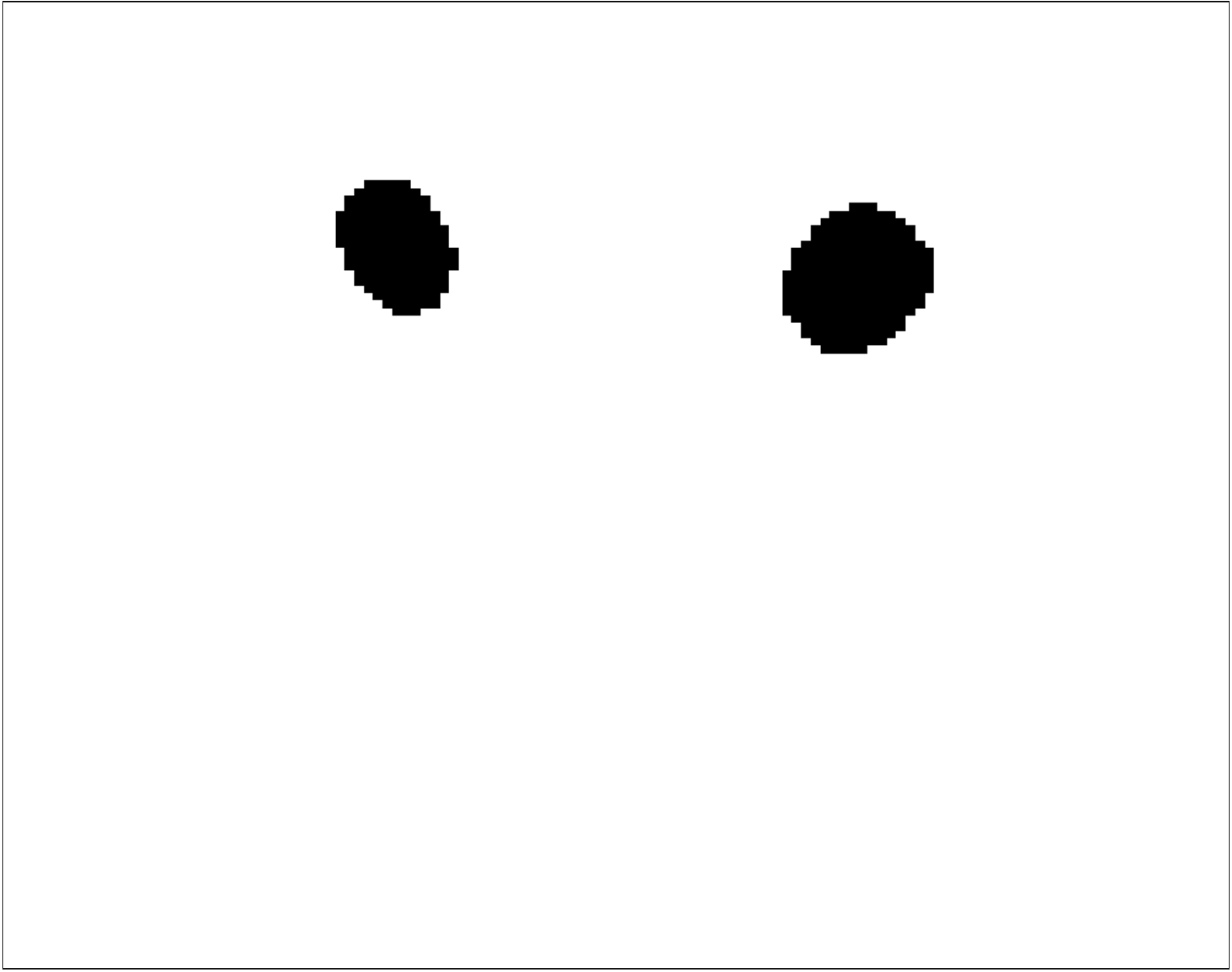}
    \caption{$\tilde{Q}_2$}
\end{subfigure} 
~ 
\begin{subfigure}[b]{0.22\textwidth}
    \includegraphics[width=\textwidth]{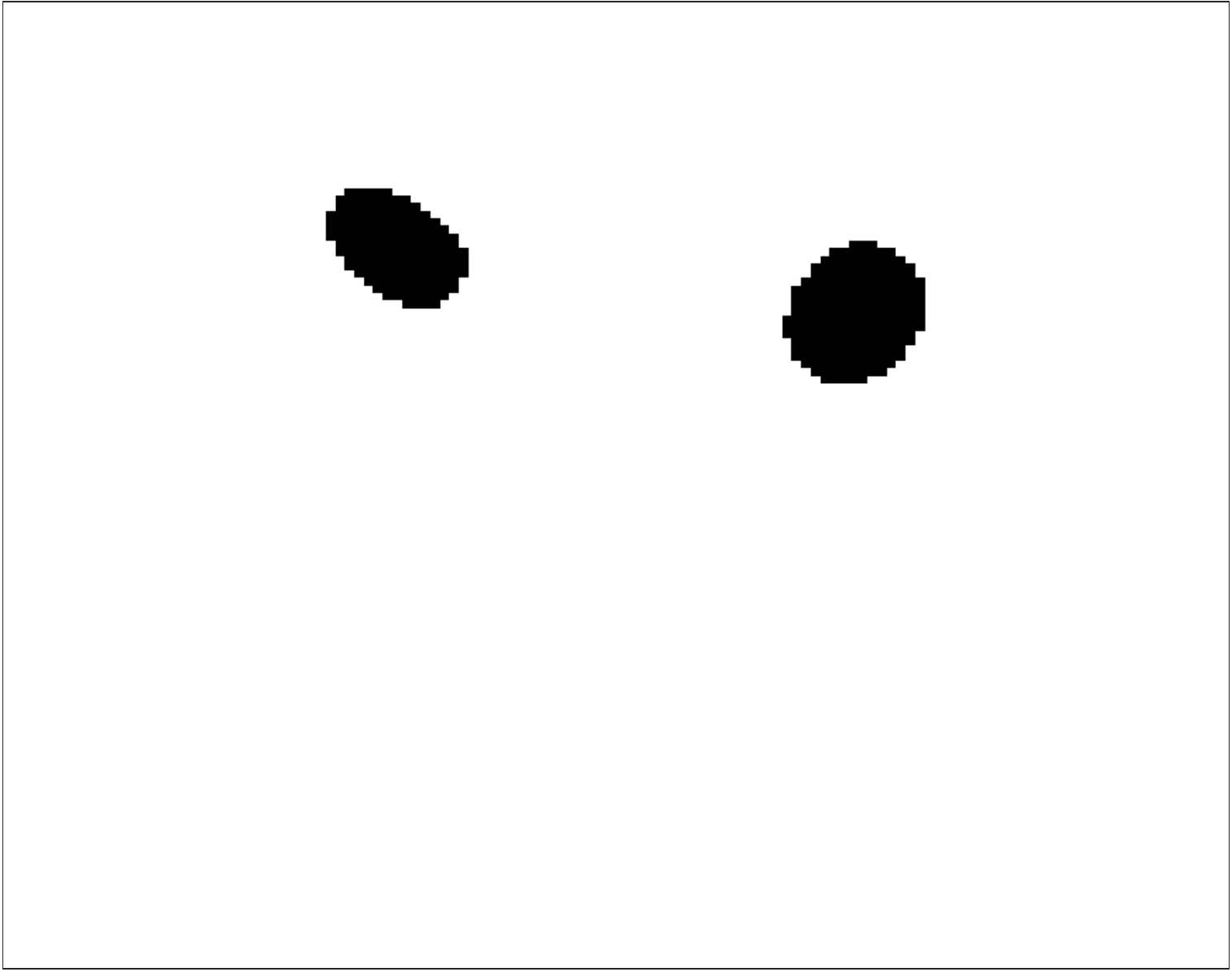}
    \caption{$\tilde{Q}_3$}
\end{subfigure} 
~
\begin{subfigure}[b]{0.22\textwidth}
    \includegraphics[width=\textwidth]{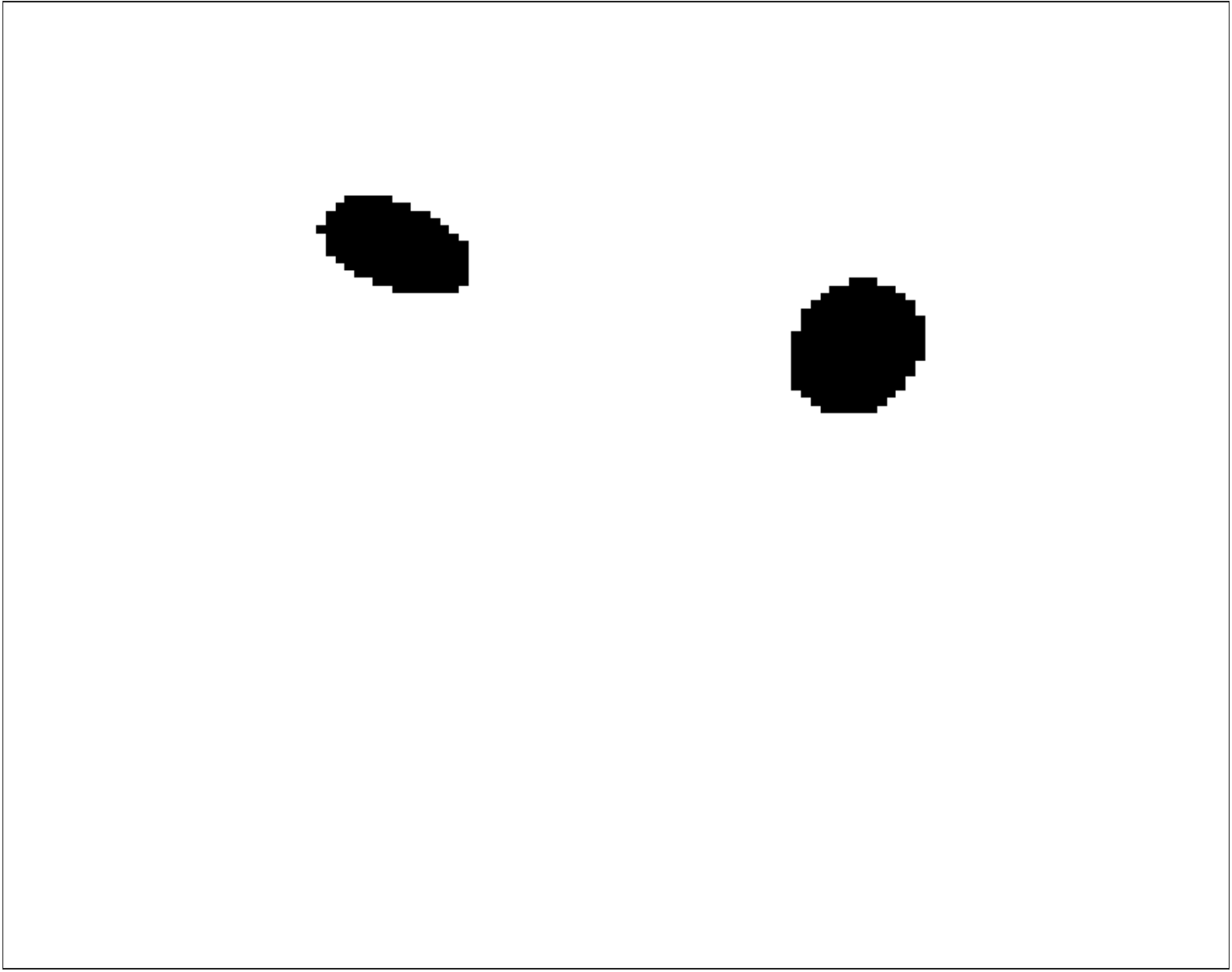}
    \caption{$\tilde{Q}_4$}
\end{subfigure} 
\\ 
\begin{subfigure}[b]{0.22\textwidth}
    \includegraphics[width=\textwidth]{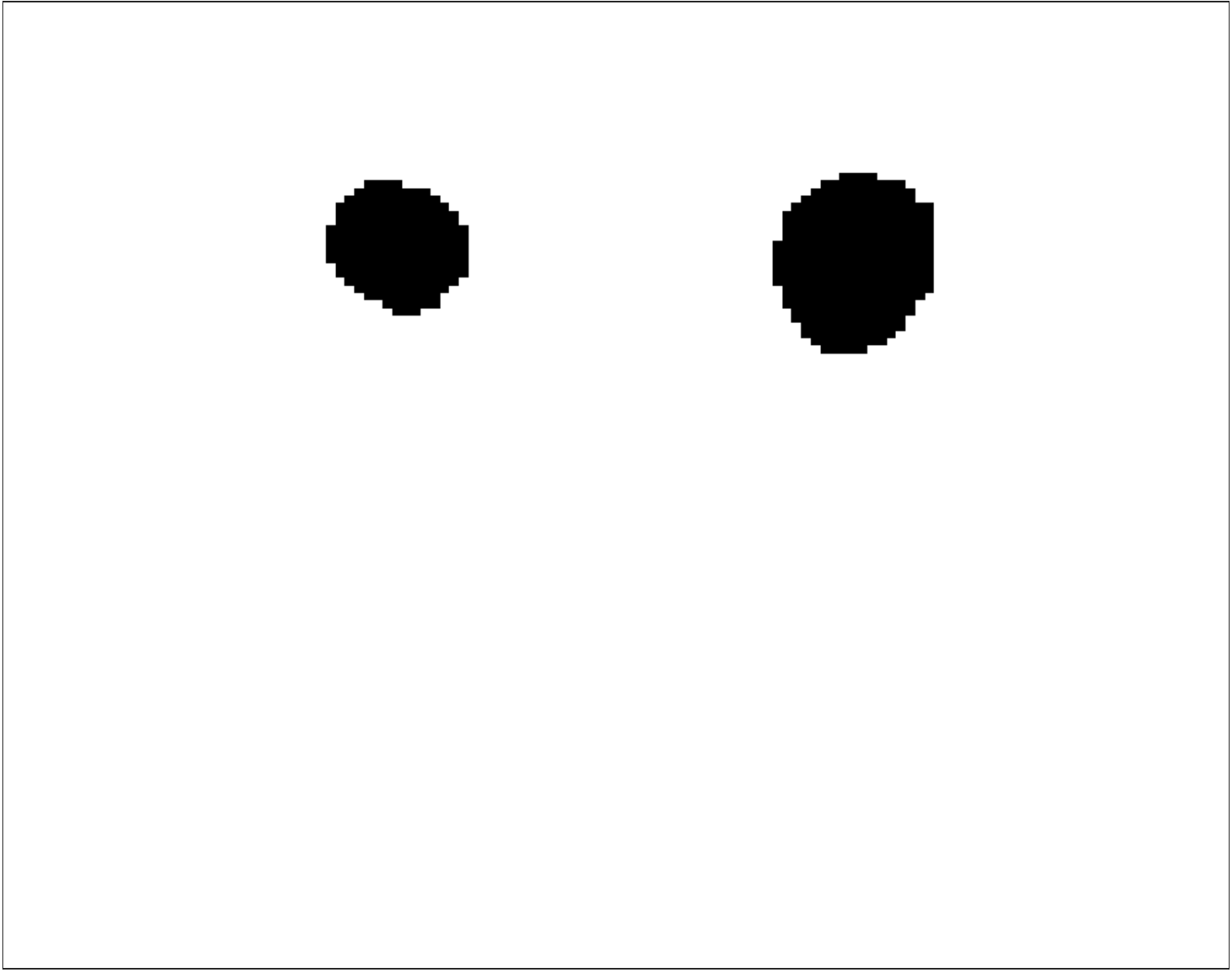}
    \caption{$\tilde{C}_1$}
\end{subfigure}
~ 
\begin{subfigure}[b]{0.22\textwidth}
    \includegraphics[width=\textwidth]{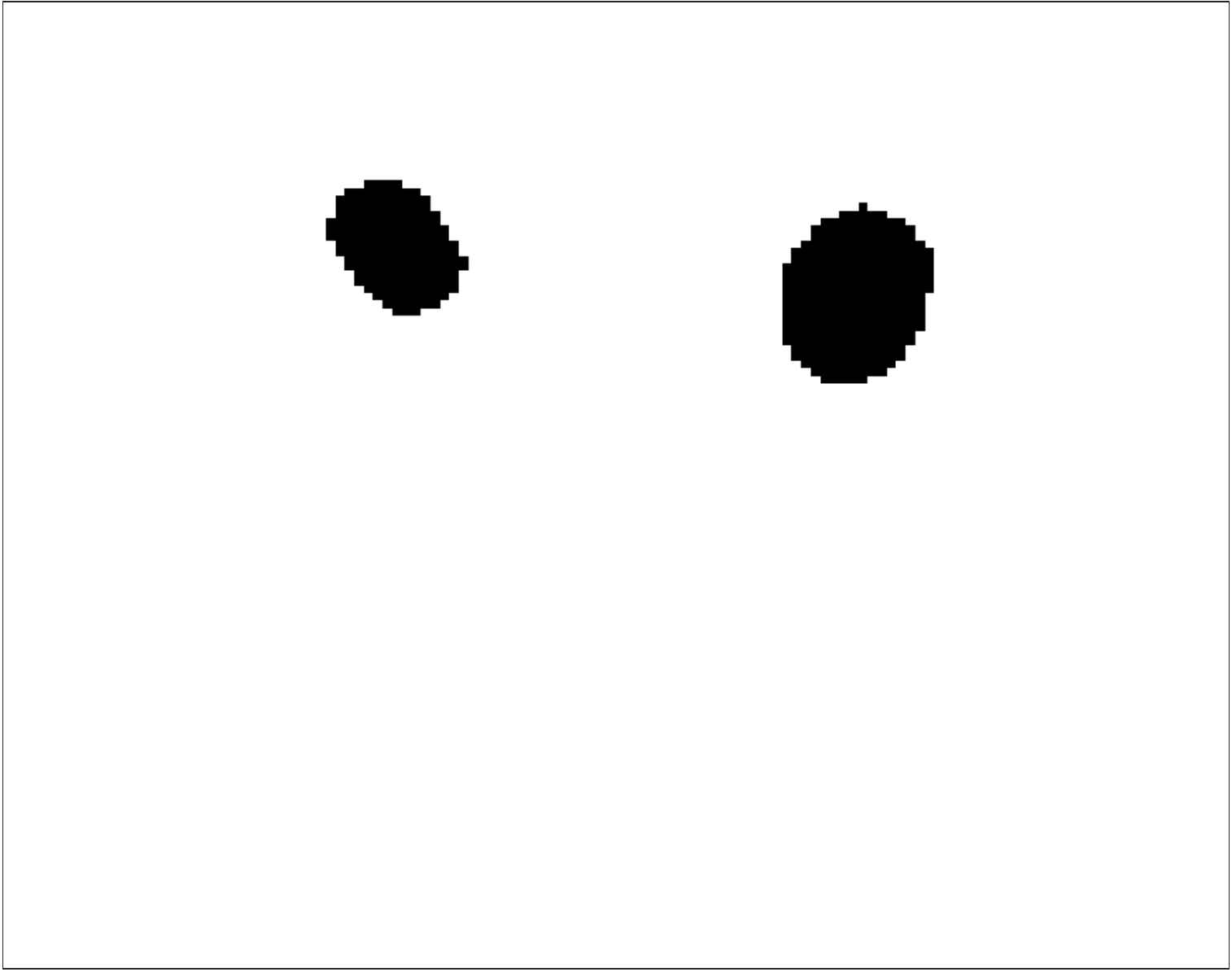}
    \caption{$\tilde{C}_2$}
\end{subfigure} 
~ 
\begin{subfigure}[b]{0.22\textwidth}
    \includegraphics[width=\textwidth]{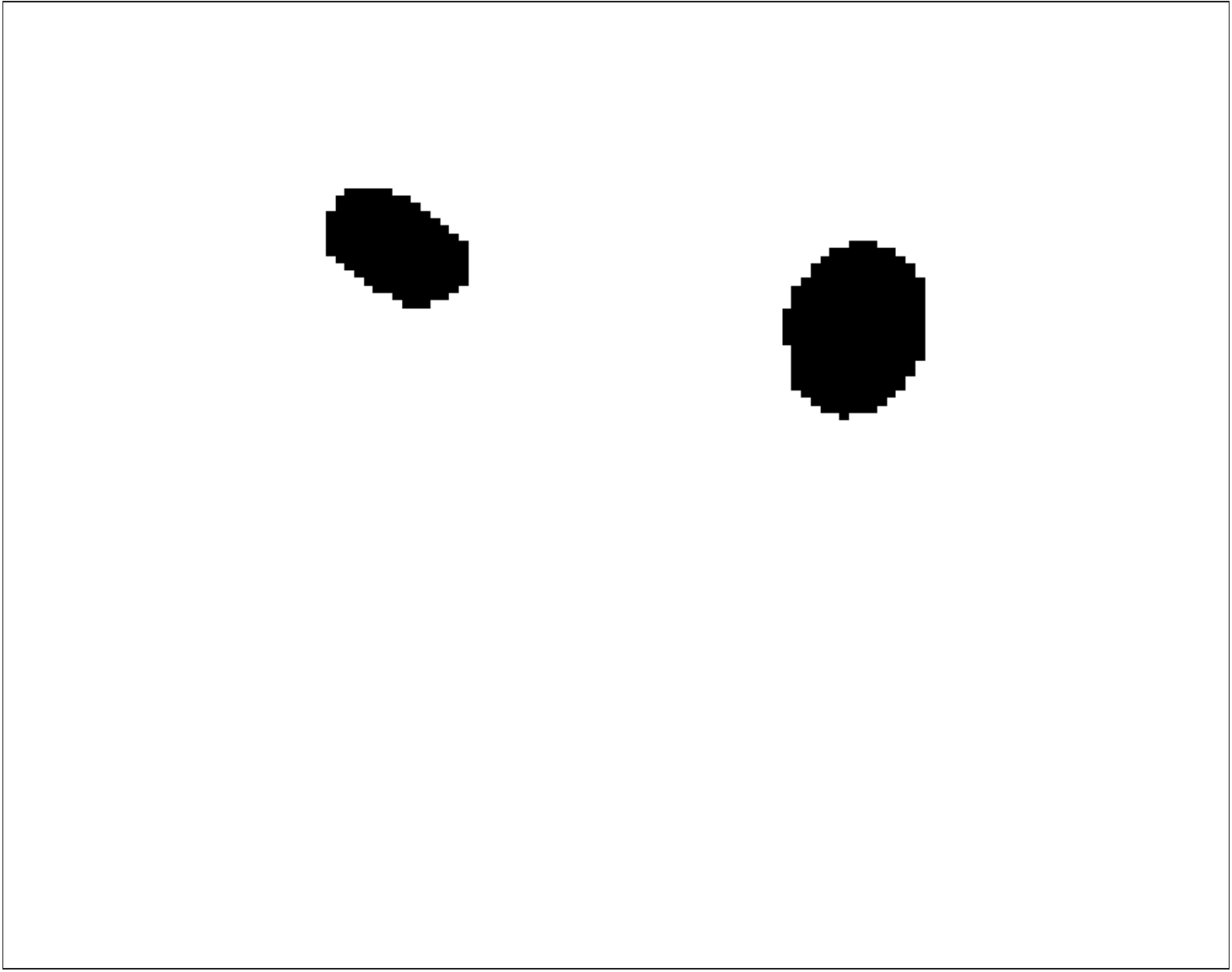}
    \caption{$\tilde{C}_3$}
\end{subfigure}
~ 
\begin{subfigure}[b]{0.22\textwidth}
    \includegraphics[width=\textwidth]{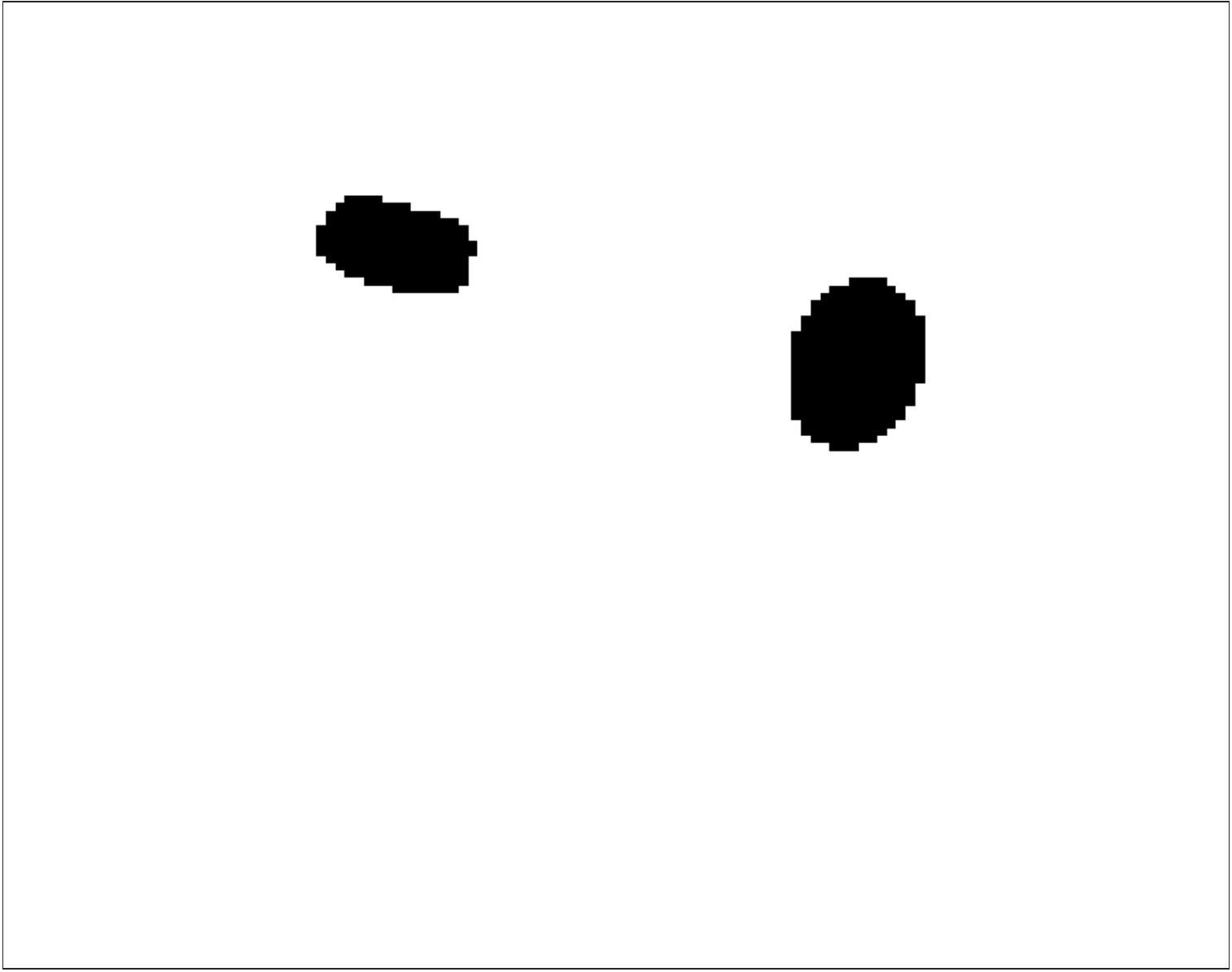}
    \caption{$\tilde{C}_4$}
\end{subfigure}
\caption{ (first row) Weighting mask $W$ calculated for each temporal image using \eqref{eq:weight_l1}.  The color bar given in the {first column} applies to all subsequent {columns}. (second {row}) {Single object} edge masks ${\tilde U}$ in  \eqref{eq:ucurve}. (third {row}) {Filled single object} masks ${\tilde Q}$ in \eqref{eq:edge-region-mask}.  (fourth {row}) Change mask ${\tilde  C}$ in \eqref{eq:changemask}. Observe from \eqref{eq:CM-edge-region-mask} that $\tilde{C}_4$ requires construction of $W_5$, $\tilde{U}_5$, and $\tilde{Q}_5$, which are not pictured.
}
\label{fig:phy_GE}
\end{figure}

\begin{figure}[h!]
\centering
\begin{subfigure}[b]{.24\textwidth}
    \includegraphics[width=\textwidth]{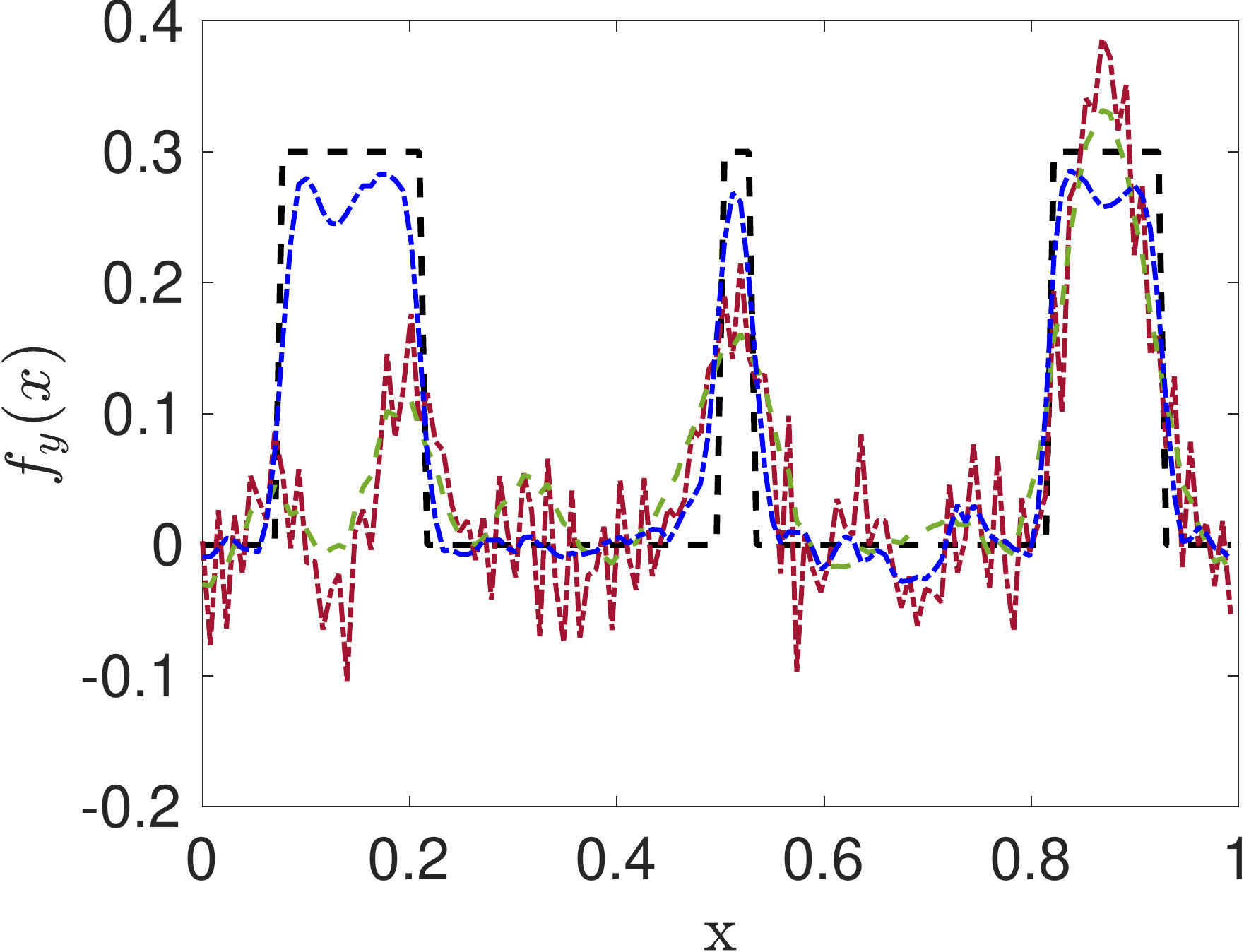}
    \caption{$\mathbf f_1$, $y=0.6512$.}
\end{subfigure}
\begin{subfigure}[b]{.24\textwidth}
    \includegraphics[width=\textwidth]{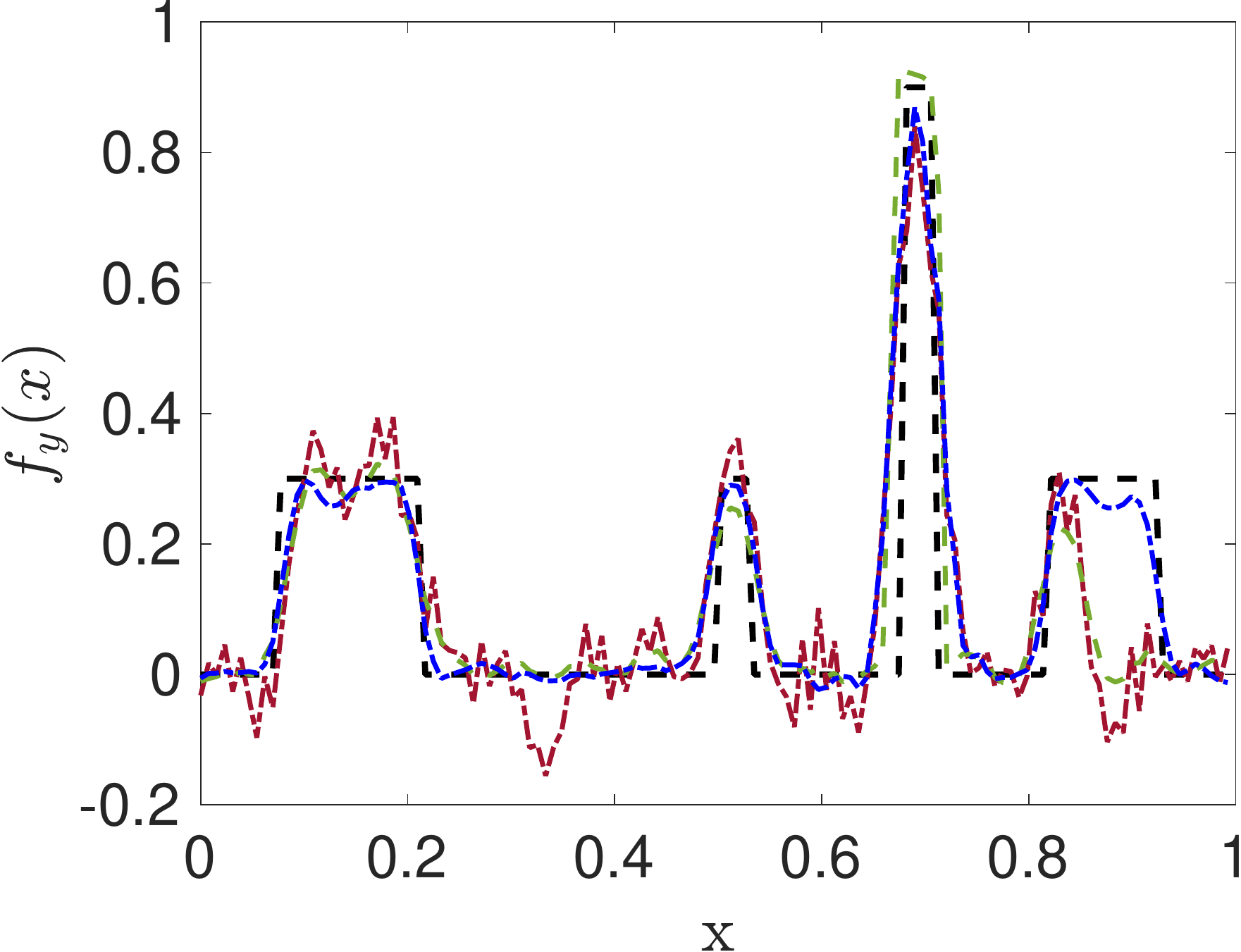}
    \caption{$\mathbf f_2$, $y=0.6589$.}
\end{subfigure}
\begin{subfigure}[b]{.24\textwidth}
    \includegraphics[width=\textwidth]{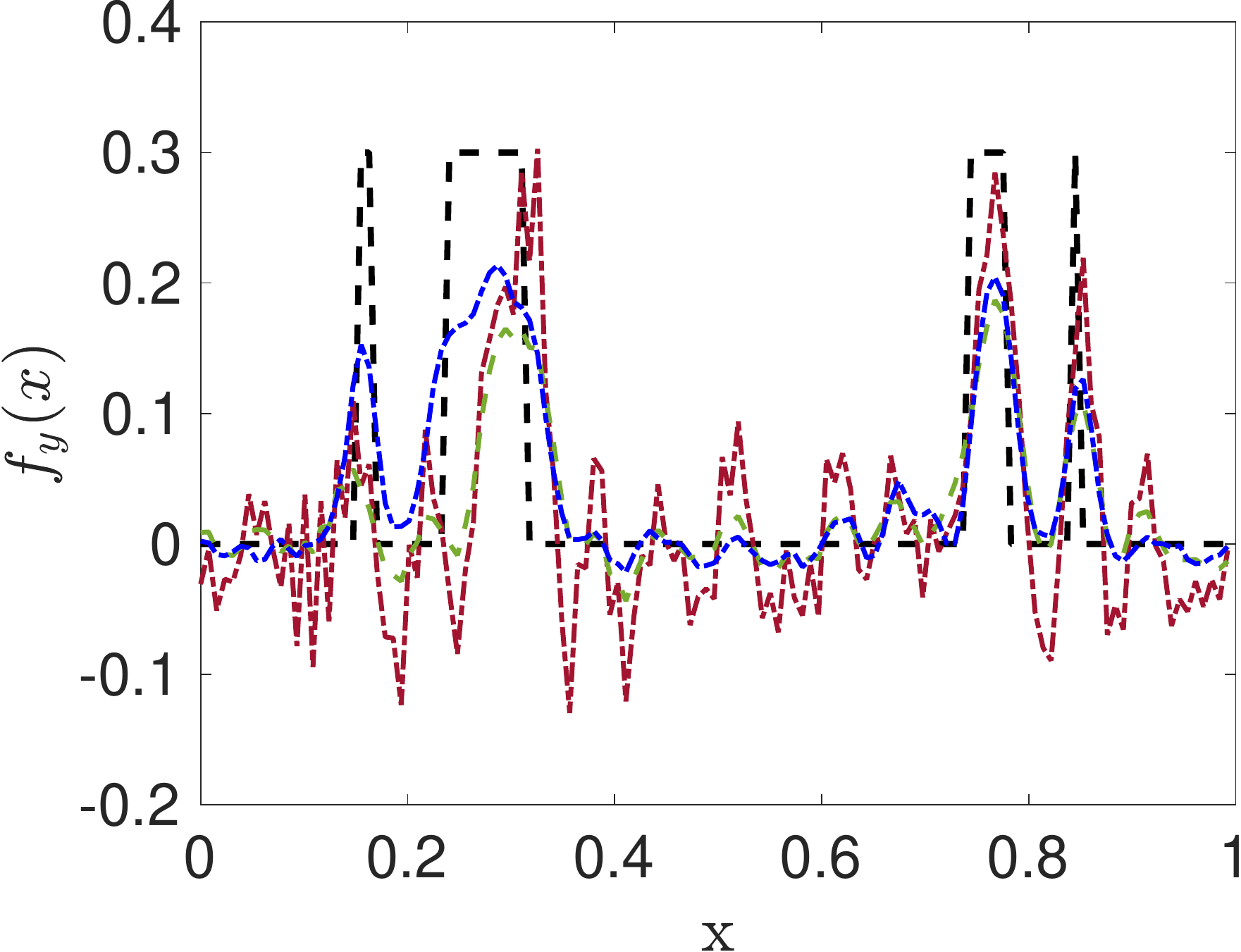}
    \caption{$\mathbf f_3$, $y=0.2016$.}
\end{subfigure}
\begin{subfigure}[b]{.24\textwidth}
    \includegraphics[width=\textwidth]{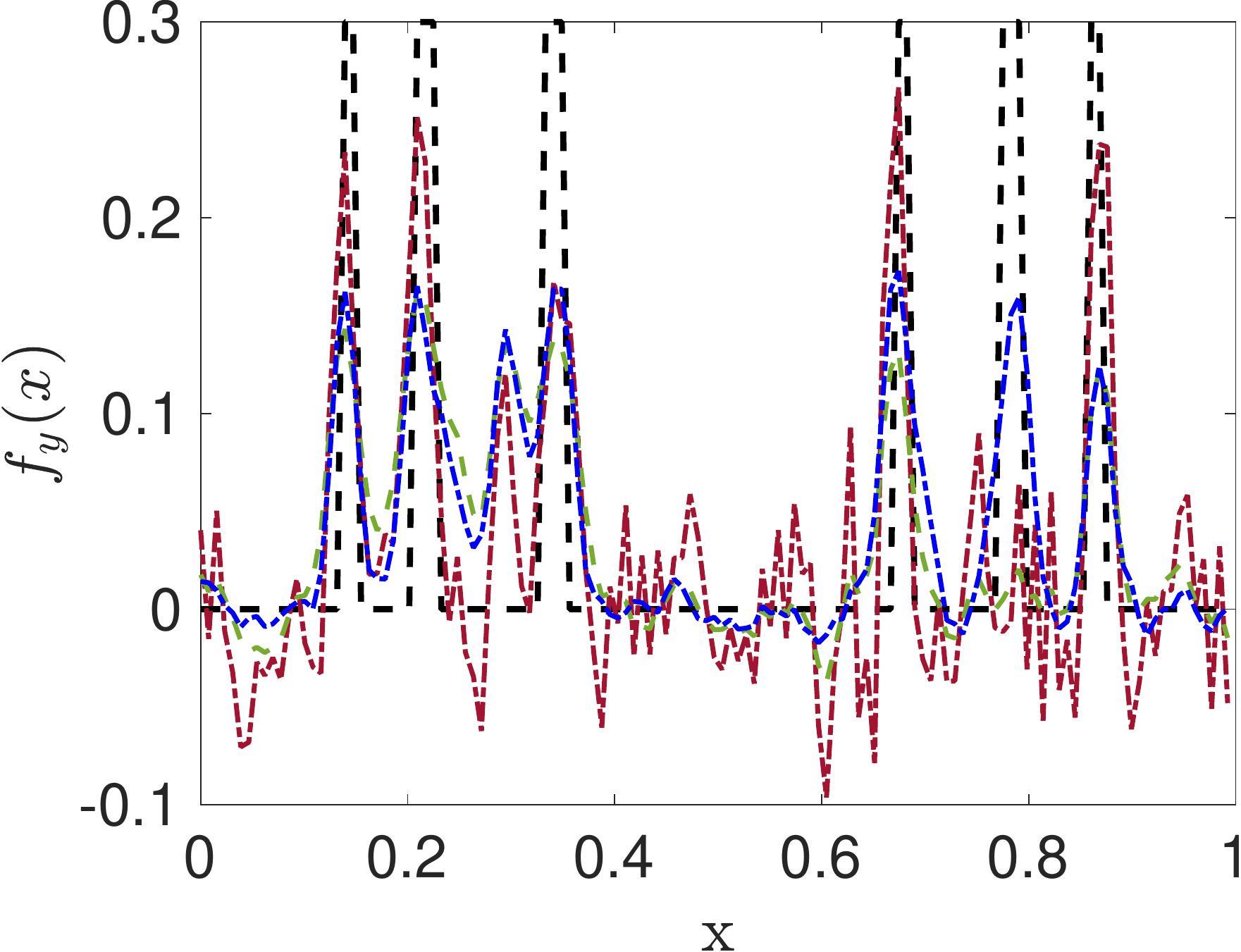}
    \caption{$\mathbf f_4$, $y=0.2248$.}
\end{subfigure}
\\
\begin{subfigure}[b]{.24\textwidth}
    \includegraphics[width=\textwidth]{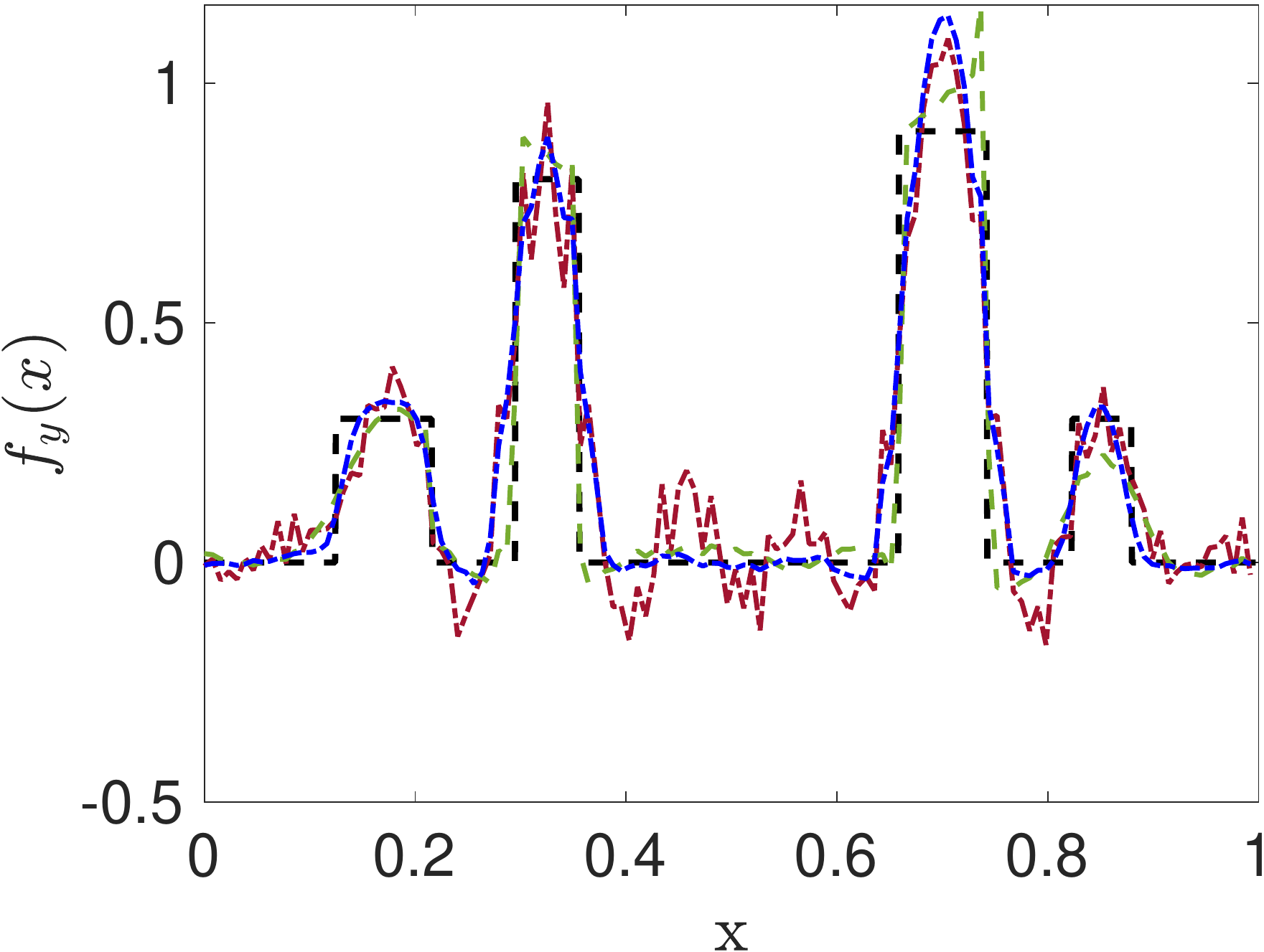}
    \caption{$\mathbf f_1$, $y = .7519$.}
\end{subfigure}
\begin{subfigure}[b]{.24\textwidth}
    \includegraphics[width=\textwidth]{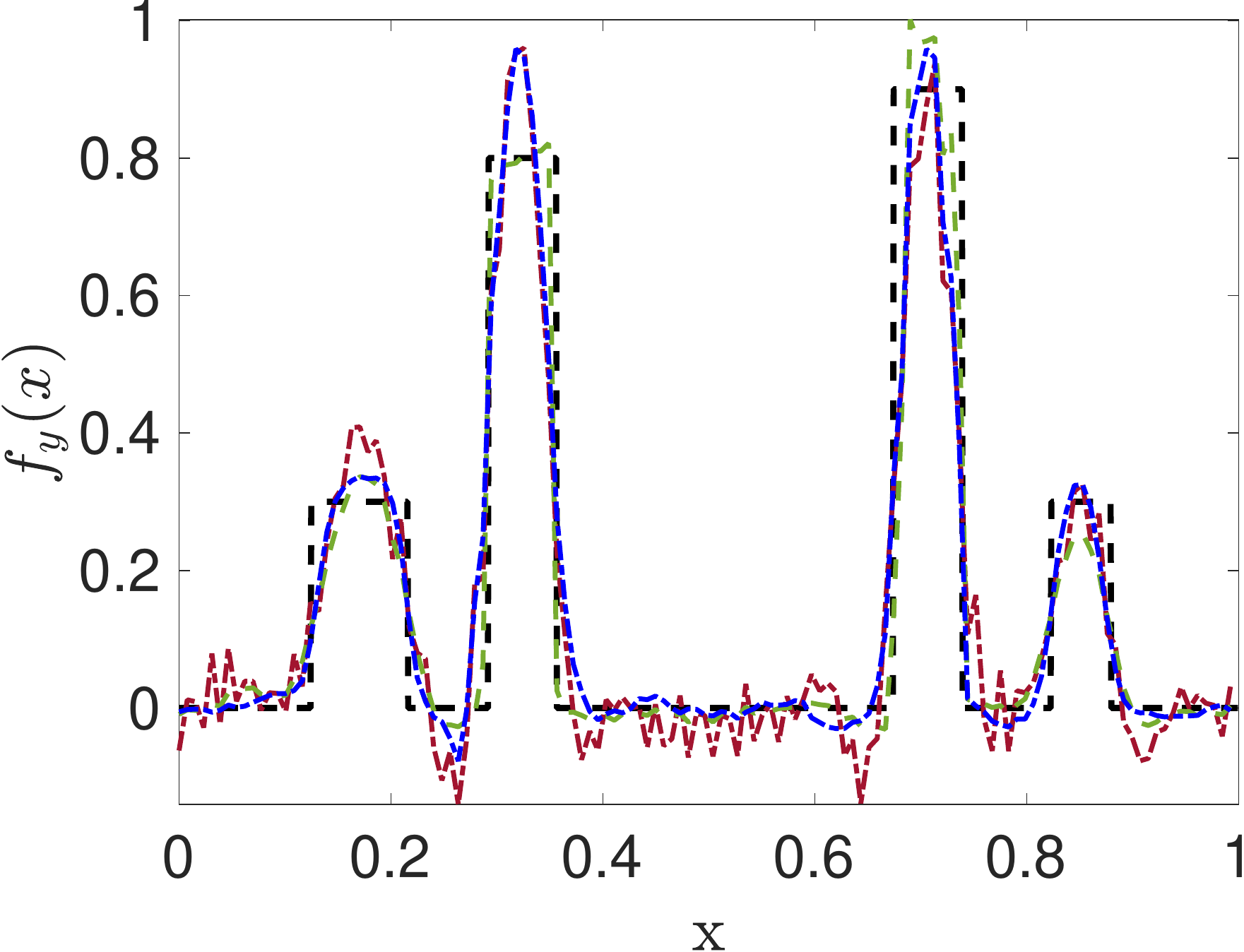}
    \caption{$\mathbf f_2$, $y = .7519$.}
\end{subfigure}
\begin{subfigure}[b]{.24\textwidth}
    \includegraphics[width=\textwidth]{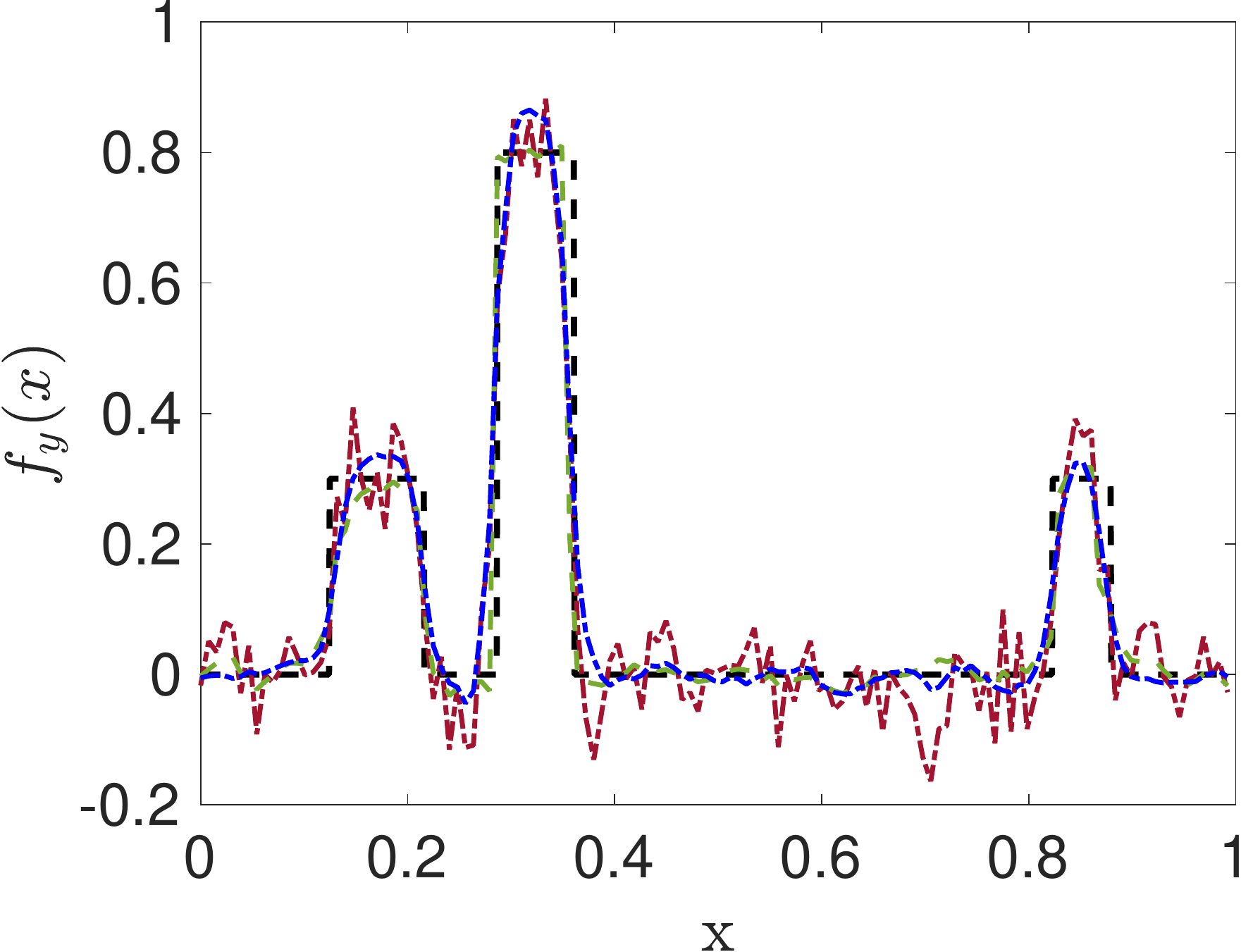}
    \caption{$\mathbf f_3$, $y = .7519$.}
\end{subfigure}
\begin{subfigure}[b]{.24\textwidth}
    \includegraphics[width=\textwidth]{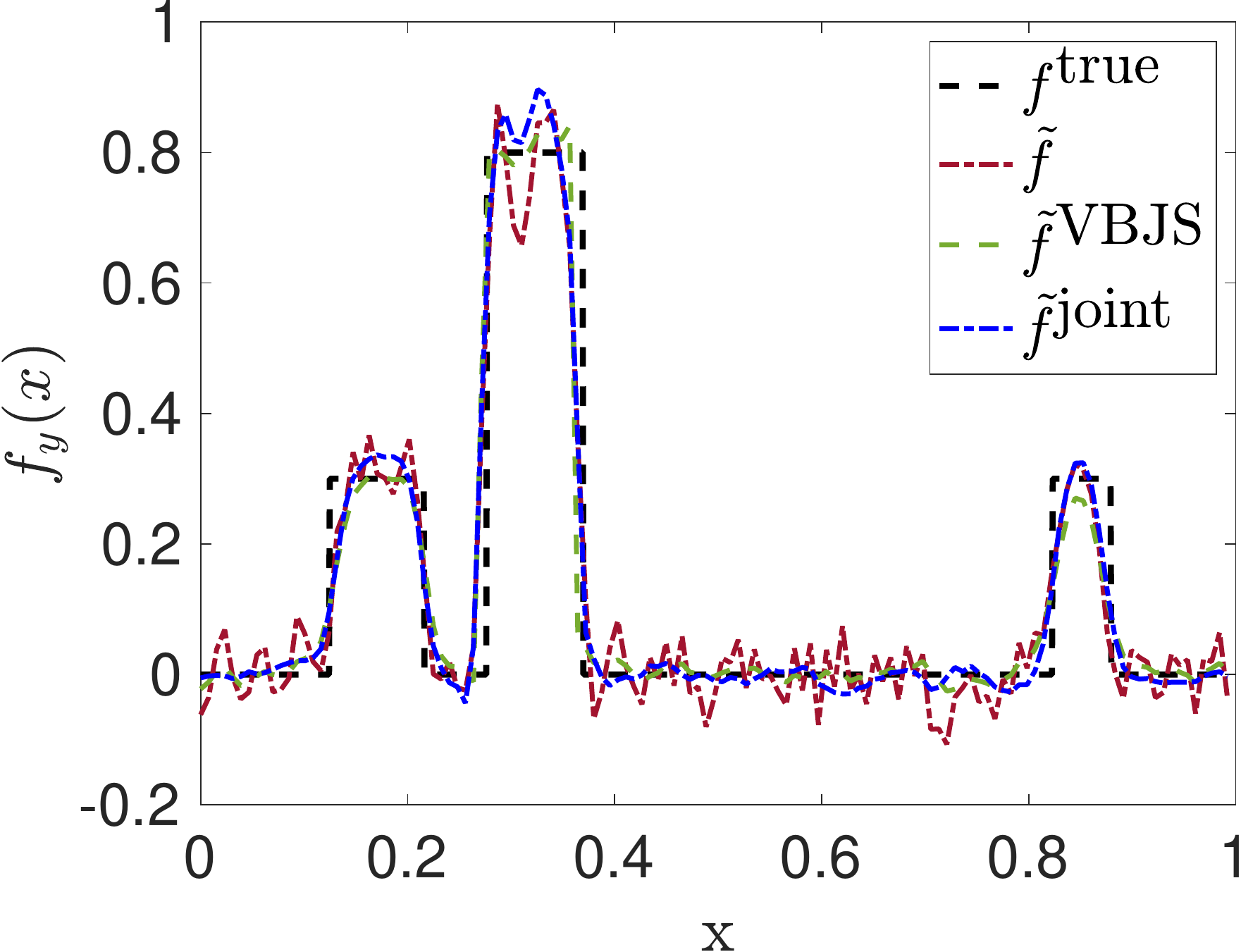}
    \caption{$\mathbf f_4$, $y = .7519$.}
\end{subfigure}

\caption{(top) Horizontal cross sections intersecting obstructed regions in the sequential regions.  (bottom) Horizontal cross sections intersecting $y = .7519$, where each of the sequential image contains a portion of the two moving ellipses. The legend is provided in the bottom right figure.
}
\label{fig:1d_cross_unmoved_GE}
\end{figure}

The image recovery given by \eqref{eq:optModel} and realized in Algorithm \ref{algo:joint_recovery} requires the construction of weight mask $W$, given in \eqref{eq:weight_l1}, and the inter-image coupling term $\Phi$. To construct $\Phi$ we need the single object edge masks, ${\tilde{U}}$ in \eqref{eq:ucurve}, followed by the {filled single object}  masks, ${\tilde Q}$ in \eqref{eq:edge-region-mask}, and finally the change masks, $\tilde{C}$ in \eqref{eq:changemask}.  As they are vital to the process, we include Figure \ref{fig:phy_GE} to show how each is respectively formed.   
We note that in order to satisfy property \ref{item:ii} for constructing \eqref{eq:ucurve} (see Remark \ref{remark:overlap}), it is necessary  to extract the skull once the edges in the image are determined.  In this particular example, since all of the background structures have the same magnitude as the skull, they are simultaneously removed.  Our second experiment in Section \ref{subsec:golf_course} illustrates the case where all of the single edge mask properties are initially satisfied so that no initial extraction is needed.
Once again we point out that $W$, and subsequently $\Phi$, are determined from the given data in \eqref{eq:forwardmodel_j}, that is, images were not pre-formed in their construction.

Figure \ref{fig:1d_cross_unmoved_GE} displays some one-dimensional cross sections of the images.  The top row shows the horizontal slices that intersect the zero-valued occlusions in the first four of the sequenced images.  
In this case we observe the benefit of {using inter-image information to recover each sequential image.} The results displayed in the bottom row correspond to a horizontal slice that intersects portions of the two moving ellipses.  Here we see the importance of accurately constructing the inter-image information matrix $\Phi$ in \eqref{eq:Phi}, specifically to avoid  sharing inconsistent  information between sequential data sets.

\subsubsection*{Numerical Convergence Comparison}

\begin{figure}[h!]
    \centering
    \begin{subfigure}[b]{.23\textwidth}
        \includegraphics[width=\textwidth]{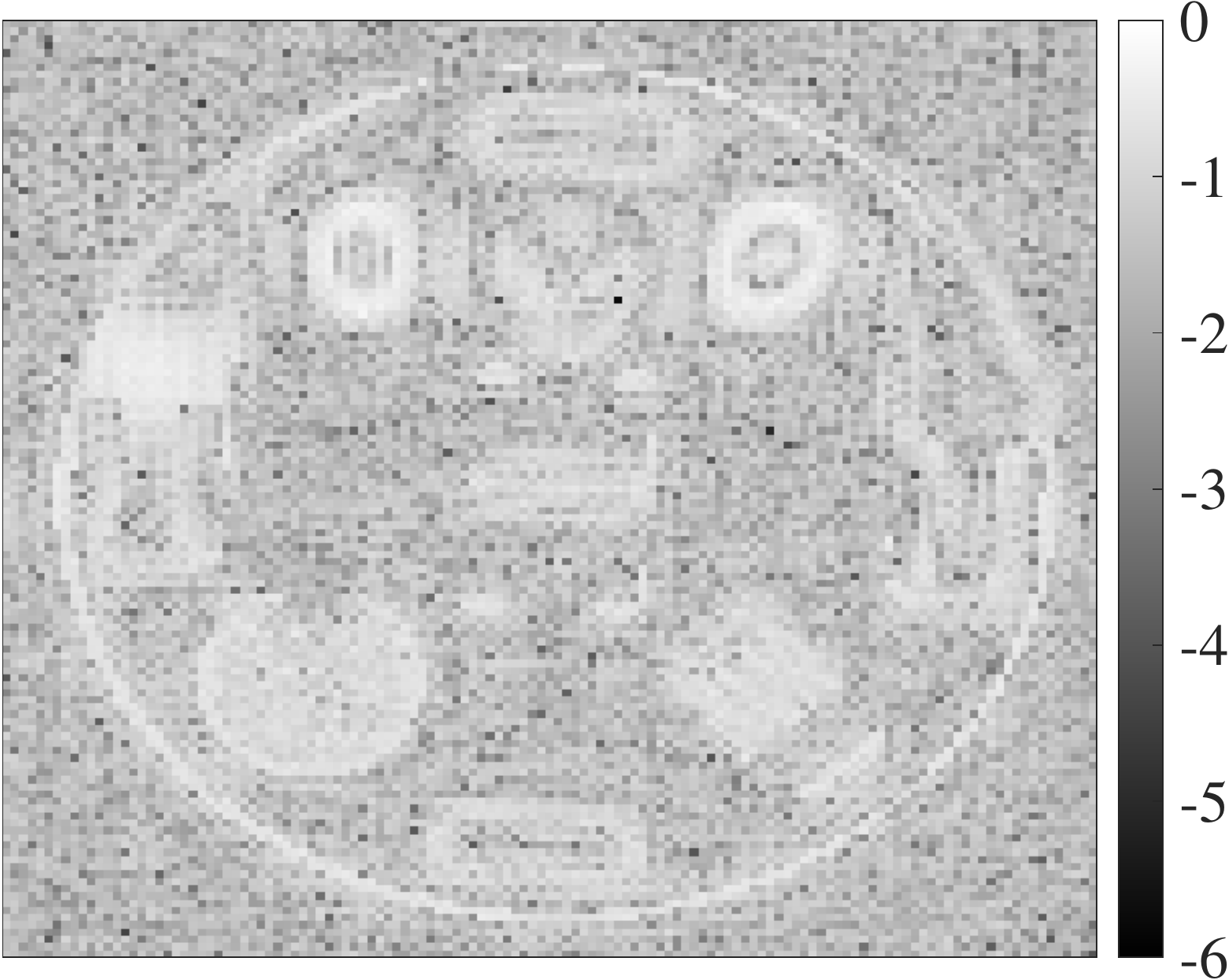}
        \caption{$\tilde f_1$}
    \end{subfigure}
    ~
    \begin{subfigure}[b]{.23\textwidth}
        \includegraphics[width=\textwidth]{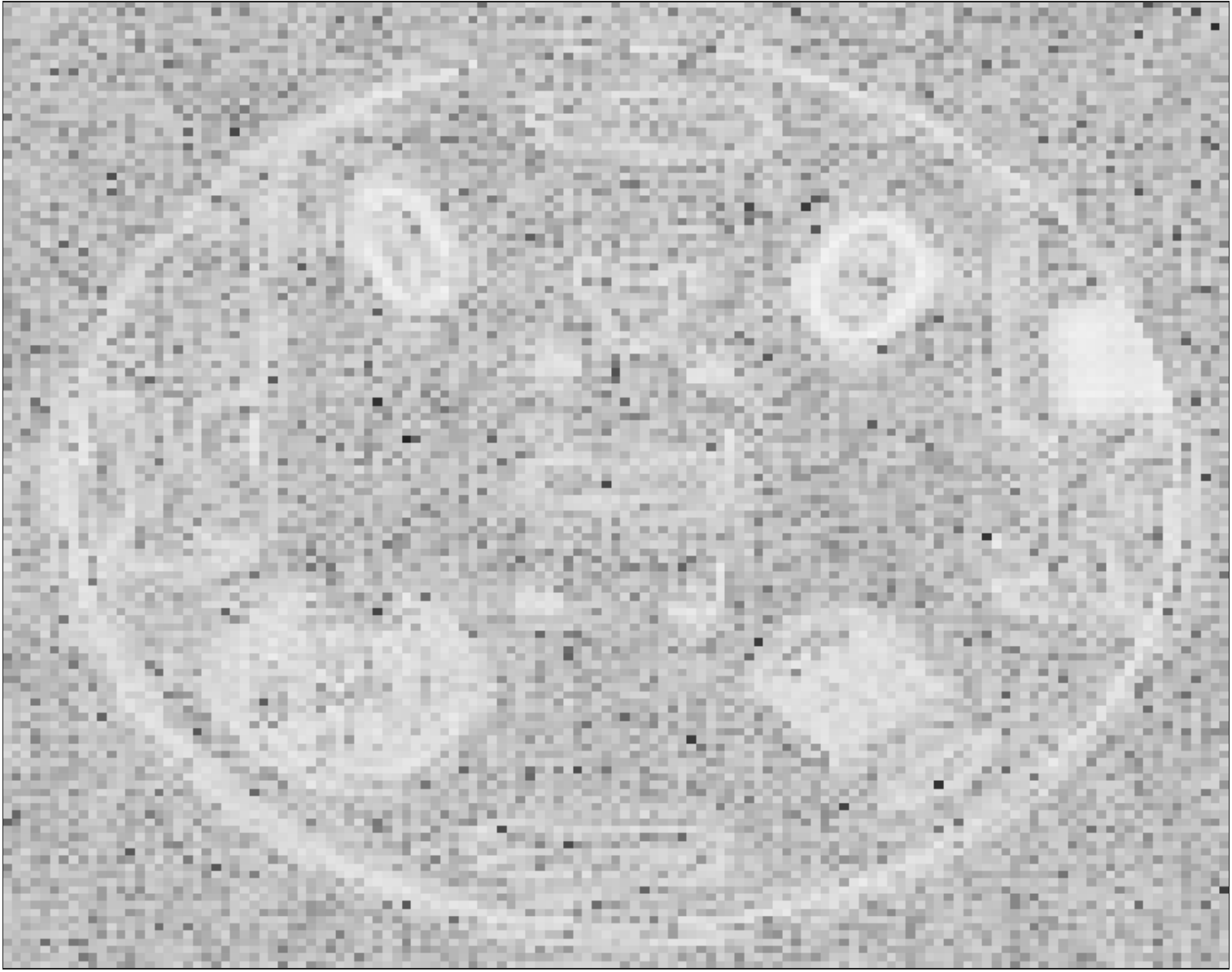}
        \caption{$\tilde f_2$}
    \end{subfigure}
    ~
    \begin{subfigure}[b]{.23\textwidth}
        \includegraphics[width=\textwidth]{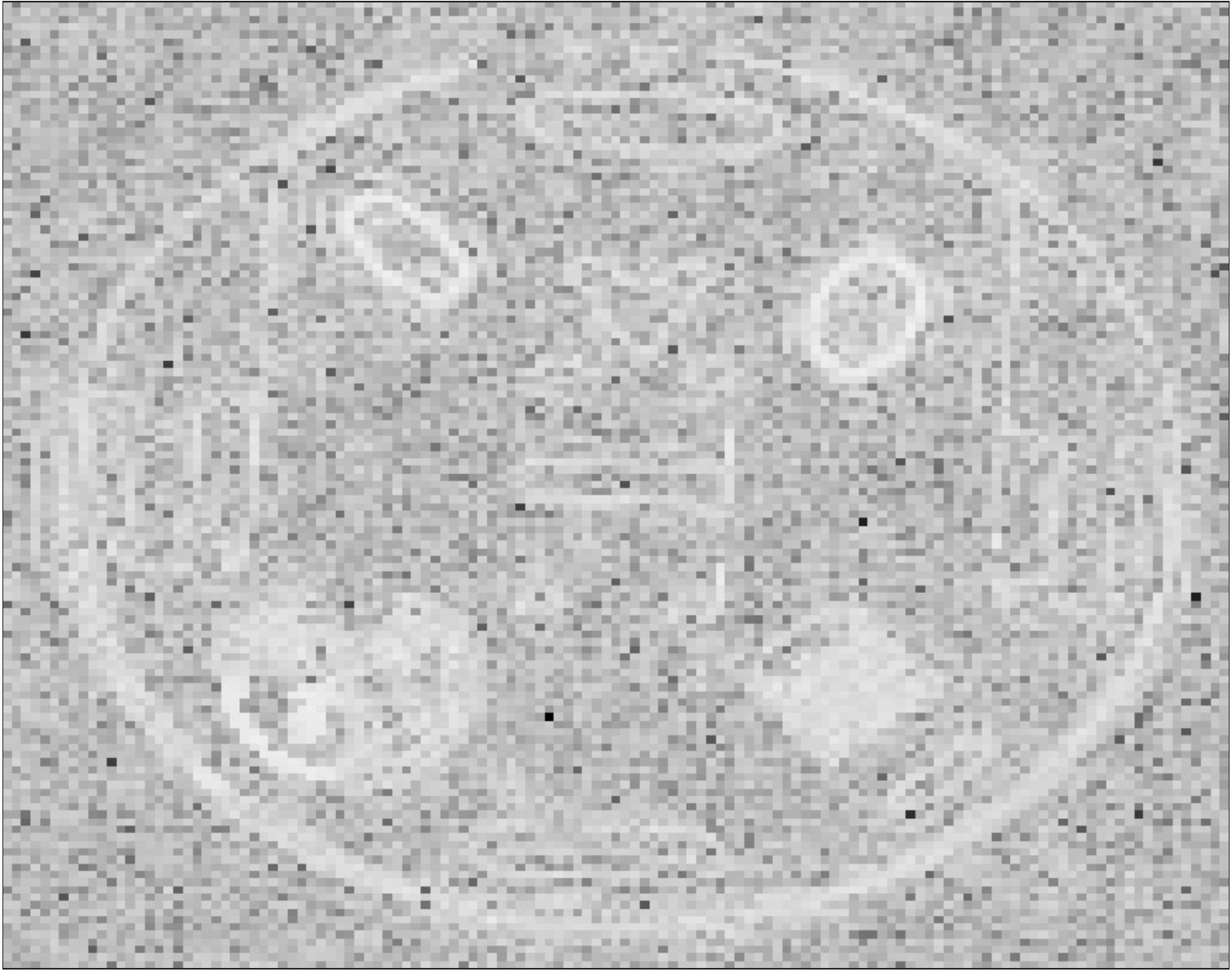}
        \caption{$\tilde f_3$}
    \end{subfigure}
    ~
    \begin{subfigure}[b]{.23\textwidth}
        \includegraphics[width=\textwidth]{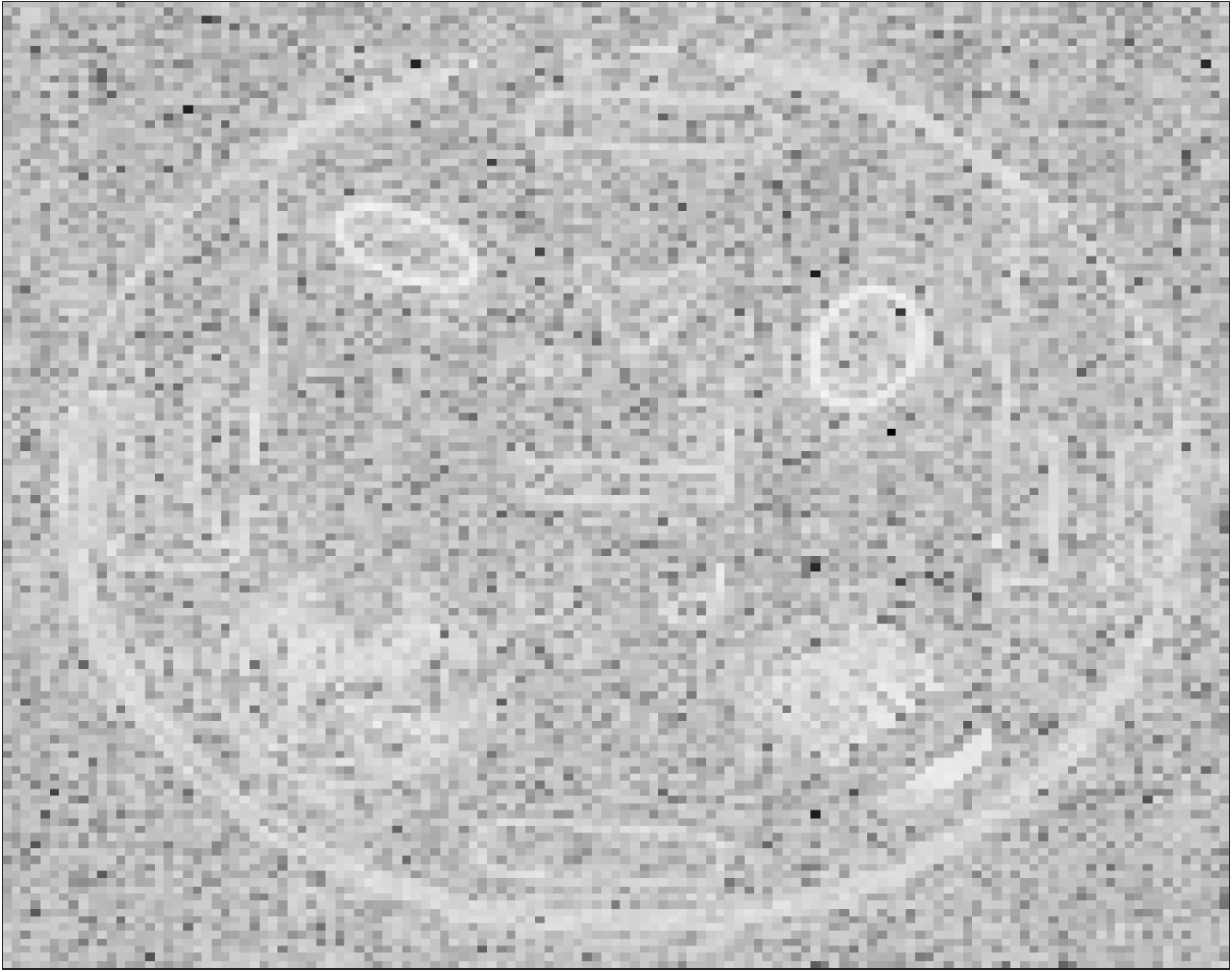}
        \caption{$\tilde f_4$}
    \end{subfigure}
    \\ 
    \begin{subfigure}[b]{.23\textwidth}
        \includegraphics[width=\textwidth]{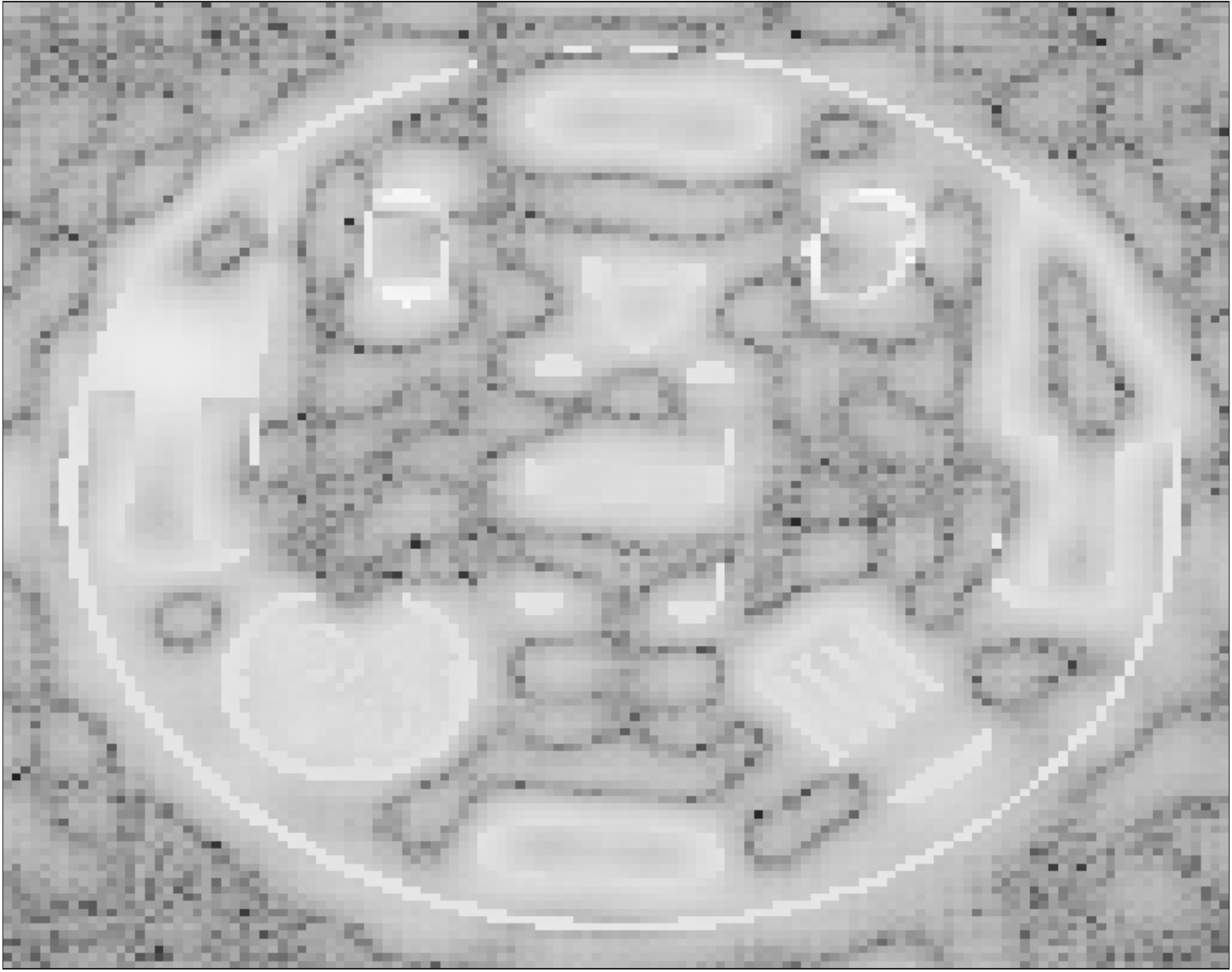}
        \caption{$\tilde f_1^\text{VBJS}$}
    \end{subfigure}
    ~
    \begin{subfigure}[b]{.23\textwidth}
        \includegraphics[width=\textwidth]{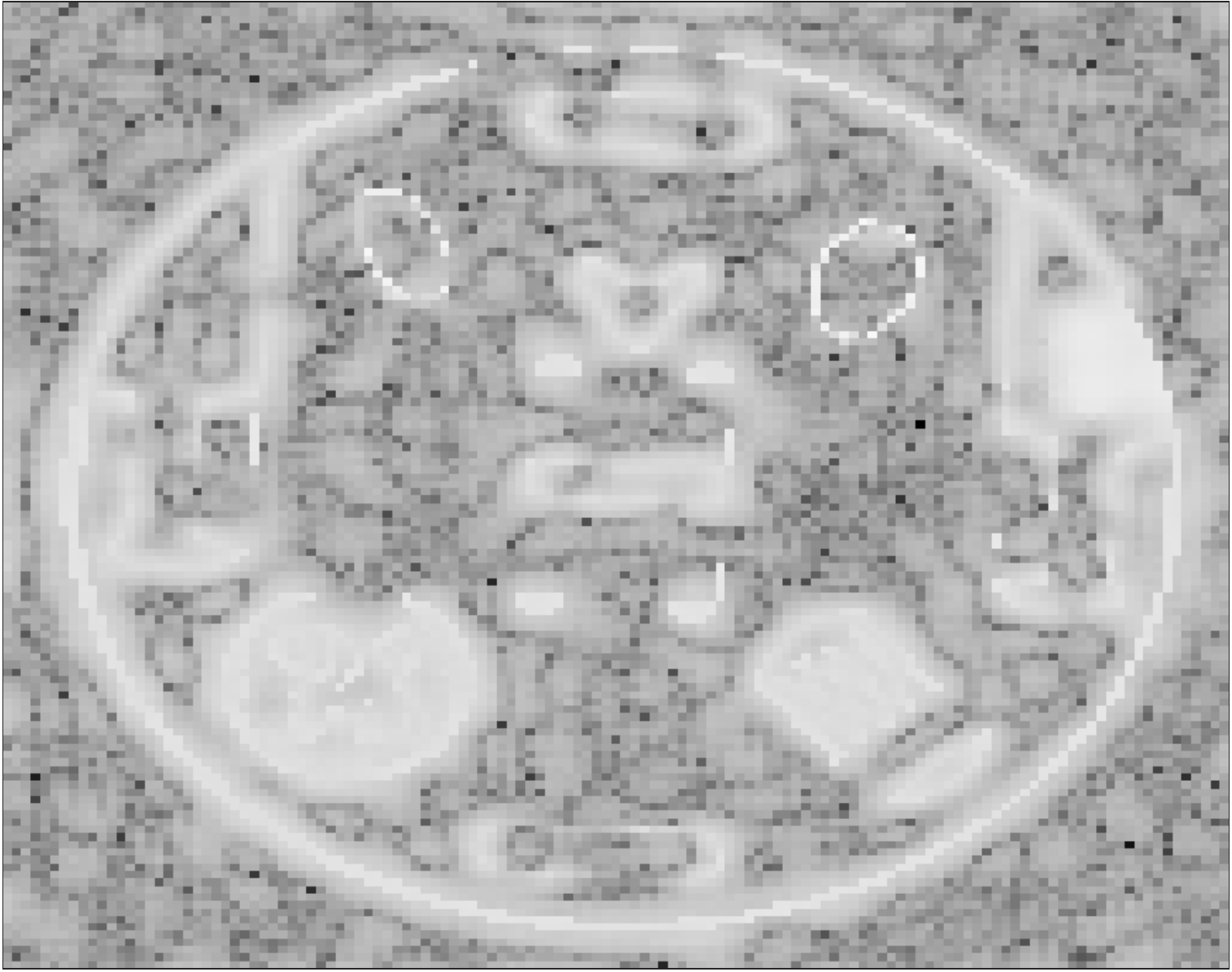}
        \caption{$\tilde f_2^\text{VBJS}$}
    \end{subfigure}
    ~
    \begin{subfigure}[b]{.23\textwidth}
        \includegraphics[width=\textwidth]{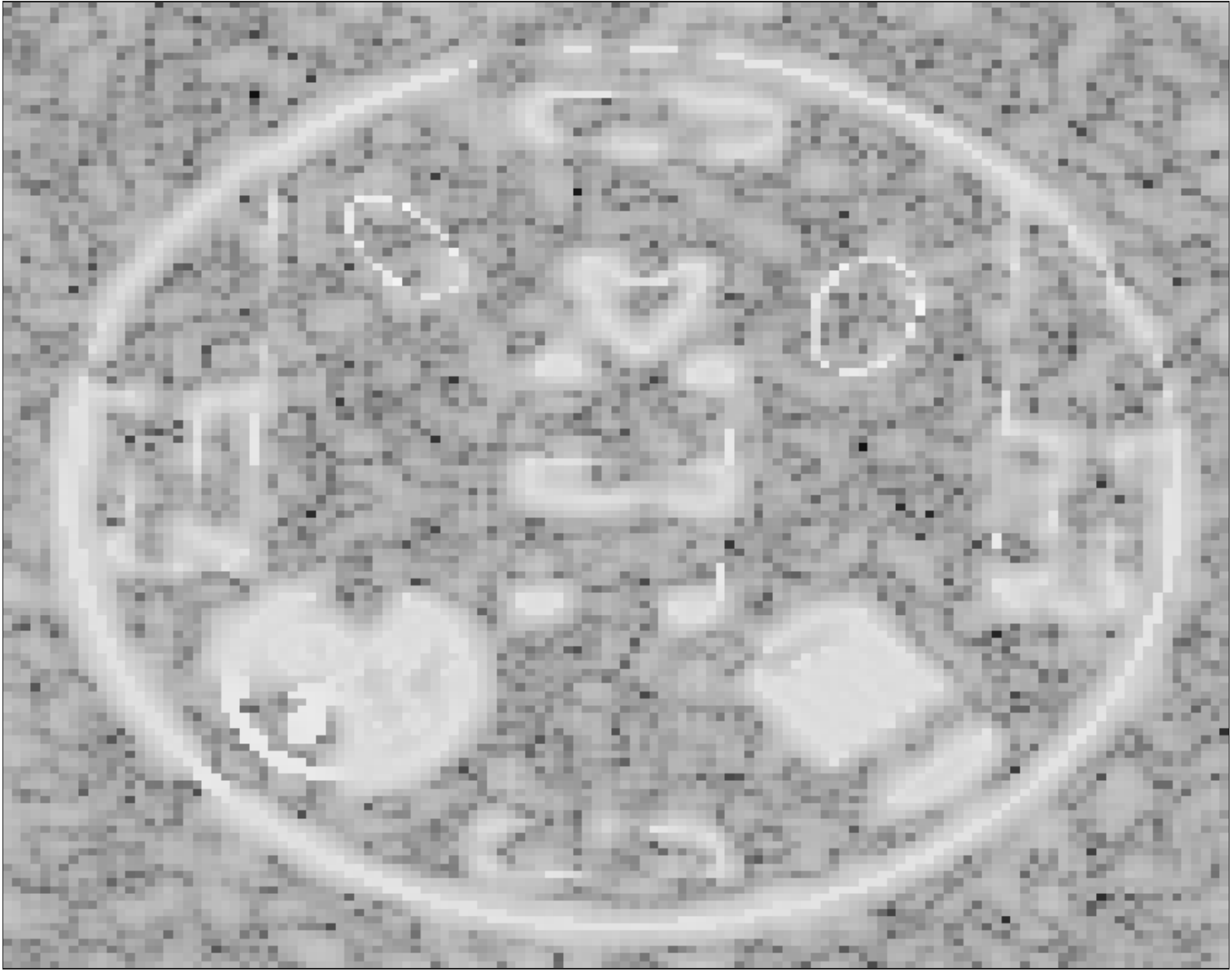}
        \caption{$\tilde f_3^\text{VBJS}$}
    \end{subfigure}
    ~
    \begin{subfigure}[b]{.23\textwidth}
        \includegraphics[width=\textwidth]{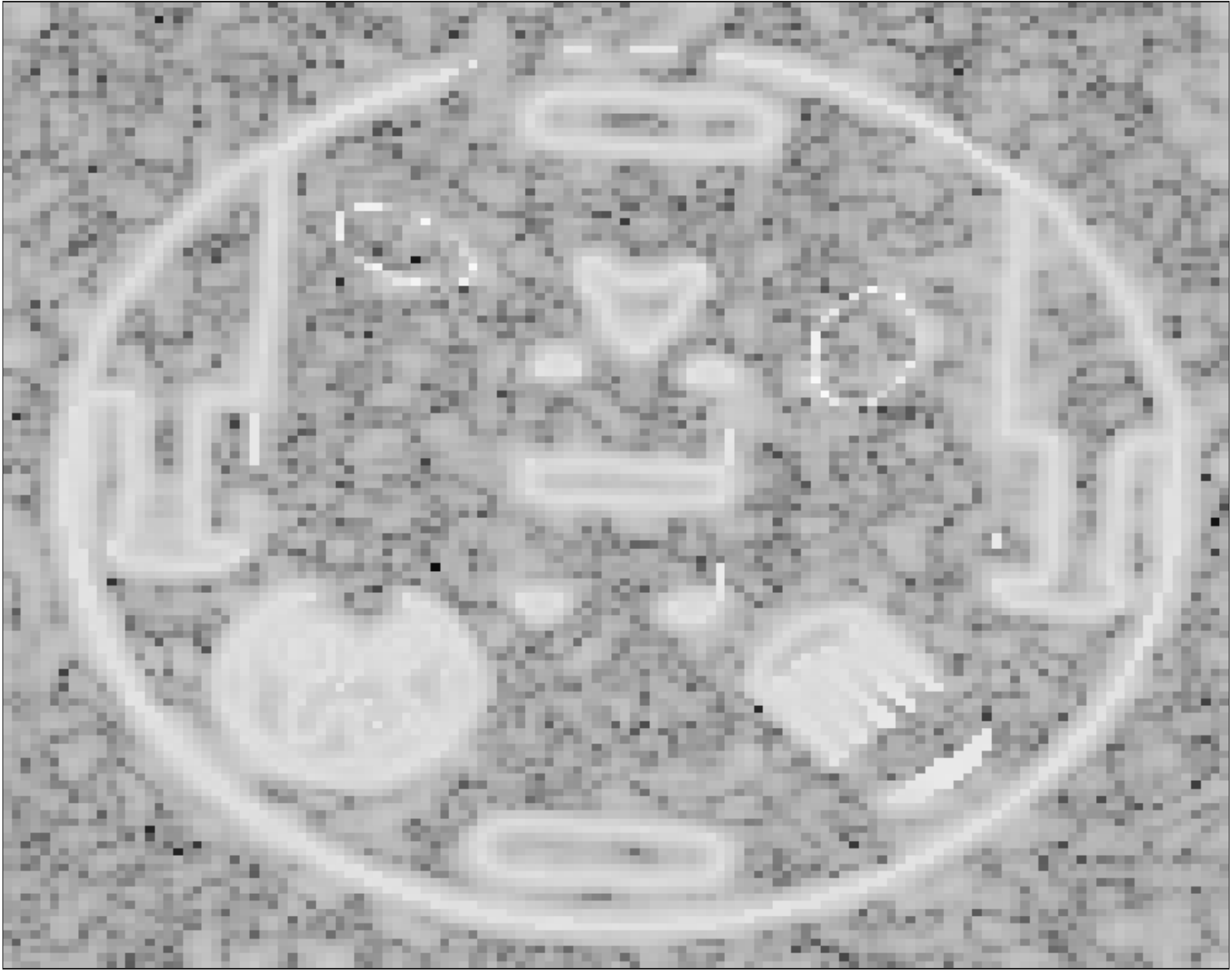}
        \caption{$\tilde f_4^\text{VBJS}$}
    \end{subfigure}
    \\
    \begin{subfigure}[b]{.23\textwidth}
        \includegraphics[width=\textwidth]{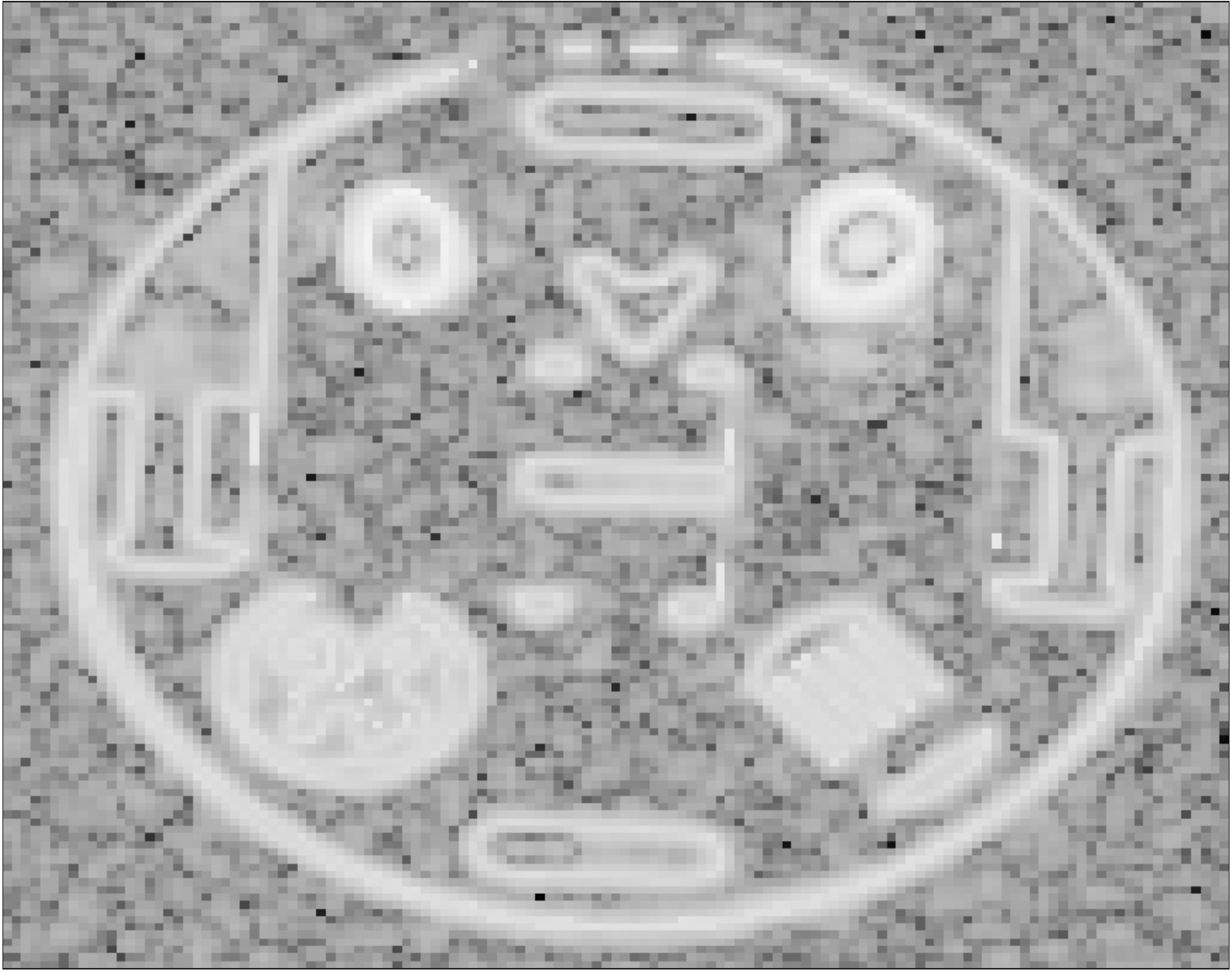}
        \caption{$\tilde f_1^\text{joint}$}
    \end{subfigure}
    ~
    \begin{subfigure}[b]{.23\textwidth}
        \includegraphics[width=\textwidth]{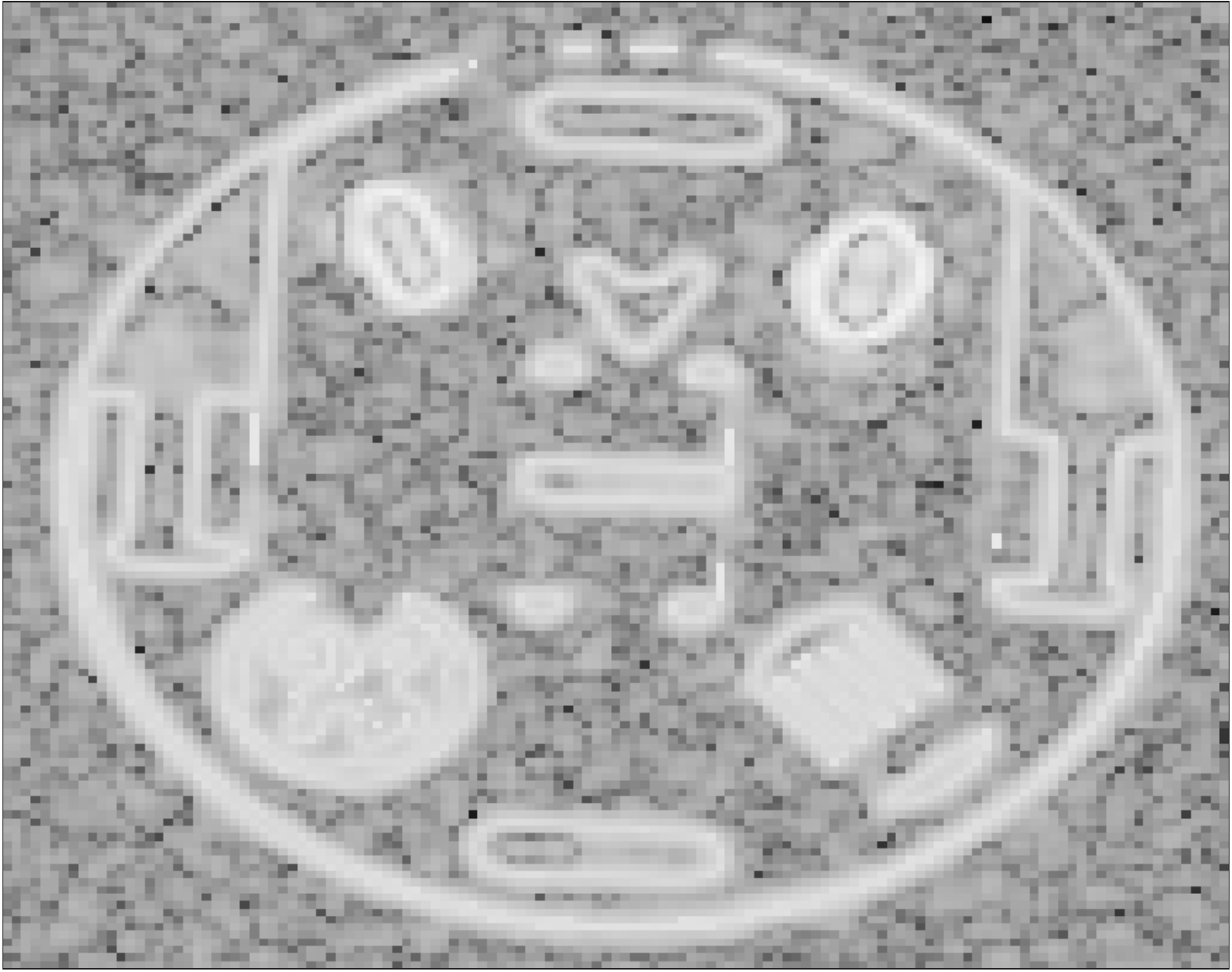}
        \caption{$\tilde f_2^\text{joint}$}
    \end{subfigure}
    ~
    \begin{subfigure}[b]{.23\textwidth}
        \includegraphics[width=\textwidth]{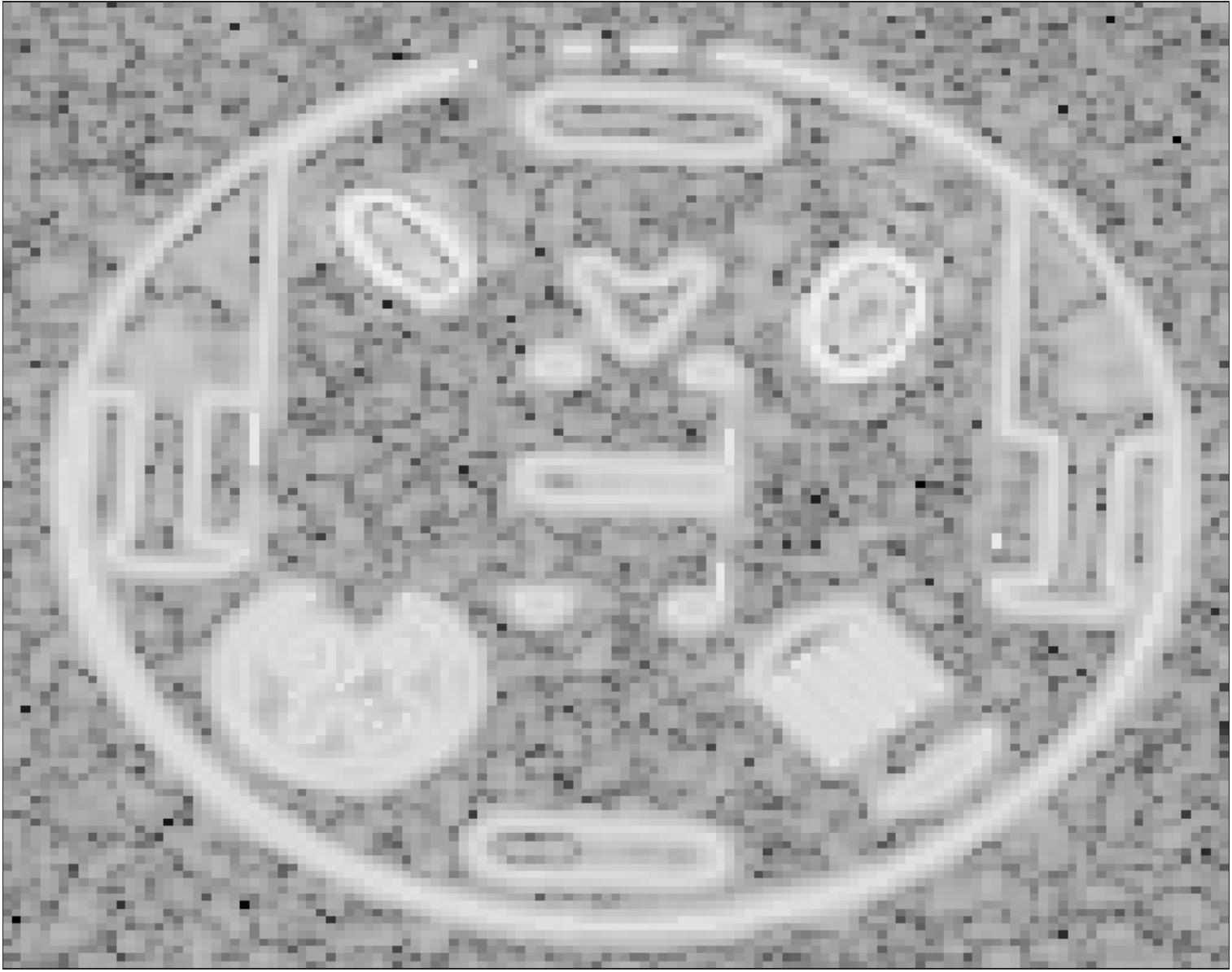}
        \caption{$\tilde f_3^\text{joint}$}
    \end{subfigure}
    ~
    \begin{subfigure}[b]{.23\textwidth}
        \includegraphics[width=\textwidth]{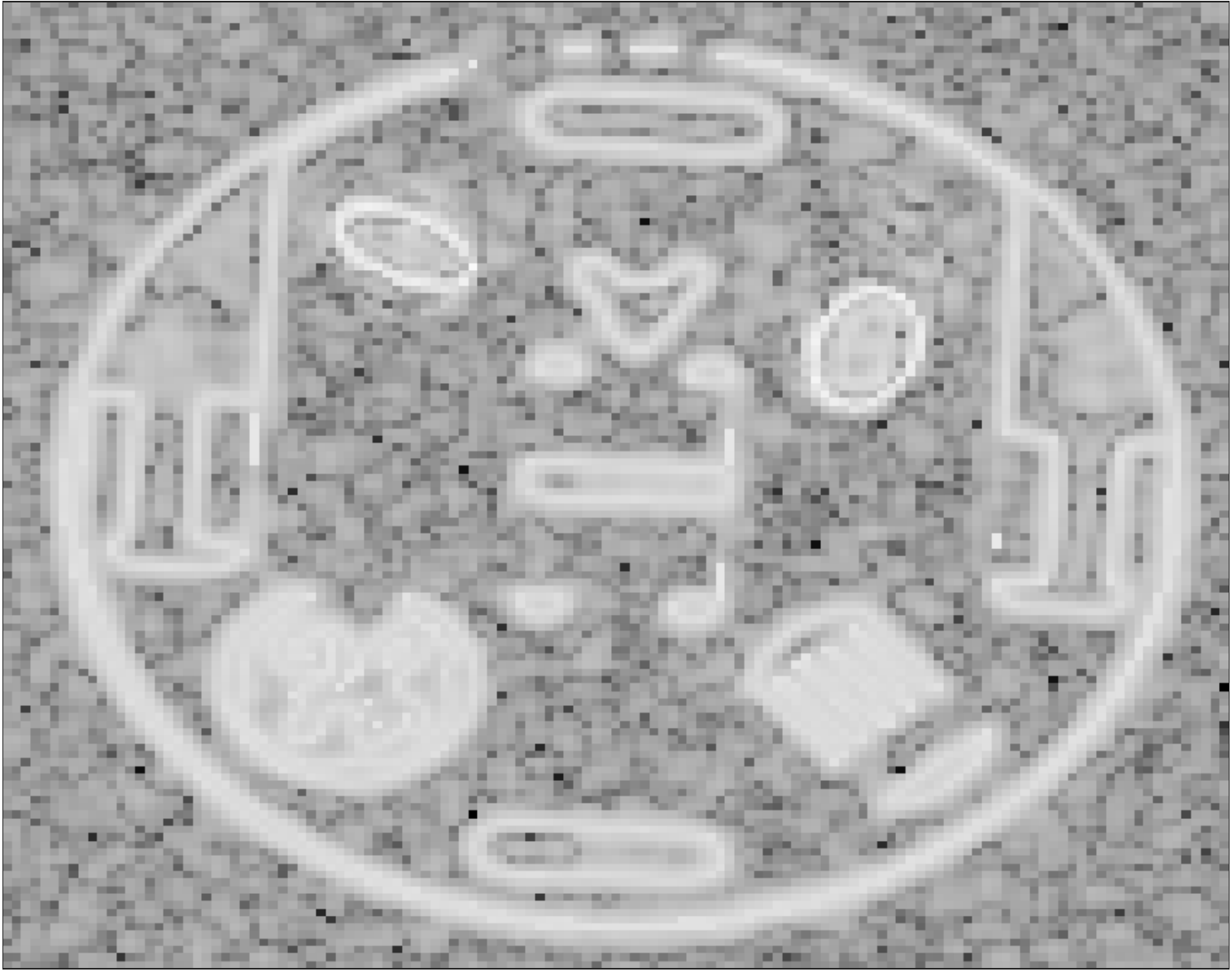}
        \caption{$\tilde f_4^\text{joint}$}
    \end{subfigure}
    \\ 
    \begin{subfigure}[b]{.23\textwidth}
        \includegraphics[width=\textwidth]{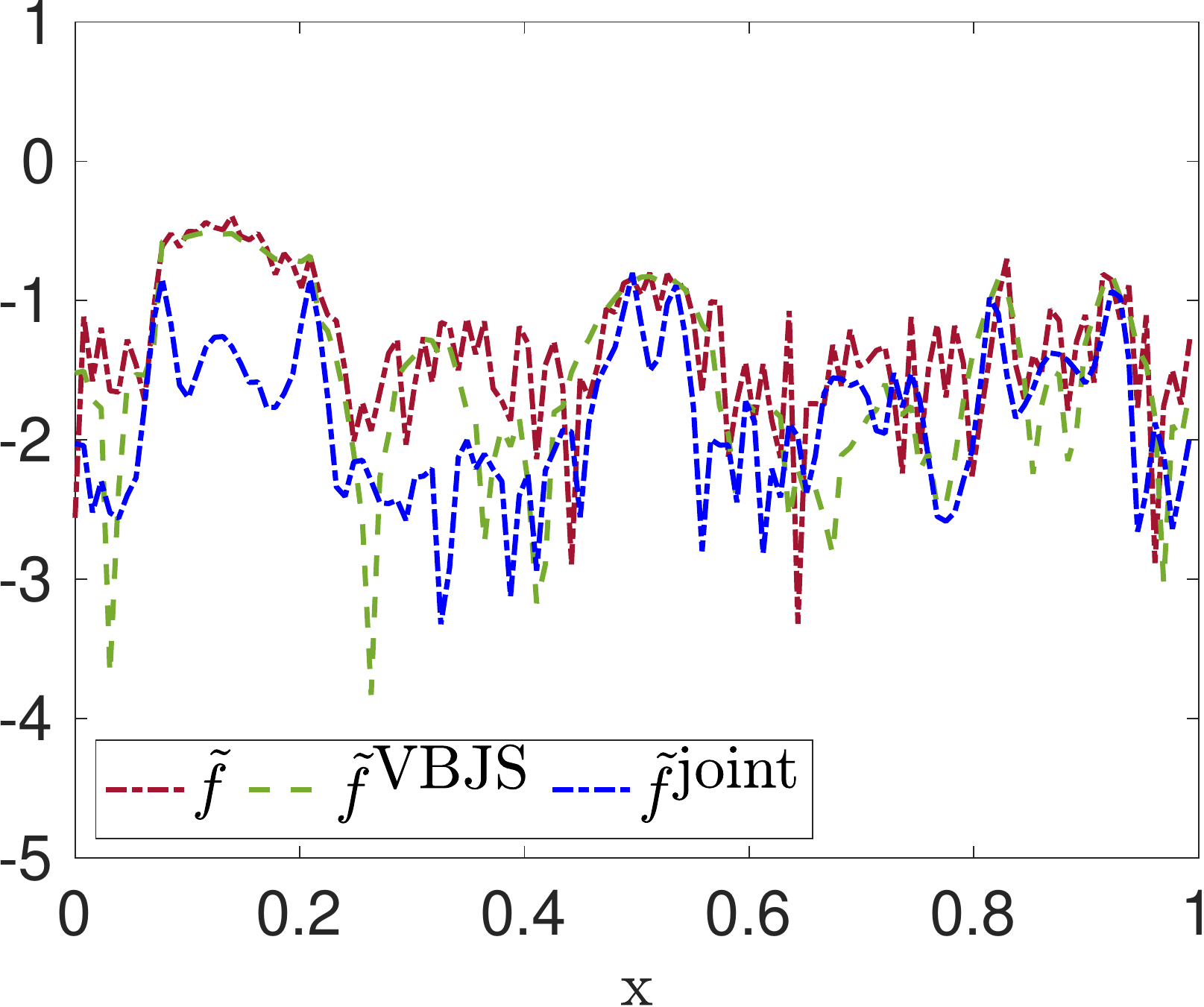}
        \caption{$\mathbf f_1$, $y=0.6512$}
    \end{subfigure}
    ~
    \begin{subfigure}[b]{.23\textwidth}
        \includegraphics[width=\textwidth]{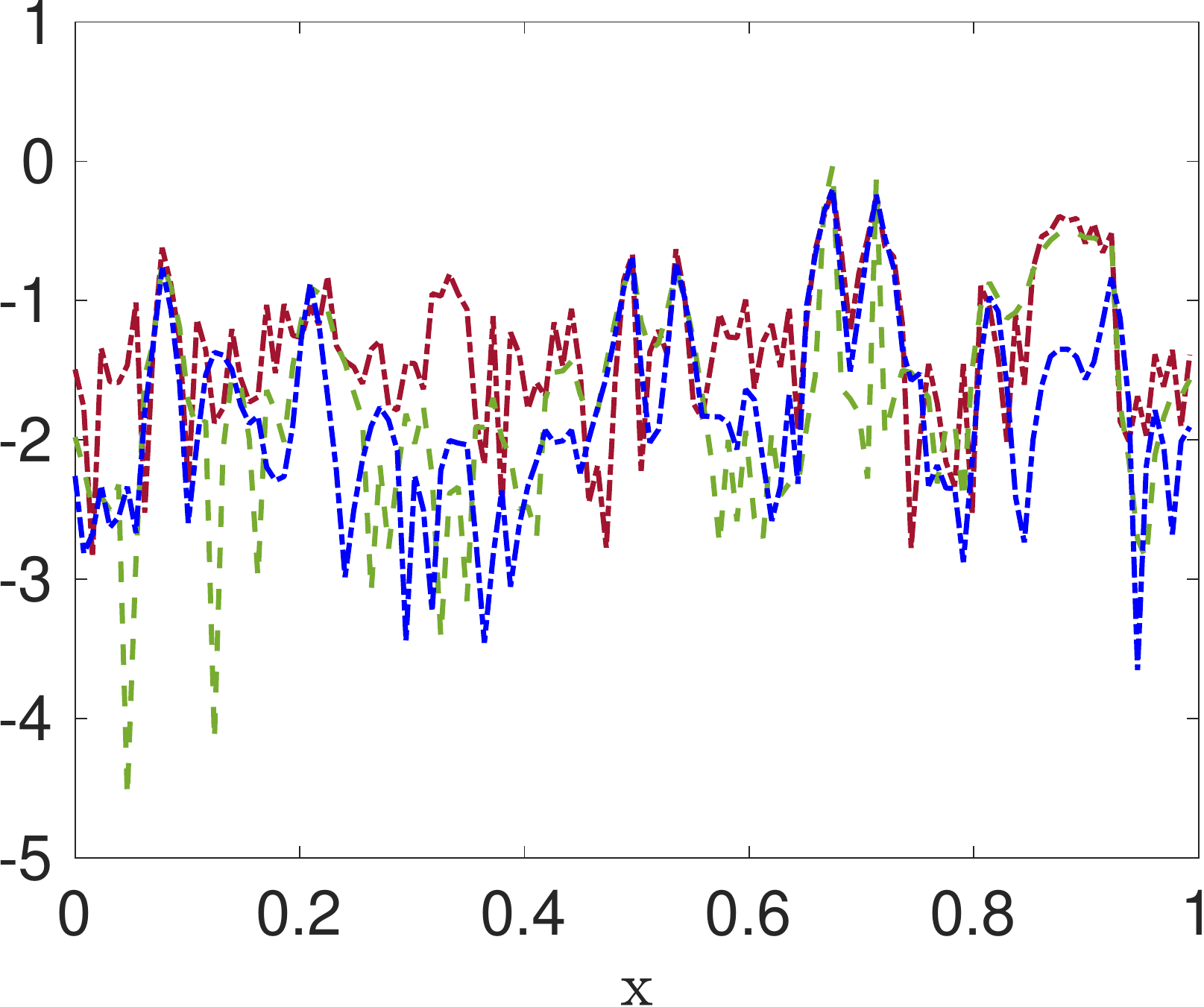}
        \caption{$\mathbf f_2$, $y=0.6512$}
    \end{subfigure}
    ~
    \begin{subfigure}[b]{.23\textwidth}
        \includegraphics[width=\textwidth]{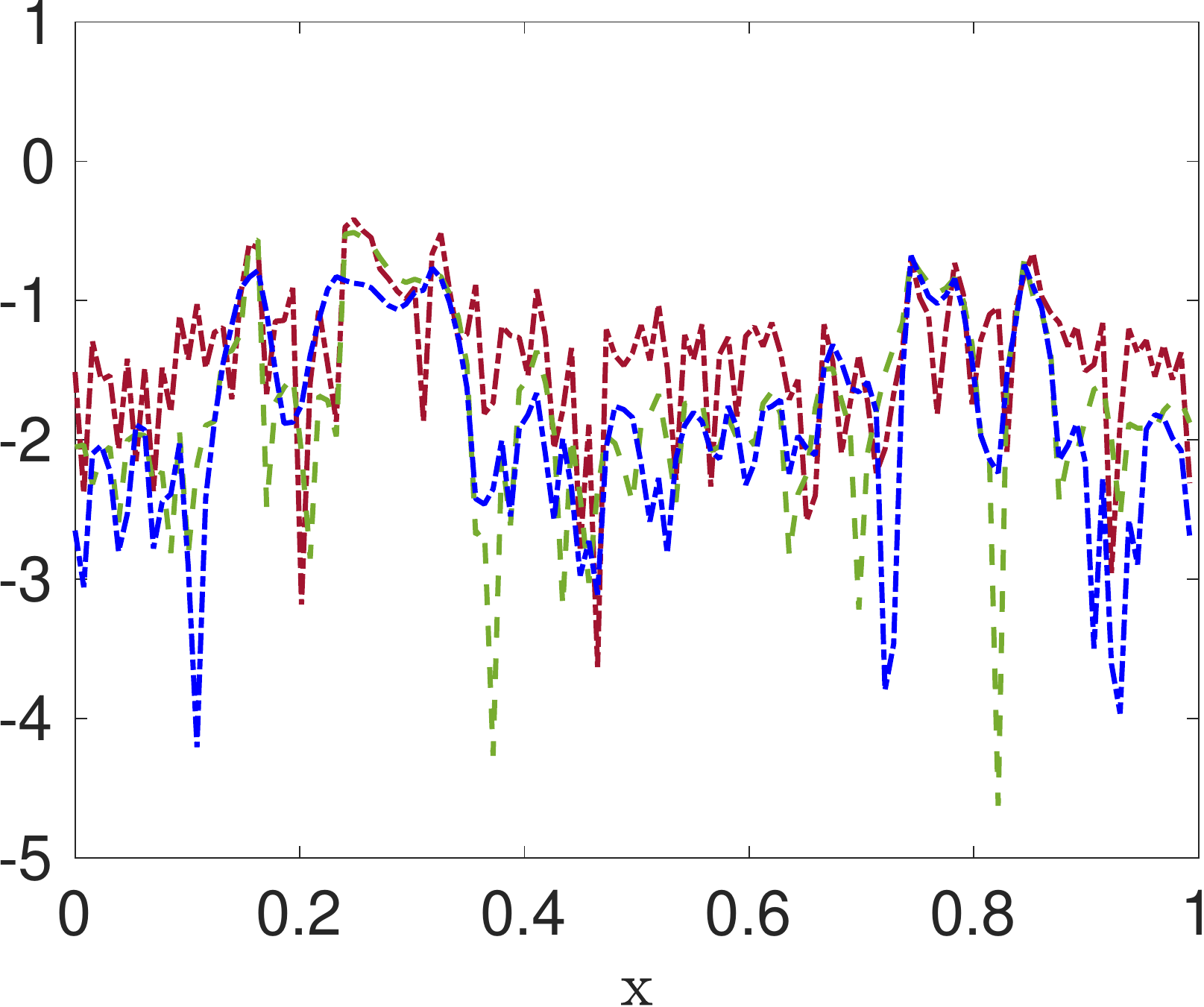}
        \caption{$\mathbf f_3$, $y=0.6512$}
    \end{subfigure}
    ~
    \begin{subfigure}[b]{.23\textwidth}
        \includegraphics[width=\textwidth]{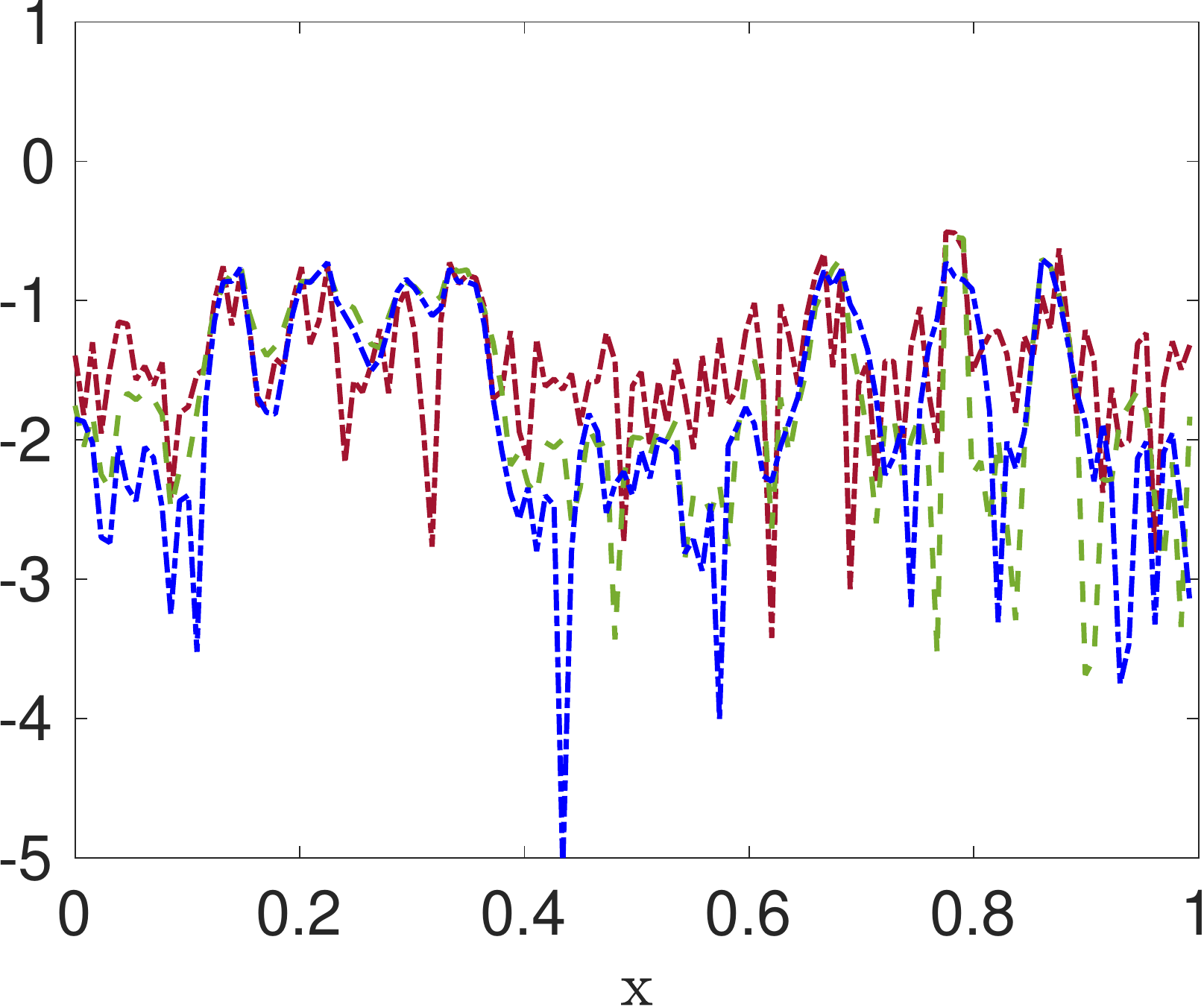}
        \caption{$\mathbf f_4$, $y=0.6512$}
    \end{subfigure}
    \caption{Error plots corresponding to image recoveries in Figure \ref{fig:rec_GE}. 
The gray scale and color legend in the top and bottom rows apply to all plots shown.}
    \label{fig:error_GE}
\end{figure}

\begin{figure}[h!]
    \centering
    \begin{subfigure}[b]{.3\textwidth}
    \includegraphics[width=\textwidth]{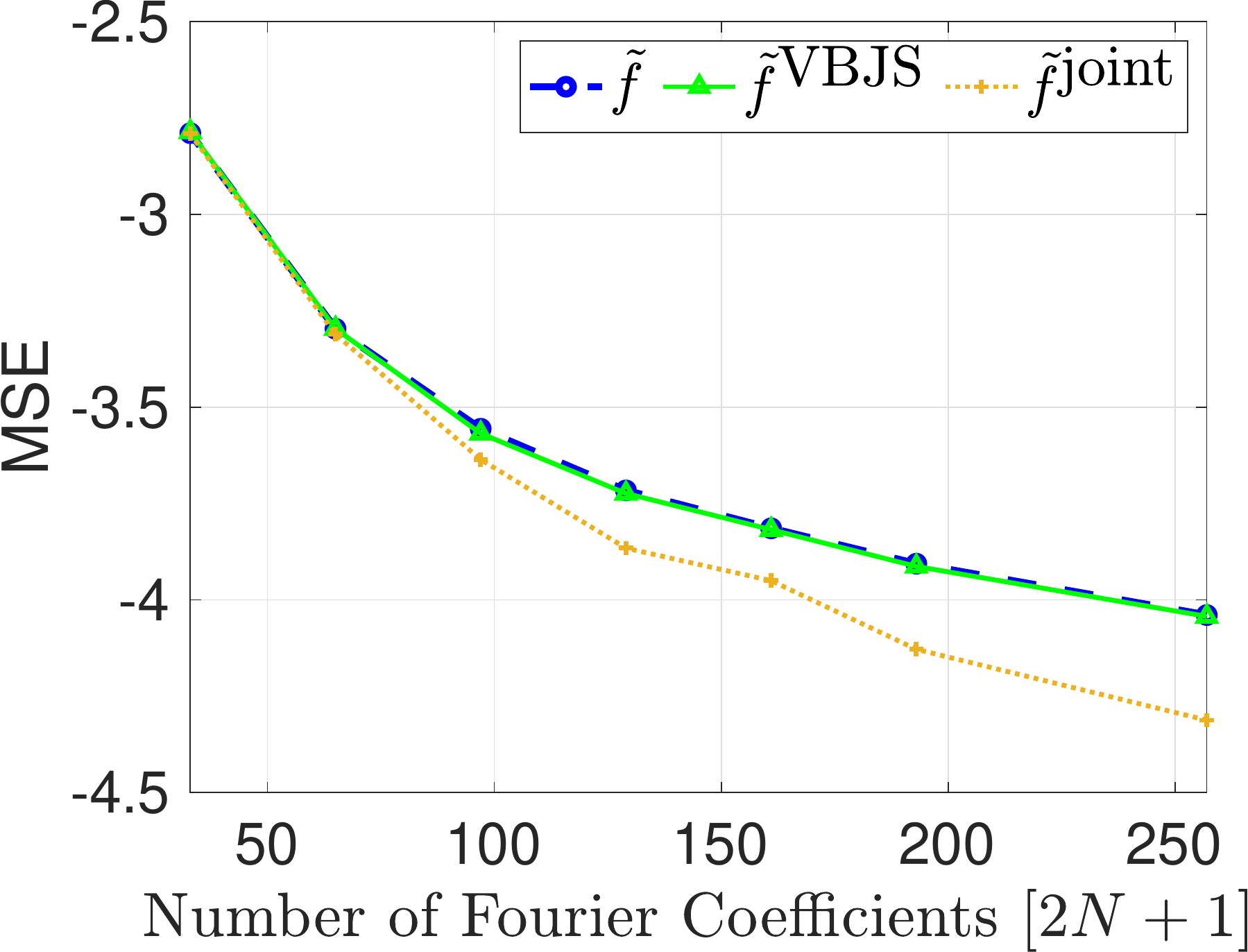}
    \caption{Whole image}
    \end{subfigure}
    \begin{subfigure}[b]{.3\textwidth}
    \includegraphics[width=\textwidth]{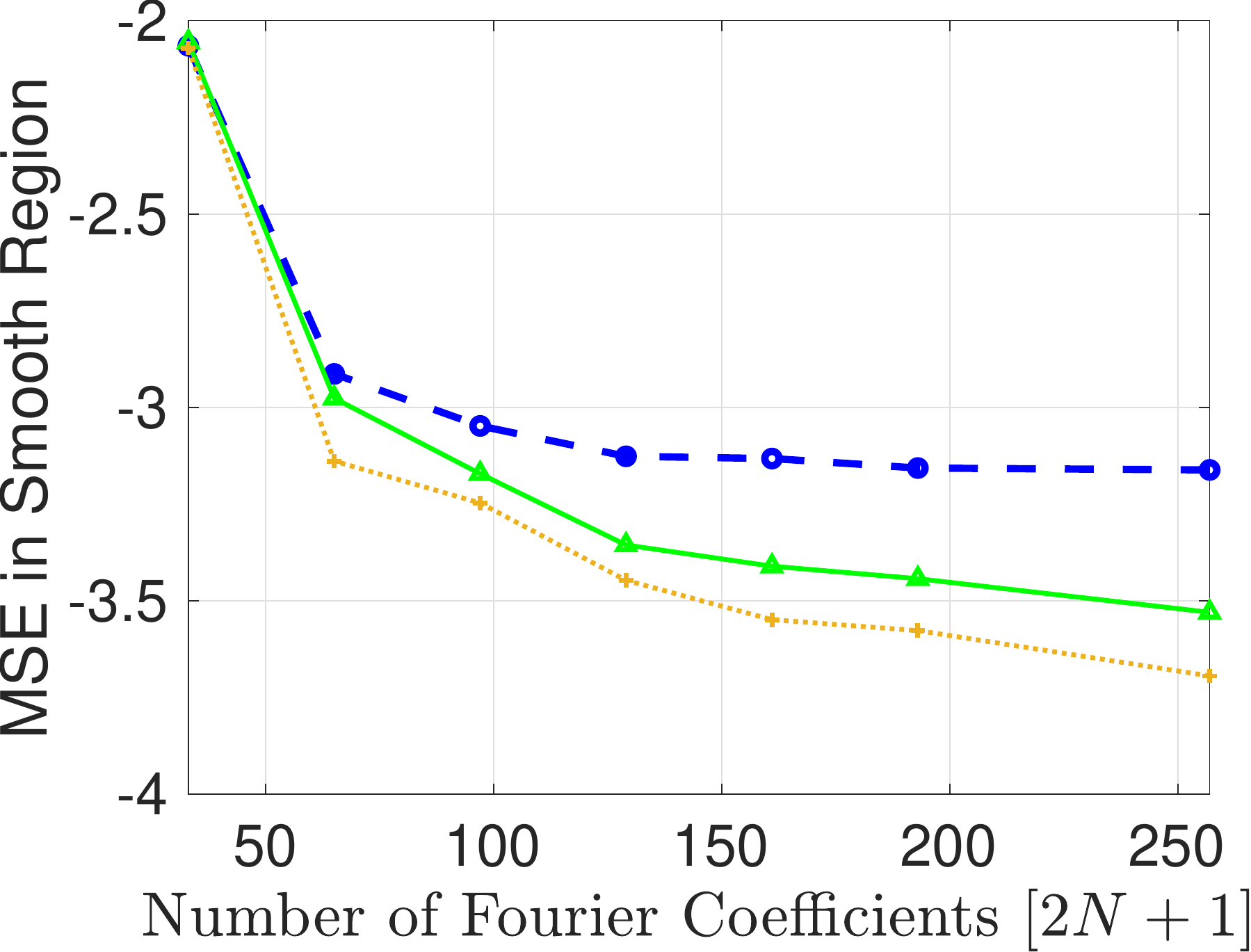}
    \caption{Smooth region}
    \end{subfigure}
    \begin{subfigure}[b]{.3\textwidth}
    \includegraphics[width=\textwidth]{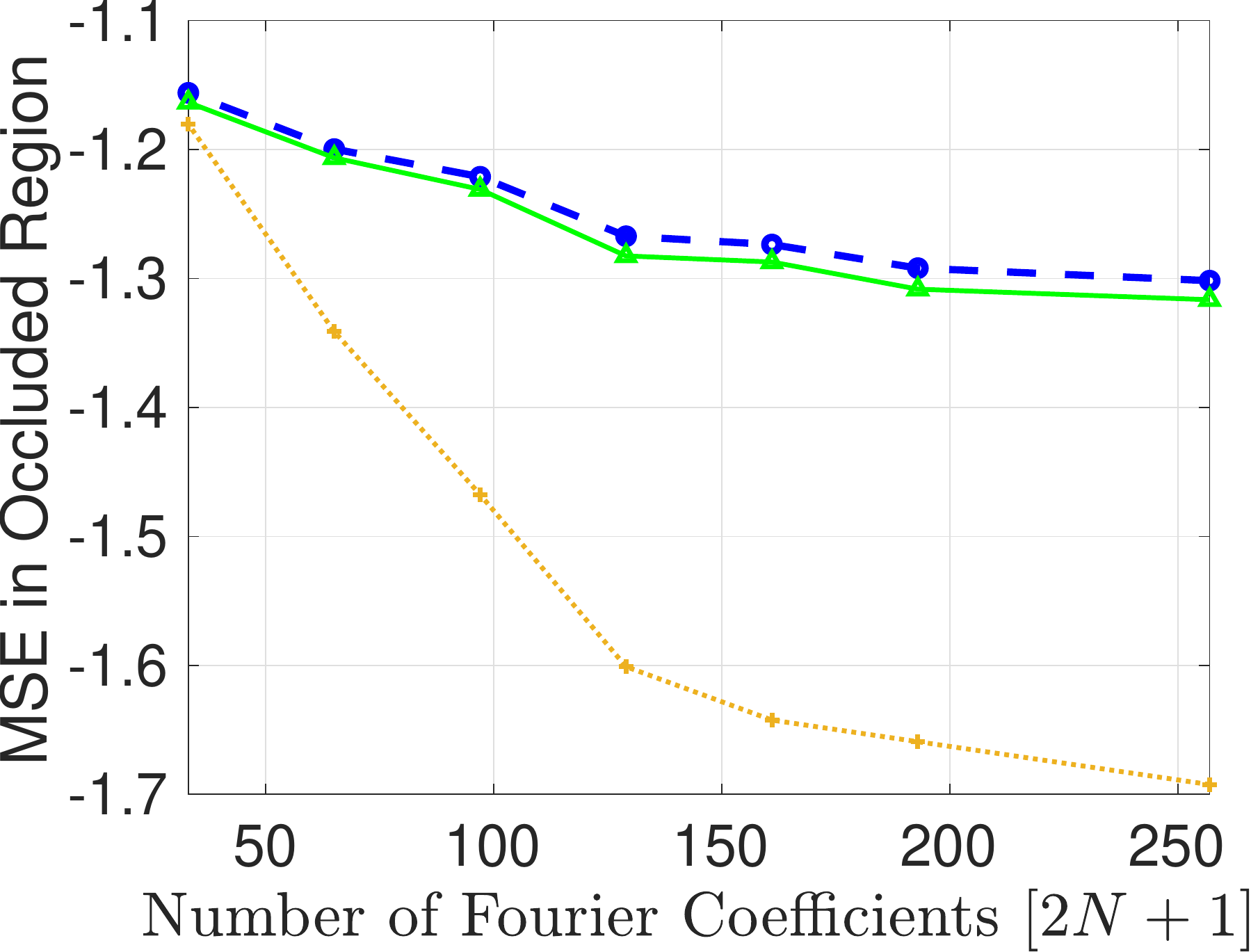}
    \caption{Occluded region}
    \end{subfigure}
    \caption{Comparison of the error in \eqref{eq:MSE} for the standard $\ell_1$ recovery (blue), VBJS (green), and joint recovery (orange) corresponding with increasing the number of Fourier samples, $2N+1$, where $N=16n$ with $n=1,\dots,7$. Error taken over (left) the entire image; (middle) the smooth regions; and (right) regions containing occlusions.}
    \label{fig:log_err}
\end{figure}

\begin{figure}[h!]
    \centering
    \begin{subfigure}[b]{.3\textwidth}
    \includegraphics[width=\textwidth]{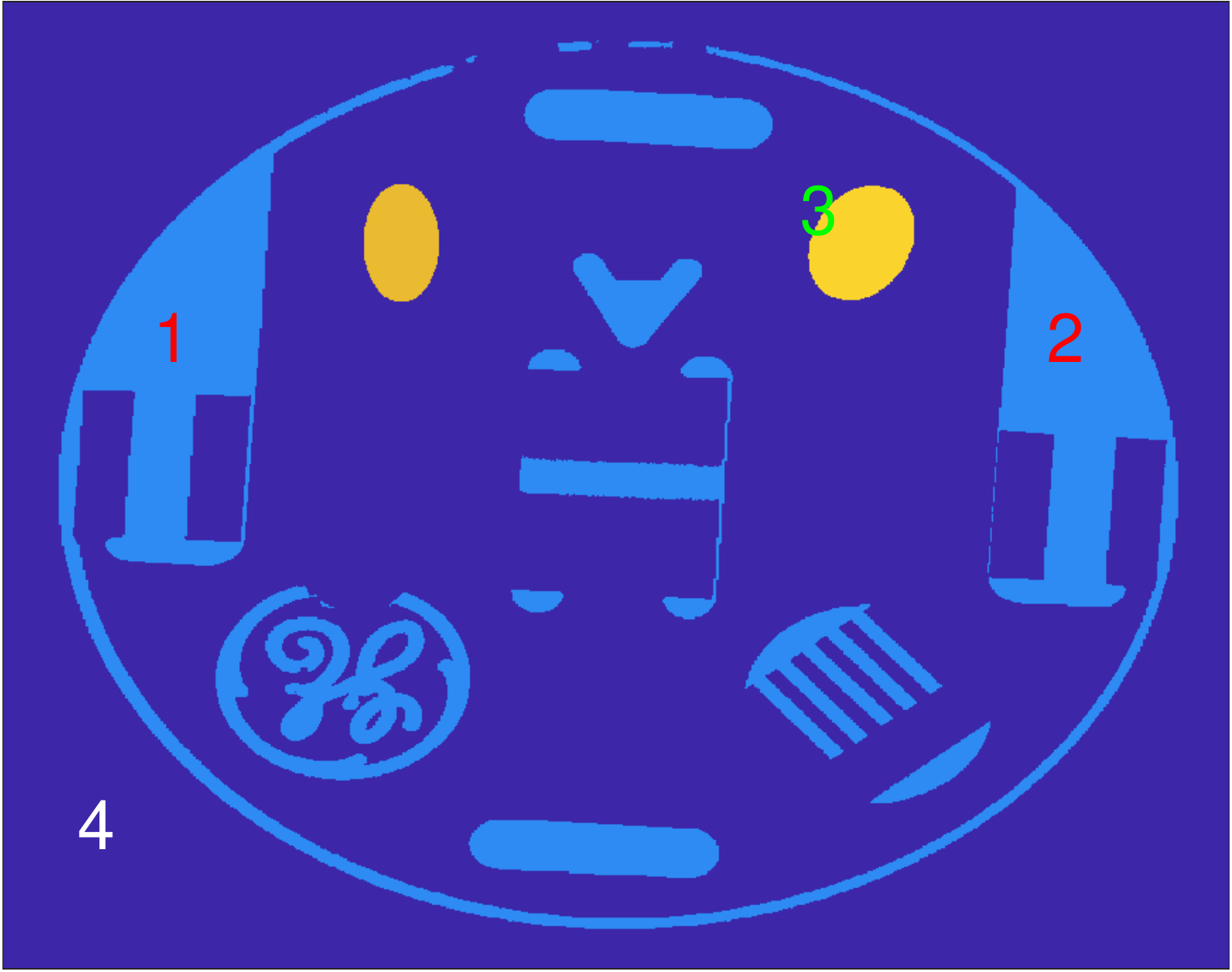}
    \caption{True image $\mathbf f_1$}
    \end{subfigure}
    \begin{subfigure}[b]{.3\textwidth}
    \includegraphics[width=\textwidth]{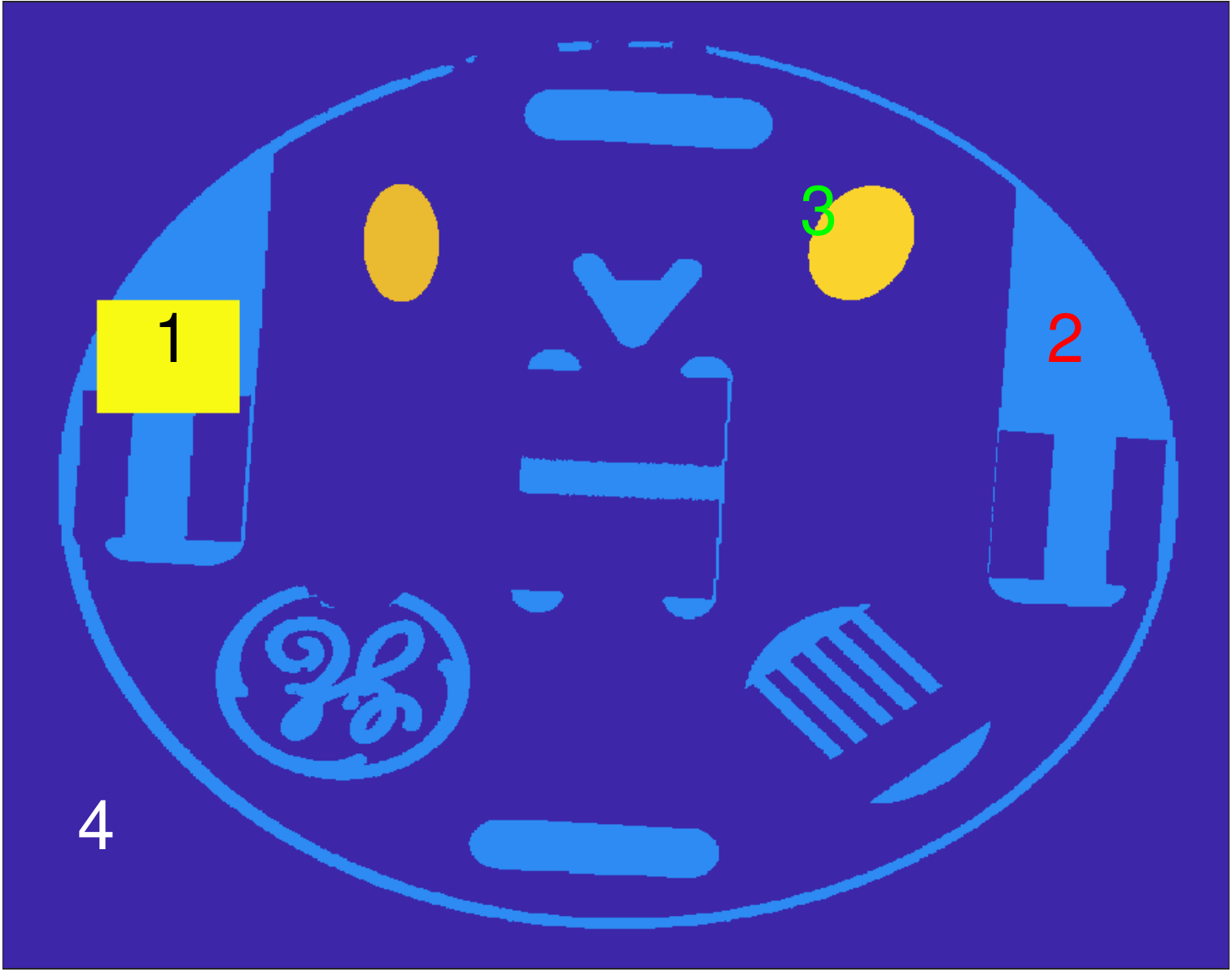}
    \caption{$\mathbf f_1$ with occlusion}
    \end{subfigure}
    \begin{subfigure}[b]{.3\textwidth}
    \includegraphics[width=\textwidth]{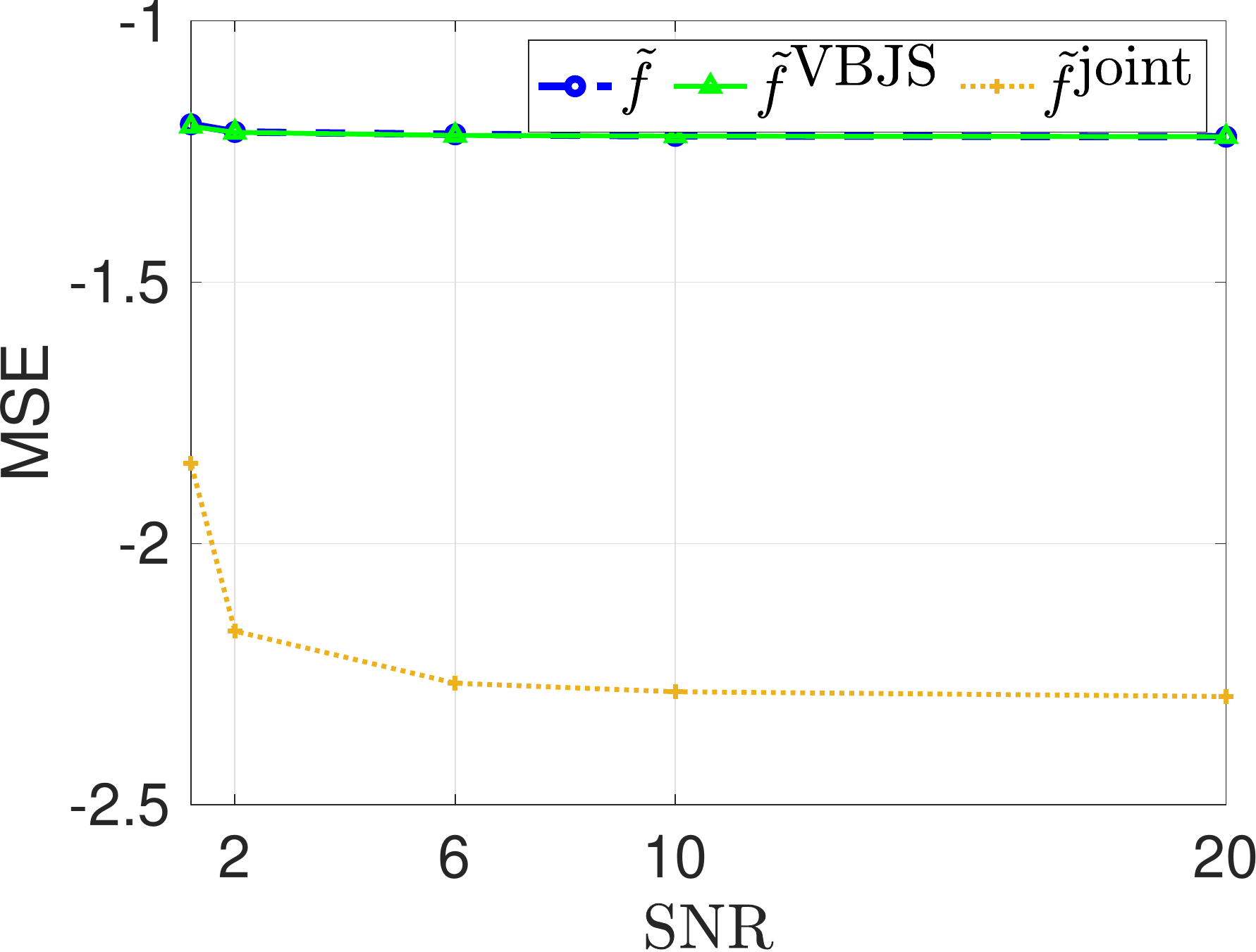}
    \caption{Point 1}
    \end{subfigure}
    \\
    \begin{subfigure}[b]{.3\textwidth}
    \includegraphics[width=\textwidth]{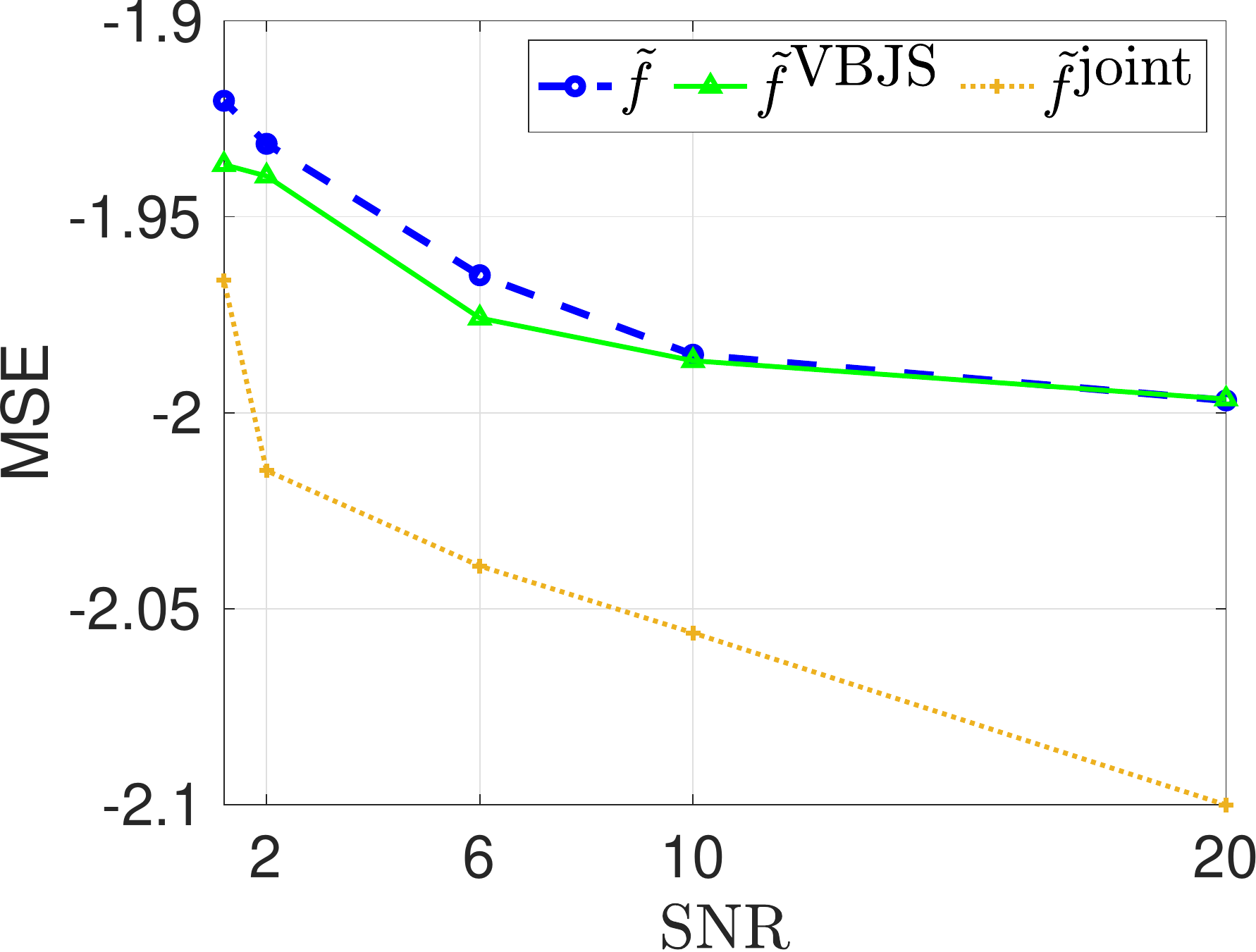}
    \caption{Point 2}
    \end{subfigure}
    \begin{subfigure}[b]{.3\textwidth}
    \includegraphics[width=\textwidth]{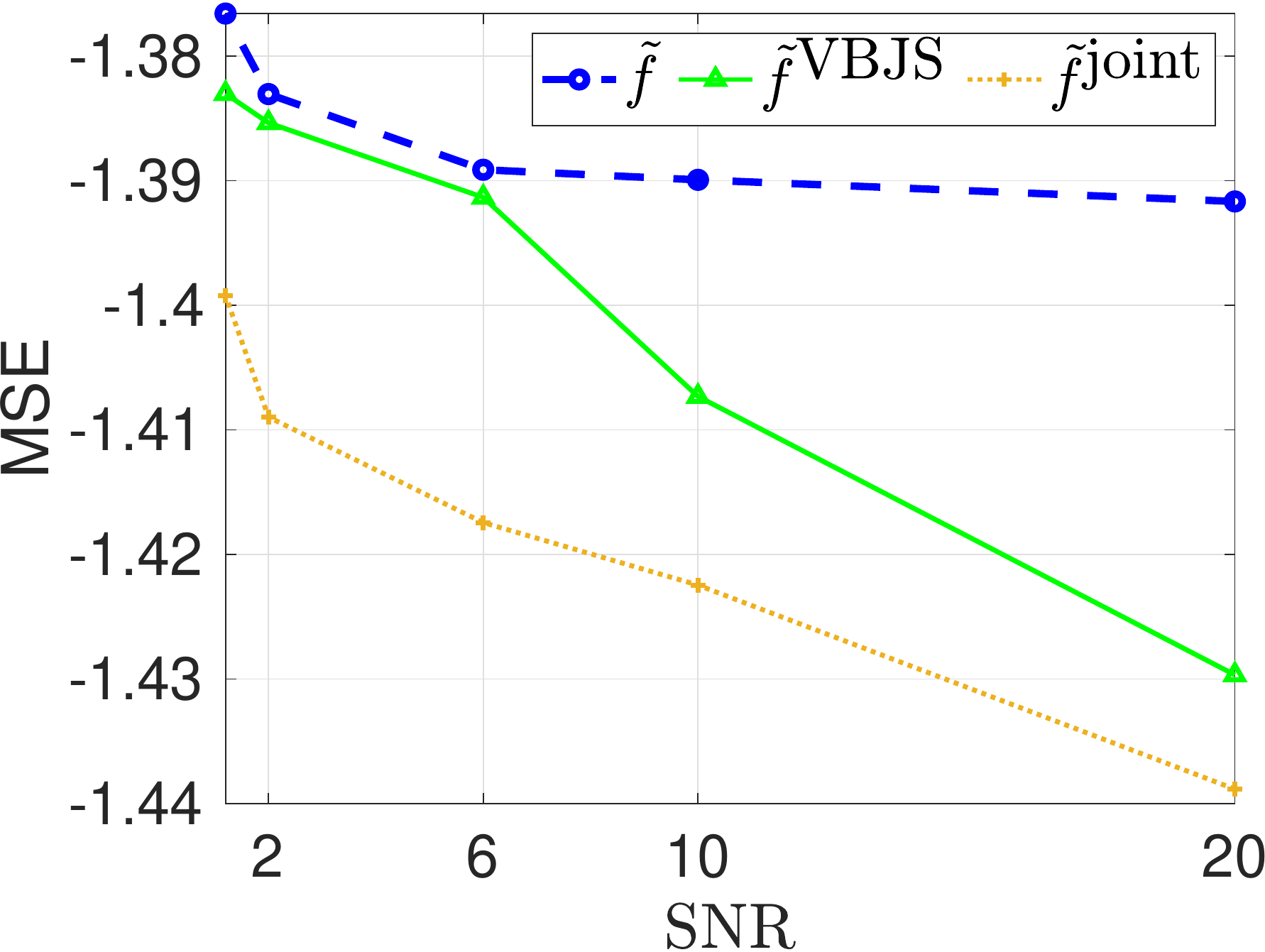}
    \caption{Point 3}
    \end{subfigure}
    \begin{subfigure}[b]{.3\textwidth}
    \includegraphics[width=\textwidth]{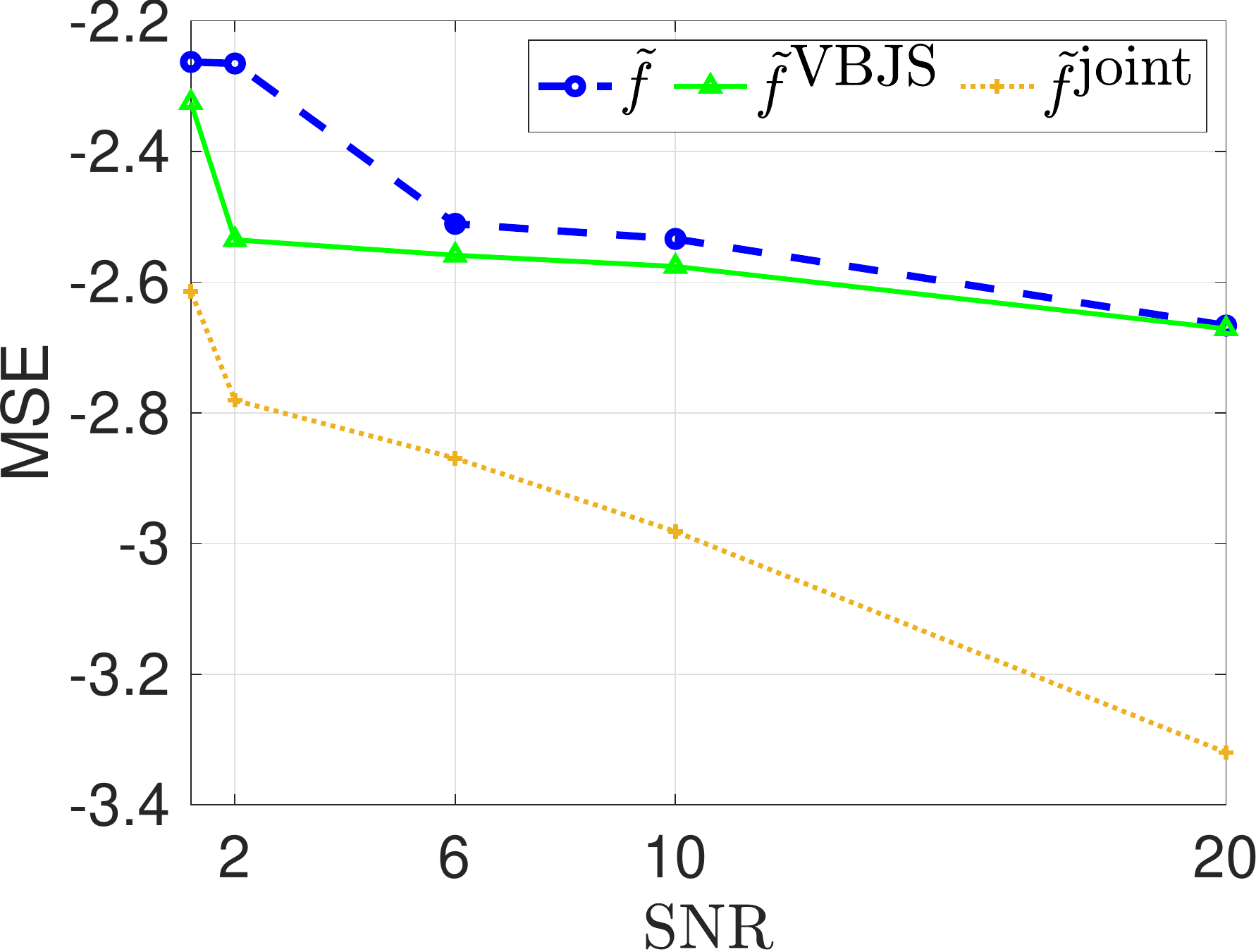}
    \caption{Point 4}
    \end{subfigure}
    
    \caption{Comparison of the error in \eqref{eq:MSE} for the standar $\ell_1$ recovery (blue), VBJS (green), and joint recovery (orange) for increasing SNR in the data (Fourier) domain with a fixed number of Fourier samples, $2N+1=129$. The points of interest, labeled as ``1'' through ``4'', in the first two plots, correspond to different aspects of image recovery.}
    \label{fig:log_err_snr}
\end{figure}
To compare the numerical convergence of the three tested methods, we consider the log-scale mean squared error (MSE)
\begin{equation}
\label{eq:MSE}
    {MSE}_{log} = \log_{10}\left( \frac{ \|\mathbf{f}-\tilde{\mathbf f} \|^2 }{ \text{size}(\mathbf{f}) }\right),
\end{equation}
where ${\bf f}$ is the pixelated underlying image and $\tilde{\mathbf f}$ is the recovered image. {Here $\text{size}{(\cdot)}$ denotes the number of entries considered in \eqref{eq:MSE}.   We are particularly interested in evaluating the error in background regions that are object free as well as in regions containing occlusions.}  We also measure the log of the pointwise error given by
\begin{equation}\label{eq:pw_log_err}
    {E}_{\log} = \log_{10} \mid\mathbf{f}-\tilde{ {\bf f}}\mid.
\end{equation}

The log error of the entire image $\mathbf{f}_j$ for $j=1,\dots, 4$ is displayed in Figure \ref{fig:error_GE}. It is evident that  only the proposed method is able to recover the obstructed regions.

Figure \ref{fig:log_err} displays the log-scale mean square error calculated in \eqref{eq:MSE} over different regions as a function of $2N+1$, the number of Fourier samples.  As shown in Figure \ref{fig:log_err}(middle), as long as information is reliable, the results for the VBJS and the joint methods are comparable in smooth regions,  and both perform better than standard $\ell_1$ regularization.  By contrast, inter-image information becomes critical once regions are occluded in any of the images in the temporal sequence, as is illustrated in Figure \ref{fig:log_err}(right).  All methods demonstrate convergence with increasing number of Fourier samples.
Next we fix the number of Fourier samples, $2N+1 = 129$, and analyze the error for varying SNR values, $1.2, 2, 6, 10$ and $20$, in the  data (Fourier) domain.   The first two plots in Figure \ref{fig:log_err_snr}  show four points of evaluation:
The  ``1'' is located inside an obstacle, the ``2'' is inside a nonzero smooth region, the ``3'' is chosen to be close to an edge, and the ``4'' is placed in the zero-valued background.  The MSE is then computed on $5\times 5$ neighborhoods centered at the respective sample points.  Once again we observe that the new joint recovery approach yields better accuracy in all cases for all levels of SNR, and is particularly effective in regions of occlusion.

\begin{figure}[h!]
    \centering
    \begin{subfigure}[b]{.19\textwidth}
    \includegraphics[width=\textwidth]{plot/true_f2_GE}
    \caption{$\mathbf f_2$}
    \end{subfigure}
    \begin{subfigure}[b]{.19\textwidth}
    \includegraphics[width=\textwidth]{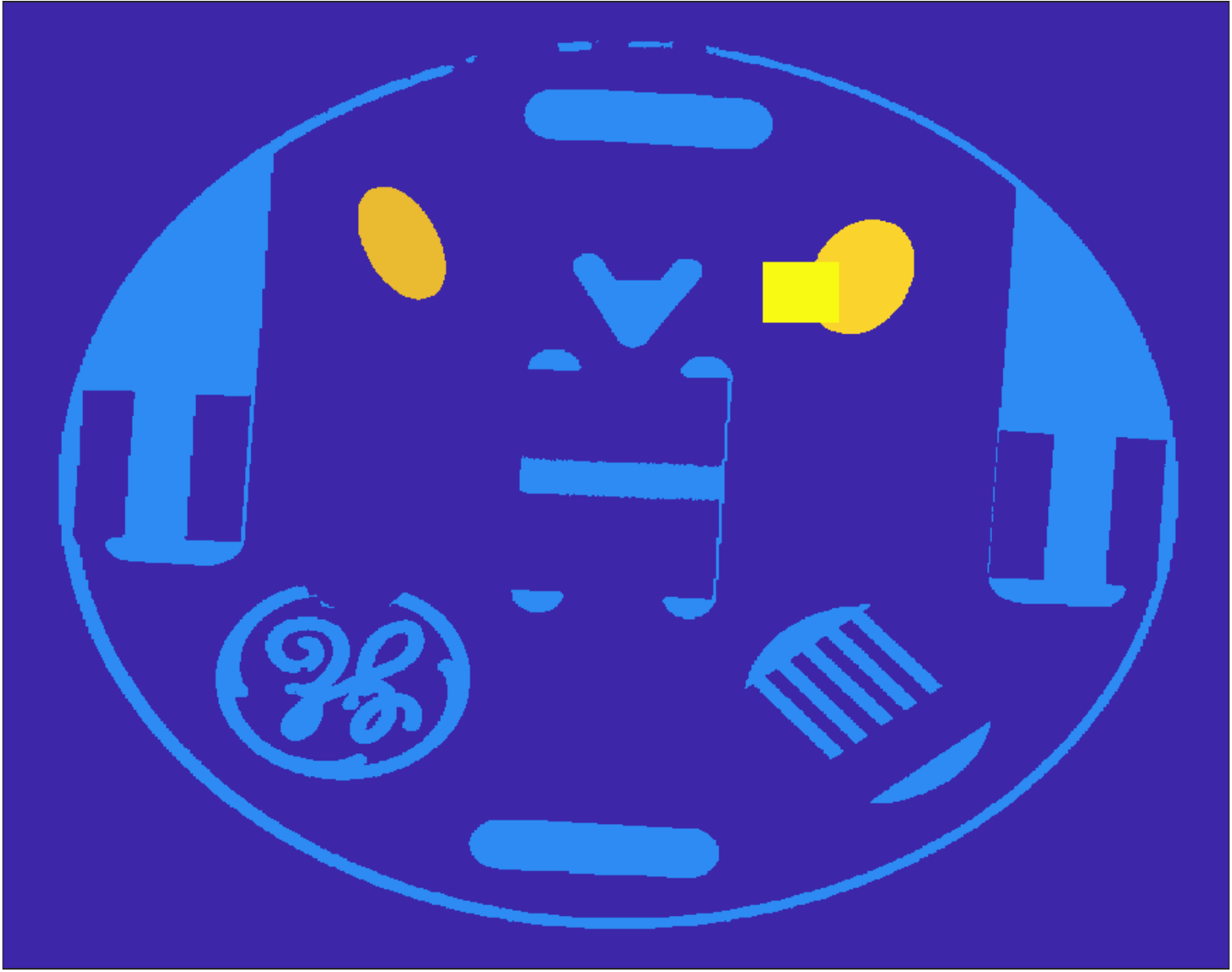}
    \caption{Obstacle$_2$}
    \end{subfigure}
    \begin{subfigure}[b]{.19\textwidth}
    \includegraphics[width=\textwidth]{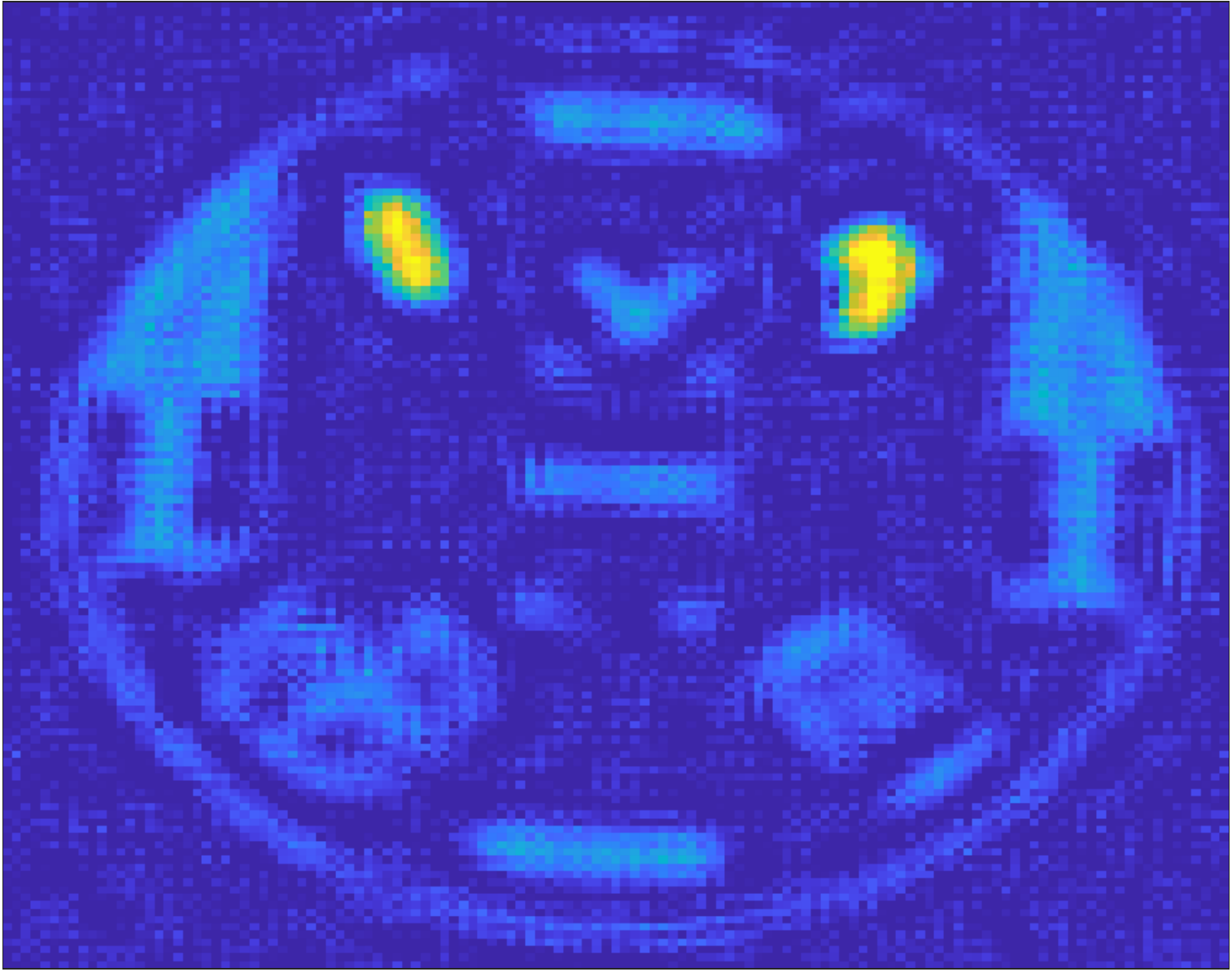}
    \caption{$\tilde{\mathbf{f}}_2$}
    \end{subfigure}
    \begin{subfigure}[b]{.19\textwidth}
    \includegraphics[width=\textwidth]{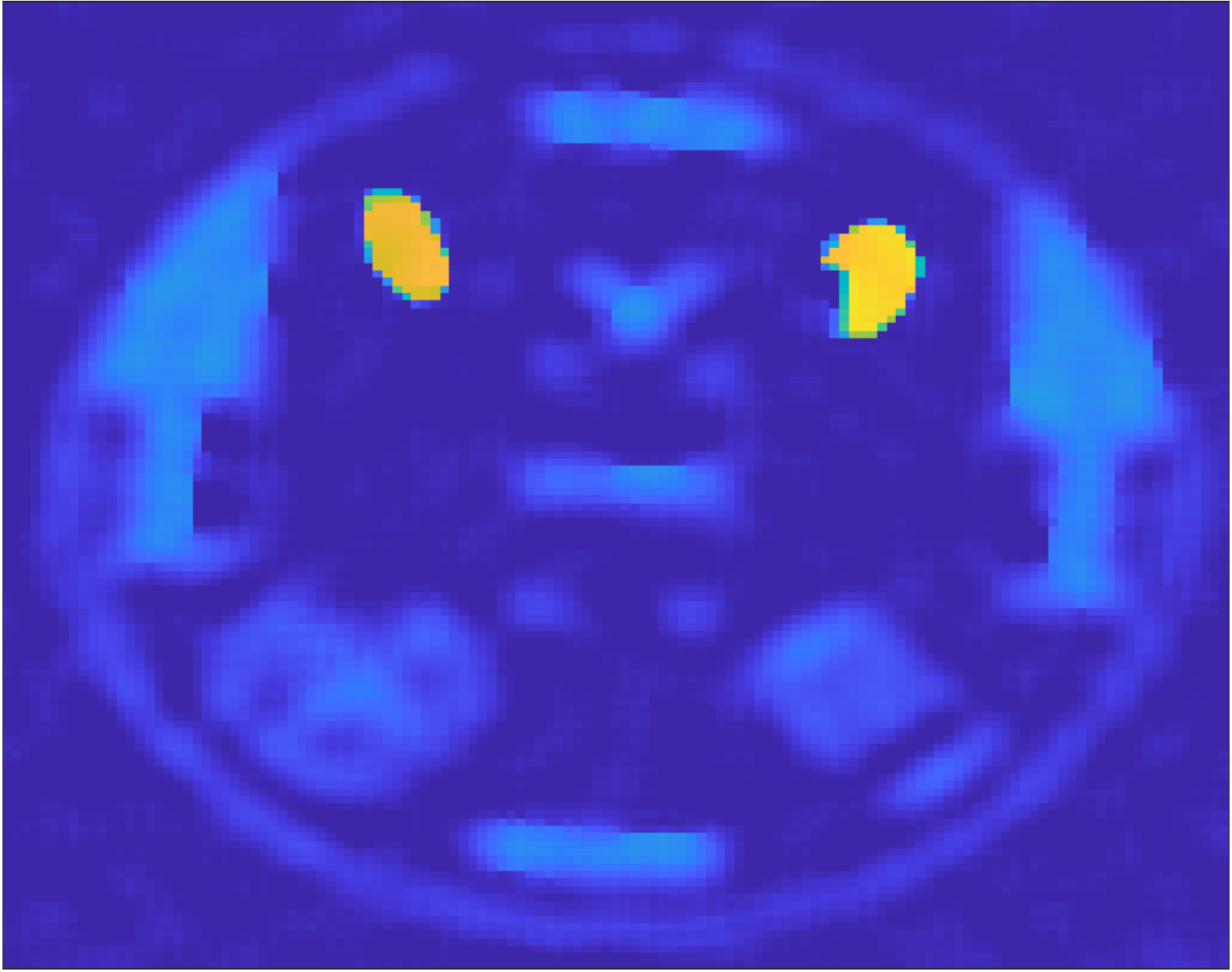}
    \caption{$\tilde{\mathbf{f}}_2^\text{VBJS}$}
    \end{subfigure}
    \begin{subfigure}[b]{.19\textwidth}
    \includegraphics[width=\textwidth]{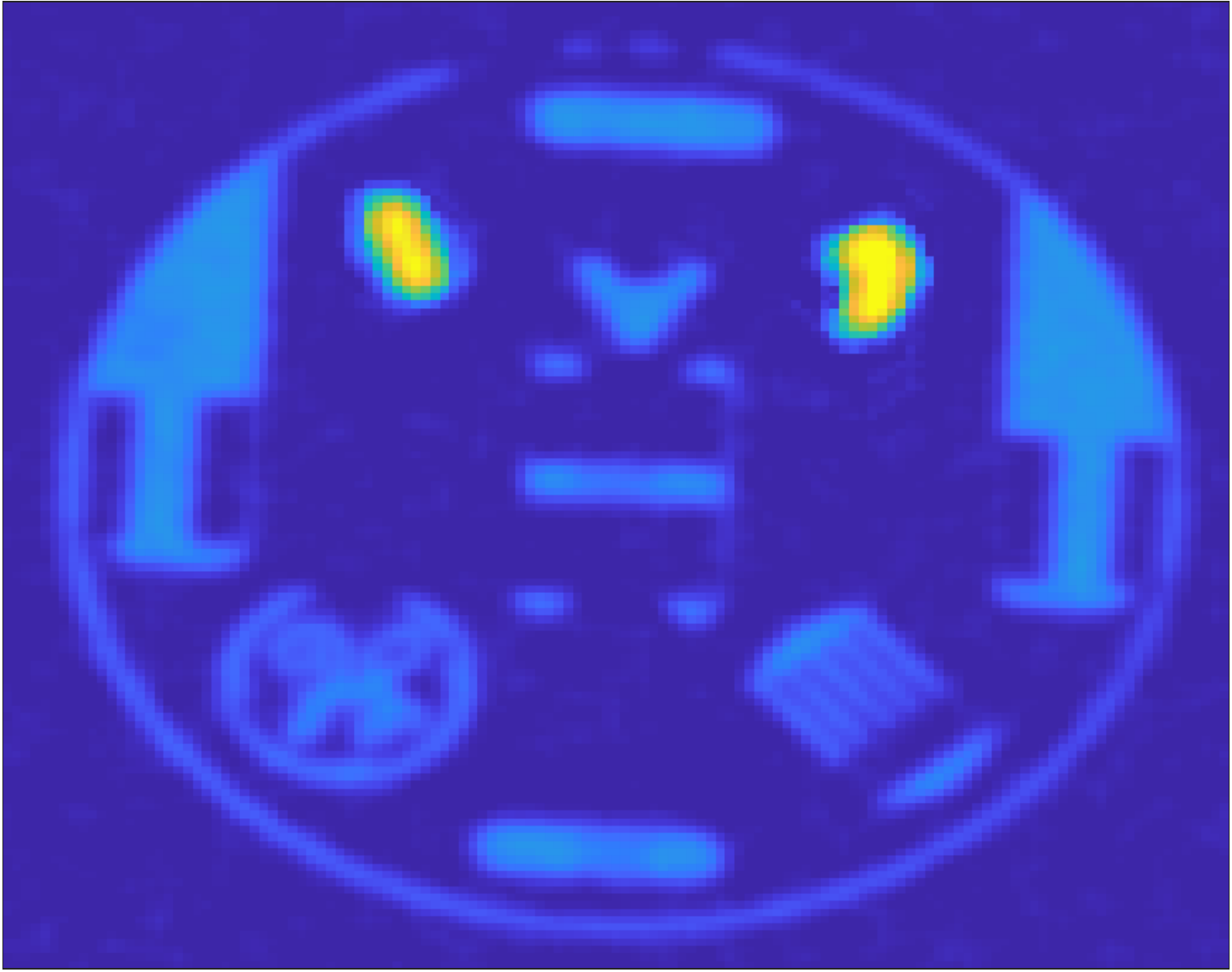}
    \caption{$\tilde{\mathbf{f}}_2^\text{joint}$}
    \end{subfigure}
    
    \caption{Reconstruction when a moving object is obscured.}
    \label{fig:block_obj_b2_GE}
\end{figure}

Finally, as our new method relies on accurate change masks to capture moving objects in adjacent images of the scene, it is important to consider the case when a moving object is obscured in one of the data acquisitions. This case is shown in Figure \ref{fig:block_obj_b2_GE}, where  we observe that even when the moving ellipse is obscured the image recovery using  our new approach is no worse than either of the  two individual recovery methods discussed.

\subsection{Experiment 2: Synthetic Aperture Radar (SAR)}
\label{subsec:golf_course}

\begin{figure}[h!]
    \centering
    \begin{subfigure}[b]{.19\textwidth}
        \includegraphics[width=\textwidth]{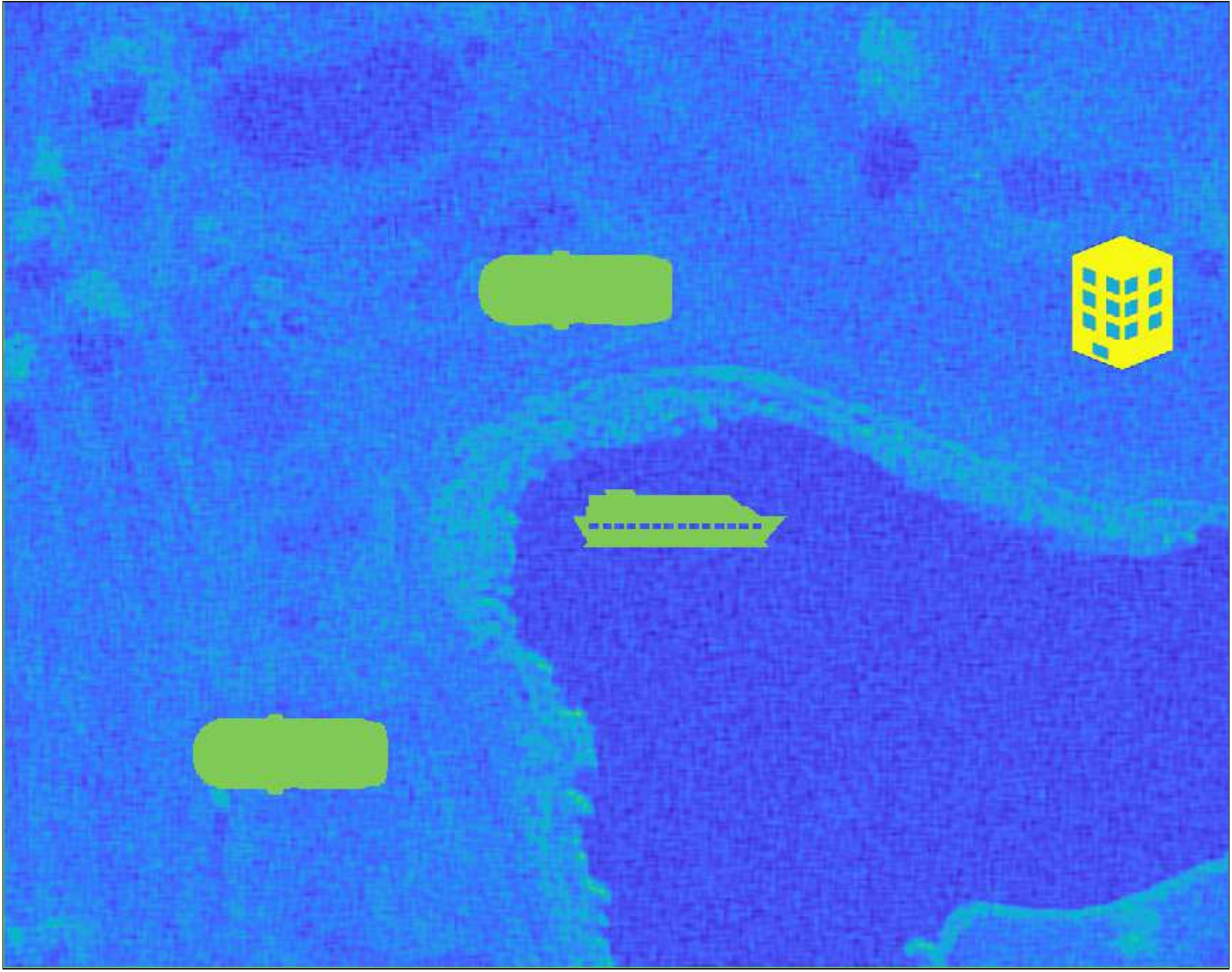}
        \caption{$\mathbf{f}_1$}
    \end{subfigure}
    ~
    \begin{subfigure}[b]{.19\textwidth}
        \includegraphics[width=\textwidth]{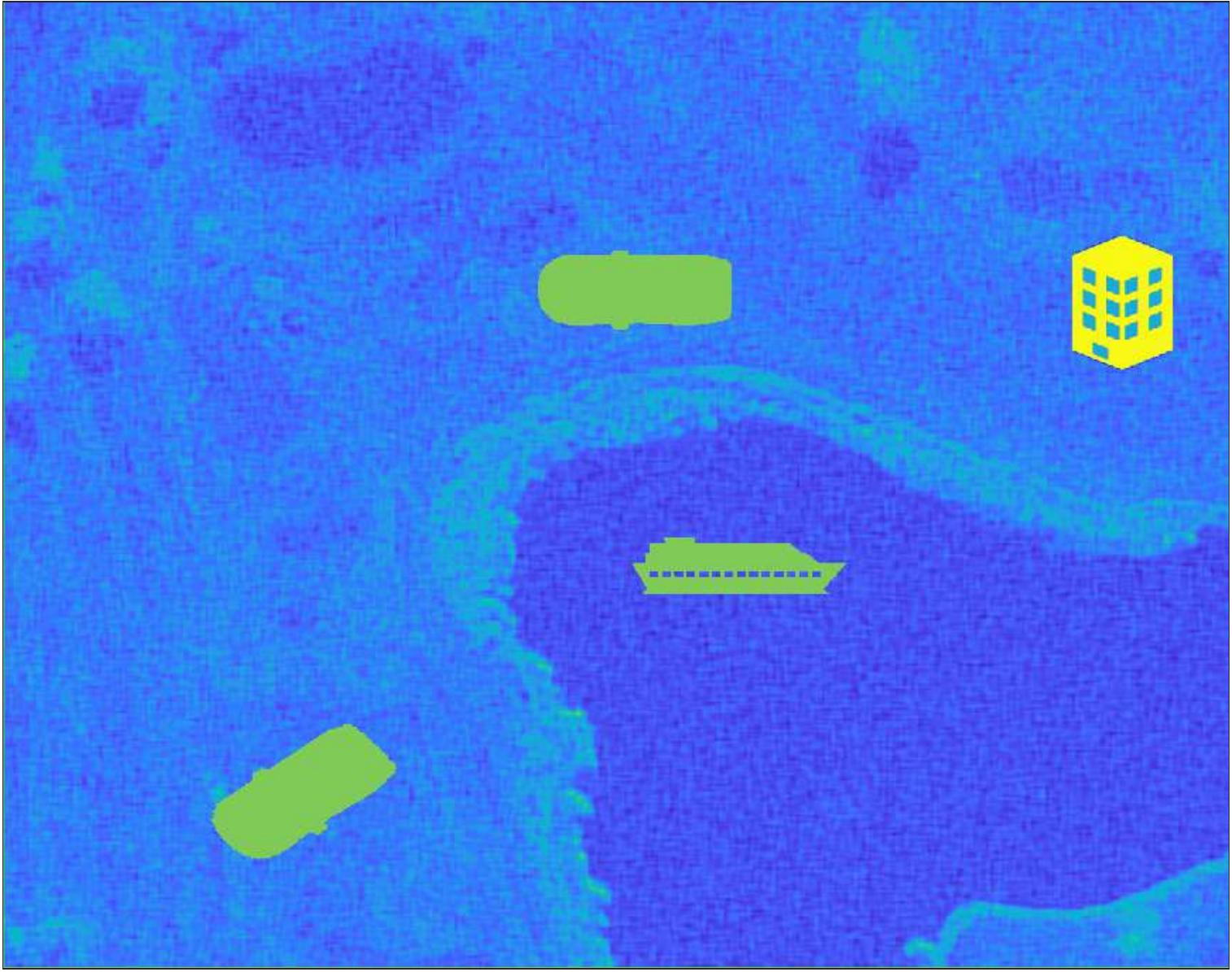}
        \caption{$\mathbf{f}_2$}
    \end{subfigure}
    ~
    \begin{subfigure}[b]{.19\textwidth}
        \includegraphics[width=\textwidth]{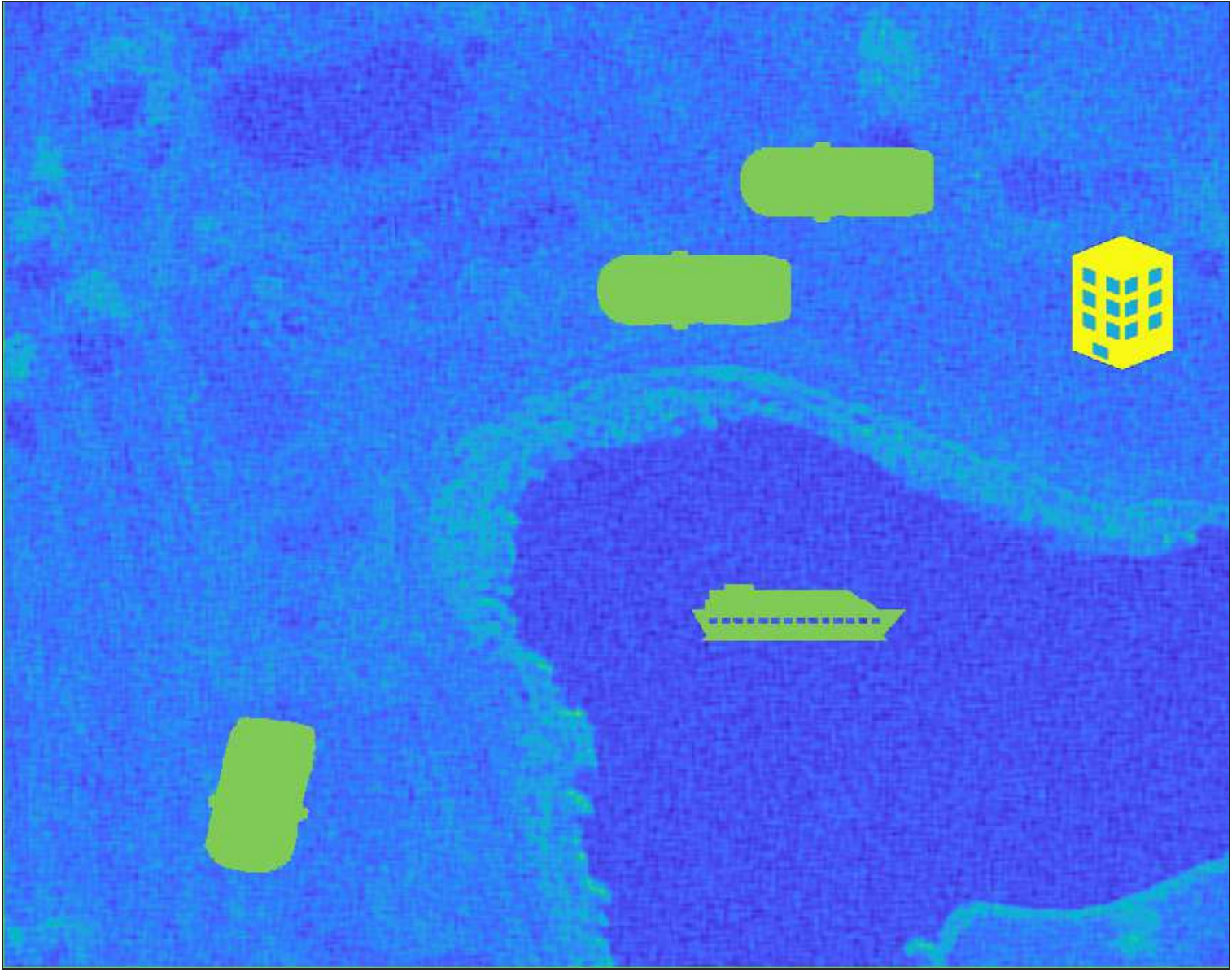}
        \caption{$\mathbf{f}_3$}
    \end{subfigure}
    ~
    \begin{subfigure}[b]{.19\textwidth}
        \includegraphics[width=\textwidth]{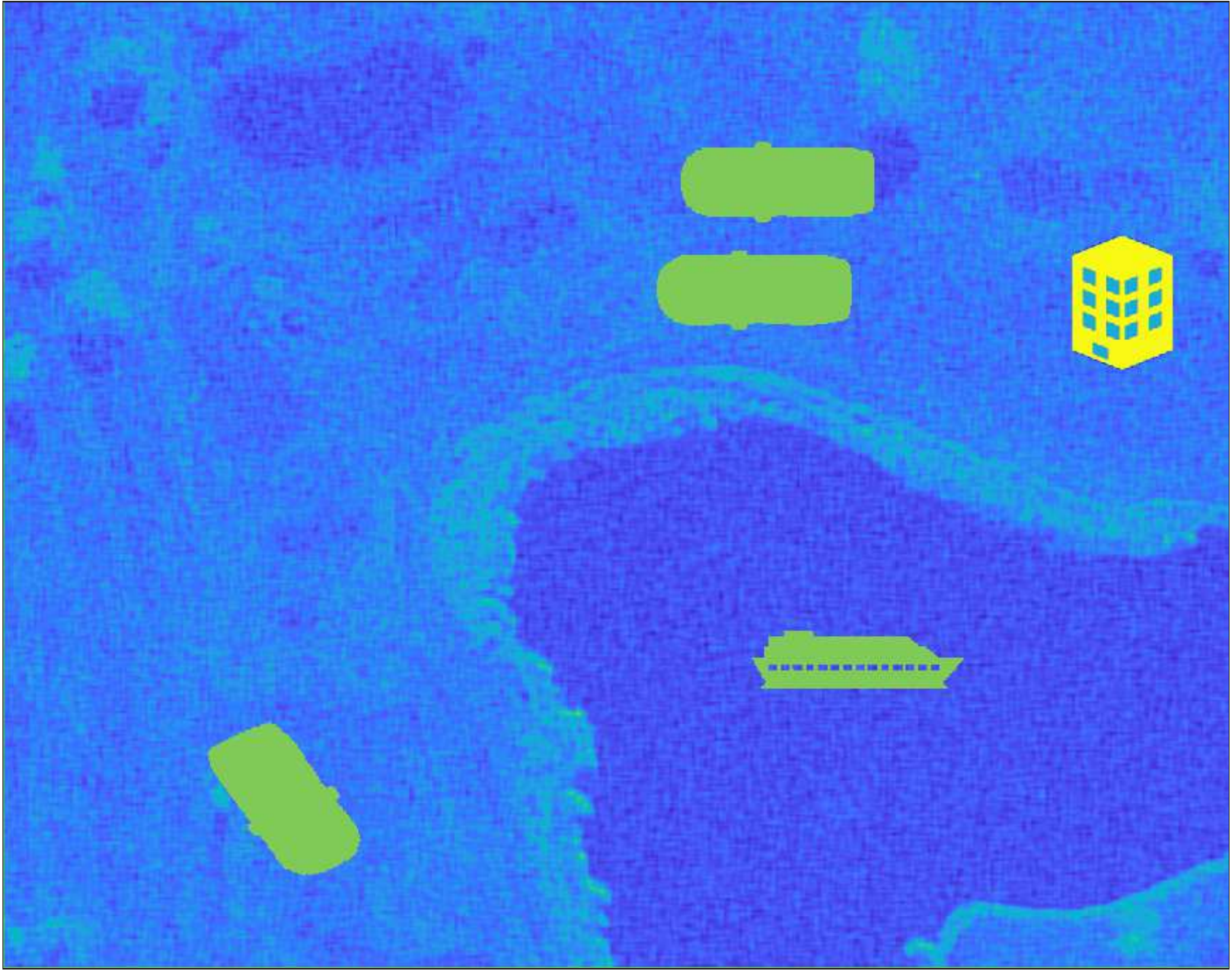}
        \caption{$\mathbf{f}_4$}
    \end{subfigure}
    \\
    \begin{subfigure}[b]{.19\textwidth}
        \includegraphics[width=\textwidth]{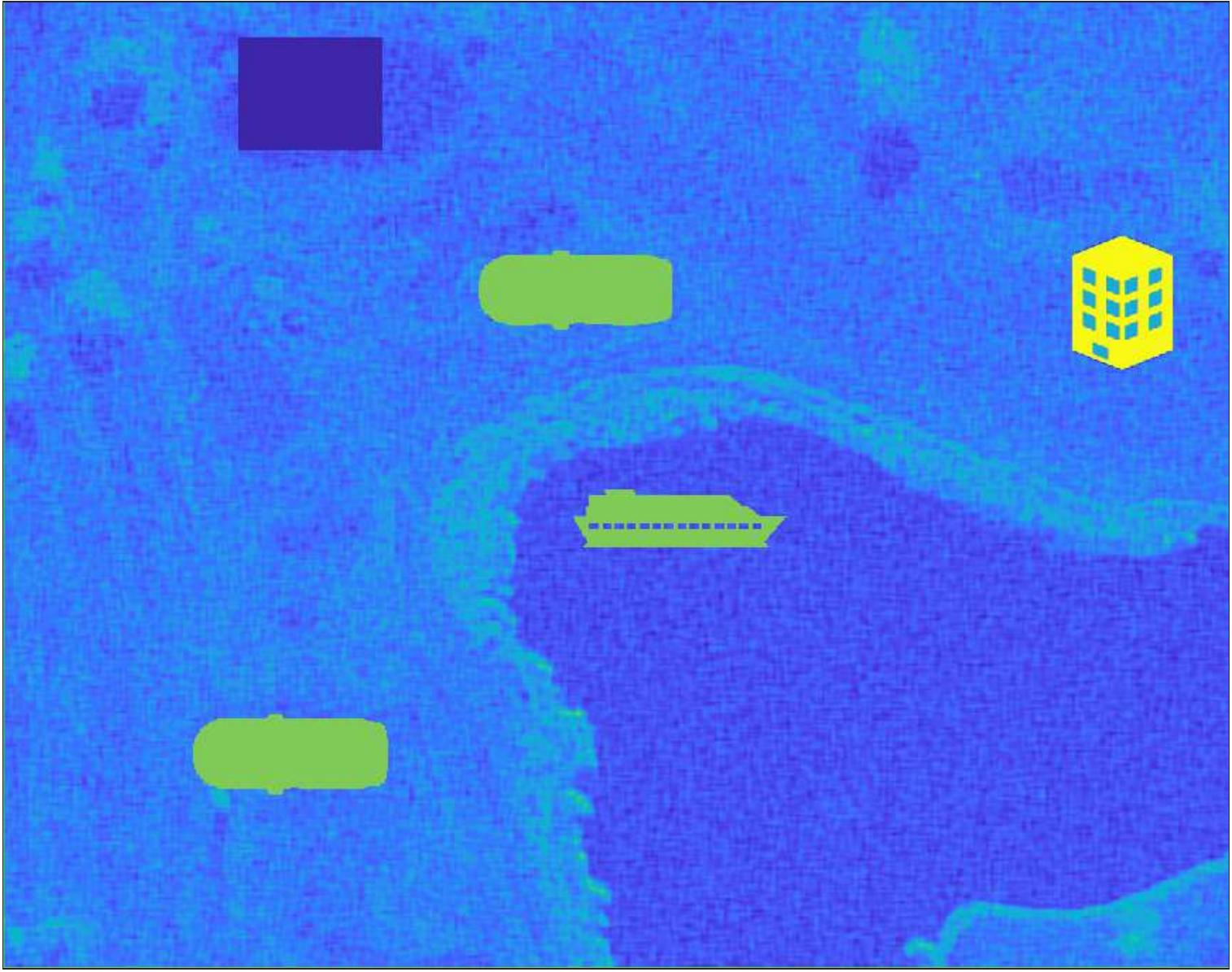}
        \caption{Obstacle$_1$}
    \end{subfigure}
    ~
    \begin{subfigure}[b]{.19\textwidth}
        \includegraphics[width=\textwidth]{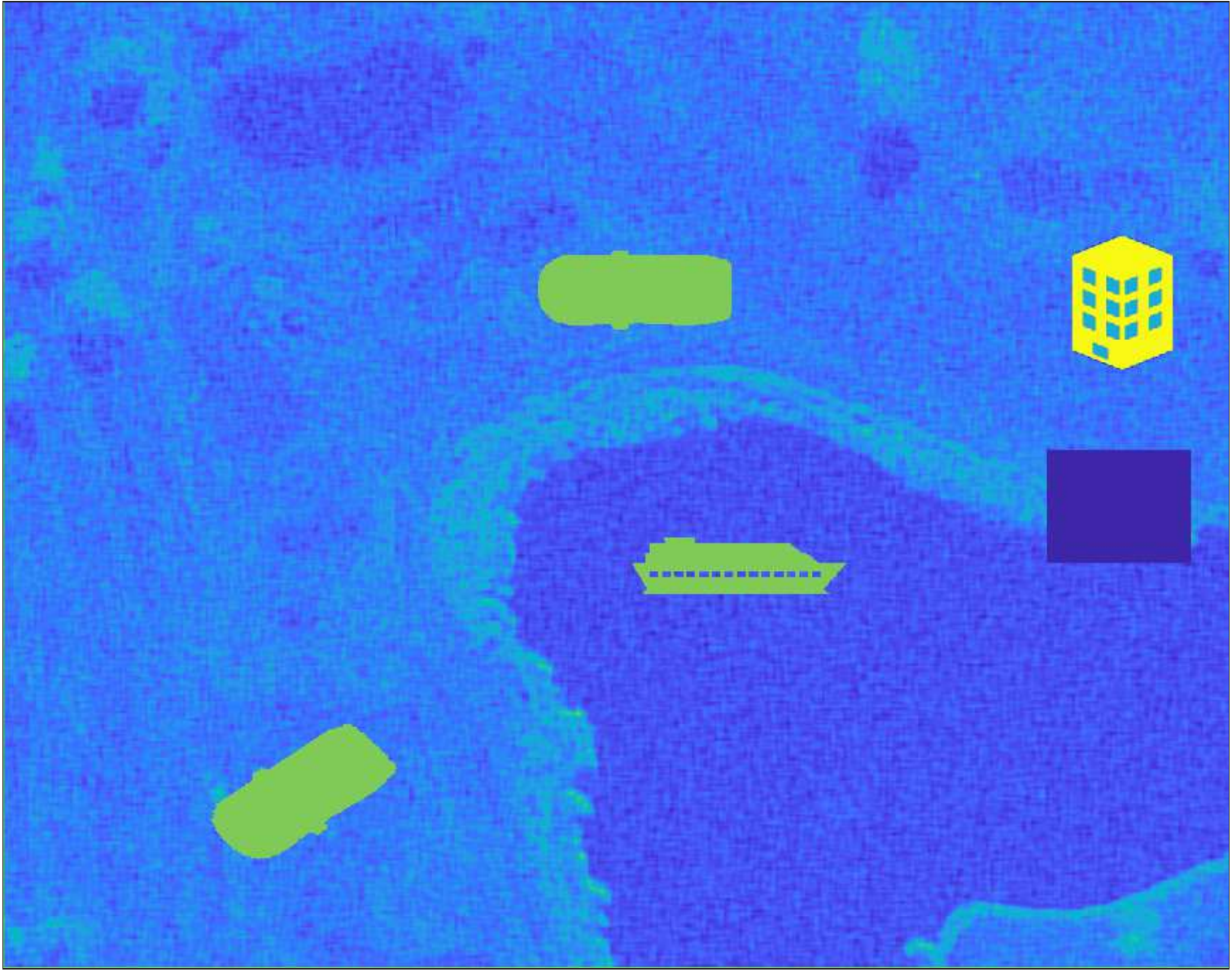}
        \caption{Obstacle$_2$}
    \end{subfigure}
    ~
    \begin{subfigure}[b]{.19\textwidth}
        \includegraphics[width=\textwidth]{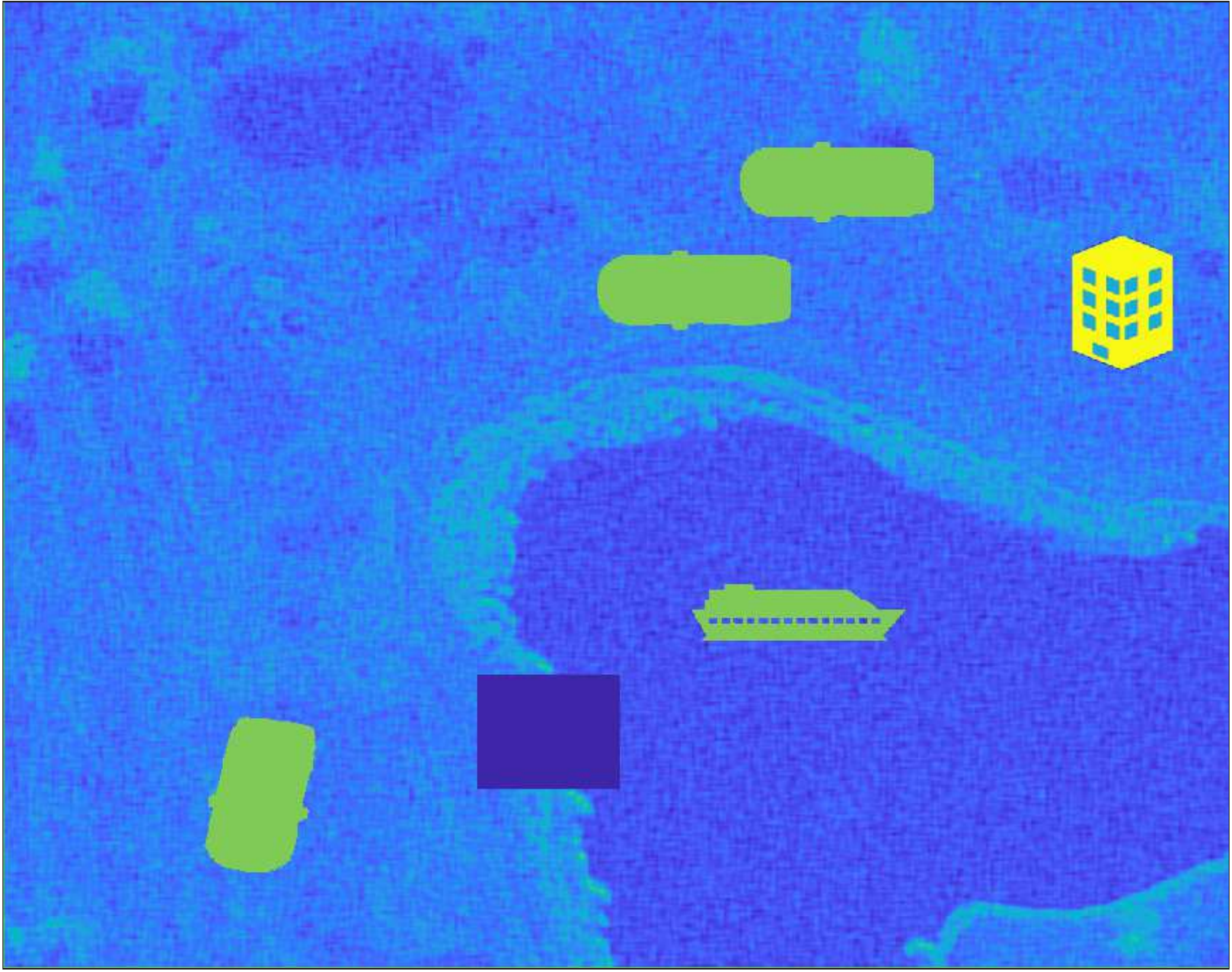}
        \caption{Obstacle$_3$}
    \end{subfigure}
    ~
    \begin{subfigure}[b]{.19\textwidth}
        \includegraphics[width=\textwidth]{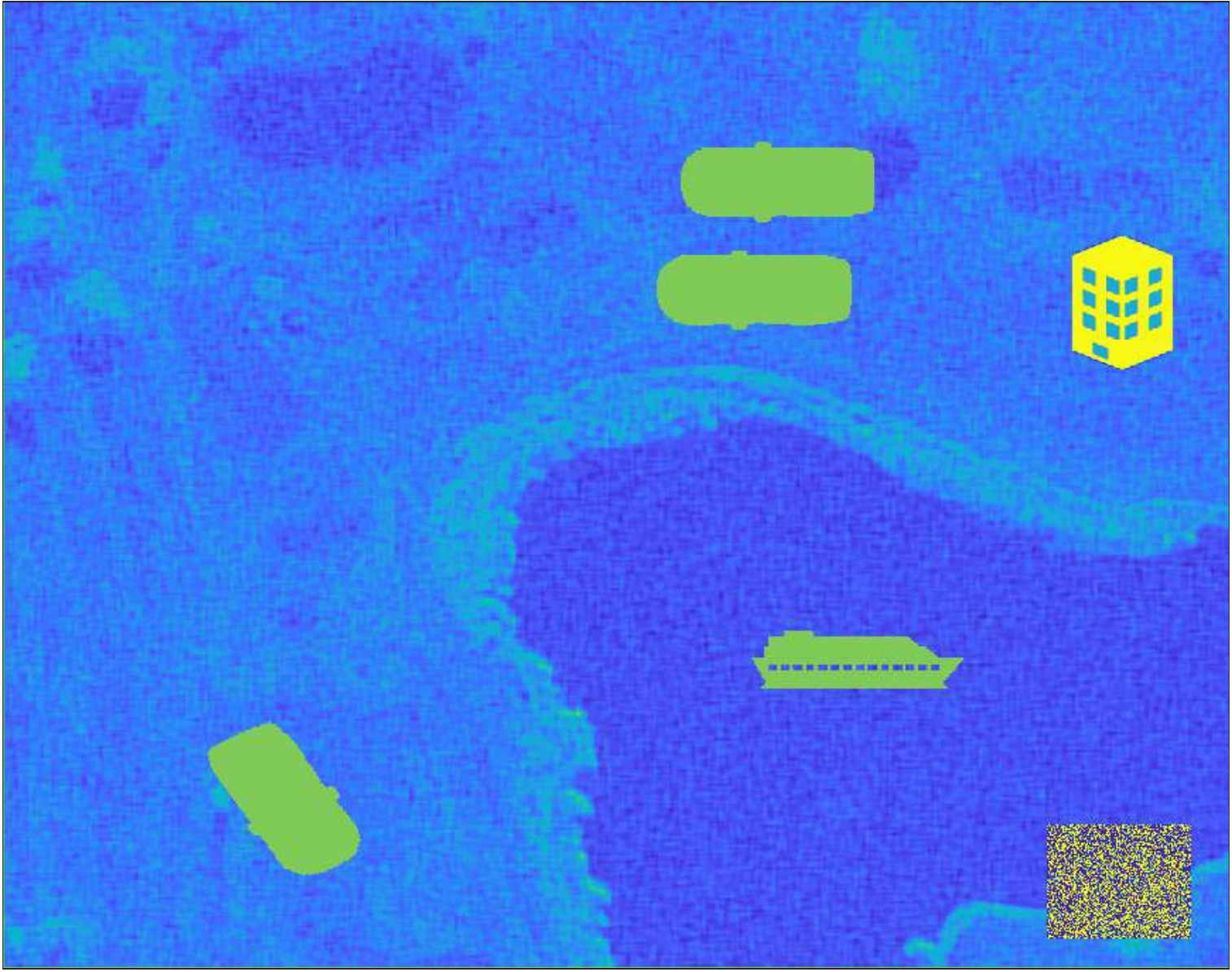}
        \caption{Obstacle$_4$}
    \end{subfigure}
    \\
    \begin{subfigure}[b]{.19\textwidth}
        \includegraphics[width=\textwidth]{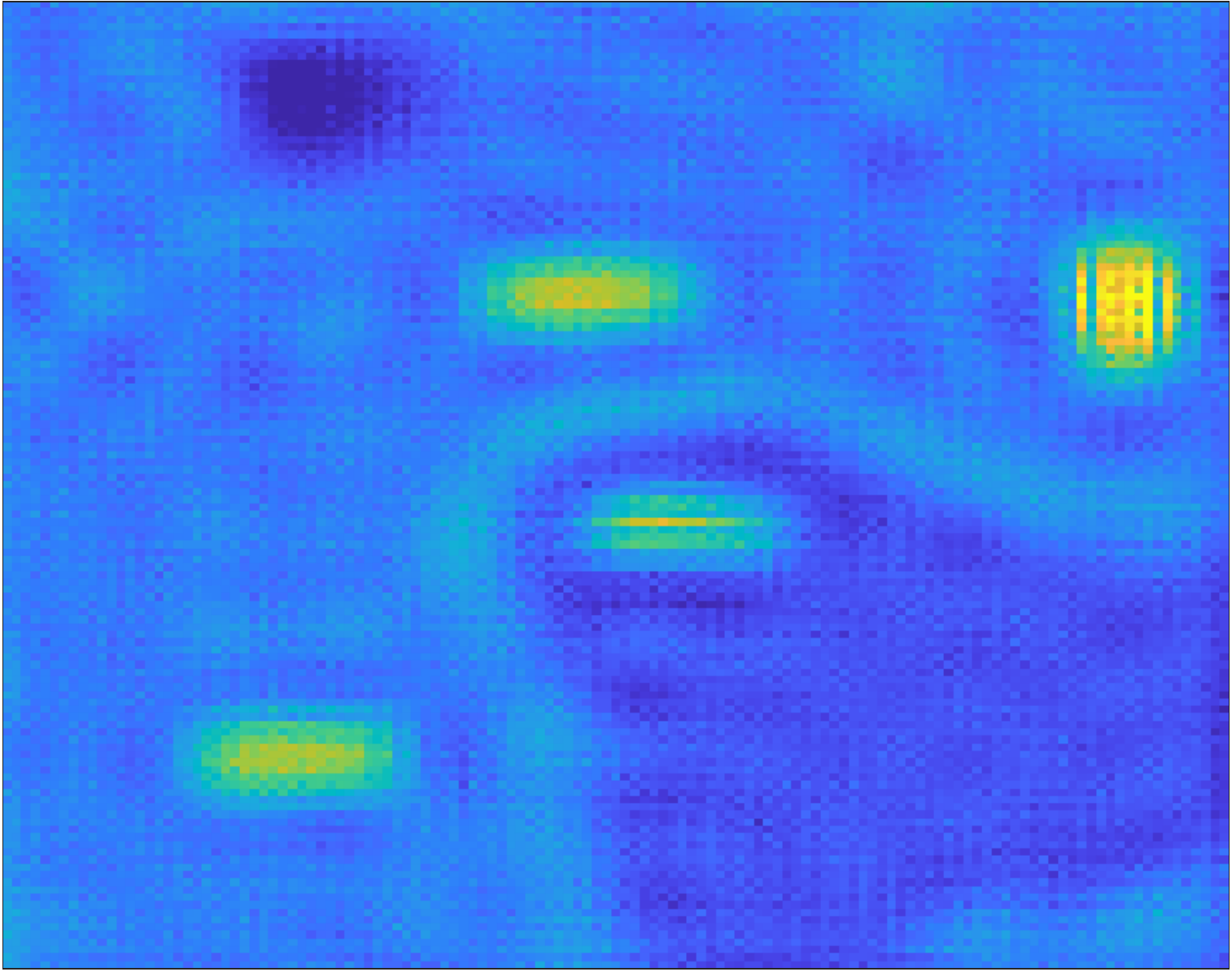}
        \caption{$\tilde{f}_{1}$}
    \end{subfigure}
    ~
    \begin{subfigure}[b]{.19\textwidth}
        \includegraphics[width=\textwidth]{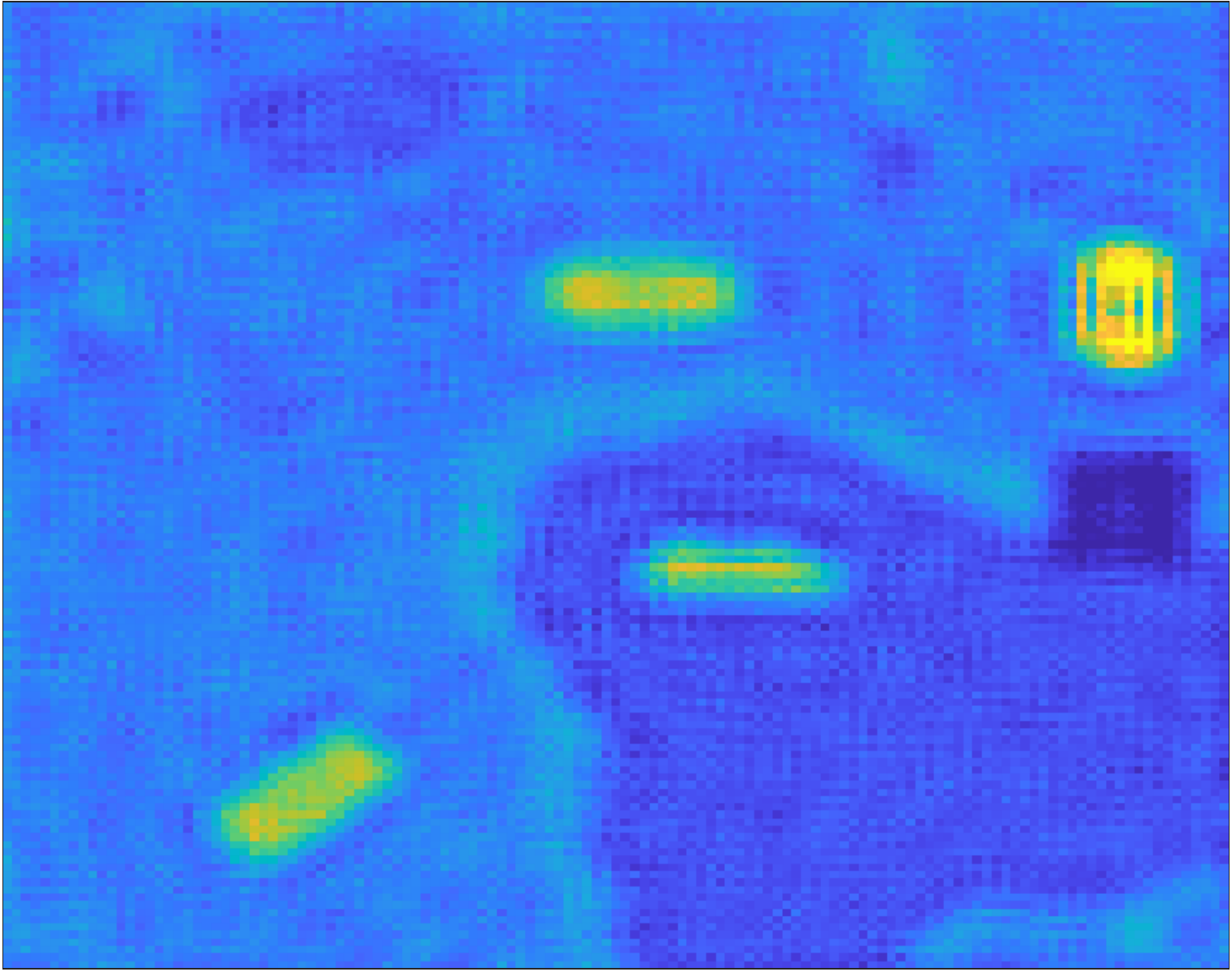}
        \caption{$\tilde{f}_{2}$}
    \end{subfigure}
    ~
    \begin{subfigure}[b]{.19\textwidth}
        \includegraphics[width=\textwidth]{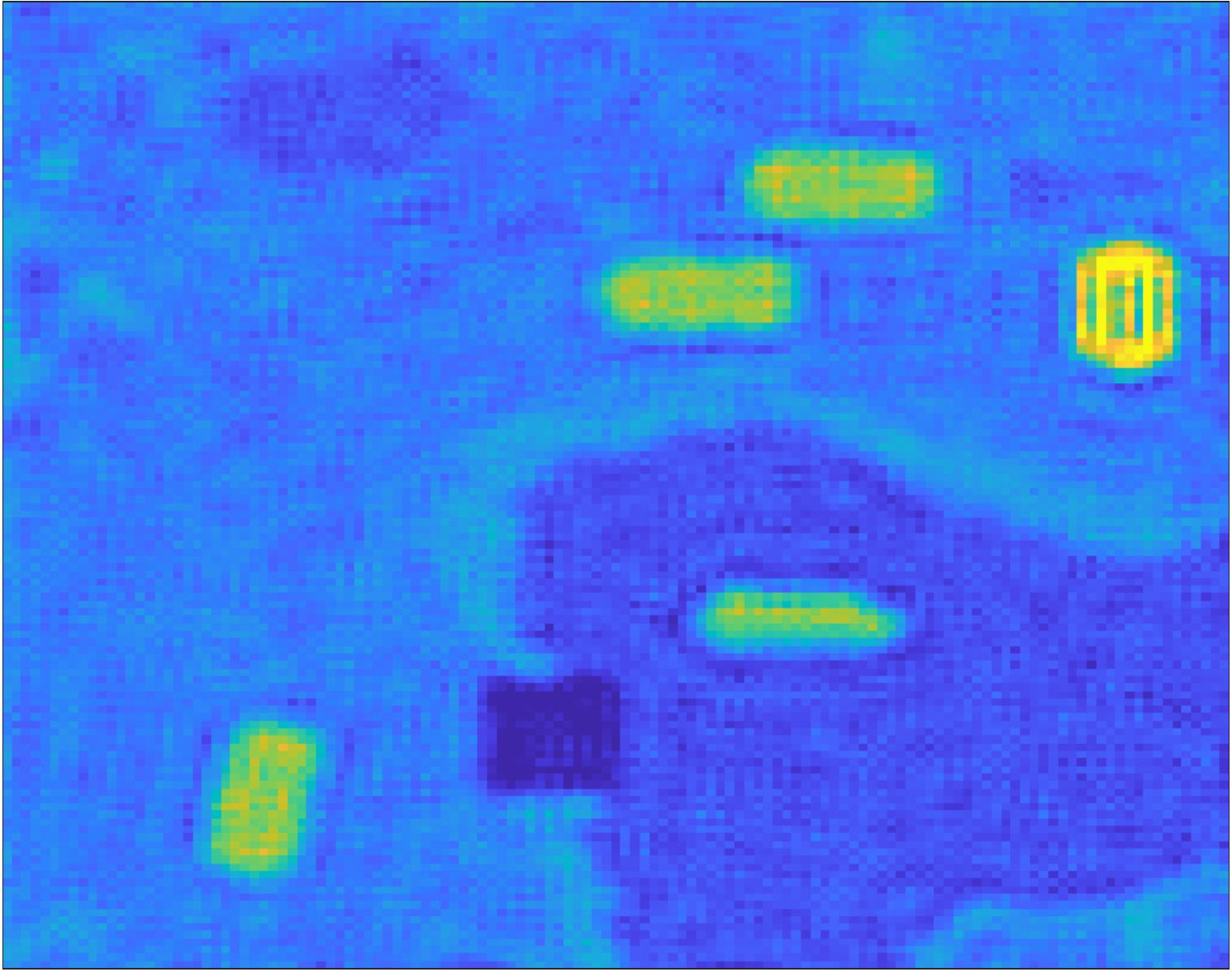}
        \caption{$\tilde{f}_{3}$}
    \end{subfigure}
    ~
    \begin{subfigure}[b]{.19\textwidth}
        \includegraphics[width=\textwidth]{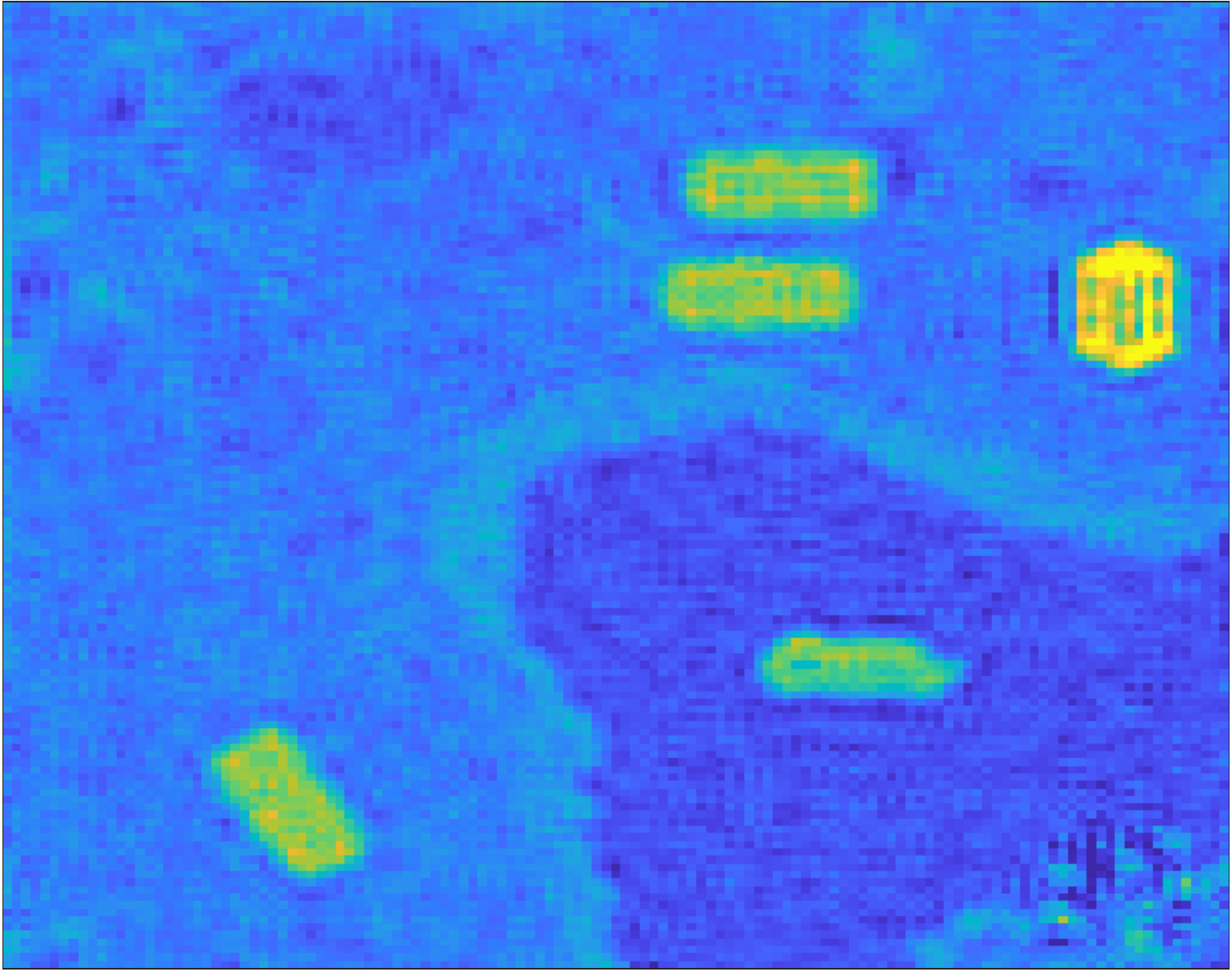}
        \caption{$\tilde{f}_{4}$}
    \end{subfigure}
    \\
    \begin{subfigure}[b]{.19\textwidth}
        \includegraphics[width=\textwidth]{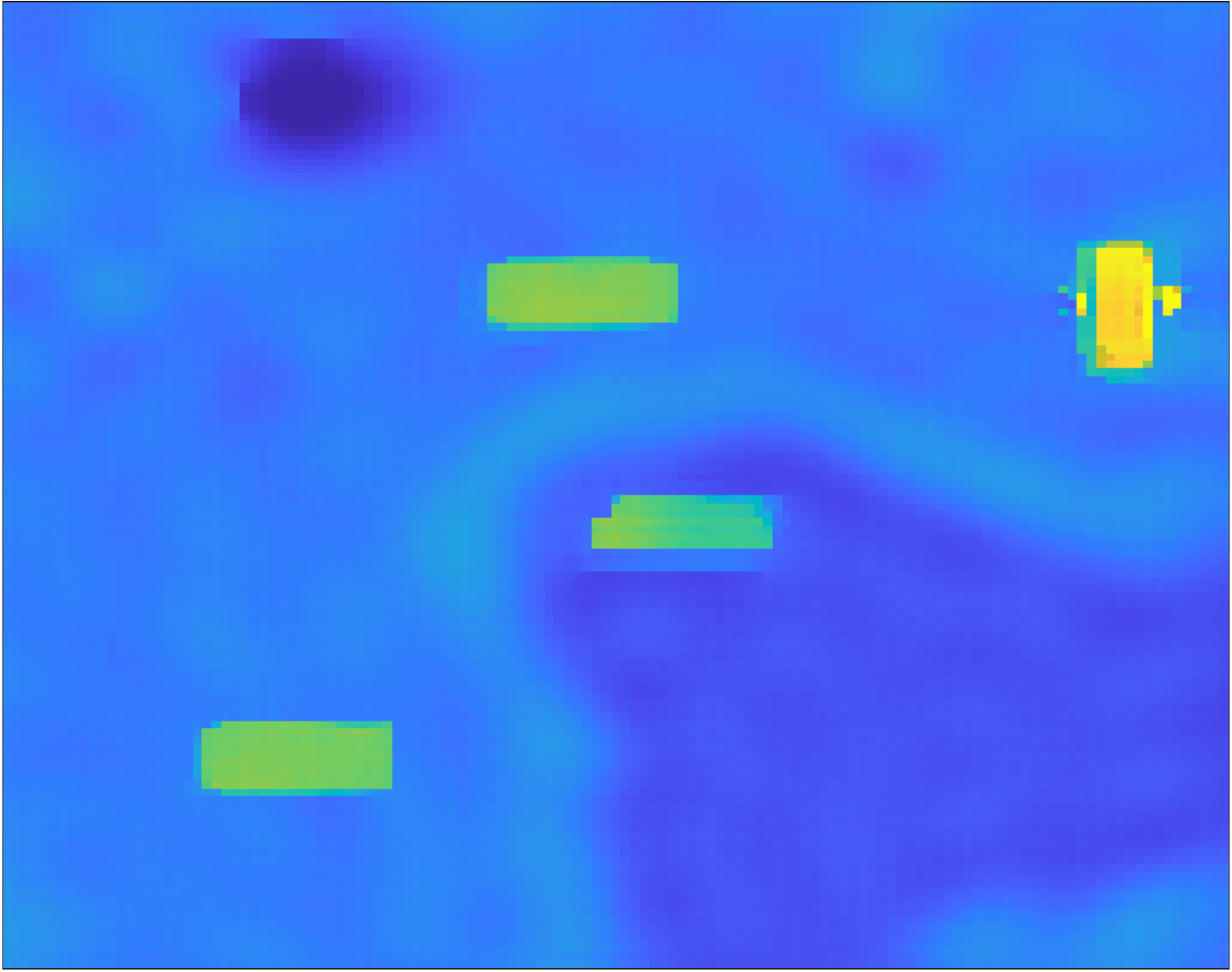}
        \caption{$\tilde{f}_1^\text{VBJS}$}
    \end{subfigure}
    ~
    \begin{subfigure}[b]{.19\textwidth}
        \includegraphics[width=\textwidth]{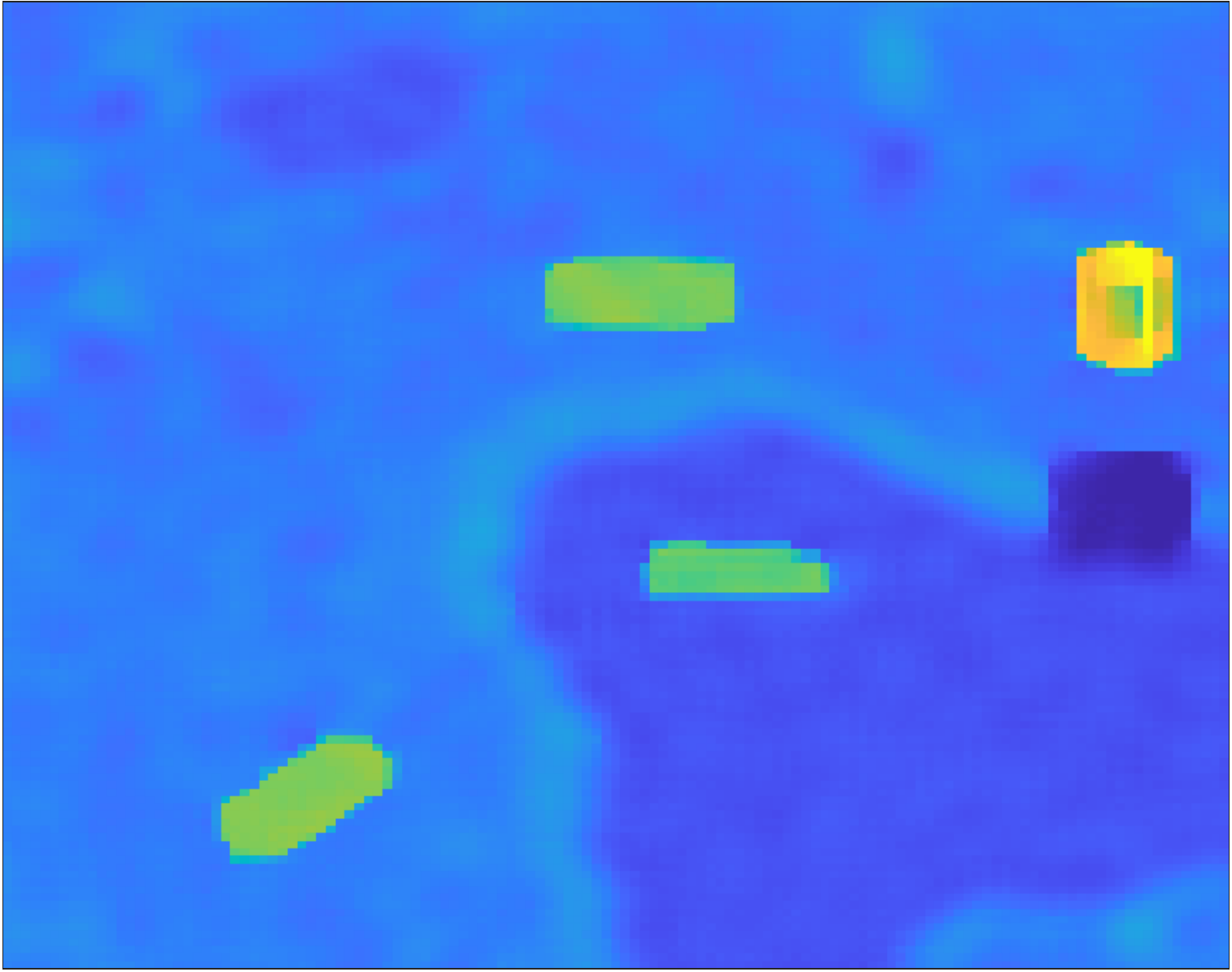}
        \caption{$\tilde{f}_2^\text{VBJS}$}
    \end{subfigure}
    ~
    \begin{subfigure}[b]{.19\textwidth}
        \includegraphics[width=\textwidth]{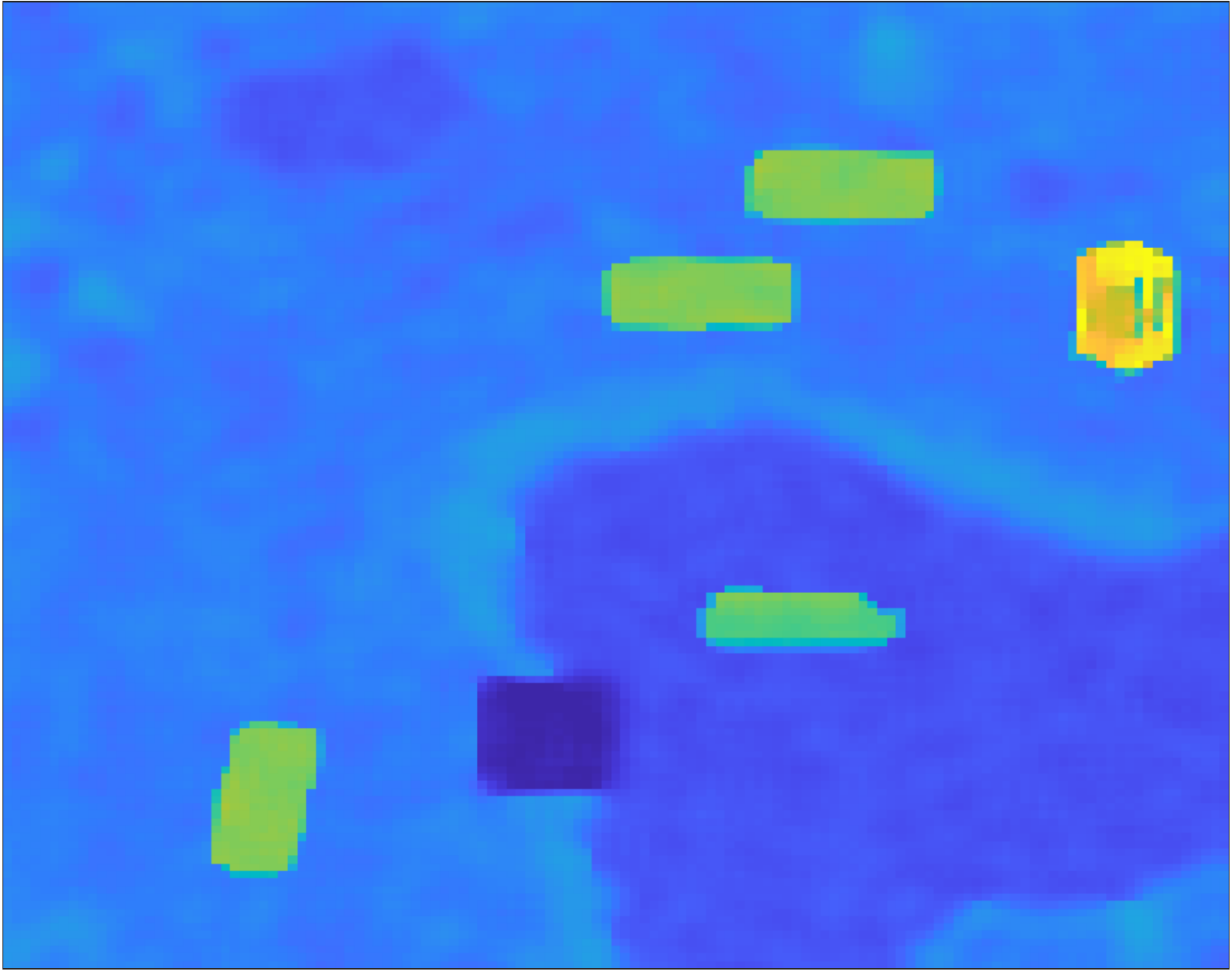}
        \caption{$\tilde{f}_3^\text{VBJS}$}
    \end{subfigure}
    ~
    \begin{subfigure}[b]{.19\textwidth}
        \includegraphics[width=\textwidth]{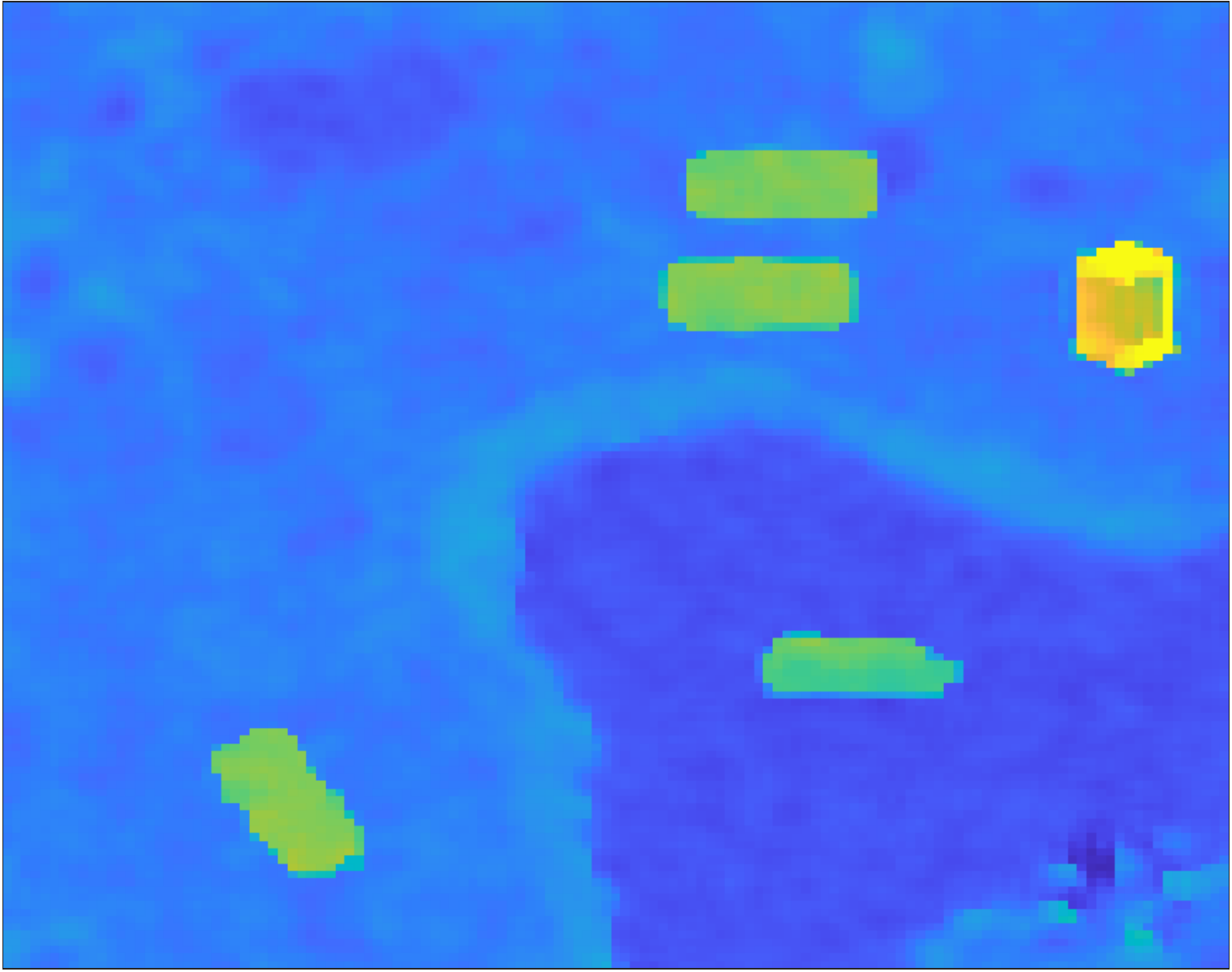}
        \caption{$\tilde{f}_4^\text{VBJS}$}
    \end{subfigure}
    \\ 
    \begin{subfigure}[b]{.19\textwidth}
        \includegraphics[width=\textwidth]{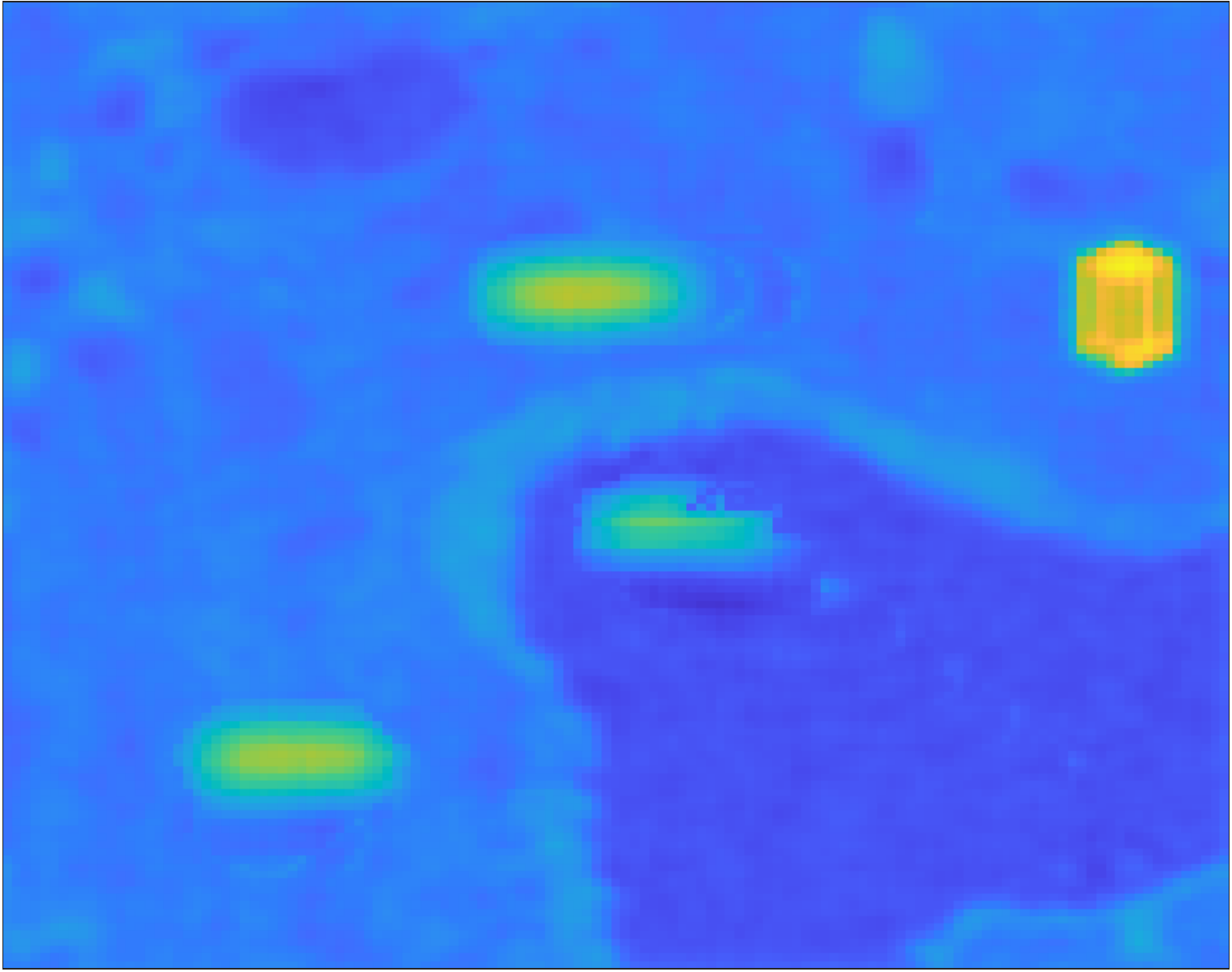}
        \caption{$\tilde{f}_1^\text{joint}$}
    \end{subfigure}
    ~
    \begin{subfigure}[b]{.19\textwidth}
        \includegraphics[width=\textwidth]{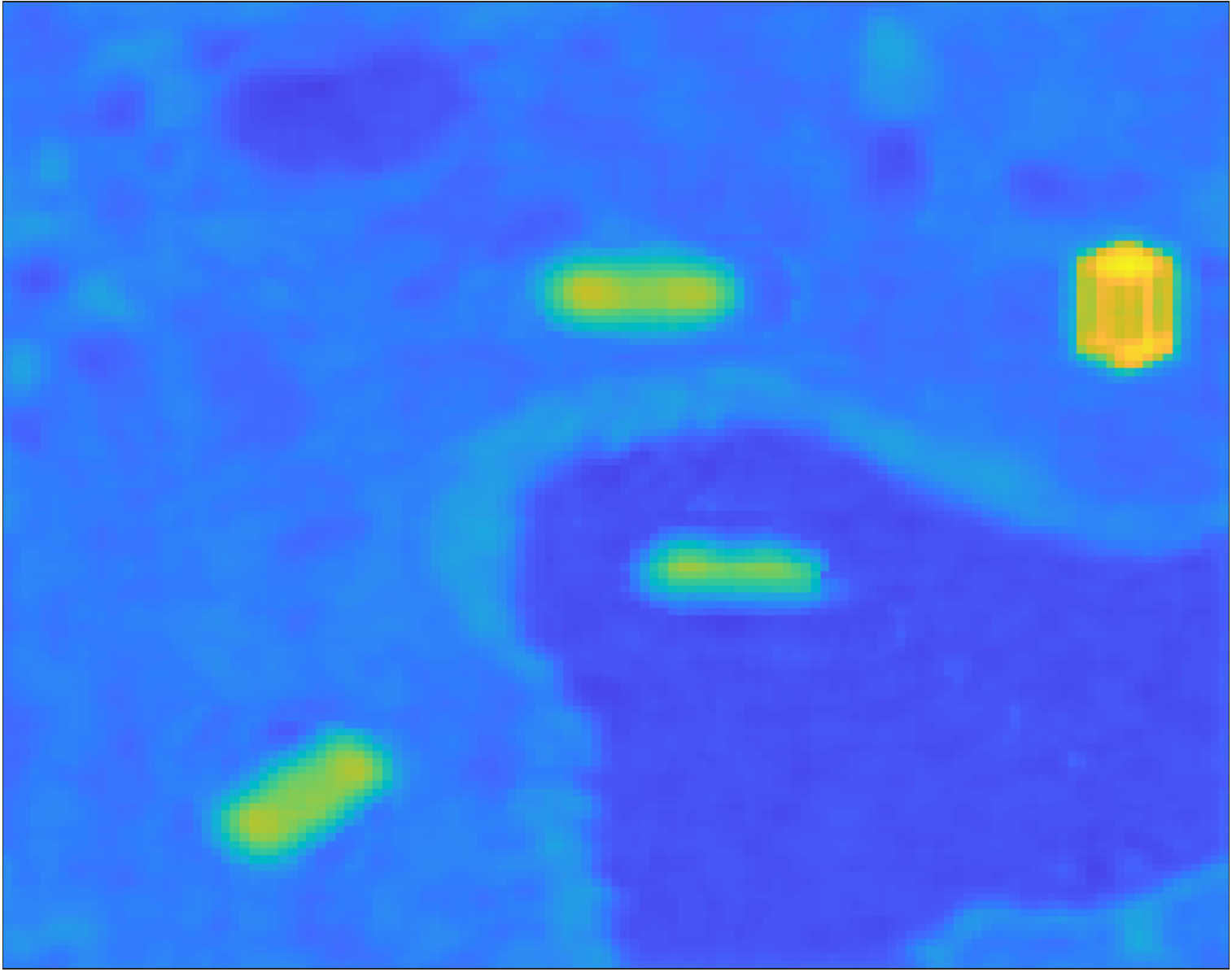}
        \caption{$\tilde{f}_2^\text{joint}$}
    \end{subfigure}
    ~
    \begin{subfigure}[b]{.19\textwidth}
        \includegraphics[width=\textwidth]{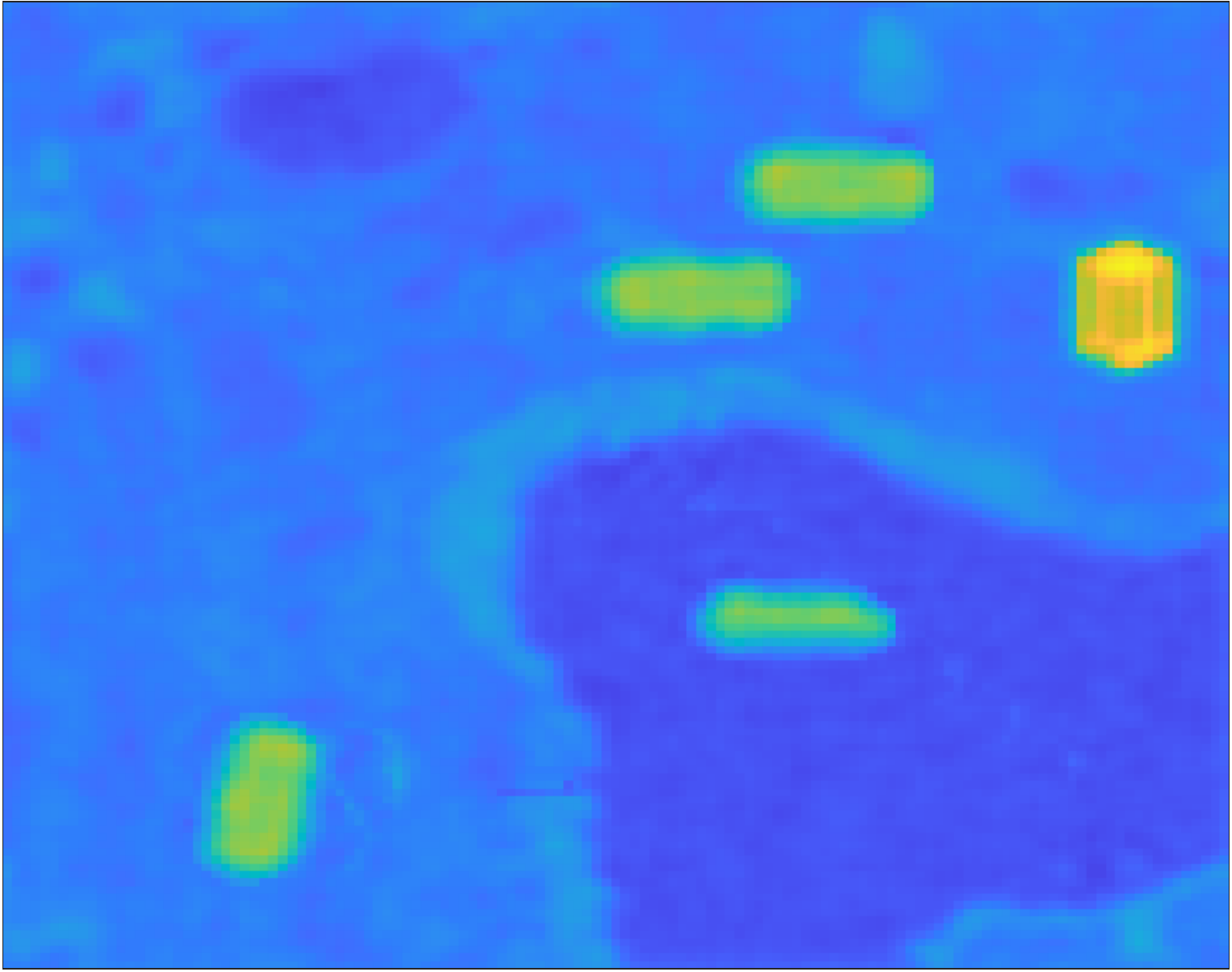}
        \caption{$\tilde{f}_3^\text{joint}$}
    \end{subfigure}
    ~
    \begin{subfigure}[b]{.19\textwidth}
        \includegraphics[width=\textwidth]{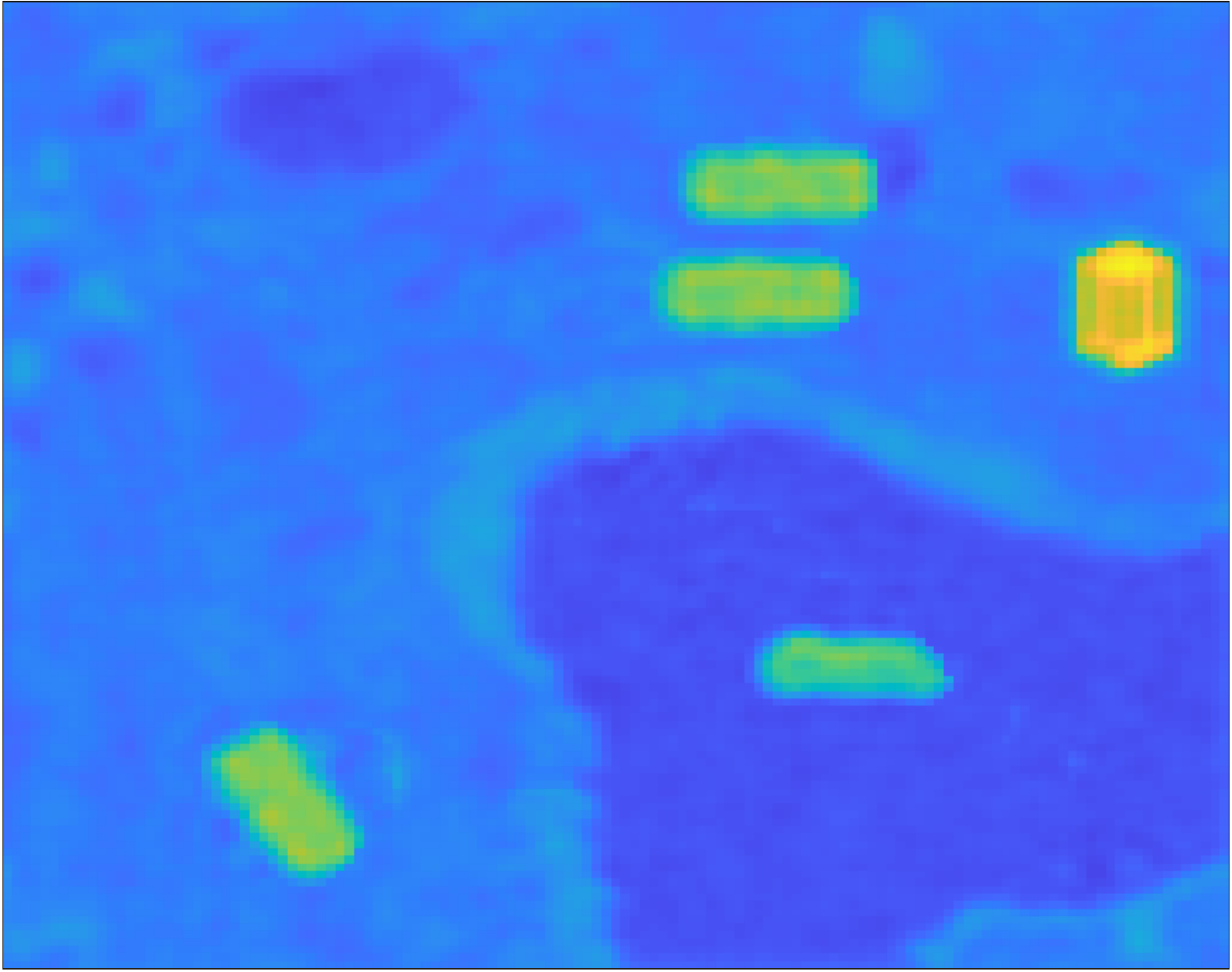}
        \caption{$\tilde{f}_4^\text{joint}$}
    \end{subfigure}
    \caption{ (first row) The underlying scene is the SAR image of a golf course, \cite{SAR_Image_ref}.  (second row) 
    Physical images in scaled color involving zero-valued blocks ({first three images}) and added zero-mean Gaussian noise with standard deviation $3$ ({last image}).  (third row)
    Separate reconstructions by standard $\ell_1$-regularization. (fourth row) 
    Recovery using VBJS. (bottom row) 
    The new joint recovery. }
\label{fig:rec_golf}
\end{figure}

Our second experiment considers a temporal sequence of six {SAR} images depicting a golf course, \cite{SAR_Image_ref},  displayed in Figure \ref{fig:rec_golf}. Observe there is no ``ground truth'' in this case.
We then add cars and boats, each of magnitude $1$, and a building of magnitude $1.5$, on top of the scene to simulate respectively moving and background objects.  In our experiment one car is inserted into the third image and a different car is deleted in the fifth image (not shown).
As displayed in the second {row} of Figure \ref{fig:rec_golf}, we also introduce several obstacles in the physical domain. Instead of using zero values in these locations across the sequence, as was done in the MRI experiment, here we use a zero-mean Gaussian distribution with standard deviation $3$ in each of the latter (4-6) part of the sequence.  For viewing purposes we have only depicted this added noise in the fourth image ({last column}, second {row}).  As in the MRI experiment, the last three {row}s in Figure \ref{fig:rec_golf} respectively compare image recovery using  the standard $\ell_1$ regularization \eqref{eq:l1_reg} using Algorithm \ref{algo:EM_l1_reg}, the VBJS technique \eqref{eq:wl1_reg} using Algorithm \ref{algo:mmv_recovery}, and our new joint recovery method \eqref{eq:optModel} as realized by Algorithm \ref{algo:joint_recovery}. It is evident that our new method yields improved results for all images, in particular in places of obstruction.  We note that the number of pixels, the noise level, and the missing frequencies are the same as in Section \ref{subsec:phantom}.

\begin{figure}[h!]
\centering 
\begin{subfigure}[b]{0.23\textwidth}
    \includegraphics[width=\textwidth]{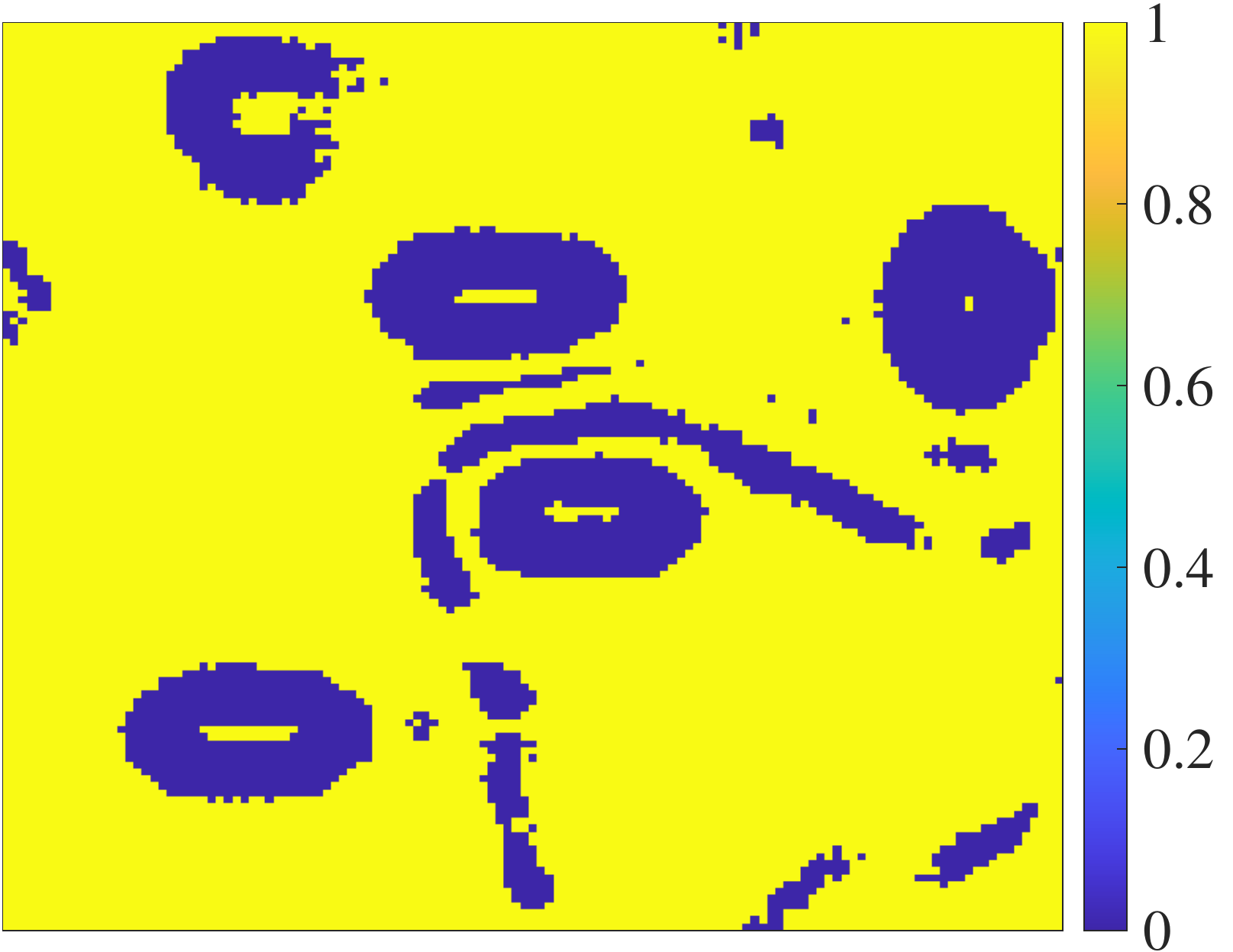}
    \caption{$W_1$}
\end{subfigure}
~
\begin{subfigure}[b]{0.23\textwidth}
    \includegraphics[width=\textwidth]{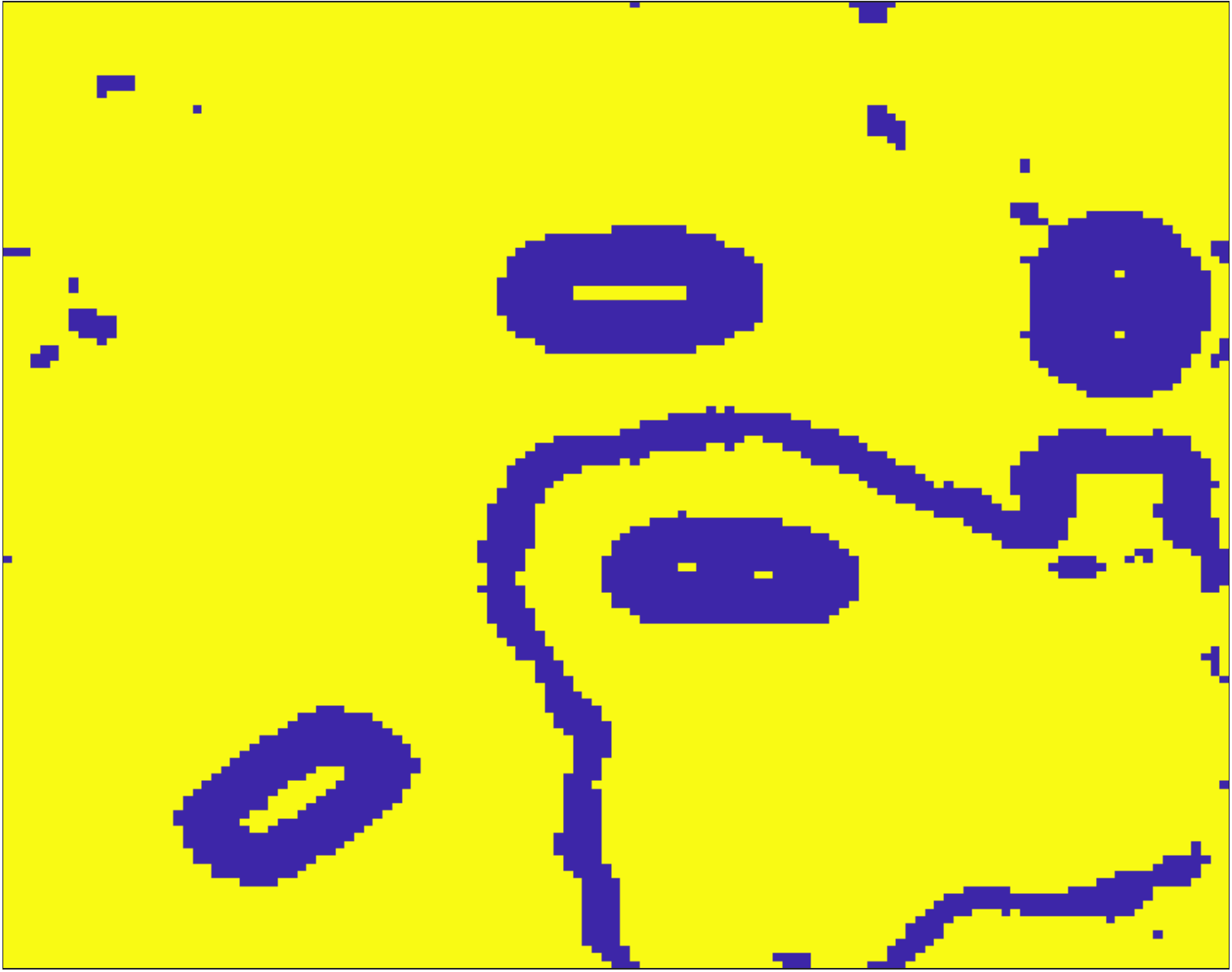}
    \caption{$W_2$}
\end{subfigure}
~
\begin{subfigure}[b]{0.23\textwidth}
    \includegraphics[width=\textwidth]{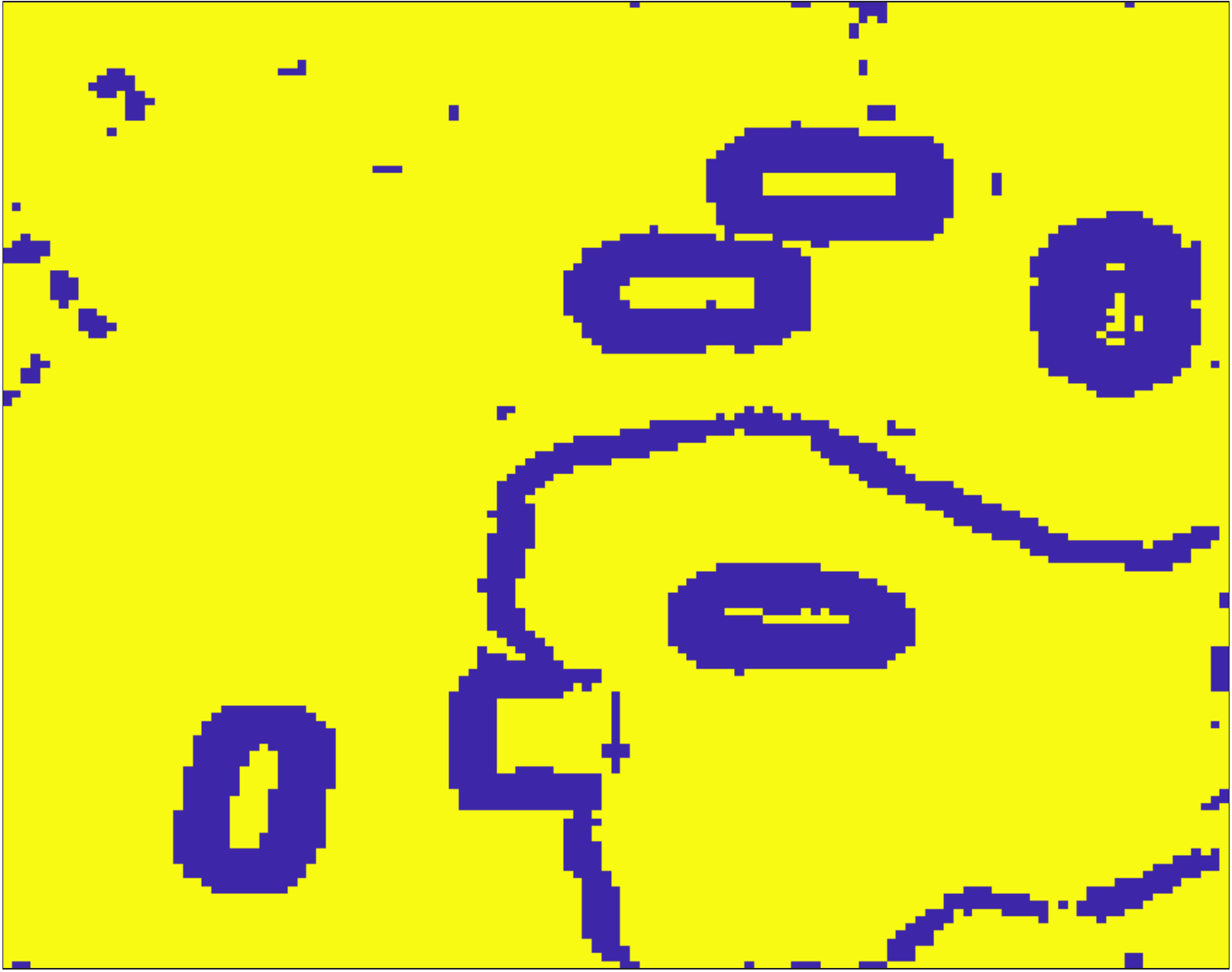}
    \caption{$W_3$}
\end{subfigure}
~
\begin{subfigure}[b]{0.23\textwidth}
    \includegraphics[width=\textwidth]{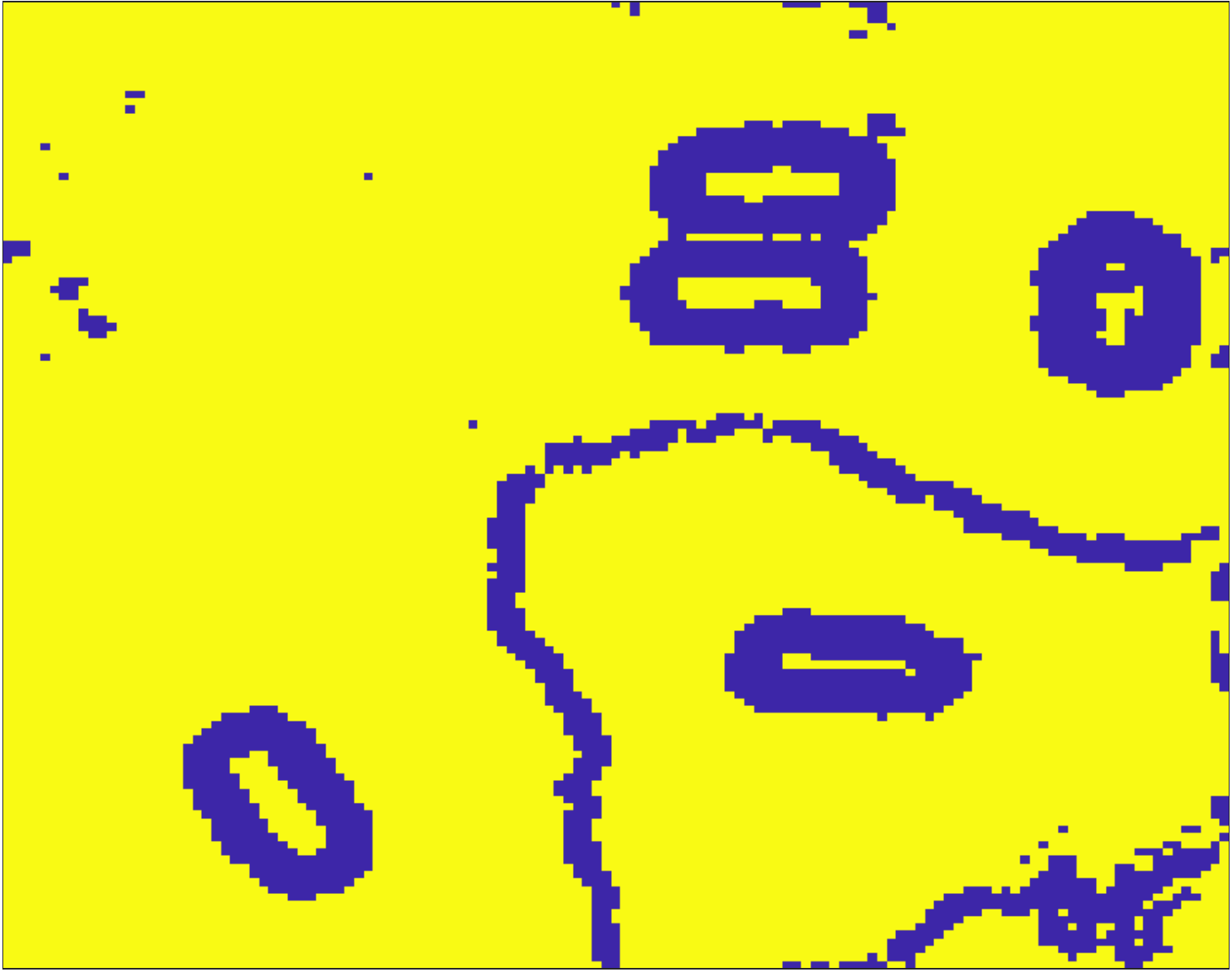}
    \caption{$W_4$}
\end{subfigure}
\\
\begin{subfigure}[b]{0.23\textwidth}
    \includegraphics[width=\textwidth]{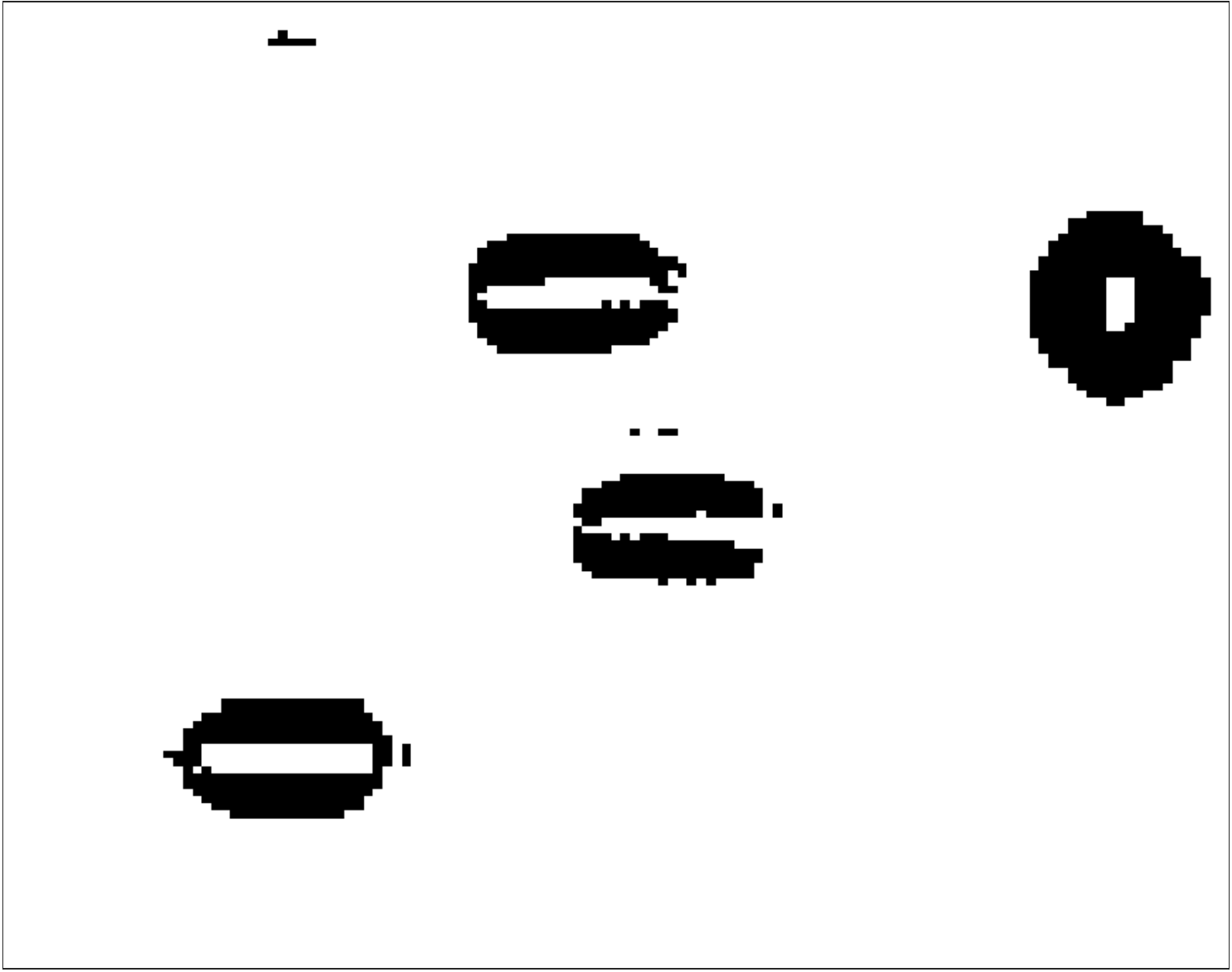}
    \caption{${\tilde U}_1$}
\end{subfigure}
~
\begin{subfigure}[b]{0.23\textwidth}
    \includegraphics[width=\textwidth]{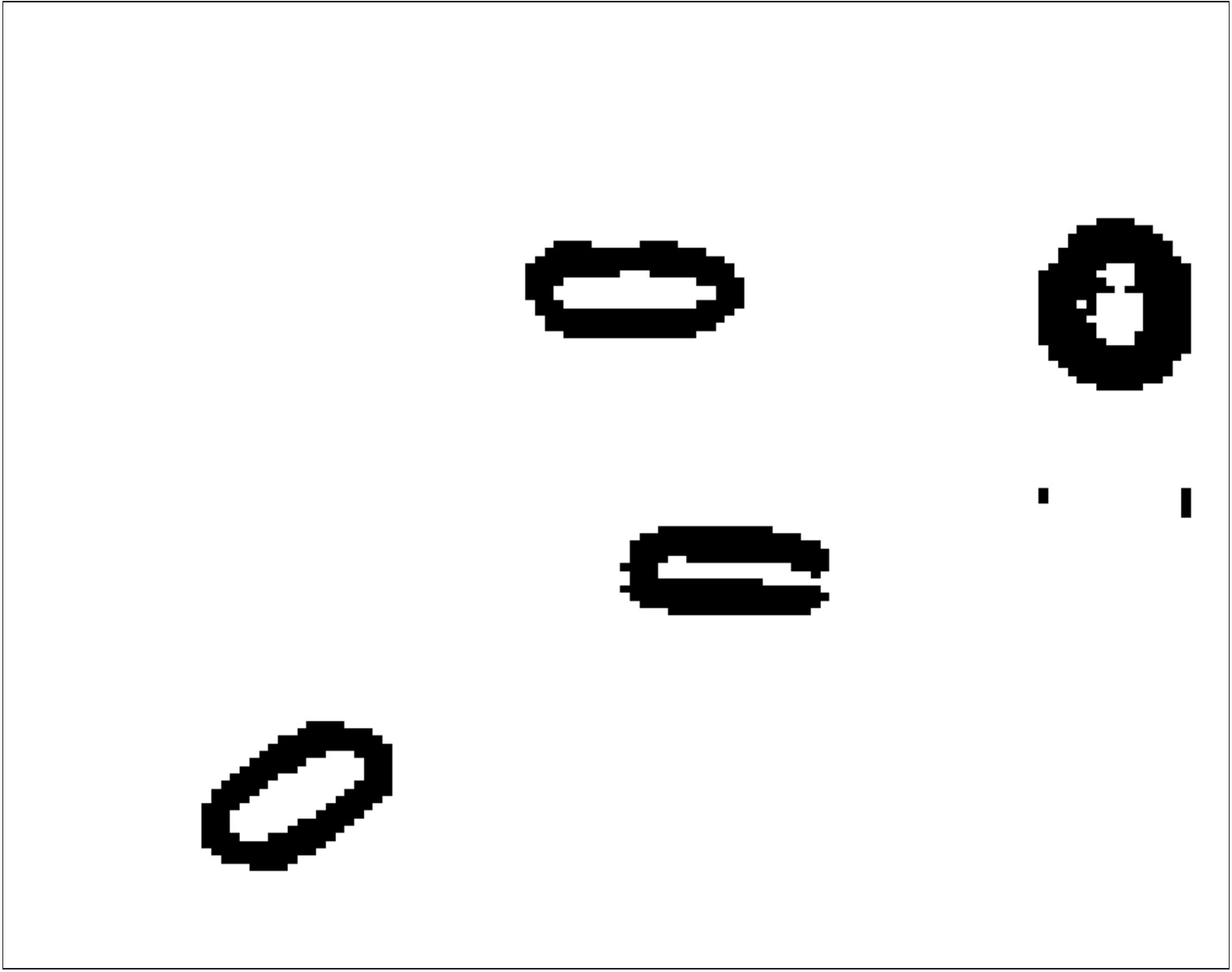}
    \caption{${\tilde U}_2$}
\end{subfigure}
~
\begin{subfigure}[b]{0.23\textwidth}
    \includegraphics[width=\textwidth]{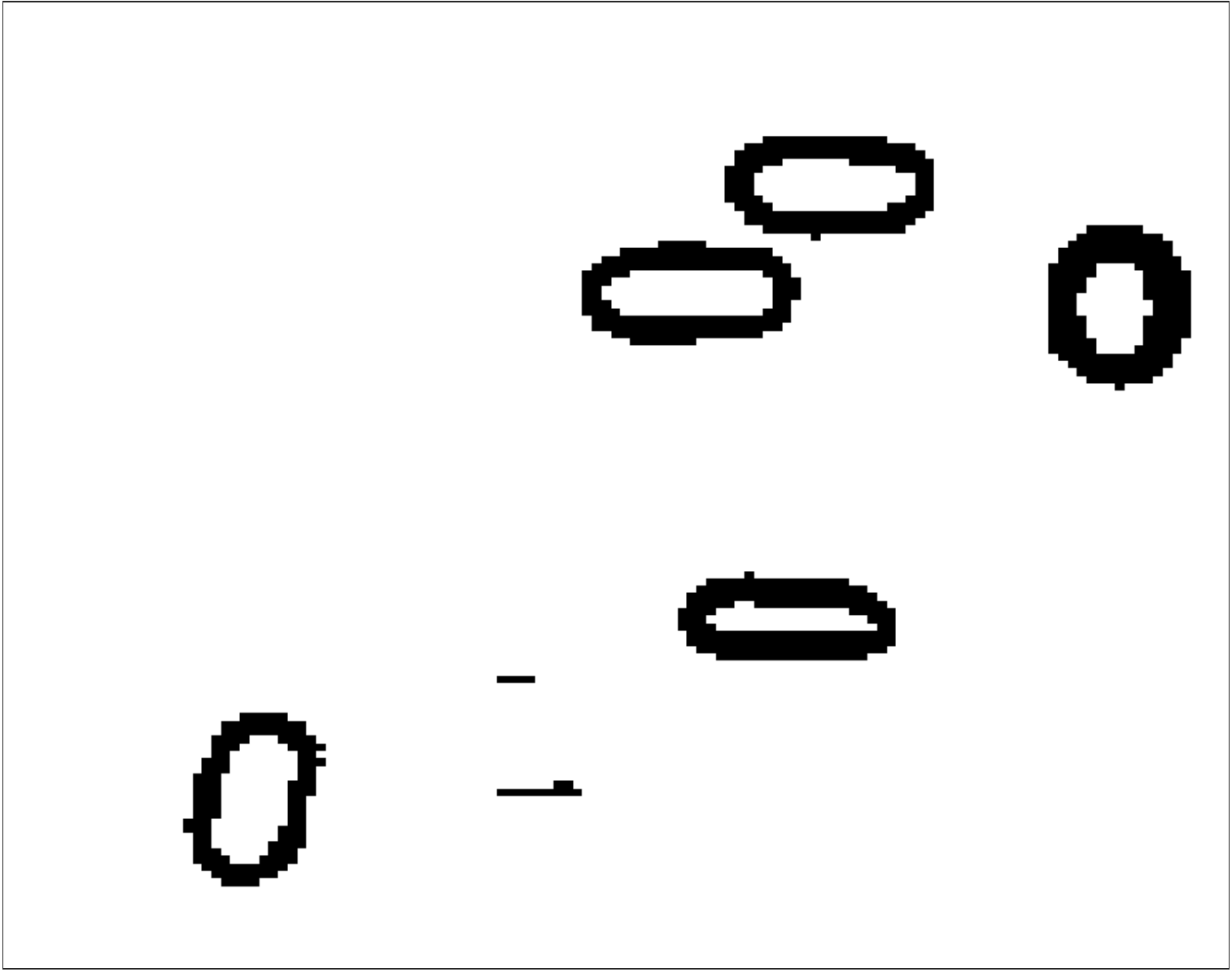}
    \caption{${\tilde U}_3$}
\end{subfigure}
~
\begin{subfigure}[b]{0.23\textwidth}
    \includegraphics[width=\textwidth]{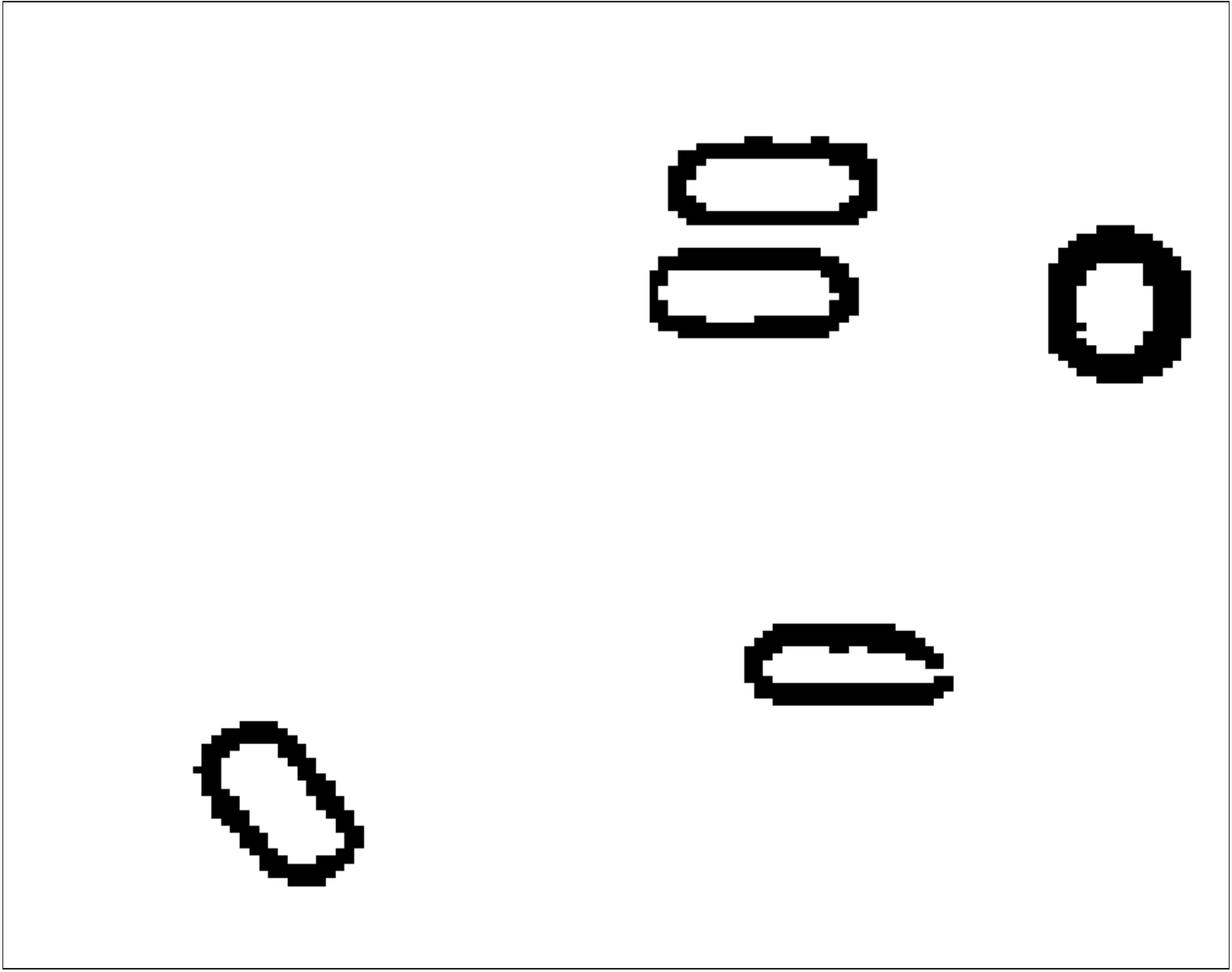}
    \caption{${\tilde U}_4$}
\end{subfigure}
\\
\begin{subfigure}[b]{0.23\textwidth}
    \includegraphics[width=\textwidth]{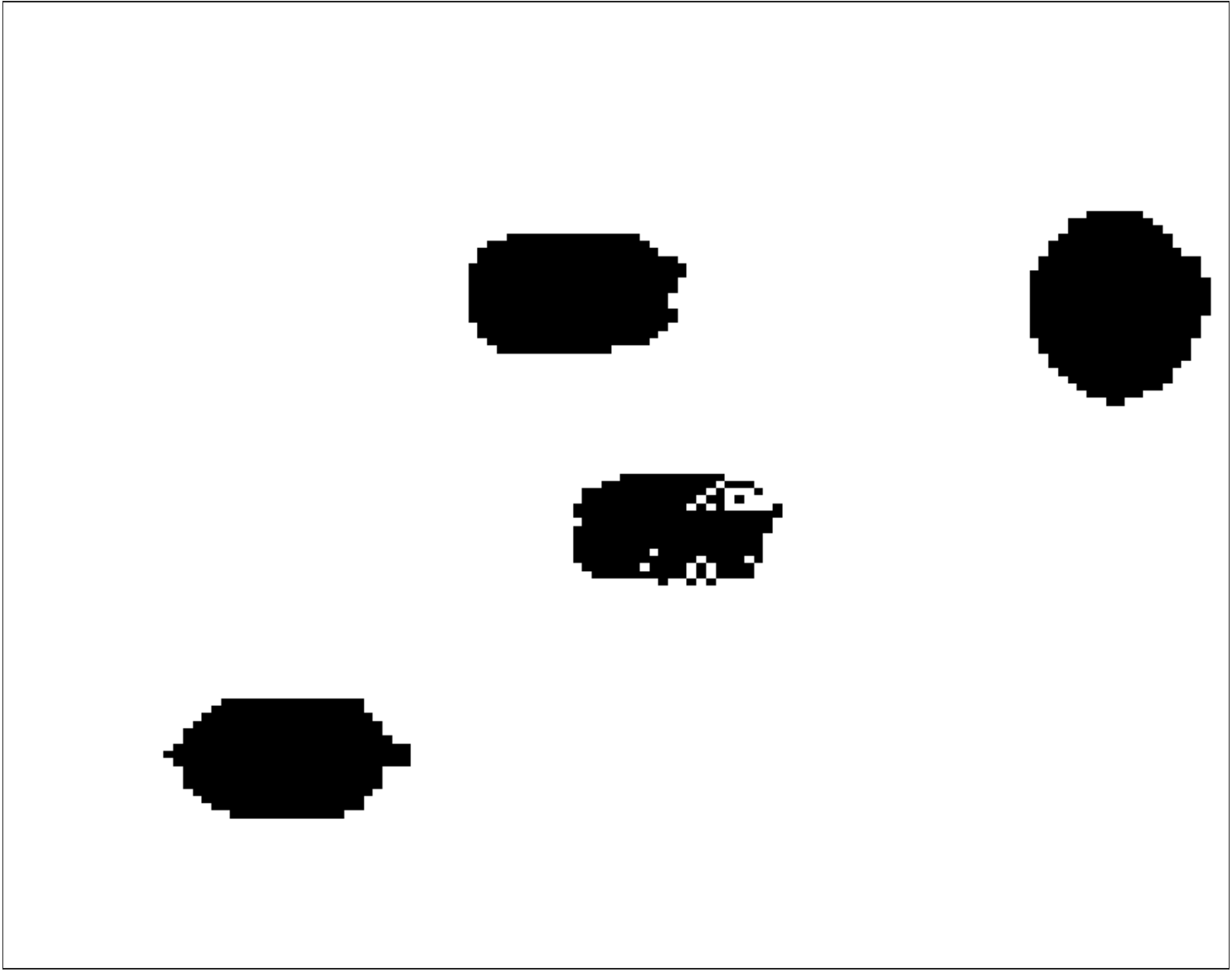}
    \caption{${\tilde Q}_1$}
\end{subfigure}
~
\begin{subfigure}[b]{0.23\textwidth}
    \includegraphics[width=\textwidth]{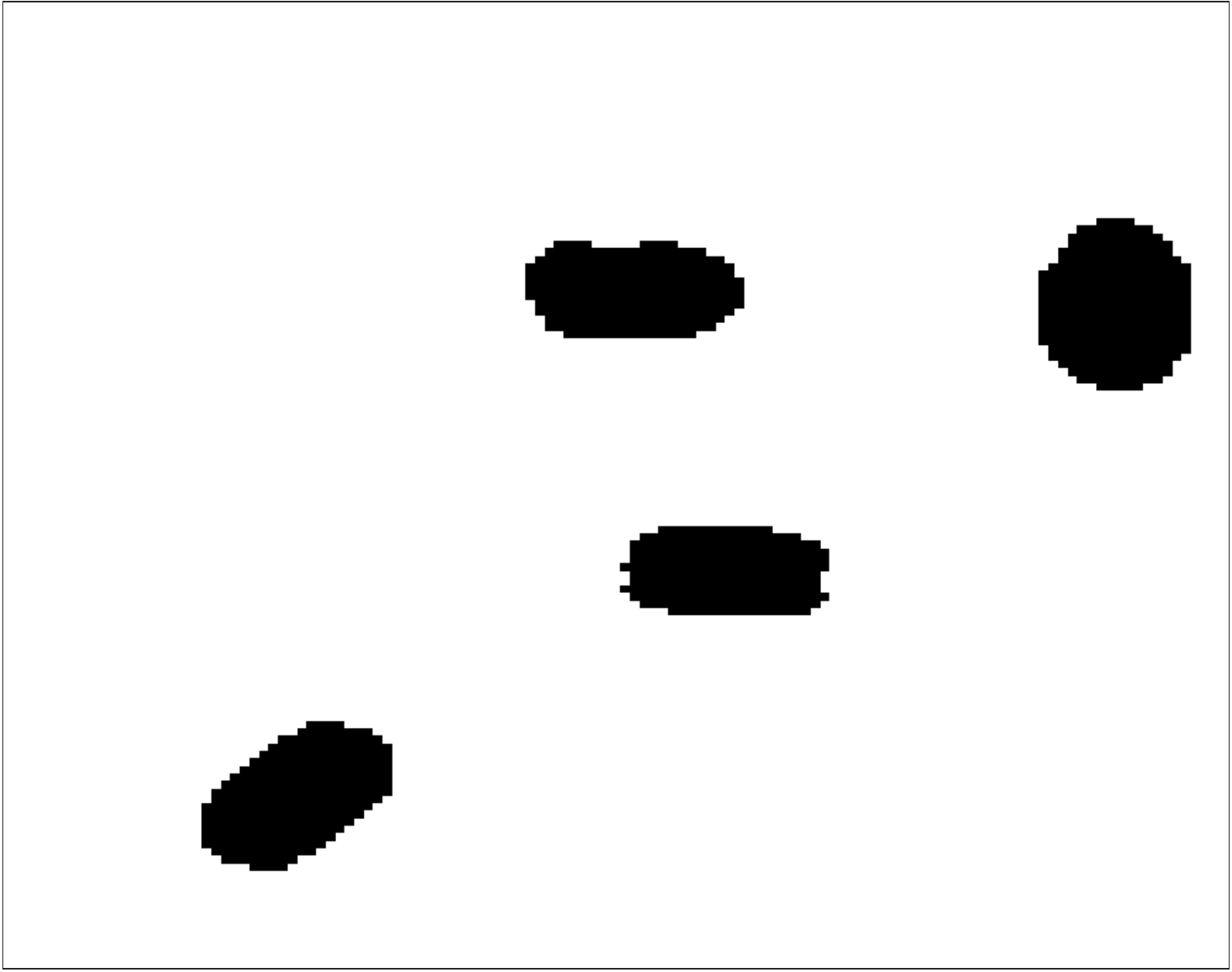}
    \caption{${\tilde Q}_2$}
\end{subfigure}
~
\begin{subfigure}[b]{0.23\textwidth}
    \includegraphics[width=\textwidth]{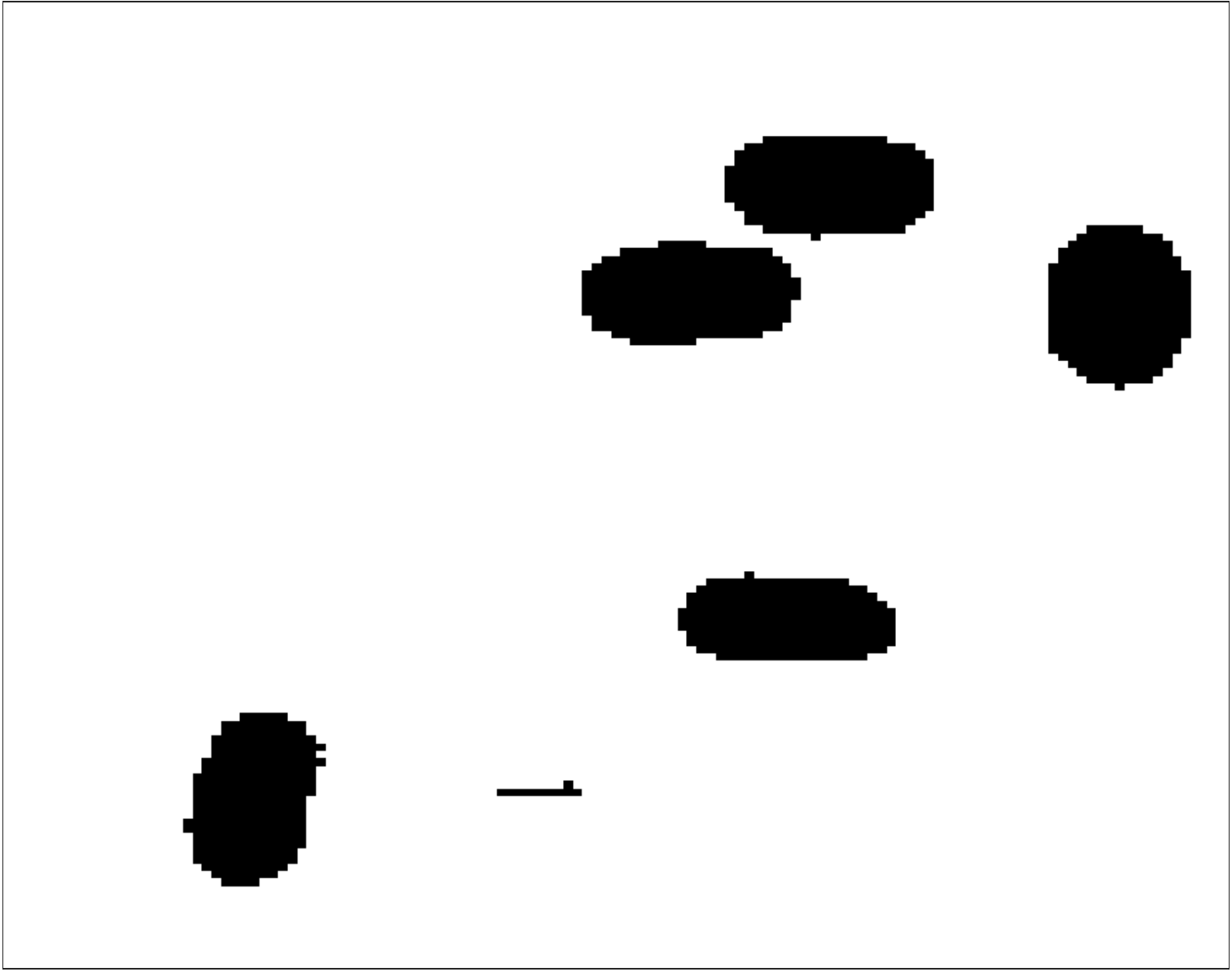}
    \caption{${\tilde Q}_3$}
\end{subfigure}
~
\begin{subfigure}[b]{0.23\textwidth}
    \includegraphics[width=\textwidth]{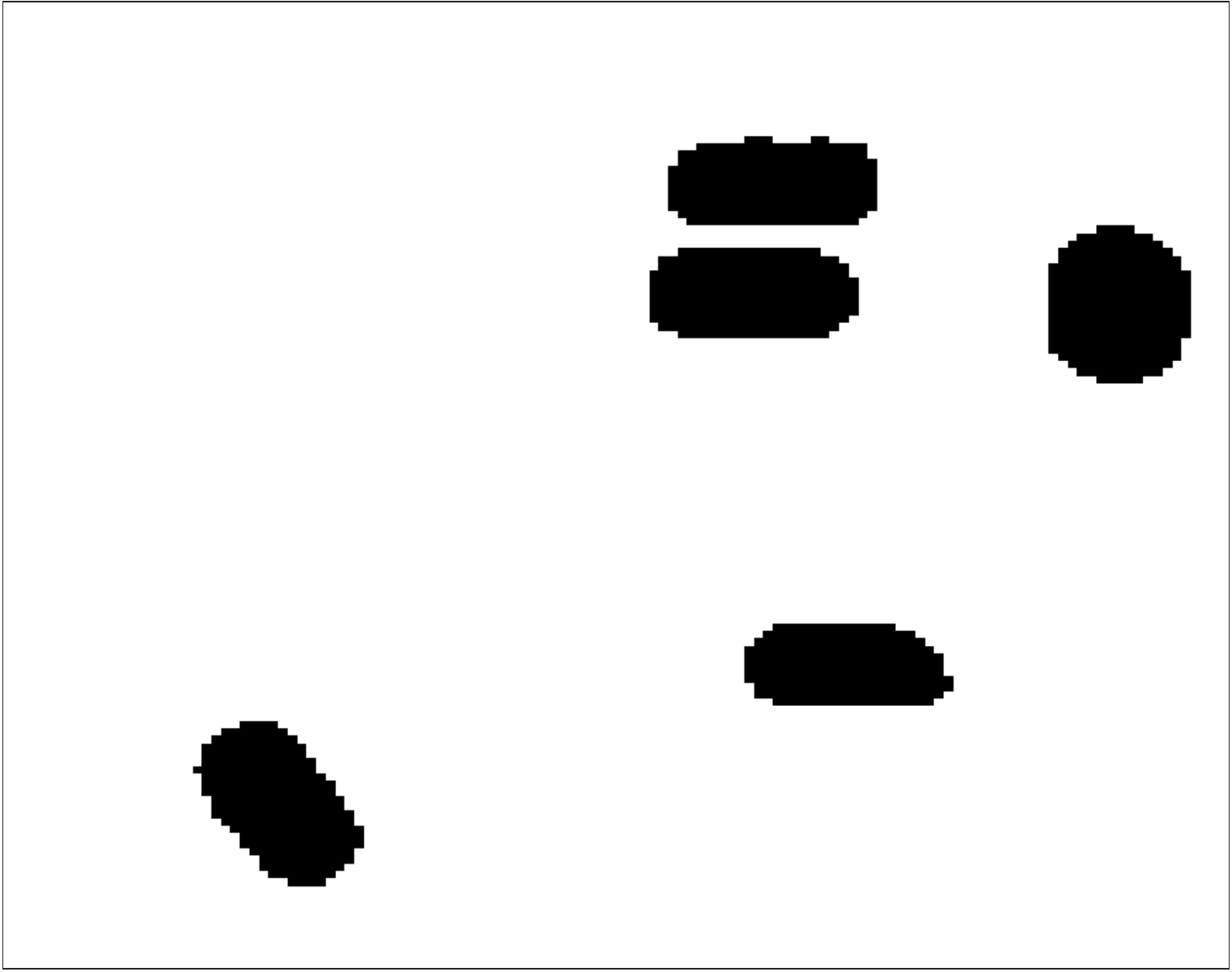}
    \caption{${\tilde Q}_4$}
\end{subfigure}
\\
\begin{subfigure}[b]{0.23\textwidth}
    \includegraphics[width=\textwidth]{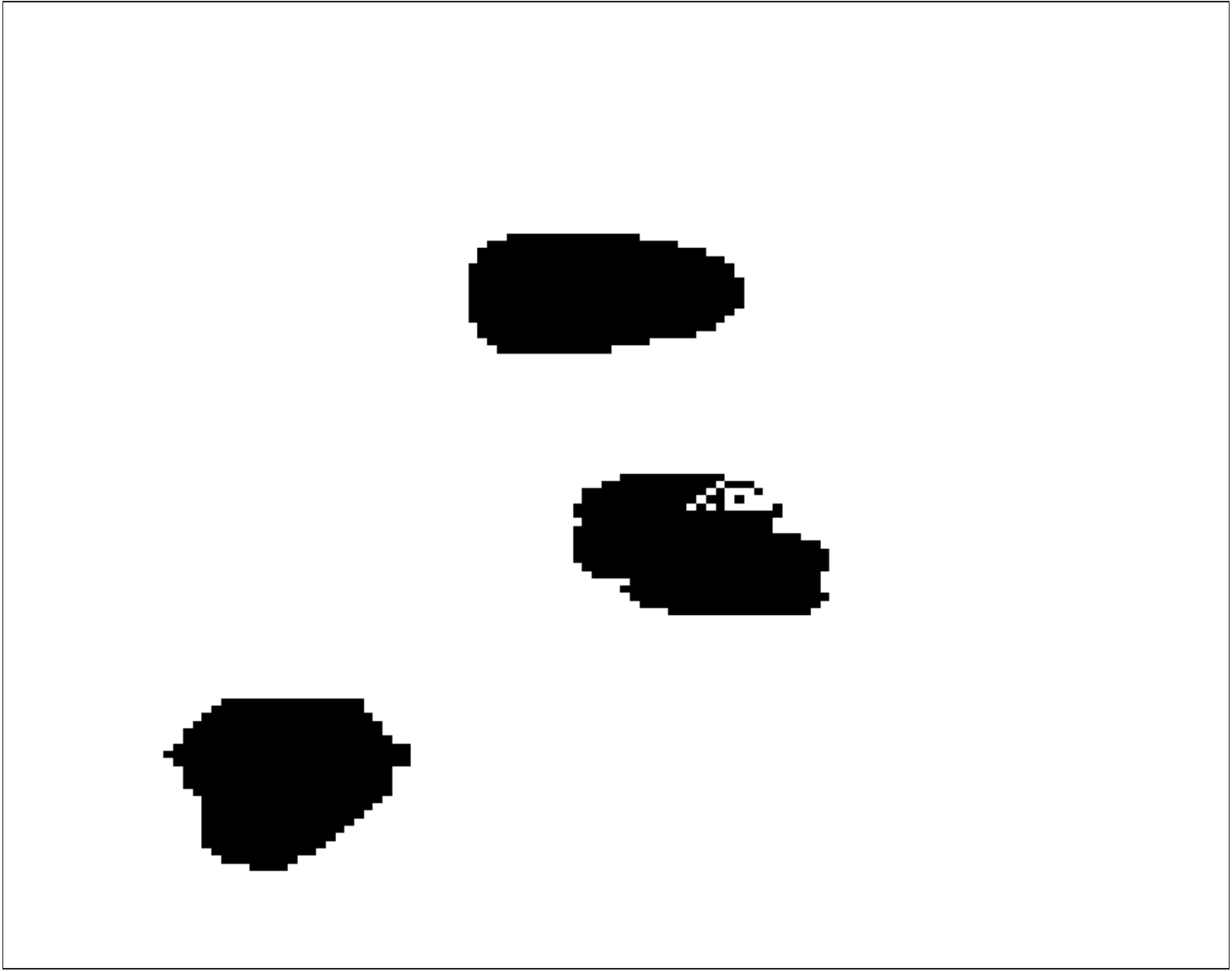}
    \caption{${\tilde{C}}_1$}
\end{subfigure}
~
\begin{subfigure}[b]{0.23\textwidth}
    \includegraphics[width=\textwidth]{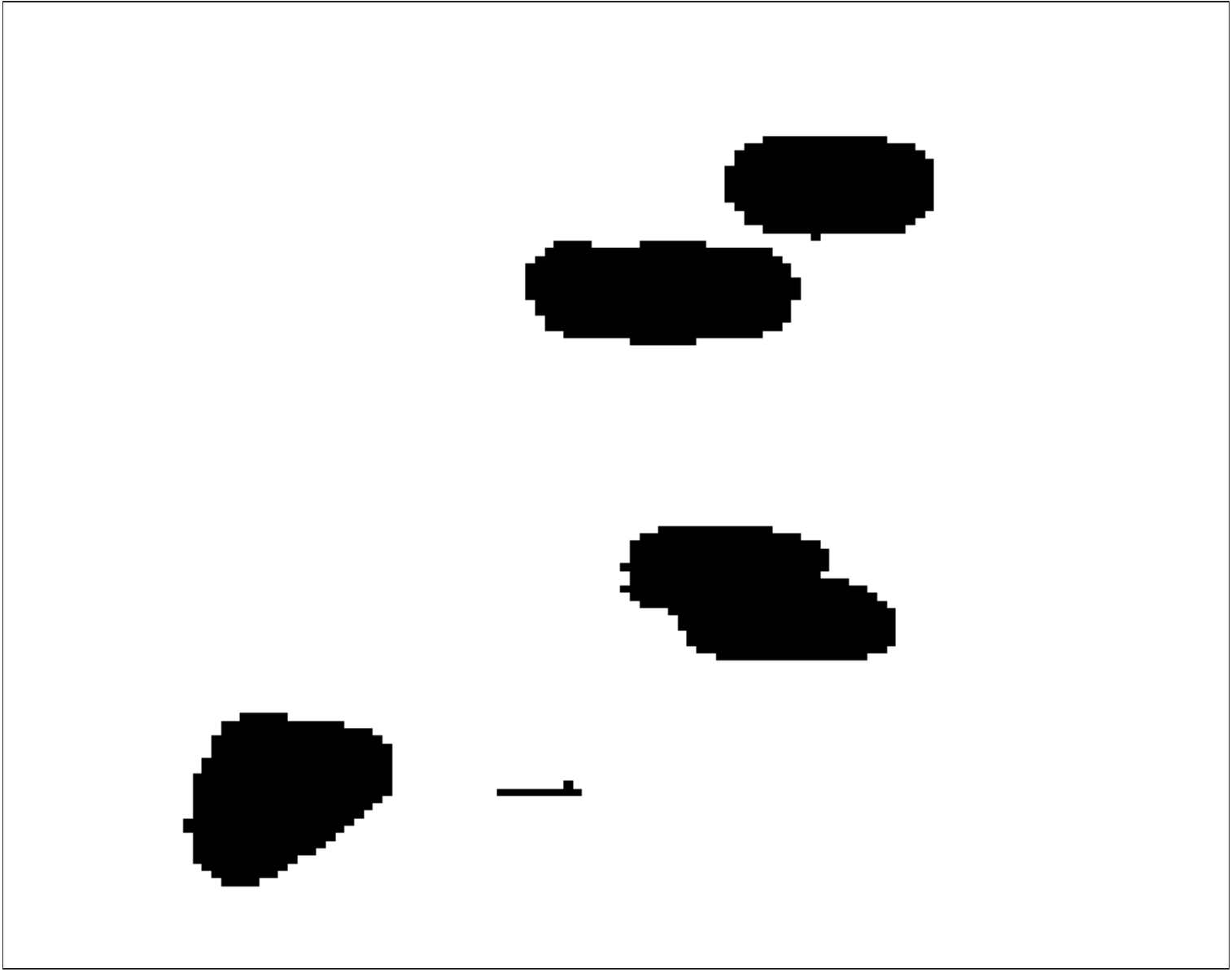}
    \caption{${\tilde{C}}_2$}
\end{subfigure}
~
\begin{subfigure}[b]{0.23\textwidth}
    \includegraphics[width=\textwidth]{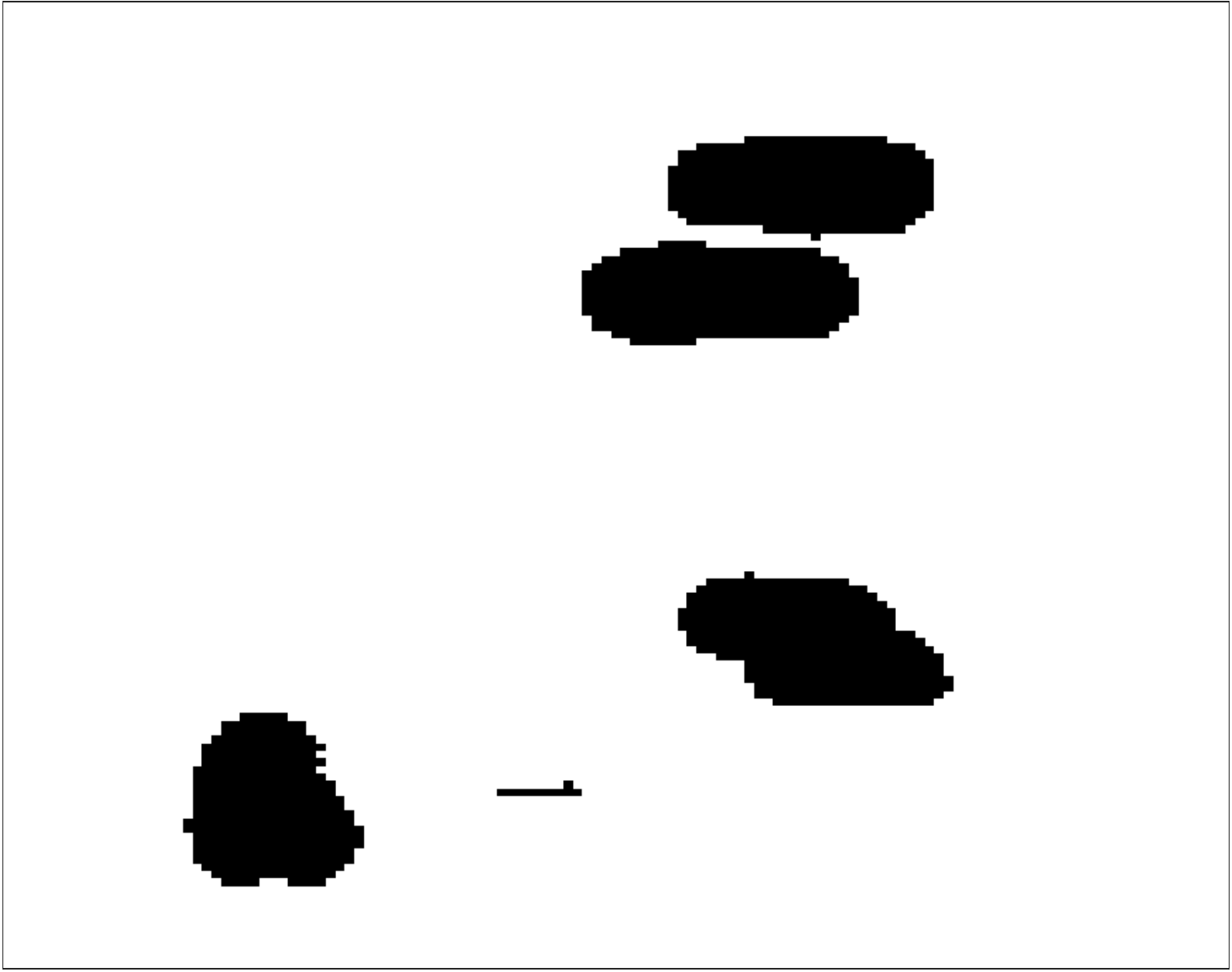}
    \caption{${\tilde{C}}_3$}
\end{subfigure}
~
\begin{subfigure}[b]{0.23\textwidth}
    \includegraphics[width=\textwidth]{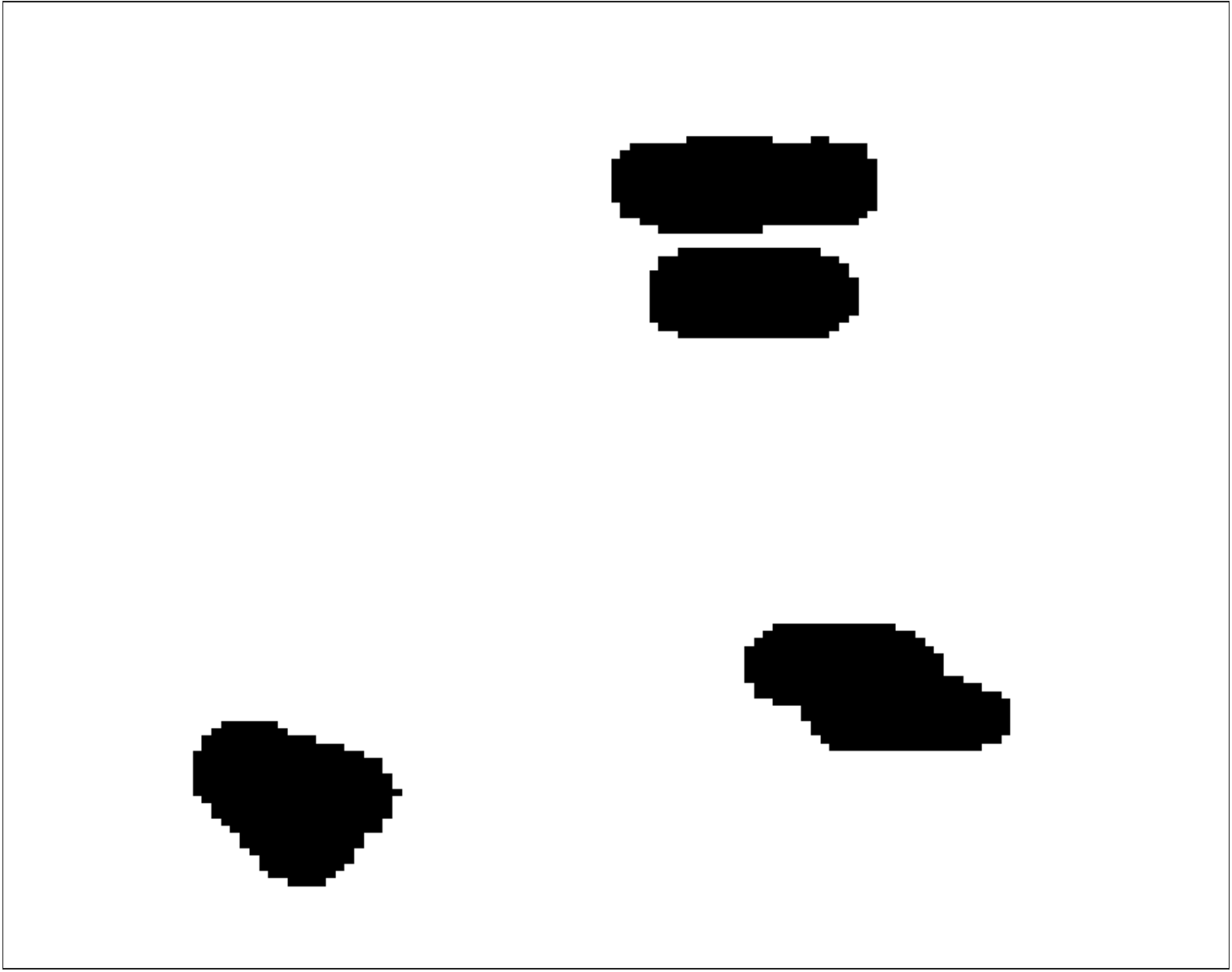}
    \caption{${\tilde{C}}_4$}
\end{subfigure}
\caption{(first {row}) Weighting mask $W$ calculated for each temporal image using \eqref{eq:weight_l1}. (second {row}) Partial edge masks ${\tilde U}$, \eqref{eq:ucurve}. (third {row}) Edge region masks ${\tilde Q}$, \eqref{eq:edge-region-mask}.  (fourth {row}) Change mask ${\tilde  C}$, \eqref{eq:CM-edge-region-mask}. Observe from \eqref{eq:CM-edge-region-mask} that $\tilde{C}_4$ requires construction of $W_5$, $\tilde{U}_5$, and $\tilde{Q}_5$, which are not pictured.}
\label{fig:phy_golf}
\end{figure}

As before, we also extract edge information {\em directly} from the given Fourier data \eqref{eq:forwardmodel_j}, although as previously noted  in this case we start from discrete data.  We then construct $M=10$ edge masks for each rotation angle $\theta$ in \eqref{eq:rotation-angles} and then follow Steps \ref{item:C1} - \ref{item:C4} to construct the inter-image matrix $\Phi$ used in \eqref{eq:optModel}.   Figure \ref{fig:phy_golf} illustrates each step of this procedure.

\begin{figure}[h!]
\centering

\begin{subfigure}[b]{.24\textwidth}
    \includegraphics[width=\textwidth]{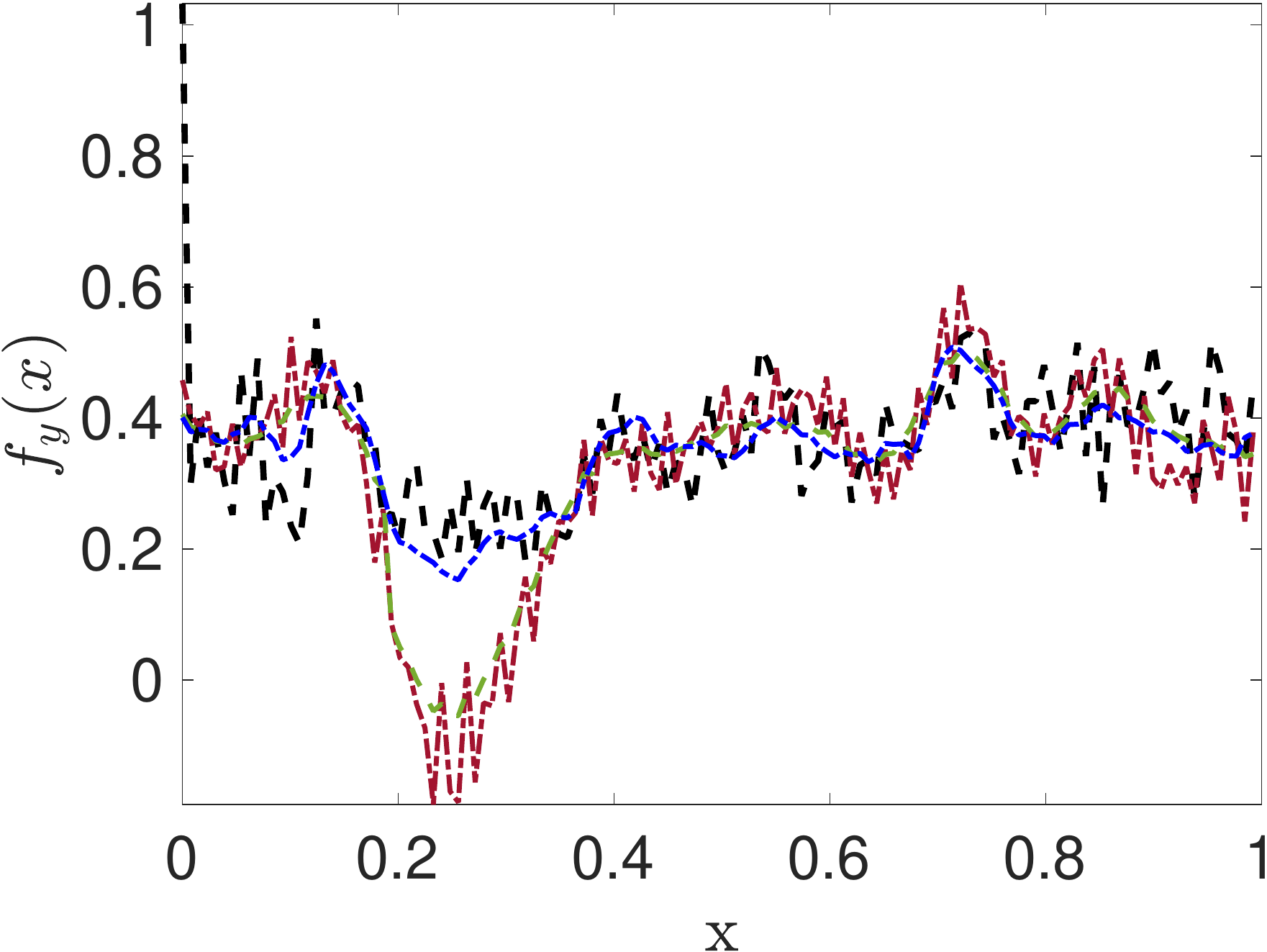}
    \caption{$\mathbf{\tilde f}_1$, $y=0.9070$}
\end{subfigure}
\begin{subfigure}[b]{.24\textwidth}
    \includegraphics[width=\textwidth]{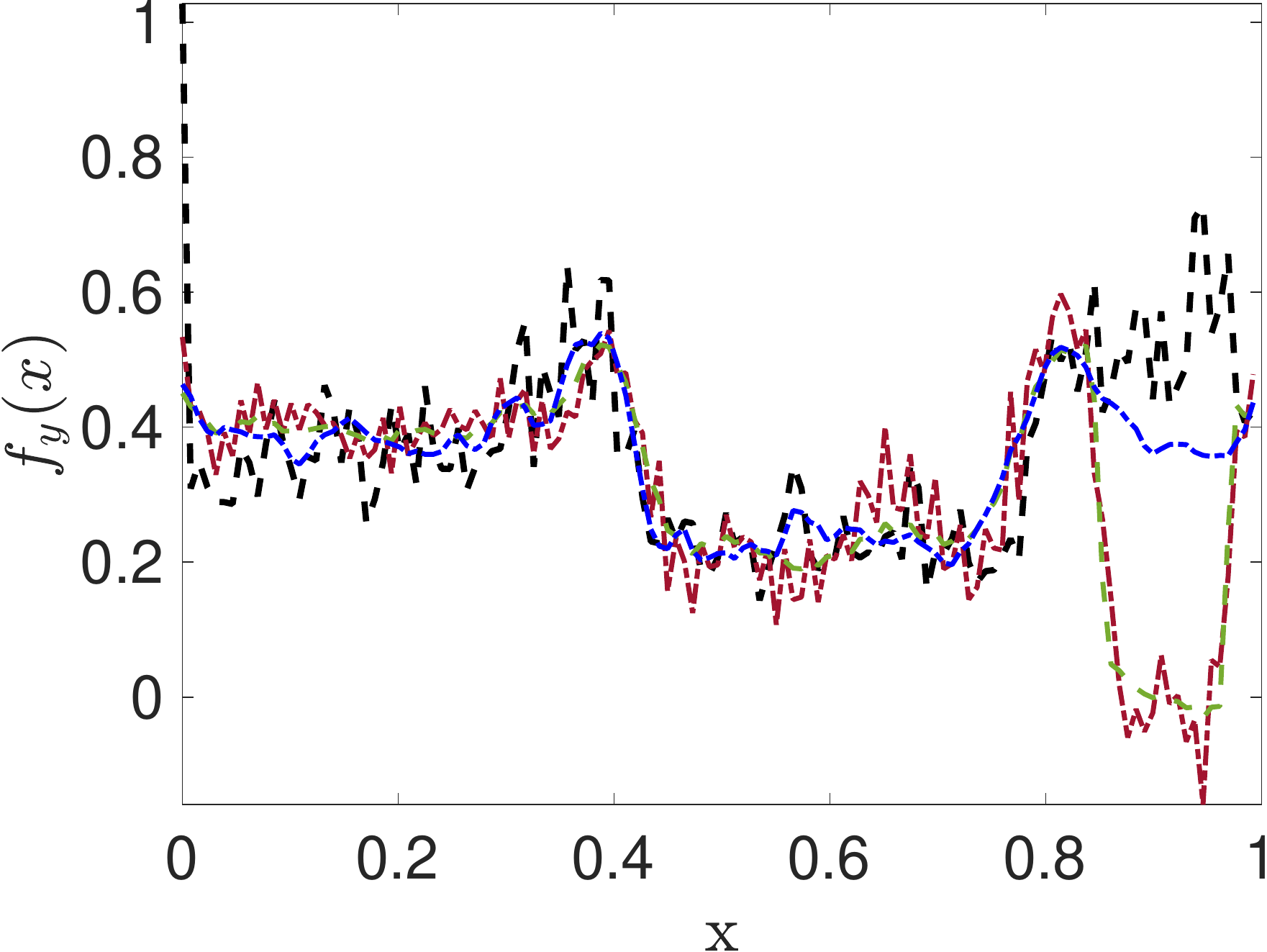}
        \caption{$\mathbf{\tilde f}_2$, $y=0.4806$}
\end{subfigure}
\begin{subfigure}[b]{.24\textwidth}
    \includegraphics[width=\textwidth]{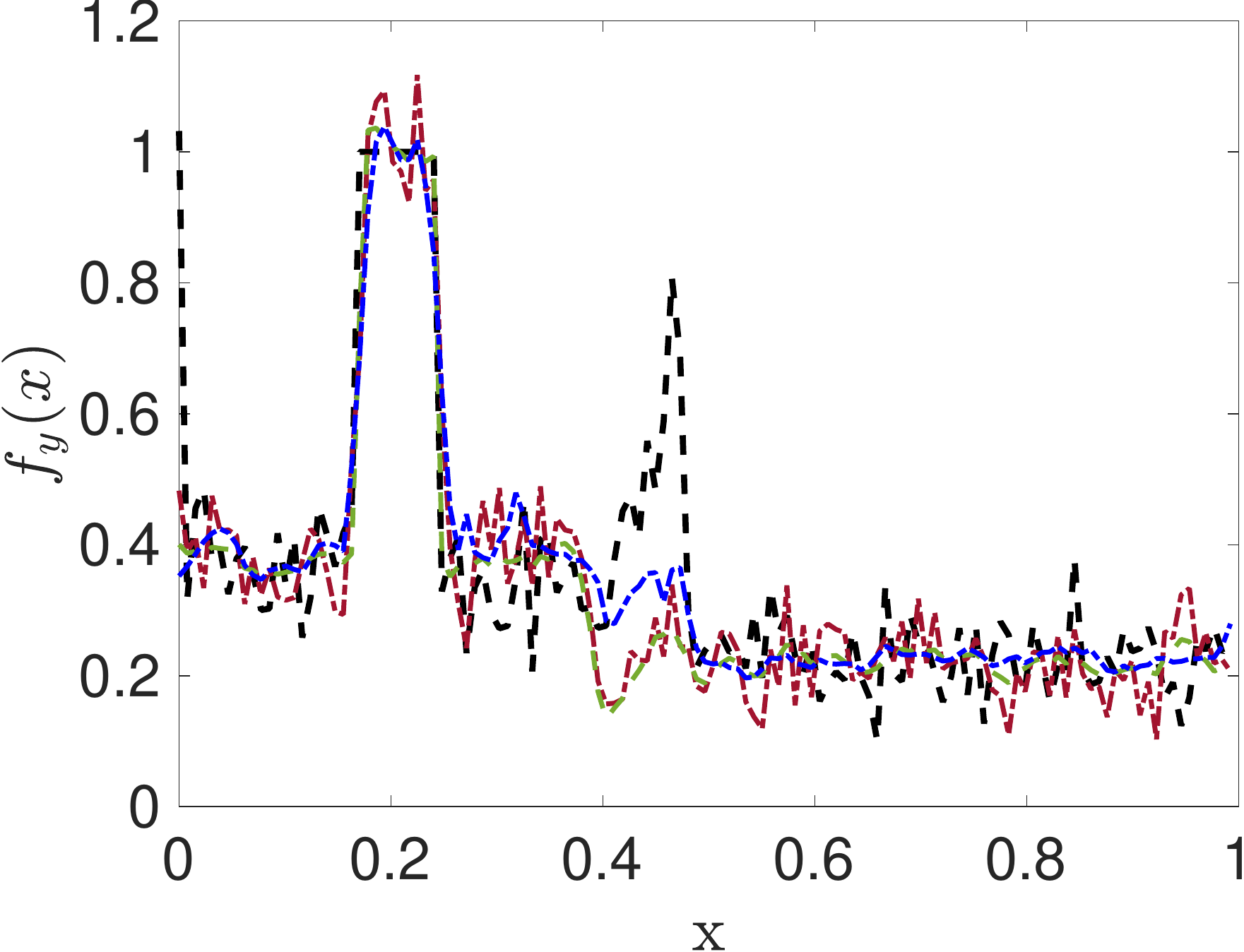}
    \caption{$\mathbf{\tilde f}_3$, $y=0.2403$}
\end{subfigure}
\begin{subfigure}[b]{.24\textwidth}
    \includegraphics[width=\textwidth]{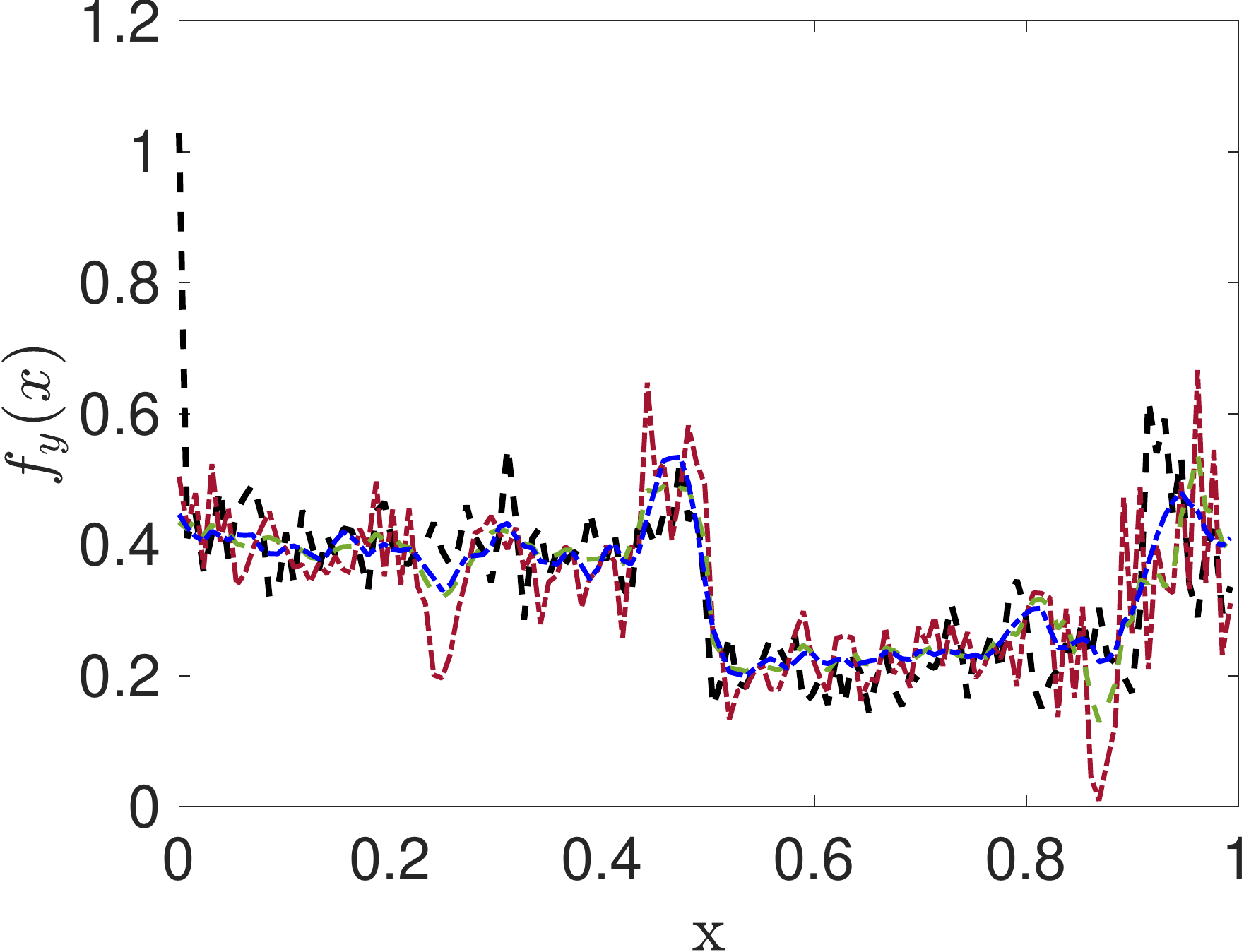}
    \caption{$\mathbf{\tilde f}_4$, $y=0.0930$}
\end{subfigure}
\begin{subfigure}[b]{.24\textwidth}
    \includegraphics[width=\textwidth]{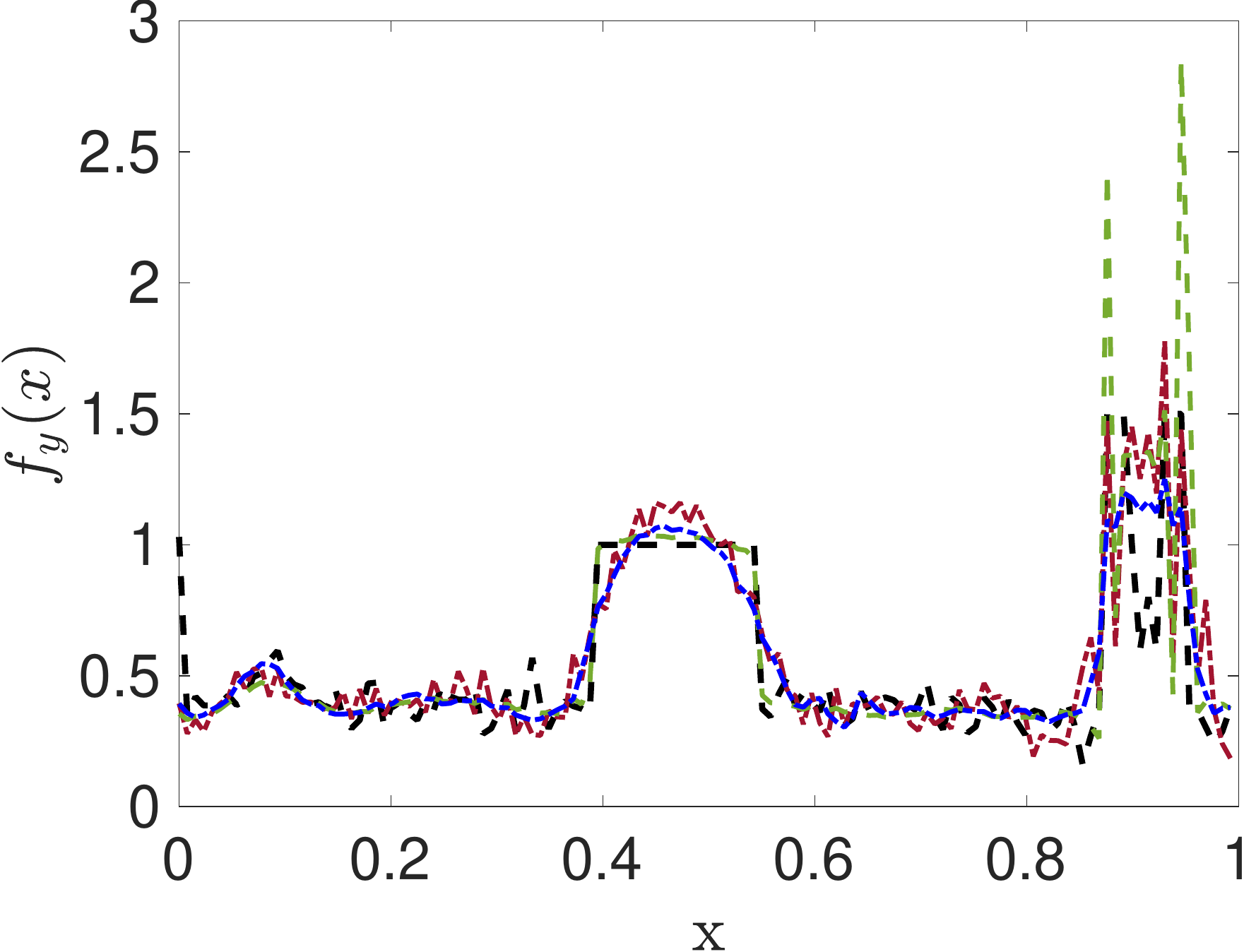}
    \caption{$\mathbf{\tilde f}_1$}
\end{subfigure}
\begin{subfigure}[b]{.24\textwidth}
    \includegraphics[width=\textwidth]{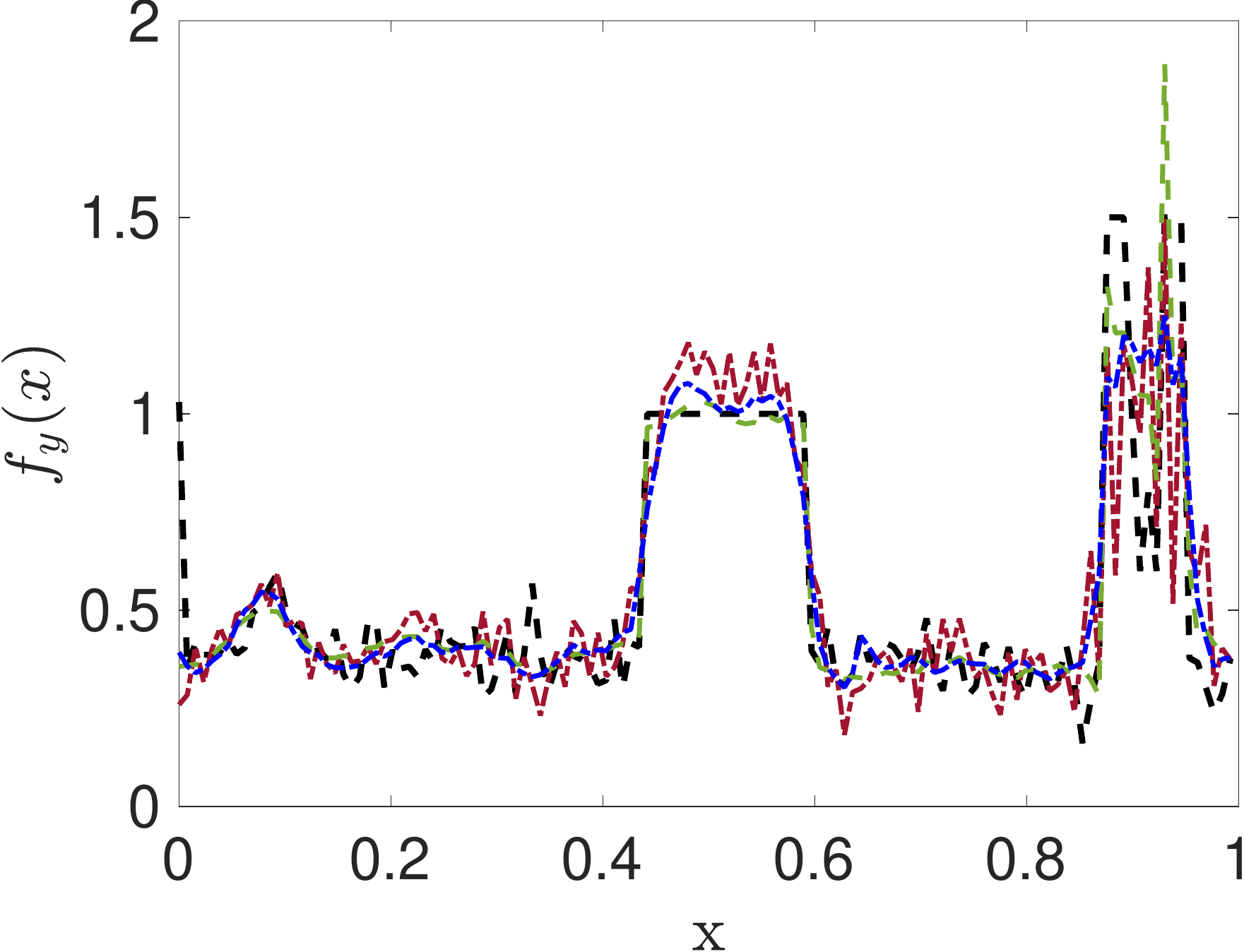}
    \caption{$\mathbf{\tilde f}_2$}
\end{subfigure}
\begin{subfigure}[b]{.24\textwidth}
    \includegraphics[width=\textwidth]{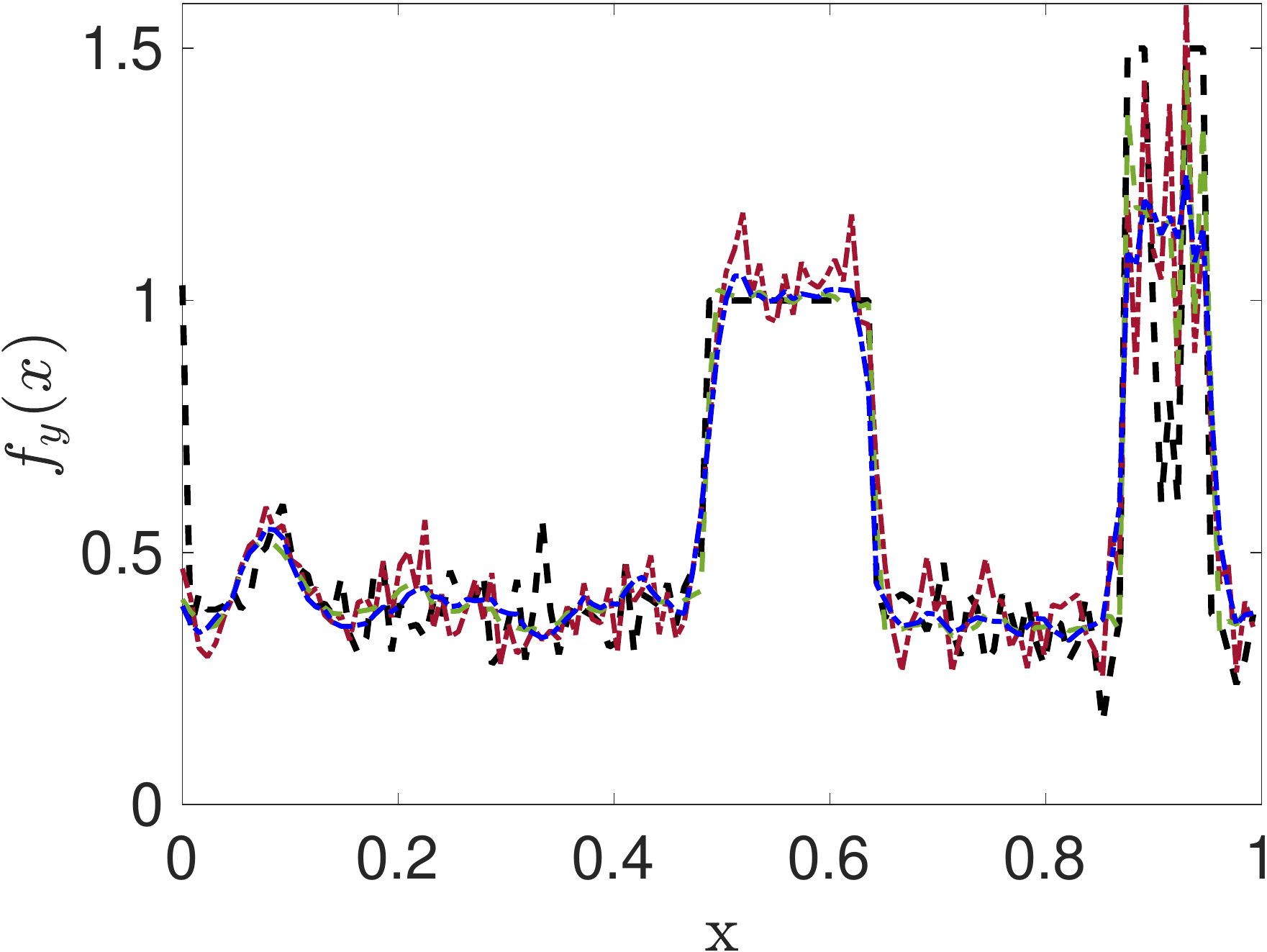}
    \caption{$\mathbf{\tilde f}_3$}
\end{subfigure}
\begin{subfigure}[b]{.24\textwidth}
    \includegraphics[width=\textwidth]{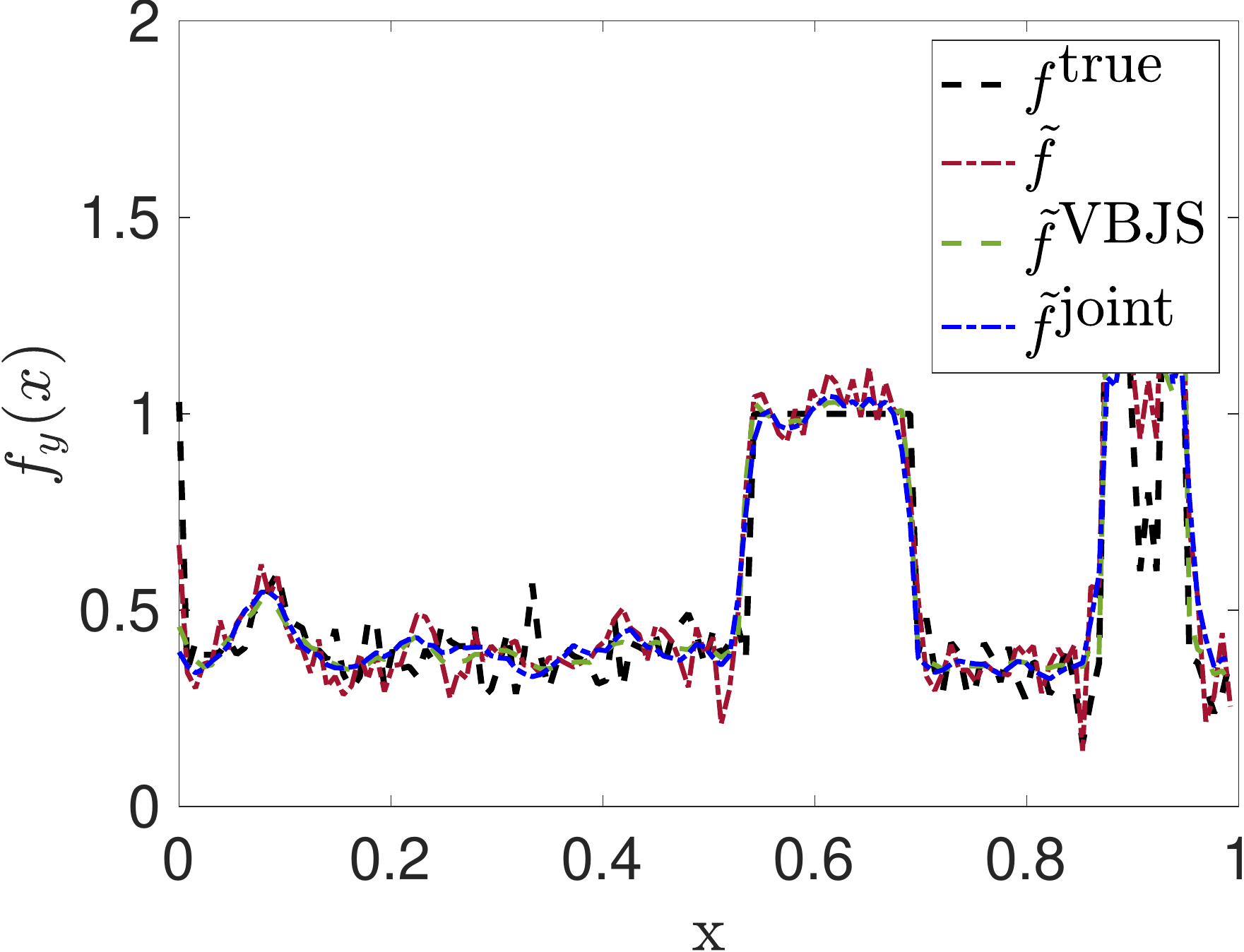}
    \caption{$\mathbf{\tilde f}_4$}
\end{subfigure}

\caption{Horizontal cross sections intersecting obstructed regions {in Figure \ref{fig:rec_golf} (top row)}. Horizontal cross sections at $y = .6899$, which  intersects the moving ellipses in {Figure \ref{fig:rec_golf} (bottom row)}. The legend in the bottom right plot is shared in other subplots. Here $f^{true}$ refers to the original SAR image with the addition of the boats and cars, since the ground truth is not available.}
\label{fig:1d_cross_unmoved_golf}
\end{figure}

One dimensional cross sections comparing the performance of each method are displayed  in Figure \ref{fig:1d_cross_unmoved_golf}. The horizontal slice is chosen to intersect an obstacle in each image of the temporal sequence (first row).  The effects of occlusions is apparently similar to what is observed for the MRI phantom in Section \ref{subsec:phantom}, with the joint recovery method being able to effectively recover information in obstructed regions.  The bottom row compares  the image recoveries at cross sections that intersect moving objects.  In this case there should be no contribution from the coupled term,  which implies that the results for VBJS and our new method should be similar.  Error comparisons cannot be made since there is no ground truth image available.

%% file: 5_summary.tex
\section{Concluding Remarks} 
\label{sec:summary}

In this work we developed a new technique that jointly recovers a set of images from a temporal sequence of under-sampled and noisy Fourier data.  There are furthermore physical obstructions that prevent the data acquisitions from ``seeing'' parts of each underlying image in the sequence.  Our method combines both intra- and inter-image information to enhance the recovery for each image.  In particular, the method builds on the  VBJS method, developed in \cite{adcock2019jointsparsity,gelb2019reducing,scarnati2019accelerated}, which effectively uses {\em intra-}image information to resolve each static image, with a coupling term that incorporates {\em inter-}image information from neighboring images.  This coupling term is important because it allows vital missing information to be borrowed from other data acquisitions in the temporal sequence.  By design, the coupling term requires an accurate depiction of shared information, which is accomplished by generating highly resolved edge masks for each image in the sequence, followed by  Steps \ref{item:C1}-\ref{item:C4}  to construct the inter-sequence ``change masks'', which are then explicitly placed in the coupling term.

Our numerical experiments demonstrate that our new joint recovery algorithm yields improved accuracy over both the standard $\ell_1$ regularization technique as well as the VBJS method, neither of which incorporates information from other data sets. A main advantage of our approach is that we are able to avoid significant data loss by directly using the given Fourier data to determine the edge masks needed both for high resolution individual recovery as well as an accurate coupling term.

A Fourier transform forward model served as a prototype in this investigation.   This is not a limitation for our method, however,  and in future investigations we will consider other models as well as data fusion problems.  In this regard our initial investigations indicate that it is still useful to use the Fourier based concentration factor method (after employing the discrete Fourier transform), although other high order edge detection methods, for example the polynomial annihilation edge detection method designed in \cite{PAGelb} and employed in \cite{adcock2019jointsparsity,gelb2019reducing,glaubitz2019high}, may also be effective at preventing data loss. 

%% file: sequential2022_article.bbl
\begin{thebibliography}{10}

\bibitem{adcock2019jointsparsity}
{\sc B.~Adcock, A.~Gelb, G.~Song, and Y.~Sui}, {\em Joint sparse recovery based
  on variances}, SIAM Journal of Scientific Computing, 41 (2019),
  pp.~A246--A268, \url{https://doi.org/10.1137/17M1155983}.

\bibitem{afaq2021analysis}
{\sc Y.~Afaq and A.~Manocha}, {\em Analysis on change detection techniques for
  remote sensing applications: A review}, Ecological Informatics,  (2021),
  p.~101310.

\bibitem{archibald2003improving}
{\sc R.~Archibald, K.~Chen, A.~Gelb, and R.~Renaut}, {\em Improving tissue
  segmentation of human brain {MRI} through preprocessing by the {G}egenbauer
  reconstruction method}, NeuroImage, 20 (2003), pp.~489--502.

\bibitem{archibald2002method}
{\sc R.~Archibald and A.~Gelb}, {\em A method to reduce the {G}ibbs ringing
  artifact in {MRI} scans while keeping tissue boundary integrity}, IEEE
  Transactions on Medical Imaging, 21 (2002), pp.~305--319.

\bibitem{ArchibaldGelb}
{\sc R.~Archibald, A.~Gelb, and R.~Platte}, {\em Image reconstruction from
  undersampled {F}ourier data using the polynomial annihilation transform},
  Journal of Scientific Computing,  (2015),
  \url{https://doi.org/10.1007/s10915-015-0088-2}.

\bibitem{PAGelb}
{\sc R.~Archibald, A.~Gelb, and J.~Yoon}, {\em Polynomial fitting for edge
  detection in irregularly sampled signals and images}, SIAM Journal on
  Numerical Analysis, 43 (2005), pp.~259--279.

\bibitem{Ash_2014}
{\sc J.~N. Ash}, {\em {A unifying perspective of coherent and non-coherent
  change detection}}, in Algorithms for Synthetic Aperture Radar Imagery XXI,
  E.~Zelnio and F.~D. Garber, eds., vol.~9093, International Society for Optics
  and Photonics, SPIE, 2014, pp.~90 -- 98,
  \url{https://doi.org/10.1117/12.2054338}.

\bibitem{boyd2011distributed}
{\sc S.~Boyd, N.~Parikh, and E.~Chu}, {\em Distributed optimization and
  statistical learning via the alternating direction method of multipliers},
  (2011), \url{https://doi.org/10.1561/2200000016}.

\bibitem{candes2006robust}
{\sc E.~J. Cand{\`e}s, J.~Romberg, and T.~Tao}, {\em Robust uncertainty
  principles: Exact signal reconstruction from highly incomplete frequency
  information}, IEEE Transactions on Information Theory, 52 (2006),
  pp.~489--509.

\bibitem{candes2006stable}
{\sc E.~J. Cand{\`e}s, J.~K. Romberg, and T.~Tao}, {\em Stable signal recovery
  from incomplete and inaccurate measurements}, Communications on Pure and
  Applied Mathematics: A Journal Issued by the Courant Institute of
  Mathematical Sciences, 59 (2006), pp.~1207--1223.

\bibitem{candes2006near}
{\sc E.~J. Cand{\`e}s and T.~Tao}, {\em Near-optimal signal recovery from
  random projections: Universal encoding strategies?}, IEEE Transactions on
  Information Theory, 52 (2006), pp.~5406--5425.

\bibitem{candes2008enhancing}
{\sc E.~J. Cand{\`e}s, M.~B. Wakin, and S.~P. Boyd}, {\em Enhancing sparsity by
  reweighted $\ell_1$ minimization}, Journal of {F}ourier Analysis and
  Applications, 14 (2008), pp.~877--905.

\bibitem{canny1986computational}
{\sc J.~Canny}, {\em A computational approach to edge detection}, IEEE
  Transactions on Pattern, Analysis and Machine Intelligence,  (1986),
  pp.~679--698.

\bibitem{chartrand2008iteratively}
{\sc R.~Chartrand and W.~Yin}, {\em Iteratively reweighted algorithms for
  compressive sensing}, in IEEE International Conference on Acoustics, Speech
  and Signal Processing, 2008, pp.~3869--3872.

\bibitem{Chen2012object}
{\sc G.~Chen, G.~J. Hay, L.~M.~T. Carvalho, and M.~A. Wulder}, {\em
  Object-based change detection}, International Journal of Remote Sensing, 33
  (2012), pp.~4434--4457, \url{https://doi.org/10.1080/01431161.2011.648285}.

\bibitem{CHEN2014assessment}
{\sc G.~Chen, K.~Zhao, and R.~Powers}, {\em Assessment of the image
  misregistration effects on object-based change detection}, ISPRS Journal of
  Photogrammetry and Remote Sensing, 87 (2014), pp.~19--27,
  \url{https://doi.org/10.1016/j.isprsjprs.2013.10.007}.

\bibitem{chen2017joint}
{\sc Y.~Chen, R.~Fang, and X.~Ye}, {\em Joint image edge reconstruction and its
  application in multi-contrast {MRI}},  (2017),
  \url{https://arxiv.org/abs/1712.02000}.

\bibitem{churchill2018edge}
{\sc V.~Churchill, R.~Archibald, and A.~Gelb}, {\em Edge-adaptive {$\ell_2$}
  regularization image reconstruction from non-uniform {F}ourier data},
  (2018), \url{https://arxiv.org/abs/1811.08487}.

\bibitem{daubechies2008iteratively}
{\sc I.~Daubechies, R.~DeVore, M.~Fornasier, and C.~S. Gunturk}, {\em
  Iteratively re-weighted least squares minimization for sparse recovery},
  Communications on Pure and Applied Mathematics: A Journal Issued by the
  Courant Institute of Mathematical Sciences, 63 (2010), pp.~1--38.

\bibitem{donoho2006compressed}
{\sc D.~L. Donoho}, {\em Compressed sensing}, IEEE Transactions on Information
  Theory, 52 (2006), pp.~1289--1306.

\bibitem{SAR_Image_ref}
{\sc M.~Ellsworth and C.~Thomas}, {\em A fast algorithm for image deblurring
  with total variation regularization}, Unmanned Tech Solutions, 4 (2014).

\bibitem{eslahi2019joint}
{\sc N.~{Eslahi} and A.~{Foi}}, {\em Joint sparse recovery of misaligned
  multimodal images via adaptive local and nonlocal cross-modal
  regularization}, in IEEE 8th International Workshop on Computational Advances
  in Multi-Sensor Adaptive Processing (CAMSAP), 2019, pp.~111--115,
  \url{https://doi.org/10.1109/CAMSAP45676.2019.9022478}.

\bibitem{gao2019extracting}
{\sc F.~Gao, M.~Wang, Y.~Cai, and S.~Lu}, {\em Extracting closed object contour
  in the image: remove, connect and fit}, Pattern Analysis and Applications, 22
  (2019), pp.~1123--1136.

\bibitem{gao2010improved}
{\sc W.~Gao, X.~Zhang, L.~Yang, and H.~Liu}, {\em An improved {S}obel edge
  detection}, in IEEE 3rd International Conference on Computer Science and
  Information Technology, vol.~5, 2010, pp.~67--71.

\bibitem{gelb2019reducing}
{\sc A.~Gelb and T.~Scarnati}, {\em Reducing effects of bad data using variance
  based joint sparsity recovery}, Journal of Scientific Computing, 78 (2019),
  pp.~94--120.

\bibitem{gelb2017detecting}
{\sc A.~Gelb and G.~Song}, {\em Detecting edges from non-uniform {F}ourier data
  using {F}ourier frames}, Journal of Scientific Computing, 71 (2017),
  pp.~737--758.

\bibitem{gelb1999detection}
{\sc A.~Gelb and E.~Tadmor}, {\em Detection of edges in spectral data}, Applied
  and Computational Harmonic Analysis, 7 (1999), pp.~101--135.

\bibitem{gelb2000detection}
{\sc A.~Gelb and E.~Tadmor}, {\em Detection of edges in spectral data {II}.
  {N}onlinear enhancement}, SIAM Journal on Numerical Analysis, 38 (2000),
  pp.~1389--1408.

\bibitem{glaubitz2019high}
{\sc J.~Glaubitz and A.~Gelb}, {\em High order edge sensors with $\ell^1$
  regularization for enhanced discontinuous {G}alerkin methods}, SIAM Journal
  on Scientific Computing, 41 (2019), pp.~A1304--A1330.

\bibitem{gong2019adaptive}
{\sc C.~Gong and L.~Zeng}, {\em Adaptive iterative reconstruction based on
  relative total variation for low-intensity computed tomography}, Signal
  Processing, 165 (2019), pp.~149--162.

\bibitem{Hsiao2005contour}
{\sc Y.-T. Hsiao, C.-L. Chuang, J.-A. Jiang, and C.-C. Chien}, {\em A contour
  based image segmentation algorithm using morphological edge detection}, in
  IEEE International Conference on Systems, Man and Cybernetics, vol.~3, 2005,
  pp.~2962--2967, \url{https://doi.org/10.1109/ICSMC.2005.1571600}.

\bibitem{HUSSAIN2013change}
{\sc M.~Hussain, D.~Chen, A.~Cheng, H.~Wei, and D.~Stanley}, {\em Change
  detection from remotely sensed images: From pixel-based to object-based
  approaches}, ISPRS Journal of Photogrammetry and Remote Sensing, 80 (2013),
  pp.~91--106, \url{https://doi.org/10.1016/j.isprsjprs.2013.03.006}.

\bibitem{inglada2007new}
{\sc J.~Inglada and G.~Mercier}, {\em A new statistical similarity measure for
  change detection in multitemporal {SAR} images and its extension to
  multiscale change analysis}, IEEE Transactions on Geoscience and Remote
  Sensing, 45 (2007), pp.~1432--1445.

\bibitem{ji2008bayesian}
{\sc S.~Ji, Y.~Xue, and L.~Carin}, {\em Bayesian compressive sensing}, IEEE
  Transactions on Signal Processing, 56 (2008), pp.~2346--2356.

\bibitem{jin2009edge}
{\sc Z.~Jin-Yu, C.~Yan, and H.~Xian-Xiang}, {\em Edge detection of images based
  on improved {S}obel operator and genetic algorithms}, in IEEE International
  Conference on Image Analysis and Signal Processing, 2009, pp.~31--35.

\bibitem{kang2019compressive}
{\sc M.-S. Kang and K.-T. Kim}, {\em Compressive sensing based {SAR} imaging
  and autofocus using improved {T}ikhonov regularization}, IEEE Sensors
  Journal, 19 (2019), pp.~5529--5540,
  \url{https://doi.org/10.1109/JSEN.2019.2904611}.

\bibitem{kazantsev2018joint}
{\sc D.~Kazantsev, J.~S. J{\o}rgensen, M.~S. Andersen, W.~R.~B. Lionheart,
  P.~D. Lee, and P.~J. Withers}, {\em Joint image reconstruction method with
  correlative multi-channel prior for {X}-ray spectral computed tomography},
  Inverse Problems, 34 (2018), p.~064001,
  \url{https://doi.org/10.1088/1361-6420/aaba86}.

\bibitem{kazantsev2014multimodal}
{\sc D.~Kazantsev, W.~R.~B. Lionheart, P.~J. Withers, and P.~D. Lee}, {\em
  Multimodal image reconstruction using supplementary structural information in
  total variation regularization}, Sensing and Imaging, 15 (2014), p.~97,
  \url{https://doi.org/10.1007/s11220-014-0097-5}.

\bibitem{Lalwanietal}
{\sc G.~Lalwani, J.~Livingston~Sundararaj, K.~Schaefer, T.~Button, and
  B.~Sitharaman}, {\em Synthesis, characterization, in vitro phantom imaging,
  and cytotoxicity of a novel graphene-based multimodal magnetic resonance
  imaging - {X}-ray computed tomography contrast agent}, Journal of Materials
  Chemistry. B, 2 (2015), p.~3519–3530,
  \url{https://doi.org/10.1039/C4TB00326H}.

\bibitem{landi2008lagrange}
{\sc G.~Landi}, {\em The {L}agrange method for the regularization of discrete
  ill-posed problems}, Computational Optimization and Applications, 39 (2008),
  pp.~347--368.

\bibitem{langer2017automated}
{\sc A.~Langer}, {\em Automated parameter selection for total variation
  minimization in image restoration}, Journal of Mathematical Imaging and
  Vision, 57 (2017), pp.~239--268.

\bibitem{LI2016a}
{\sc H.~Li, M.~Gong, Q.~Wang, J.~Liu, and L.~Su}, {\em A multiobjective fuzzy
  clustering method for change detection in {SAR} images}, Applied Soft
  Computing, 46 (2016), pp.~767--777,
  \url{https://doi.org/10.1016/j.asoc.2015.10.044}.

\bibitem{Liu_2012}
{\sc Y.~Liu, J.~Ma, Y.~Fan, and Z.~Liang}, {\em Adaptive-weighted total
  variation minimization for sparse data toward low-dose {X}-ray computed
  tomography image reconstruction}, Physics in Medicine and Biology, 57 (2012),
  pp.~7923--7956, \url{https://doi.org/10.1088/0031-9155/57/23/7923}.

\bibitem{martinez2014edge}
{\sc A.~Martinez, A.~Gelb, and A.~Gutierrez}, {\em Edge detection from
  non-uniform {F}ourier data using the convolutional gridding algorithm},
  Journal of Scientific Computing, 61 (2014), pp.~490--512.

\bibitem{mcdermid2008object}
{\sc G.~J. McDermid, J.~Linke, A.~D. Pape, D.~N. Laskin, A.~J. McLane, and
  S.~E. Franklin}, {\em Object-based approaches to change analysis and thematic
  map update: Challenges and limitations}, Canadian Journal of Remote Sensing,
  34 (2008), pp.~462--466.

\bibitem{papari2008adaptive}
{\sc G.~Papari and N.~Petkov}, {\em Adaptive pseudo dilation for gestalt edge
  grouping and contour detection}, IEEE Transactions on Image Processing, 17
  (2008), pp.~1950--1962.

\bibitem{papari2011edge}
{\sc G.~Papari and N.~Petkov}, {\em Edge and line oriented contour detection:
  State of the art}, Image and Vision Computing, 29 (2011), pp.~79--103.

\bibitem{Rigie2015joint}
{\sc D.~S. Rigie and P.~J.~L. Rivi{\`{e}}re}, {\em Joint reconstruction of
  multi-channel, spectral {CT} data via constrained total nuclear variation
  minimization}, Physics in Medicine and Biology, 60 (2015), pp.~1741--1762,
  \url{https://doi.org/10.1088/0031-9155/60/5/1741}.

\bibitem{sanders2020effective}
{\sc T.~Sanders, R.~B. Platte, and R.~D. Skeel}, {\em Effective new methods for
  automated parameter selection in regularized inverse problems}, Applied
  Numerical Mathematics, 152 (2020), pp.~29--48.

\bibitem{scarnati2019accelerated}
{\sc T.~Scarnati and A.~Gelb}, {\em Accurate and efficient image reconstruction
  from multiple measurements of {F}ourier samples}, 2020.

\bibitem{shchukina2017pitfalls}
{\sc A.~Shchukina, P.~Kasprzak, R.~Dass, M.~Nowakowski, and K.~Kazimierczuk},
  {\em Pitfalls in compressed sensing reconstruction and how to avoid them},
  Journal of Biomolecular NMR, 68 (2017), pp.~79--98,
  \url{https://doi.org/10.1007/s10858-016-0068-3}.

\bibitem{sobel19683x3}
{\sc I.~Sobel and G.~Feldman}, {\em A $3\times 3$ isotropic gradient operator
  for image processing}, A Talk at the Stanford Artificial Project,  (1968),
  pp.~271--272.

\bibitem{song2020multimodal}
{\sc P.~{Song}, X.~{Deng}, J.~F.~C. {Mota}, N.~{Deligiannis}, P.~L. {Dragotti},
  and M.~R.~D. {Rodrigues}}, {\em Multimodal image super-resolution via joint
  sparse representations induced by coupled dictionaries}, IEEE Transactions on
  Computational Imaging, 6 (2020), pp.~57--72,
  \url{https://doi.org/10.1109/TCI.2019.2916502}.

\bibitem{song2016coupled}
{\sc P.~{Song}, J.~F.~C. {Mota}, N.~{Deligiannis}, and M.~R.~D. {Rodrigues}},
  {\em Coupled dictionary learning for multimodal image super-resolution}, in
  IEEE Global Conference on Signal and Information Processing (GlobalSIP),
  2016, pp.~162--166, \url{https://doi.org/10.1109/GlobalSIP.2016.7905824}.

\bibitem{TEWKESBURY2015a}
{\sc A.~P. Tewkesbury, A.~J. Comber, N.~J. Tate, A.~Lamb, and P.~F. Fisher},
  {\em A critical synthesis of remotely sensed optical image change detection
  techniques}, Remote Sensing of Environment, 160 (2015), pp.~1--14,
  \url{https://doi.org/10.1016/j.rse.2015.01.006}.

\bibitem{thonfeld2016robust}
{\sc F.~Thonfeld, H.~Feilhauer, M.~Braun, and G.~Menz}, {\em Robust {C}hange
  {V}ector {A}nalysis {(RCVA)} for multi-sensor very high resolution optical
  satellite data}, International Journal of Applied Earth Observation and
  Geoinformation, 50 (2016), pp.~131--140.

\bibitem{tipping2001sparse}
{\sc M.~E. Tipping}, {\em Sparse {B}ayesian learning and the relevance vector
  machine}, Journal of Machine Learning Research, 1 (2001), pp.~211--244.

\bibitem{viswanathan2012iterative}
{\sc A.~Viswanathan, A.~Gelb, and D.~Cochran}, {\em Iterative design of
  concentration factors for jump detection}, Journal of Scientific Computing,
  51 (2012), pp.~631--649.

\bibitem{vogel2002Computational}
{\sc C.~R. Vogel}, {\em Regularization parameter selection methods}, Society
  for Industrial and Applied Mathematics, 2002, ch.~7, pp.~97--127,
  \url{https://doi.org/10.1137/1.9780898717570.ch7}.

\bibitem{wen2011parameter}
{\sc Y.-W. Wen and R.~H. Chan}, {\em Parameter selection for
  total-variation-based image restoration using discrepancy principle}, IEEE
  Transactions on Image Processing, 21 (2011), pp.~1770--1781.

\bibitem{wipf2004sparse}
{\sc D.~P. Wipf and B.~D. Rao}, {\em Sparse {B}ayesian learning for basis
  selection}, IEEE Transactions on Signal Processing, 52 (2004),
  pp.~2153--2164.

\bibitem{xie2014reweighted}
{\sc W.~Xie, Y.~Deng, K.~Wang, X.~Yang, and Q.~Luo}, {\em Reweighted $\ell_1$
  regularization for restraining artifacts in {FMT} reconstruction images with
  limited measurements}, Optics Letters, 39 (2014), pp.~4148--4151.

\bibitem{yang2015tv}
{\sc X.~Yang, R.~Hofmann, R.~Dapp, T.~Van~de Kamp, T.~dos Santos~Rolo, X.~Xiao,
  J.~Moosmann, J.~Kashef, and R.~Stotzka}, {\em {TV}-based conjugate gradient
  method and discrete {L}-curve for few-view {CT} reconstruction of {X}-ray in
  vivo data}, Optics Express, 23 (2015), pp.~5368--5387.

\bibitem{YE2016a}
{\sc S.~Ye, D.~Chen, and J.~Yu}, {\em A targeted change-detection procedure by
  combining change vector analysis and post-classification approach}, ISPRS,
  114 (2016), pp.~115--124,
  \url{https://doi.org/10.1016/j.isprsjprs.2016.01.018}.

\end{thebibliography}
